\documentclass[11pt]{amsart}

\usepackage[a4paper,margin=2.5cm]{geometry}
\usepackage{amsmath,amssymb,amsfonts,amsthm,amscd,amsxtra,mathtools}
\usepackage{mathrsfs}
\usepackage{extarrows}
\usepackage{bbm}
\usepackage{graphicx}
\usepackage{xcolor}
\usepackage{float}
\usepackage{enumitem}
\usepackage{microtype}
\usepackage{verbatim}
\usepackage{booktabs}
\usepackage{tabularx}
\usepackage{array}
\usepackage{tikz}
\usetikzlibrary{patterns,arrows.meta,decorations.pathreplacing}
\usetikzlibrary{fit}
\usepackage[colorlinks=true,linkcolor=blue,citecolor=blue,urlcolor=blue]{hyperref}
\allowdisplaybreaks
\setlist{nosep,leftmargin=2.2em}
\numberwithin{equation}{section}

\newcommand{\R}{\mathbb R}
\newcommand{\Z}{\mathbb Z}
\newcommand{\N}{\mathbb N}
\newcommand{\T}{\mathbb T}

\newcommand{\E}{\mathbb E}
\renewcommand{\P}{\mathbb P}
\newcommand{\ind}{\mathbbm{1}}
\newcommand{\slitN}{C^{N,\lambda}}
\newcommand{\pto}{\xrightarrow{\P}}
\newcommand{\dto}{\Rightarrow}

\newcommand{\im}{\operatorname{Im}}
\newcommand{\rea}{\operatorname{Re}}

\newcommand{\dist}{\operatorname{dist}}
\newcommand{\eps}{\varepsilon}
\newcommand{\id}{\operatorname{id}}

\newcommand{\chl}{\mathcal C}

\theoremstyle{plain}
\newtheorem{thm}{Theorem}[section]
\newtheorem{proposition}[thm]{Proposition}
\newtheorem{lemma}[thm]{Lemma}
\newtheorem{corollary}[thm]{Corollary}

\theoremstyle{definition}
\newtheorem{definition}[thm]{Definition}
\newtheorem{assumption}[thm]{Assumption}

\theoremstyle{remark}
\newtheorem{remark}[thm]{Remark}

\title{A Brownian web invariance principle:\\ Localized disturbance flows and Hastings--Levitov(0)}

\author{Guanyi Chen$^1$}
\address[Guanyi Chen]{School of Mathematical Sciences, Peking University, Beijing, China}
\email{cgy7698@stu.pku.edu.cn}
\thanks{$^1$School of Mathematical Sciences, Peking University, Beijing, China}

\author{Xinyi Li$^2$}
\address[Xinyi Li]{BICMR, Peking University, Beijing, China}
\email{xinyili@bicmr.pku.edu.cn}
\thanks{$^2$BICMR, Peking University, Beijing, China}

\author{Eviatar B. Procaccia$^3$}
\address[Eviatar B. Procaccia]{Faculty of Data and Decision Sciences, Technion--Israel Institute of Technology, Haifa, Israel}
\email{eviatarp@technion.ac.il}
\thanks{$^3$Faculty of Data and Decision Sciences, Technion--Israel Institute of Technology, Haifa, Israel}

\author{Yuxuan Zong$^1$}
\address[Yuxuan Zong]{School of Mathematical Sciences, Peking University, Beijing, China}

\email{yxzong25@stu.pku.edu.cn}

\begin{document}

    \begin{abstract}
    We establish a Brownian web invariance principle for localized disturbance flows,  providing a criterion that strengthens the weak-flow convergence of~\cite{norris2015weak} to path-Hausdorff convergence of their full trajectory collections.  We apply this approach to the radial, stationary, and cylindrical versions of the Hastings--Levitov$(0)$ aggregation model, proving path-Hausdorff convergence of the full harmonic measure path collections to the (Cylindrical) Brownian web, together with joint convergence of the full gap--finger collections to geometric embeddings of the corresponding double Brownian webs. This extends the radial HL$(0)$ small-particle limits of~\cite{norris2012hastings} for the harmonic measure flow and finite gap--finger collections. In the radial and cylindrical HL$(0)$, we also show that the rescaled one-arm domination time converges to a multiple of the global coalescence time of coalescing Brownian motion on $\mathbb T$, with convergence of expectations. The main new ingredient is a three-point survival estimate, proved via a multiscale stopping-time decomposition.

    \end{abstract}

    \maketitle
    \tableofcontents

    \section{Introduction and main results}\label{S.introduction}

    The Brownian web (BW), introduced by Arratia~\cite{arratia1979coalescing}, is a random compact set of paths consisting of coalescing standard Brownian motions started from every point in the space--time plane $\R^2$.  Following the metric-space characterization and the path-Hausdorff convergence criteria of~\cite{toth1998true,fontes2004brownian,newman2005convergence}, the BW has become a universal diffusive scaling limit of various (1+1)-dimensional coalescing systems~\cite{berestycki2015coalescing,ferrari2005two,coletti2009scaling,coletti2014convergence,coupier2021directedspanningforest}; see the survey~\cite{schertzer2017brownian} for an overview. The Cylindrical Brownian web (CBW) is its periodic analogue, in which space is a circle~\cite{coupier2019directed}.

    The starting point is the framework of localized disturbance flows of Norris and Turner~\cite{norris2015weak}. A localized disturbance flow is obtained by composing i.i.d.\ order-preserving maps, each of which is a small, localized perturbation of the identity.    Norris and Turner showed that, as the disturbances become small and spatially localized, finite collections of trajectories converge to coalescing Brownian motions, and the flows converge to the coalescing Brownian flow in the weak-flow topology, both on the circle and locally on the line. These results leave open path-Hausdorff convergence of the full trajectory collection, which requires ruling out extra paths in the limit, leading to our first question:

    \par\smallskip
    \textbf{Question 1:} \emph{What conditions ensure convergence of the full trajectory collections of localized disturbance flows to the BW or CBW in the path-Hausdorff topology?}
    \par\smallskip

    The aggregation models considered below provide a natural setting for this question, as their boundary evolution can be described in harmonic-measure coordinates.  Diffusion-limited aggregation (DLA), introduced by Witten and Sander~\cite{witten1981diffusion}, attaches particles according to harmonic measure and produces highly ramified clusters. Fundamental questions about the growth rate,  scaling limits, the number of arms, and domination times remain difficult even in flat or cylindrical geometries~\cite{kesten1987hitting,kesten1987long,kol1998solution,benjamini2008diffusion,marchetti2012stationary}.  For DLA grown from a long line and its stationary version, see~\cite{procaccia2019stationary,procaccia2020stationary,procaccia2021dimension,mu2022scaling}.  Hastings and Levitov \cite{hastings1998laplacian}  introduced an off-lattice analogue of this growth process, in which the cluster is encoded by a composition of conformal maps.  This construction is also closely related to Loewner--Kufarev evolution \cite{carleson2001aggregation,Johansson2012scaling,norris2024stability}. In radial HL$(0)$, the aggregate is obtained by composing i.i.d.\ rotations of a fixed slit map. Equivalently, particles attach at i.i.d.\ positions in harmonic-measure coordinates, and the evolution of the boundary in these coordinates is a localized disturbance flow.

    We consider three HL$(0)$ geometries: radial, stationary, and cylindrical HL(0), denoted by HL(0), SHL(0), and CHL(0), respectively~\cite{hastings1998laplacian, norris2012hastings, berger2022growth,berger2025logarithmic, procaccia2023cylindrical,chen2025one}.  In all three models, each new particle induces a one-dimensional disturbance in harmonic-measure coordinates. We study three families of paths (see Figure~\ref{fig:hl-observables}): the \emph{harmonic measure paths},  which track the harmonic-measure coordinates of boundary points forward in aggregation time, with harmonic measure describing the hitting distribution of incoming particles on the current tree; the \emph{fingers}, which run backward along the particle tree from the tip of each particle to the initial boundary; and the \emph{gaps}, which start at the tip of each particle and run to infinity through the complement of the complete particle tree.

    For radial HL$(0)$, Norris and Turner~\cite{norris2012hastings} established the small-particle limit, including weak-flow convergence of the harmonic measure flow and joint convergence of finite collections of gaps and fingers. We now turn to the full path collections generated by all particles of the random aggregate. This leads to our second question:

    \par\smallskip
    \textbf{Question 2:} \emph{What are the scaling limits of the full collections of harmonic measure paths, fingers, and gaps in the three HL$(0)$ geometries?}
    \par\smallskip

    In this paper, we give sufficient conditions for the path-Hausdorff convergence in \textbf{Question 1} and identify the joint scaling limits in \textbf{Question 2} as the corresponding Brownian webs and their geometric embeddings. These results concern the full path collections, including paths starting at random points determined by the aggregate. We also show that, in radial and cylindrical HL$(0)$, the rescaled one-arm domination time converges in distribution to an explicit multiple of the global coalescence time of coalescing Brownian motions on the circle, with convergence of expectations. Beyond HL$(0)$, our approach provides a framework for studying BW convergence in other correlated models.

    \subsection{Main results I: localized disturbance flows}\label{SS.BW-invariance-principle}
    In this section, we define localized disturbance flows and give sufficient conditions for their full collections of forward trajectories to converge to  the (Cylindrical) BW in the path-Hausdorff topology.

    Let $\Xi\in\{\R,\T\}$, where $\T=\R/(2\pi\Z)$, and let $d_\Xi$ be the Euclidean distance if $\Xi=\R$ and the circular distance if $\Xi=\T$.  For each integer $n\ge1$, let $g_n:\R\to\R$ be a non-decreasing right-continuous map; when $\Xi=\T$, we require in addition that $g_n(x+2\pi)=g_n(x)+2\pi$ for all $x$, so that $g_n$ induces a map of $\T$. Define the deviation  $\Delta_n(x)=g_n(x)-x$, which describes how $g_n$ moves each point $x$.  For $\theta\in\R$, define the disturbance translated by $\theta$ and its deviation by
    \begin{align}\label{map-deviation}
        g_{n,\theta}(x):=g_n(x-\theta)+\theta, \qquad \Delta_{n,\theta}(x):=g_{n,\theta}(x)-x=\Delta_n(x-\theta).
    \end{align}
    For each $n$, let $\rho_n>0$ be an intensity parameter, and let $\mathcal P_n$ be a Poisson point process (PPP) of disturbance times and locations $(t,\theta)$. When $\Xi=\T$, $\mathcal P_n$ has intensity $\rho_n\,(2\pi)^{-1}dt\,d\theta$ on $\R\times[0,2\pi)$, so disturbances occur at total rate $\rho_n$. For fixed times $s,t\in\R$ with $s\le t$, if the points of $\mathcal P_n\cap((s,t]\times[0,2\pi))$ are $(t_i,\theta_i)_{1\le i\le m}$ with $t_1<\cdots<t_m$, define
    \[
        \Phi_{s,t}^{n}:=g_{n,\theta_m}\circ\cdots\circ g_{n,\theta_1}, \qquad \Phi_{s,s}^{n}:=\id.
    \]
    When $\Xi=\mathbb R$, $\mathcal P_n$ has intensity $\rho_n\,dt\,d\theta$ on $\mathbb R\times\mathbb R$, so infinitely many disturbances occur in every time interval. We assume that, as $R\uparrow\infty$, the ordered compositions over $\mathcal P_n\cap((s,t]\times[-R,R])$ converge in probability uniformly on compact time intervals for each deterministic starting point to a jointly measurable flow $\Phi^n=(\Phi^n_{s,t}:s\le t)$ of non-decreasing right-continuous maps, with adapted c\`adl\`ag trajectories, stationary independent increments and the cocycle property $\Phi^n_{s,u}=\Phi^n_{t,u}\circ\Phi^n_{s,t}$ for $s\le t\le u$. We impose the zero-drift condition (the first identity below) and define the covariance rate $K_n$ by
    \begin{align}\label{drift-covariance}
        \int_\Xi \Delta_n(y)\,\mu_\Xi(dy)=0, \qquad K_n(h):=\rho_n\int_\Xi \Delta_n(y)\Delta_n(y+h)\,\mu_\Xi(dy),
    \end{align}
    where $\mu_{\mathbb R}(dy)=dy$ and $\mu_{\mathbb T}(dy)=(2\pi)^{-1}dy$ on $[0,2\pi)$. When $\Xi=\mathbb R$, the first integral in \eqref{drift-covariance} is understood as a symmetric principal value.  In words, a single trajectory has no drift, and $K_n(h)$ is the covariance rate of two trajectories at distance $h$; in particular, $K_n(0)$ is the variance rate of a single trajectory. Now we define localized disturbance flows as follows:

    \begin{definition}[Localized disturbance flows]\label{def:abstract-localized}
    The flows $\Phi^n$ satisfying the construction and zero-drift assumptions above are called localized disturbance flows with limiting variance $\sigma^2\in(0,\infty)$ if the following conditions hold.
    \begin{enumerate}[label=\textup{(L\arabic*)},leftmargin=2.7em]
    \item \textbf{(Small jumps)} There exists $j_n\downarrow0$ such that $\|\Delta_n\|_\infty\le j_n$.
    \item  \textbf{(Limiting variance)}  $K_n(0)\rightarrow\sigma^2.$
    \item \textbf{(Localization)} For every $c_n\downarrow0$ with $j_n=o(c_n)$,
    \[
        \kappa_n(c_n):= \sup_{d_\Xi(h,0)\ge c_n}|K_n(h)|\rightarrow0.
    \]
    \end{enumerate}
    \end{definition}

    For $d>0$, define the quadratic-variation rate $Q_n(d)$ of the gap between two trajectories at distance $d$ and their single-disturbance coalescence rate $\Lambda_n(d)$ by
    \begin{align}
        Q_n(d)&:=\rho_n\int_\Xi\bigl(\Delta_n(z+d)-\Delta_n(z)\bigr)^2\,\mu_\Xi(dz)=2\bigl(K_n(0)-K_n(d)\bigr),\label{eq:abstract-Q-covariance-identity} \\
        \Lambda_n(d)&:=\rho_n\,\mu_\Xi\{z\in\Xi:g_n(z+d)=g_n(z)\}.\label{eq:abstract-Q-coalescence-rates}
    \end{align}
    For $t_0\in\R$, $t>0$, and $A\subset\Xi$, define two sets of positions at time $t_0+t$. The first consists of positions reached by trajectories started at grid times $s\le t_0$ and passing through $A$ at time $t_0$; the second consists of positions reached from $A$ at time $t_0$:
    \[
        \begin{aligned}
        \mathsf D_{\mathrm{traj}}^{\Phi^n}(t_0,t;A) &:= \bigl\{\Phi^n_{s,t_0+t}(x): s\in n^{-1}\mathbb Z,\ s\le t_0,\ x\in\Xi,\ \Phi^n_{s,t_0}(x)\in A\bigr\}, \\
        \mathsf D^{\Phi^n}(t_0,t;A) &:= \{\Phi^n_{t_0,t_0+t}(y):y\in A\}.
        \end{aligned}
    \]
    Let $\mathcal F_t$ denote the natural filtration generated by the driving disturbances up to time $t$.     For $A\subset\Xi$ and $q>0$, write $A^q:=\{x\in\Xi:d_\Xi(x,A)\le q\}$ for the closed $q$-neighbourhood of $A$. For an interval or arc $A$, let $|A|$ denote its length.  For subsets $E,F\subset\Xi$, define the one-sided Hausdorff distance
    \begin{align*}
        d_H^\to(E,F):=\sup_{x\in E}\inf_{y\in F}d_\Xi(x,y),
    \end{align*}
    with $d_\Xi(x,\varnothing)=+\infty$, $d_H^\to(\varnothing,F)=0$, and $d_H^\to(E,\varnothing)=+\infty$ for $E\ne\varnothing$, so that $d_H^\to(E,F)\le q$ if and only if $E\subset F^q$.

    Our main theorem requires two further conditions.  The first is a uniform positive lower bound on a combination of the gap quadratic variation and the coalescence rate; it ensures that small gaps between adjacent paths close quickly. The second is a type of finite covering condition: for each fixed $t>0$, the positions at time $t_0+t$ of all paths passing through an interval or arc $A$ at time $t_0$ can, with high probability, be approximated by the images of a finite set of points chosen at time $t_0$, with a vanishing error and a controlled number of points.
    \begin{assumption}\label{ass:abstract-no-extra-path}
    There exist positive sequences $a_n,r_n,q_n\downarrow0$, and $M_n\uparrow\infty$ such that
    \begin{equation}\label{parameters-control}
        j_n=o(a_n),\qquad r_n=O(\sqrt{a_n}), \qquad \frac{a_n\log M_n}{r_n}\to0,
    \end{equation}
    and the following conditions hold.
    \begin{enumerate}[leftmargin=2.9em]
    \item[\textup{(I)}] \textbf{(Small-gap  bound)} There exists a constant $c_0>0$ such that for all sufficiently large $n$,
    \begin{equation}\label{eq:abstract-micro-coalescence}
        \inf_{0<d\le a_n} \{Q_n(d)+j_nr_n\Lambda_n(d)\}\ge c_0.
    \end{equation}

    \item[\textup{(II)}] \textbf{(Finite covering)} For every $t_0\in\mathbb R$ and every interval or arc $A\subset\Xi$ with $|A|\le1$, there is an $\mathcal F_{t_0}$-measurable finite set $\mathcal R_n(t_0,A)\subset A^{q_n}$ such that
    \begin{equation}\label{eq:abstract-shadow-size}
        \sup_{t_0}\sup_{|A|\le1} \mathbb P(|\mathcal R_n(t_0,A)|>M_n)\rightarrow0.
    \end{equation}
    For every $t>0$,
    \begin{equation}\label{eq:abstract-shadow-cover}
        \sup_{t_0}\sup_{|A|\le1} \mathbb P\left\{ d_H^\to\left( \mathsf D_{\mathrm{traj}}^{\Phi^n}(t_0,t;A), \mathsf D^{\Phi^n}(t_0,t;\mathcal R_n(t_0,A)) \right)>q_n \right\}\rightarrow0.
    \end{equation}
    \end{enumerate}
    \end{assumption}

    To state the path-Hausdorff convergence, we interpolate the flow trajectories at grid times $n^{-1}\Z$. Define $p_\Xi:\R\to\Xi$ by $p_\R(x)=x$ and $p_\T(x)=x+2\pi\Z$.  For $s\in n^{-1}\mathbb Z$, $x\in\Xi$, and $t\ge s$, let $k$ be the unique integer such that $t\in[k/n,(k+1)/n)$.  When $\Xi=\T$, the interpolation is performed on $\mathbb R$, using a real representative of $x$ and the compositions of the maps $g_{n,\theta}$. Define the grid interpolation of $t\mapsto\Phi^n_{s,t}(x)$ by
    \[
        I_n\Phi^n_{s,t}(x):= p_\Xi\!\left((k+1-nt)\Phi^n_{s,k/n}(x) +(nt-k)\Phi^n_{s,(k+1)/n}(x)\right).
    \]
    Thus $I_n\Phi^n_{s,\cdot}(x)$ equals $\Phi^n_{s,\cdot}(x)$ at grid times $n^{-1}\Z$, with linear interpolation between adjacent grid times. Define the interpolated spatial trajectory collection
    \[
        \mathcal X_n^\Phi:= \overline{\bigl\{t\mapsto I_n\Phi^n_{s,t}(x): x\in\Xi,\ s\in n^{-1}\mathbb Z,\ t\ge s\bigr\}}.
    \]

    Let $\mathcal W_{\Xi}$ denote the standard BW $\mathcal W$ when $\Xi=\R$, and the standard CBW $\mathcal W^{\uparrow}$ when $\Xi=\T$. Their precise constructions are recalled in Section~\ref{S.BW}. For $\sigma>0$, $\sigma\mathcal W_\Xi$ denotes the BW or CBW on the same spatial domain $\Xi$ whose individual paths have variance $\sigma^2$.  Our first main result is the following BW invariance principle.

    \begin{thm}\label{thm:abstract-boundary}
    Let $\Phi^n$ be localized disturbance flows with limiting variance $\sigma^2$ such that Assumption~\ref{ass:abstract-no-extra-path} holds. Then as $n\to\infty$,
    \[
        \mathcal X_n^\Phi\dto\sigma\mathcal W_{\Xi}
    \]
    in the path-Hausdorff topology over $\Xi$.
    \end{thm}

    \subsection{Main results II: \texorpdfstring{Hastings--Levitov}{Hastings--Levitov}\texorpdfstring{$(0)$}{(0)}}\label{SS.HL-invariance-principle}

    We now state the resulting BW invariance principles for the three HL$(0)$ geometries.  The models and observables are defined precisely in Subsection~\ref{S.HL-models}.

    Fix $\eta,\lambda>0$. We use $\mathcal A\in\{\mathcal H,\mathcal S,\mathcal C\}$ to index radial HL$(0)$, SHL$(0)$, and CHL$(0)$, respectively. The boundary spaces are given before rescaling, with $\mathbb T_N:=\mathbb R/(2\pi N\mathbb Z)$. The table below lists the model parameters and space-time rescalings at scale $N$.
    \begin{center}
    \small
    \setlength{\tabcolsep}{3.2pt}
    \renewcommand{\arraystretch}{1.05}
    \begin{tabularx}{\textwidth}
    {@{} >{\raggedright\arraybackslash}p{0.07\textwidth} >{\centering\arraybackslash}p{0.05\textwidth} >{\centering\arraybackslash}p{0.12\textwidth} >{\raggedright\arraybackslash}X >{\centering\arraybackslash}p{0.14\textwidth} >{\raggedright\arraybackslash}p{0.25\textwidth} @{}}
    \toprule
    model & $\mathcal A$ & boundary & model parameters & time scale & rescaled coordinates\\
    \midrule
    $\mathrm{HL}(0)$ & $\mathcal H$ & $\mathbb T$ & particle size $N^{-1}\eta$ & $N^3t$ & $\arg z$, $N^{-1}\log |z|$\\
    $\mathrm{SHL}(0)$ & $\mathcal S$ & $\mathbb R$ & particle height $\lambda$ & $N^2t$ & $\operatorname{Re} z/N$, $\operatorname{Im} z/N^2$\\
    $\mathrm{CHL}(0)$ & $\mathcal C$ & $\mathbb T_{N}$ & particle height $\lambda$, period $2\pi N$ & $N^2t$ & $\operatorname{Re} z/N \bmod 2\pi$, $\operatorname{Im} z/N^2$\\
    \bottomrule
    \end{tabularx}
    \end{center}
    The different time scales reflect the Poisson arrival rates: HL$(0)$ has unit total rate, whereas SHL$(0)$ and CHL$(0)$ have unit rate per unit boundary length. 
    
\begin{figure}[htbp]
      \centering
      \input{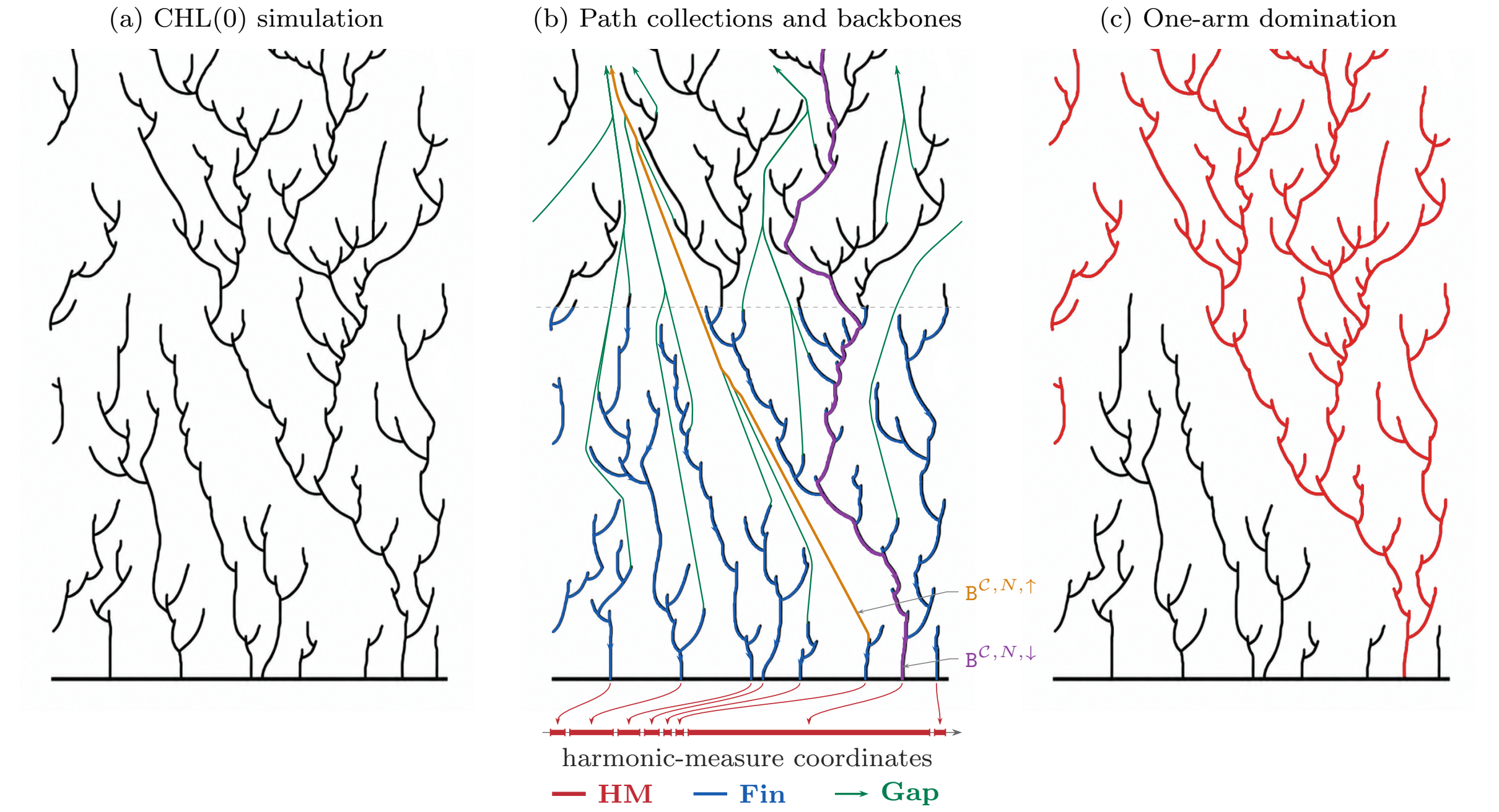}
      \caption{A CHL$(0)$ simulation and illustration of its observables. In (b), the red arrows indicate the harmonic-measure intervals corresponding to the arms, while the blue fingers and green gaps are geometric paths. The backbone and dual backbone are highlighted in purple and orange. In (c), the infinite arm is highlighted in red.}
      \label{fig:hl-observables}
    \end{figure}

    For each $\mathcal A\in\{\mathcal H,\mathcal S,\mathcal C\}$ and $N$, let $\mathbf{HM}^{\mathcal A,N}$, $\mathbf{Fin}^{\mathcal A,N}$, and $\mathbf{Gap}^{\mathcal A,N}$ denote the rescaled full path collections of Definition~\ref{def:process-generated-observables}. Here ``full'' means that paths from every particle in the growing aggregate are included. For each particle, two forward harmonic measure paths parametrized by time start from the two sides of its attachment on the boundary, and the interval between them tracks the harmonic measure of the subtree rooted at that particle as the aggregate grows.  $\mathbf{Fin}^{\mathcal A,N}$ is generated by the fingers, that is, backward paths along the particle tree from the particle tips to the initial boundary; and $\mathbf{Gap}^{\mathcal A,N}$ is generated by the gaps, that is, paths that leave particle tips and escape to infinity through the complement of the complete particle tree. Figure~\ref{fig:hl-observables}(a)--(b) shows a CHL$(0)$ simulation and an illustration of the corresponding path collections.

    Define the growth rates $v_{\mathcal A}$ and the variances     $\sigma_{\mathcal A}^2$ by
    \begin{align}\label{growth-rates-variance}
        v_{\mathcal A}:=
        \begin{cases}
        \dfrac{\eta^2}{4}, & \mathcal A=\mathcal H, \\[2mm]
        \dfrac{\pi\lambda^2}{2}, & \mathcal A\in\{\mathcal S,\mathcal C\},
        \end{cases}
        \qquad \sigma_{\mathcal A}^2:=
        \begin{cases}
        \dfrac{2\eta^3}{3\pi}, & \mathcal A=\mathcal H, \\[2mm]
        \dfrac{4\lambda^3}{3}, & \mathcal A\in\{\mathcal S,\mathcal C\}.
        \end{cases}
    \end{align}
    Let $(\mathcal W_{\mathcal A,\ge0}^{\uparrow}, \mathcal W_{\mathcal A,\ge0}^{\downarrow})$ denote the restrictions to nonnegative times of the standard double CBW for $\mathcal A\in\{\mathcal H,\mathcal C\}$ and of the standard double BW for $\mathcal A=\mathcal S$, as defined in Subsection~\ref{SS.convergence criteria}.  Our first result for HL$(0)$ concerns the full collection of harmonic measure paths.
    \begin{thm}\label{thm:harmonic-measure-webs}
    For every $\mathcal A\in\{\mathcal H,\mathcal S,\mathcal C\}$, as $N\to\infty$,
    \[
        \mathbf{HM}^{\mathcal A,N} \dto \sigma_{\mathcal A}\mathcal W_{\mathcal A,\ge0}^{\uparrow}
    \]
    in the path-Hausdorff topology defined in Subsection~\ref{SS.convergence criteria}.
    \end{thm}

    The fingers and gaps take values in the rescaled geometric spaces $\mathbb G_{\mathcal H}=\{z\in\mathbb C:|z|\ge1\}$, $\mathbb G_{\mathcal S}=\{z\in\mathbb C:\operatorname{Im}z\ge0\}$, and $\mathbb G_{\mathcal C}=\mathbb T\times[0,\infty)$, respectively; see Subsection~\ref{S.HL-models}. To describe their limits, we embed the harmonic-measure paths in $\mathbb G_{\mathcal A}$ using the growth rate $v_{\mathcal A}>0$, sending time $t$ to height $v_{\mathcal A}t$ in SHL$(0)$ and CHL$(0)$, and to radius $e^{v_{\mathcal A}t}$ in HL$(0)$.  Write $E^{\mathcal H}=E^{\mathcal C}=\mathbb T$ and $E^{\mathcal S}=\mathbb R$. For $v, t\ge0$, and $x\in E^{\mathcal A}$, define the geometric embedding by
    \begin{equation*}
        \mathfrak e_{\mathcal A}^{v}(t,x)=
        \begin{cases}
        \exp(vt+\mathrm ix),&\mathcal A=\mathcal H, \\
        x+\mathrm ivt,&\mathcal A\in\{\mathcal S,\mathcal C\}.
        \end{cases}
    \end{equation*}
    For CHL(0), $x+\mathrm ih$ denotes $(x,h)\in\mathbb T\times[0,\infty)$. For a forward or backward path collection $K$, let $\star\in\{\uparrow,\downarrow\}$ denote its time direction and define
    \[
        \mathfrak E_v^{\mathcal A,\star}(K) :=\overline{\left\{ u\mapsto\mathfrak e_{\mathcal A}^{v}(u,\pi(u)): \pi\in K \right\}},
    \]
    where each path is embedded on its time domain and the closure is taken in the corresponding  path topology of Subsection~\ref{S.HL-models}.  The following theorem establishes joint convergence of the full finger and gap collections to geometric embeddings of the corresponding double BW or CBW.
    \begin{thm}\label{thm:joint-geometric-webs}
    For every $\mathcal A\in\{\mathcal H,\mathcal S,\mathcal C\}$, as $N\to\infty$,
    \[
        (\mathbf{Fin}^{\mathcal A,N},\mathbf{Gap}^{\mathcal A,N}) \dto \left( \mathfrak E^{\mathcal A,\downarrow}_{v_{\mathcal A}} (\sigma_{\mathcal A}\mathcal W_{\mathcal A,\ge0}^{\downarrow}), \mathfrak E^{\mathcal A,\uparrow}_{v_{\mathcal A}} (\sigma_{\mathcal A}\mathcal W_{\mathcal A,\ge0}^{\uparrow}) \right)
    \]
    in the product topology induced by $D_{\mathcal A,H}$ and $D_{\mathcal A,H}^{\uparrow}$ of Subsection~\ref{S.HL-models}.
    \end{thm}
    The convergences in Theorems~\ref{thm:harmonic-measure-webs} and~\ref{thm:joint-geometric-webs} hold jointly, with the same double Brownian web in the limits.

    For HL$(0)$ and CHL$(0)$, let $C_{\mathcal A}^{\downarrow}$ and $C_{\mathcal A}^{\uparrow}$ be the unique  bi-infinite paths of the double CBW on $\mathbb T\times\mathbb R$ with variance $\sigma_{\mathcal A}^{2}$; their joint law is described in \cite[Theorems~3.1 and~3.6]{coupier2019directed}.  The backbone $\mathtt{B}^{\mathcal A,N,\downarrow}$ is the unique infinite path in the particle tree, while the selected dual backbone $\mathtt{B}^{\mathcal A,N,\uparrow}$ is a gap chosen using an independent uniform boundary point. Their precise definitions are given in Subsection~\ref{S.HL-models} (see Figure~\ref{fig:hl-observables}(b)). Their joint scaling limit is given by the following corollary.
    \begin{corollary}\label{cor:backbone-dual-backbone}
    For every $\mathcal A\in\{\mathcal H,\mathcal C\}$, as $N\to\infty$,
    \[
        \bigl(\mathtt{B}^{\mathcal A,N,\downarrow},\mathtt{B}^{\mathcal A,N,\uparrow}\bigr) \dto \bigl( u\mapsto\mathfrak e_{\mathcal A}^{v_{\mathcal A}}(u,C_{\mathcal A}^{\downarrow}(u)), u\mapsto\mathfrak e_{\mathcal A}^{v_{\mathcal A}}(u,C_{\mathcal A}^{\uparrow}(u)) \bigr),\qquad u\ge0,
    \]
    in the product topology induced by $D_{\mathcal A}$ and $D_{\mathcal A}^{\uparrow}$ of Subsection~\ref{S.HL-models}.
    \end{corollary}

    We next apply these scaling limits to the one-arm domination time in HL$(0)$ and CHL$(0)$. An \emph{arm} is a connected component of the particle tree after removing the initial boundary; see Definition~\ref{def:one-arm-domination-time}. By \cite{rohde2005some,norris2012hastings,chen2025one}, almost surely exactly one arm receives infinitely many particles, and the one-arm domination time is the last time at which a particle attaches to an arm other than this infinite arm.  See Figure~\ref{fig:hl-observables}(c) for an illustration, where the infinite arm is highlighted in red. Let $\upsilon^{\mathcal H}_{N,\eta}$ and $\upsilon^{\mathcal C}_{N,\lambda}$ denote the one-arm domination times of HL$(0)$ and CHL$(0)$,  and set
    \[
        \tau_{\mathcal H}^N:=N^{-3}\upsilon^{\mathcal H}_{N,\eta}, \qquad \tau_{\mathcal C}^N:=N^{-2}\upsilon^{\mathcal C}_{N,\lambda},
    \]
    for radial HL(0) with particle size $N^{-1}\eta$. Let $\upsilon_\infty$ be the global coalescence time of coalescing standard Brownian motions started from every point of $\mathbb T$ at time $0$. By \cite[Theorem~3.5 and Remark~3.8]{zhou2008stepping},    $\mathbb E[\upsilon_\infty]=2\pi^2/3$.  The following theorem identifies the scaling limit of the one-arm domination time and its limiting expectation. This sharpens \cite[Theorem~1.3]{chen2025one}, which gives matching upper and lower bounds for the expected one-arm domination time up to constants.

    \begin{thm}\label{thm:one-arm-domination}
    The one-arm domination times in HL(0) and CHL(0) satisfy
    \[
        \tau_{\mathcal H}^N\dto \frac{3\pi}{2\eta^3}\,\upsilon_\infty, \qquad \tau_{\mathcal C}^N\dto \frac{3}{4\lambda^3}\,\upsilon_\infty,\qquad \lim_{N\to\infty}\mathbb E[\tau_{\mathcal H}^N]=\frac{\pi^3}{\eta^3}, \qquad \lim_{N\to\infty}\mathbb E[\tau_{\mathcal C}^N]=\frac{\pi^2}{2\lambda^3}.
    \]

    \end{thm}

    \subsection{Proof strategy}
    
    We outline the two parts of the proof. The first establishes path-Hausdorff convergence for localized disturbance flows.  The second builds on these estimates to prove the scaling limits for the three HL$(0)$ geometries.

    \medskip
    \noindent\textbf{Localized disturbance flows.} By covariance localization, finite collections of trajectories converge to coalescing Brownian motions, as in     \cite[Propositions~2.3 and~8.2]{norris2015weak}. Together with the non-crossing property, this gives tightness. It remains to rule out additional paths in subsequential path-Hausdorff limits by verifying condition (B2) of Theorem~\ref{thm:BW-CBW-convergence-criterion} for each such limit. For each fixed $t>0$, this requires that the probability that paths passing through an interval or arc of length $\varepsilon$ at time $t_0$ occupy at least three distinct positions at time $t_0+t$ is $o(\varepsilon)$, uniformly in $t_0$ and spatial location.

    To verify (B2), since the one-gap survival bound gives only an $O(\varepsilon)$ estimate, we prove a three-point survival estimate (Theorem~\ref{thm:abstract-three-point}) using Assumption~\ref{ass:abstract-no-extra-path}\textup{(I)} and covariance localization. Consider three ordered paths and the two gaps between adjacent paths. Once a relevant gap is at most $3a_n$, the one-gap survival bound is $O(a_n)$. This already gives the required bound when the sum of the two adjacent gaps is at least $r_n$, since $(a_n/r_n)r_n=a_n$. Otherwise, we use a three-step multiscale stopping-time decomposition. First, we stop when either adjacent gap closes or one of the gaps reaches $r_n$.  On the latter event, both gaps are at most $r_n+O(j_n)$. Second, as long as all relevant gaps remain larger than the small scale $3a_n$, covariance localization makes the drift of the product of the gaps uniformly negative; this controls the time until one gap returns to scale $3a_n$. Third, once a gap is again at most $3a_n$, the one-gap estimate bounds the probability that all three paths stay distinct. For a triple of initial diameter $\ell$, the period before the first gap becomes small is handled by the two-gap estimate and contributes the term $C(t)\ell^2$; the three-step multiscale decomposition controls the remaining contribution. Together these estimates give, for some $c>0$ and each fixed $t>0$,
    \[
        \mathbb P\{\text{three paths remain distinct at time }t\} \le C(t)\left\{\ell^2+\frac{a_n}{r_n}\ell+\frac{j_n}{r_n}\exp\!\left(-\frac{ct}{r_n^2}\right)\right\}.
    \]

    To pass from a single triple to a controlled finite collection of starting points, we adapt the dyadic cutting argument of \cite[Lemma~9.4]{toth1998true}:  repeatedly halving a short interval lets us apply the three-point estimate on smaller and smaller intervals (see Proposition~\ref{prop:abstract-binary-reduction}). Combined with Assumption~\ref{ass:abstract-no-extra-path}\textup{(II)}, which reduces all paths through the interval to finitely many labels, the three-point survival estimate shows that the probability that paths in any subsequential limit passing through an interval of length $\varepsilon$ occupy at least three distinct positions at a fixed later time $t$ is at most $C(t)\varepsilon^2$, uniformly in the starting time and location (see Subsection~\ref{SS.abstract-boundary-proof}). This verifies (B2), and the BW convergence criterion then     identifies every subsequential limit with     $\sigma\mathcal W_\Xi$, proving     Theorem~\ref{thm:abstract-boundary}.

    \medskip
    \noindent\textbf{HL$(0)$ geometries.} At finite $N$, harmonic measure paths started from distinct points at the same time do not coalesce, whereas Root paths can. Therefore, we study the Root path collections $\mathbf R^{\mathcal A,N}$ in $E^{\mathcal A}$ for each $\mathcal A\in\{\mathcal H,\mathcal S,\mathcal C\}$. At each earlier grid time, a particle's Root path gives the harmonic-measure coordinate of the first point already present on its tree path traced from the tip toward the initial boundary.  After time reversal, the Root paths are driven by localized disturbances  and reusing the proof of Theorem~\ref{thm:abstract-boundary} (see Subsection~\ref{S.three-point-estimate-I}),  as in Proposition~\ref{prop:process-root-webs}, gives
    \[
        \mathbf R^{\mathcal A,N} \dto \sigma_{\mathcal A}\mathcal W^\downarrow_{\mathcal A,\ge0}.
    \]

    After restricting the paths to grid-time endpoints as in Lemma~\ref{lem:particle-preparation}, the Root and harmonic measure paths do not cross, by the generalized inverse relation. The restriction errors vanish as $N\to\infty$.  Therefore, the limits of finitely many such paths, together with the meeting times of the backward Root paths, identify the dual coupling in the limit.  Applying Lemma~\ref{lem:half-axis-dual-wedge-criterion} yields
    \[
        (\mathbf R^{\mathcal A,N},\mathbf{HM}^{\mathcal A,N}) \dto \bigl( \sigma_{\mathcal A}\mathcal W^\downarrow_{\mathcal A,\ge0}, \sigma_{\mathcal A}\mathcal W^\uparrow_{\mathcal A,\ge0} \bigr)
    \]
    and proves Theorem~\ref{thm:harmonic-measure-webs}.

    We then return to the geometric path collections. Our geometric estimates in Subsections~\ref{SS.geometric-tightness} and~\ref{SS.binary-decomposition-section6} show that the Fingers are close to the embedded Root paths, and the Gaps are close to the embedded harmonic measure paths:
    \[
        \mathbf{Fin}^{\mathcal A,N} \approx \mathfrak e_{\mathcal A}^{v_{\mathcal A}}(\mathbf R^{\mathcal A,N}), \qquad \mathbf{Gap}^{\mathcal A,N} \approx \mathfrak e_{\mathcal A}^{v_{\mathcal A}}(\mathbf{HM}^{\mathcal A,N}).
    \]
    Together with the joint limit in harmonic-measure coordinates, these approximations give Theorem~\ref{thm:joint-geometric-webs}.  In HL(0) and CHL(0), the same boundary evolution also gives the harmonic measure of each arm. By controlling later attachments to arms with small harmonic measure, we compare the one-arm domination time with the coalescence of the boundary paths and identify its limit in terms of the global coalescence time of the limiting coalescing Brownian motions, proving Theorem~\ref{thm:one-arm-domination}; see Subsection~\ref{SS.one-arm-domination-section6}.

    \subsection{Related work and key contributions}\label{relate-work}

    We make three main contributions. First, for localized disturbance flows, Norris and Turner~\cite{norris2015weak} proved that finite collections of trajectories converge to coalescing Brownian motions and that the flows converge in the weak-flow topology. That topology retains only the leftmost and rightmost paths from each space--time point and does not require ruling out extra paths in the limit~\cite[Section~6.7.1.1]{schertzer2017brownian}. Under the small-gap bound and finite-covering conditions of Assumption~\ref{ass:abstract-no-extra-path}, Theorem~\ref{thm:abstract-boundary} upgrades this to path-Hausdorff convergence of the full trajectory collection to the BW or CBW. 

    Second, earlier scaling limits for HL$(0)$ control finitely many paths, or the flow in the weak-flow sense. For radial HL$(0)$, Norris and Turner~\cite{norris2012hastings} proved weak-flow convergence of the harmonic measure flow and joint convergence of finite samples of fingers and gaps, started from fixed space--time points, to backward and forward coalescing Brownian motions. For SHL$(0)$, the known limits concern a single boundary point or a single finger~\cite{berger2022growth}. We prove that the harmonic measure paths, fingers and gaps generated by every particle of the aggregate converge jointly, in the path-Hausdorff topology, to the double BW or CBW and its geometric embeddings (Theorems~\ref{thm:harmonic-measure-webs} and~\ref{thm:joint-geometric-webs}). We do so for radial HL$(0)$, for SHL$(0)$, and for CHL$(0)$, which is conformally equivalent to radial HL$(0)$ (Remark~\ref{rmk:radial-from-cylindrical}). These limits control all paths simultaneously, including those that start at random points determined by the aggregate. For instance, the backbone, which is not the path from any fixed point, converges together with the selected dual backbone to the geometric embeddings of the bi-infinite paths of the double CBW (Corollary~\ref{cor:backbone-dual-backbone}).
    
    Third, in radial HL$(0)$ and CHL$(0)$, the rescaled one-arm domination times converge in distribution to explicit multiples of the global coalescence time of coalescing Brownian motions on $\mathbb T$, with convergence of expectations (Theorem~\ref{thm:one-arm-domination}). For CHL$(0)$, this sharpens the bounds up to constants for the expected domination time obtained by Chen, Procaccia and Zong~\cite{chen2025one}.

    Our convergence proofs use the criteria of Fontes, Isopi, Newman and Ravishankar~\cite{fontes2004brownian}, with the principal task being the verification of $(B2)$. In several lattice and drainage-network models, this condition is verified through FKG-type inequalities and two-path coalescence estimates. Coupier, Saha, Sarkar and Tran~\cite{coupier2020collision} developed an alternative approach, using Lyapunov functions to bound the expected first collision time among three ordered paths. Our argument is related to this collision-time approach. 

    The main technical ingredient is a three-point survival estimate, proved via a multiscale stopping-time decomposition. The argument distinguishes two regimes: nearby paths may be strongly correlated, while the covariance rate of well-separated paths becomes asymptotically negligible. The decomposition combines one-gap estimates near coalescence with a two-gap estimate at larger separations, without requiring independence of the paths at any finite scale. For all three HL$(0)$ geometries, we use this estimate to establish Brownian web convergence of the Root path collections. The covariance and gap conditions underlying our argument can also be verified for other correlated tree models, including stretched IDLA~\cite{berger2014stretched}. Finite-dimensional convergence of its rescaled paths to coalescing Brownian motions will be established in forthcoming work by the third author with Elboim and Yadin.

    \subsection{Structure of the paper}
    Section~\ref{S.BW} recalls the BW and CBW and the convergence criteria used throughout the paper. Section~\ref{S.abstract-flow} proves Theorem~\ref{thm:abstract-boundary} for localized disturbance flows, with the three-point survival estimate providing the main step in verifying $(B2)$. Section~\ref{S.HL-application} introduces the three HL$(0)$ models, their observables, and their path collections, and then establishes the geometric estimates needed later. Section~\ref{S.HL-observables-BW} proves Theorems~\ref{thm:harmonic-measure-webs} and~\ref{thm:joint-geometric-webs}, Corollary~\ref{cor:backbone-dual-backbone} and then derives Theorem~\ref{thm:one-arm-domination}. 

    \medskip
    \noindent{\textbf{Convention.}} Throughout the paper, $Y_n\dto Y$ and $Y_n\pto Y$ denote convergence in distribution (in the topology under discussion) and in probability, respectively. For $a_N>0$, $Z_N=O_{\P}(a_N)$ means that $(Z_N/a_N)_{N\ge1}$ is tight. In $O_\lambda(a_N)$, the constant may depend on $\lambda$ and other fixed parameters specified in the statement, but not on $N$. If $\lambda$ ranges over a compact subset of $(0,\infty)$, the bound holds uniformly in $\lambda$.  On circular  boundaries, order, intervals, crossings, and differences are interpreted using real-valued representatives.
    
    \medskip
    \noindent{\textbf{Acknowledgements.}}
    X.L. and Y.Z. thank Rongfeng Sun and Jinjiong Yu for helpful discussions and references. Y.Z. also thanks Ruixuan Li for suggesting the use of interpolation. X.L. and Y.Z. are supported by the National Key R\&D Program of China (No.\ 2021YFA1002700), and the Beijing Natural Science Foundation (JQ26001). 

    \medskip
    
    \noindent{\textbf{Declaration of AI Use.}} All  mathematical results, research ideas, proof strategies, and the overall proof framework were independently developed by the authors. AI tools were used only to assist with language polishing, figure generation, proofreading,  and checking and refining mathematical derivations. The authors take full responsibility for the content of the paper.

    \section{The (Cylindrical) Brownian web}\label{S.BW}
    In this section, we recall the BW on $\mathbb R$ and the CBW on $\mathbb T=\R/(2\pi\Z)$ as compact collections of coalescing paths, together with the convergence criteria used below. We mainly follow \cite{fontes2004brownian,schertzer2017brownian,coupier2019directed}.

    \subsection{Definitions}
    We first define the path spaces for the line and for the circle, and the corresponding spaces of compact sets of paths. Let $\overline\R=[-\infty,\infty]$. For $(x,t)\in\R^2$, define
    \[
        \Phi(x,t):=\frac{\tanh x}{1+|t|}, \qquad \Psi(t):=\tanh t,
    \]
    with $\Psi(\pm\infty)=\pm1$ and $\Phi(\pm\infty,t)=\pm(1+|t|)^{-1}$.  For the line, let $R_c^2$ be the completion of $\R^2$ under
    \[
        \rho_{\R}\bigl((x_1,t_1),(x_2,t_2)\bigr) := |\Phi(x_1,t_1)-\Phi(x_2,t_2)| \vee |\Psi(t_1)-\Psi(t_2)|.
    \]
    For the periodic case, let $d_{\T}$ denote the circular distance on $\T$, and let $\mathrm{Cyl}_c$ be the completion of $\T\times\R$ under
    \begin{equation}\label{eq:cyl-metric}
        \rho_{\T}\bigl((x_1,t_1),(x_2,t_2)\bigr) := |\Psi(t_1)-\Psi(t_2)| + \frac{1-\max\{|\Psi(t_1)|,|\Psi(t_2)|\}}{\pi} d_{\T}(x_1,x_2).
    \end{equation}
    In both completions, the lines $t=\pm\infty$ are collapsed to single points $(*,\pm\infty)$. For $\Xi\in\{\R,\T\}$, write $(\mathcal E_\Xi,\rho_\Xi)=(R_c^2,\rho_{\R})$ or $(\mathrm{Cyl}_c,\rho_{\T})$, respectively.  Set $\mathsf X_{\R}=\overline\R$ and $\mathsf X_{\T}=\T$. A path with starting time $\sigma_\pi\in\overline\R$ is a map $\pi:[\sigma_\pi,\infty]\to\mathsf X_\Xi\cup\{*\}$, taking values in $\mathsf X_\Xi$ at finite times, with $\pi(\infty)=*$ and $\pi(\sigma_\pi)=*$ if $\sigma_\pi=-\infty$, such that $t\mapsto(\pi(t),t)$ is continuous into $(\mathcal E_{\Xi},\rho_{\Xi})$. Let $\Pi_{\Xi}$ be the set of such paths, equipped with
    \begin{equation}\label{eq:path-metric}
        d_{\Pi_\Xi}(\pi_1,\pi_2) := |\Psi(\sigma_{\pi_1})-\Psi(\sigma_{\pi_2})| \vee
        \begin{cases}
        \displaystyle \sup_{t\in\R} \left| \Phi(\pi_1(t\vee\sigma_{\pi_1}),t) - \Phi(\pi_2(t\vee\sigma_{\pi_2}),t) \right|, & \Xi=\R, \\[3mm]
        \displaystyle \sup_{t\in\R} \rho_{\T}\Bigl( (\pi_1(t\vee\sigma_{\pi_1}),t\vee\sigma_{\pi_1}), (\pi_2(t\vee\sigma_{\pi_2}),t\vee\sigma_{\pi_2}) \Bigr), & \Xi=\T.
        \end{cases}
    \end{equation}
    For $\Xi=\R$, a path starting at $+\infty$ carries a label $x\in\overline\R$, used as its constant extension in \eqref{eq:path-metric}. For $\Xi=\T$, on every bounded time interval, the metric \eqref{eq:cyl-metric} induces the same topology as that used in \cite{coupier2019directed}. Moreover, each of the regions $t\ge T$ and $t\le -T$ has $\rho_{\T}$-diameter $O(1-\tanh T)$. Hence the results of \cite[Theorems~2.3--2.5]{coupier2019directed} can first be applied on bounded time intervals and then extended to the whole time axis by letting $T\to\infty$. Let $\mathfrak K_{\Xi}$ be the space of nonempty compact subsets of $(\Pi_{\Xi},d_{\Pi_\Xi})$, with the Hausdorff metric     $d_{\mathfrak K_{\Xi}}$ and Borel $\sigma$-field  $\mathcal F_{\mathfrak K_{\Xi}}$. For a compact path collection $K$, write $K(z)$ for the paths in $K$ started from $z$.

    The following result characterizes the BW and CBW used below as compact path collections.
    \begin{thm}[{\cite[Theorem~2.1]{fontes2004brownian}, \cite[Theorem~2.3]{coupier2019directed}}]
    \label{thm:BW-CBW}
    There are unique laws for the standard BW $\mathcal W$ and the CBW $\mathcal W^\uparrow$, valued in $(\mathfrak K_{\R},\mathcal F_{\mathfrak K_{\R}})$ and $(\mathfrak K_{\T},\mathcal F_{\mathfrak K_{\T}})$. Use the notation
    \[
        (E,K,\omega_z)=
        \begin{cases}
        (\R^2,\mathcal W,\pi_z), & \Xi=\R, \\
        (\T\times\R,\mathcal W^\uparrow,\pi_z^\uparrow), & \Xi=\T.
        \end{cases}
    \]
    Then
    \begin{enumerate}[label=\textup{(\roman*)}]
    \item For any  $z\in E$, $K(z)$ almost surely contains a unique path $\omega_z$ starting from $z$;
    \item For any $z_1,\ldots,z_m \in E$, $(\omega_{z_1},\ldots,\omega_{z_m})$ are coalescing standard Brownian motions started from $z_1,\ldots,z_m$;
    \item For any dense countable subset $D\subset E$, $K=\overline{\{\omega_z:z\in D\}}$ in $(\Pi_\Xi,d_{\Pi_\Xi})$.
    \end{enumerate}
    \end{thm}

    We also recall the corresponding backward BWs. Define the time reversal map by $(x,t)\mapsto(-x,-t)$, applied pointwise to paths ($\pi\mapsto-\pi$) and to path collections.  For $\Xi\in\{\R,\T\}$, set  $\Pi_\Xi^\downarrow:=-\Pi_\Xi$, $d_{\Pi_\Xi}^\downarrow(-\pi_1,-\pi_2):=d_{\Pi_\Xi}(\pi_1,\pi_2),$ with corresponding Hausdorff space $\mathfrak K_\Xi^\downarrow$. Then   $-\widehat{\mathcal W}\stackrel d=\mathcal W\text{ and } -\mathcal W^\downarrow\stackrel d=\mathcal W^\uparrow$, where $\widehat{\mathcal W}$ and $\mathcal W^\downarrow$ are the dual (backward) BW and CBW of Theorem~\ref{thm:double BWs} below.

    We use the same crossing convention in both settings.  A forward  path $\pi$ crosses a backward  path $\hat\pi$ if there exist $s_1<s_2$ in the intersection of their time domains such that
    \begin{equation}\label{eq:forward-backward-crossing}
        \bigl(\pi(s_1)-\hat\pi(s_1)\bigr) \bigl(\pi(s_2)-\hat\pi(s_2)\bigr)<0.
    \end{equation}
    On $\T$, the inequality is required to hold for some continuous real-valued representations of the two paths. Touching without a change of spatial order is not regarded as crossing.

    The following result characterizes the joint laws of the BW and CBW with their duals through non-crossing.
    \begin{thm}[{\cite[Theorem~6.2.4]{schertzer2017brownian}, \cite[Theorem~2.4]{coupier2019directed}}]\label{thm:double BWs}
    There are unique laws of the pairs $(\mathcal W,\widehat{\mathcal W})$ and $(\mathcal W^\uparrow,\mathcal W^\downarrow)$, valued in $\mathfrak K_{\R}\times\mathfrak K_{\R}^{\downarrow}$ and $\mathfrak K_{\T}\times\mathfrak K_{\T}^{\downarrow}$, such that $\mathcal W$ and $-\widehat{\mathcal W}$ are standard BWs, $\mathcal W^\uparrow$ and $-\mathcal W^\downarrow$ are standard CBWs, and almost surely no forward path crosses any backward path, in the sense of \eqref{eq:forward-backward-crossing}. These pairs are called the double BW and the double CBW.

    \end{thm}

    \subsection{Convergence criteria}\label{SS.convergence criteria}
    We next recall the criterion used to identify path-Hausdorff limits.  The main difficulty is to rule out extra paths in subsequential limits.   To this end, given an interval or arc $A$ at time $t_0$, the criterion counts the distinct positions at time $t_0+t$ of the paths that start no later than $t_0$ and pass through $A$ at time $t_0$.

    For $K\in\mathfrak K_{\R}$, $t_0\in\R$, $t>0$, and $A=[a,b]$ with $a<b$, define
    \begin{equation}\label{eq:eta-BW}
        \eta_K^{\mathbb R}(t_0,t;A):=\#\{\pi(t_0+t):\pi\in K,\ \sigma_\pi\le t_0,\ \pi(t_0)\in A\}.
    \end{equation}
    For the cylinder, let $A=[a\to b]$ be the counterclockwise arc from $a$ to $b$, of length $|a\to b|=(b-a)\bmod 2\pi$.  For $K\in\mathfrak K_{\T}$, set
    \begin{equation}\label{eq:eta-cyl}
        \eta_K^{\mathbb T}(t_0,t;A):=\#\{\pi(t_0+t):\pi\in K,\ \sigma_\pi\le t_0,\ \pi(t_0)\in A\}.
    \end{equation}

    The following theorem is the non-crossing convergence criterion used throughout the paper.

    \begin{thm}[{\cite[Theorem~2.2]{fontes2004brownian}, \cite[Theorem~6.6.5]{schertzer2017brownian}, \cite[Theorem~2.5]{coupier2019directed}}]\label{thm:BW-CBW-convergence-criterion}
    Use the notation
    \[
        (E,\mathscr K,\mathscr W,A_\delta(a),\eta_K)=
        \begin{cases}
        (\R^2,\mathfrak K_{\R},\mathcal W,[a,a+\delta],\eta_K^{\mathbb R}), & \Xi=\R, \\
        (\T\times\R,\mathfrak K_{\T},\mathcal W^\uparrow, [a\to a+\delta],\eta_K^{\mathbb T}), & \Xi=\T.
        \end{cases}
    \]
    Let $\chi^n$ be random non-crossing $\mathscr K$-valued compact path collections, and suppose that the following hold.
    \begin{enumerate}[label=\textup{(\Alph*)}]
    \item[(I)] There exist a countable dense subset  $D\subset E$ and paths $\omega_z^n\in\chi^n$, $z\in D$, such that for any $z_1,\ldots,z_k\in D$, the paths $(\omega_{z_i}^n)_{1\le i\le k}$ converge in distribution to coalescing standard Brownian motions on $\Xi$ started at $z_1,\ldots,z_k$;
    \item[(B2)] For every $t>0$,
    \[
        \lim_{\delta\downarrow0}\delta^{-1}\limsup_{n\to\infty}\sup_{(a,t_0)\in \Xi\times\R} \P\{\eta_{\chi^n}(t_0,t;A_\delta(a))\ge3\}=0.
    \]
    \end{enumerate}
    Then $\chi^n\dto\mathscr W$.
    \end{thm}
    \begin{remark}
    In the non-crossing setting, the $(B1)$ condition in \cite[Theorem~2.2]{fontes2004brownian} is redundant given condition $(I)$; see \cite[Section~6.6.2]{schertzer2017brownian}. So we state the criterion without it.
    \end{remark}
    Condition $(I)$ gives the convergence of finitely many paths from deterministic starting points to coalescing Brownian motions, while $(B2)$ rules out extra paths in any path-Hausdorff limit.
    Convergence to the backward BW and CBW is obtained by applying the same criterion after time reversal.

    For later use in the HL(0) convergence argument of Section~\ref{S.HL-observables-BW}, we record a version of the dual wedge criterion of \cite{schertzer2017brownian} that involves only countably many backward paths. For backward paths $\hat\pi_1,\hat\pi_2$ with finite terminal times, let $r$ be the minimum of their terminal times and define their last meeting time
    \[
        T(\widehat\pi_1,\widehat\pi_2) :=\sup\{u\le r:\widehat\pi_1(u)=\widehat\pi_2(u)\}, \qquad \sup\varnothing=-\infty.
    \]
    If $\widehat\pi_1(r)\neq \widehat\pi_2(r)$, relabel if necessary so that $\widehat\pi_1(r)<\widehat\pi_2(r)$, and define the wedge
    \[
        V(\widehat\pi_1,\widehat\pi_2) :=\{(x,u):T(\widehat\pi_1,\widehat\pi_2)<u<r, \ \widehat\pi_1(u)<x<\widehat\pi_2(u)\};
    \]
    otherwise there is no wedge. A forward path $\pi$ \emph{enters $V$ from outside} if, for some finite $\sigma_\pi\le s<t$,
    \[
        (\pi(s),s)\notin\overline V, \qquad (\pi(t),t)\in V,
    \]
    as in \cite[Definition~6.3.8]{schertzer2017brownian}.

    Fix a deterministic countable dense set $D_0\subset\mathbb R^2$. For each $w\in D_0$, let $w_n\to w$ and let $\pi_n^w\in\Pi_{\R}^{\downarrow}$ be a backward path terminating at $w_n$.

    \begin{lemma}\label{lem:half-axis-dual-wedge-criterion}
    Let $\chi_n$ be non-crossing $\mathfrak K_{\R}$-valued random variables satisfying \textup{(I)} of Theorem~\ref{thm:BW-CBW-convergence-criterion}. Suppose that every finite collection of the paths $\pi_n^w$, together with their pairwise meeting times (viewed in $[-1,1]$ through $\Psi$), converges jointly in distribution to coalescing standard backward Brownian motions terminating at the corresponding points $w\in D_0$, together with their meeting times. If almost surely no path in $\chi_n$ enters, from outside, any wedge generated by these backward paths, then $\chi_n\dto\mathcal W$.
    \end{lemma}

    \begin{proof}
    By non-crossing and condition~\textup{(I)}, tightness follows from \cite[Proposition~B.2]{fontes2004brownian} and \cite[Proposition~6.6.4]{schertzer2017brownian}. Let $D$ and $(\omega_z^n)_{z\in D}$ be as in condition~\textup{(I)}. Since $D$ and $D_0$ are countable, by taking a subsequence and applying Skorohod representation, we may assume that
    \[
        \chi_n\to K,\qquad \omega_z^n\to\omega^z,\qquad \pi_n^w\to \pi^w,
    \]
    jointly for $z\in D$ and $w\in D_0$, together with all pairwise backward meeting times. Path-Hausdorff convergence and $\omega_z^n\in\chi_n$ imply that $\omega^z\in K$ for every $z\in D$. Hence
    \[
        \mathcal W_0:=\overline{\{\omega^z:z\in D\}}\subset K.
    \]
    By Theorem~\ref{thm:BW-CBW}, $\mathcal W_0$ and $\widehat{\mathcal W}_0:=\overline{\{\pi^w:w\in D_0\}}$ are respectively a forward and backward standard BW.

    The convergence of the backward paths and their meeting times implies that wedge exclusion (no path enters a wedge from outside) passes to the limit, including entry through the bottom point; see \cite[Remark~6.6.7]{schertzer2017brownian}. Indeed, if a path in $K$  entered a limiting wedge from outside, then path-Hausdorff convergence,  together with convergence of the wedge boundaries and their meeting time,     would give a path in $\chi_n$ entering the corresponding prelimit wedge for all sufficiently large $n$, a contradiction.

    Let $\mathcal W_*$ be the forward BW dual to $\widehat{\mathcal W}_0$. Since no path in $K$ enters a wedge generated by the backward collection $\{\pi^w:w\in D_0\}$, the wedge characterization \cite[Remark~6.3.10 and Section~6.6.2]{schertzer2017brownian} gives $\mathcal W_0\subset K\subset\mathcal W_*$. Both $\mathcal W_0$ and $\mathcal W_*$ are standard BWs, so the above inclusion implies that they are equal almost surely. Hence $K=\mathcal W_0$. Thus every subsequential limit is a standard BW, and therefore $\chi_n\dto\mathcal W$.
    \end{proof}


    On $\mathbb T$, let two backward paths be distinct at time $r$ and first meet at $\tau<r$ when traced backwards. Represent them by continuous real-valued functions $L,R$, exchanging the paths and shifting by multiples of $2\pi$ if necessary, such that $L(\tau)=R(\tau),\ 0<R(u)-L(u)<2\pi, \tau<u\le r$.  Their wedge is
    \[
        (p_{\mathbb T}\times\mathrm{id})\{(x,u):\tau<u<r,\ L(u)<x<R(u)\}.
    \]
    The argument of \cite[Section~4.2, proof of Theorem~1.9]{sun2008brownian} also applies on $\mathbb T$, showing that any path with finite starting time avoiding such wedges belongs to the forward CBW dual to the limiting backward CBW.

    For $t\ge0$, extend the limiting backward BW or CBW to negative times. A forward path cannot enter a wedge with $\tau<0$ without crossing a side, while wedges with $\tau>0$ are excluded by prelimit wedge exclusion and convergence of meeting times. For the chosen countable backward paths, almost surely no wedge has $\tau=0$. Thus all paths with finite starting times are identified, and closure and the time-tail bound recover the whole collection.

    The HL$(0)$ observables are only defined at nonnegative times, and for them we use the following compactification. For path collections on $t\ge0$, write
    \[
        \iota_\Xi(x,t):=\left(\frac{c(x)}{1+t},\frac{t}{1+t}\right),\qquad c(x)=
        \begin{cases}
        \tanh x,&\Xi=\R, \\
        (\cos x,\sin x),&\Xi=\T.
        \end{cases}
    \]
    Extend each forward path to earlier times by its starting value and each backward path to later times by its terminal value. For forward paths, set
    \begin{align}\label{path-distance}
        d_{\rm path}(\pi,\pi'):=\left|\frac{\sigma_\pi}{1+\sigma_\pi}-\frac{\sigma_{\pi'}}{1+\sigma_{\pi'}}\right|\vee \sup_{t\ge0}\|\iota_\Xi(\pi(t\vee\sigma_\pi),t)-\iota_\Xi(\pi'(t\vee\sigma_{\pi'}),t)\|_\infty,
    \end{align}
    and use the analogous formula for backward paths with terminal times. Take the metric completions and the induced Hausdorff metrics on compact collections. The limits of forward paths with starting times tending to infinity form
    \[
        \mathcal I=\left\{ \left(1,\left(\frac{q}{1+t},\frac{t}{1+t}\right)_{t\ge0}\right) :q\in\overline{c(\Xi)}\right\},\qquad      d_{\rm path}(\pi,\mathcal I)\le \frac{2}{1+\sigma_\pi}.
    \]
    On bounded time intervals, this topology is uniform convergence in the compactified spatial coordinate together with convergence of endpoint times. Hence, for bounded endpoint times, convergence in \eqref{eq:path-metric}, together with the preceding tail bound, implies convergence in this topology. Let $\mathcal W^\uparrow_{\Xi,\ge0}$ be the closure in this completion of the BW/CBW paths with finite starting times $s\ge0$, and let $\mathcal W^\downarrow_{\Xi,\ge0}$ be the closure of the corresponding dual backward paths with finite terminal times $t\ge0$, restricted to $[0,t]$, including $t=0$.

    \section{BW convergence of localized disturbance flows}\label{S.abstract-flow}

    In this section, we prove Theorem~\ref{thm:abstract-boundary}.  We first verify finite-dimensional convergence, that is, condition (I) of Theorem~\ref{thm:BW-CBW-convergence-criterion}. Subsection~\ref{SS.three-point survival} then derives the three-point survival estimate from Assumption~\ref{ass:abstract-no-extra-path}(I). Finally, in Subsection~\ref{SS.abstract-boundary-proof}, we combine this estimate with Assumption~\ref{ass:abstract-no-extra-path}(II) to verify (B2) and complete the proof.

    \medskip
    \paragraph{\textbf{Finite-dimensional convergence.}} We first establish finite-dimensional convergence for the flow trajectories and then transfer it to the grid interpolation.  For $m\ge1$, let  $z_i^n=(x_i^n,s_i^n)\to z_i=(x_i,s_i)$, $1\le i\le m$,     be deterministic points in $\Xi\times\R$. For $t<s$, each trajectory started at time $s$ is set equal to its value at $s$, and the starting time is included in the path topology. We use the local Skorokhod topology for these trajectories.

    \begin{proposition}\label{prop:abstract-fdd}
    The paths $\bigl(\Phi^n_{s_i^n,\cdot}(x_i^n)\bigr)_{i=1}^m$  converge jointly in distribution to coalescing Brownian motions on $\Xi$ with variance rate $\sigma^2$, started from $z_1,\ldots,z_m$.
    \end{proposition}

    \begin{proof}
    We follow the martingale argument in  \cite[proofs of Propositions~2.3 and~8.2]{norris2015weak}. For $s\in\R$ and $x\in\R$, the process $ \Phi^n_{s,t}(x)-x, t\ge s $, is a square-integrable martingale with jumps bounded by $j_n$ and predictable quadratic variation $K_n(0)(t-s)$; for two trajectories, their predictable covariation grows at rate $K_n$ evaluated at their current separation. On $\R$, the same martingale and covariation identities follow from the spatially truncated construction. By the zero-drift condition in \eqref{drift-covariance} and $\Delta_n\in L^2(\R)$, the truncated drifts vanish uniformly on bounded spatial intervals, and the identities pass to the limit after stopping the paths on such intervals. The maximal inequality \cite[Theorem~3.2, $\mathsf{BJ}_{2,4}$]{marinelli2014maximal} gives, for $h>0$,
    \begin{equation}\label{eq:clock-cell-fourth-moment}
        \E\!\left[ \sup_{0\le v\le h} |\Phi^n_{s,s+v}(x)-x|^4 \,\middle|\,\mathcal F_s \right] \le C\{K_n(0)^2h^2+j_n^2K_n(0)h\}.
    \end{equation}
    Together with $K_n(0)\to\sigma^2$ and $j_n\to0$, this yields tightness with continuous subsequential limits, to which the martingale characterization applies. By condition \textup{(L3)} in Definition \ref{def:abstract-localized}, the  covariation rates of distinct trajectories vanish uniformly as long as the trajectories stay separated. Hence the martingale characterization in    \cite[proofs of Propositions~2.3 and~8.2]{norris2015weak}    identifies the limiting trajectories, up to collision, as Brownian motions with variance $\sigma^2$ and zero covariation. For an ordered pair, the limiting gap is a non-negative continuous local martingale and therefore is absorbed at zero. Thus the limiting trajectories are coalescing Brownian motions. The same argument applies when $z_i^n\to z_i$, since the initial positions converge and $ K_n(0)|s_i^n-s_i|\rightarrow0.$
    \end{proof}

    To compare the trajectories with their grid interpolations, we control the largest deviation within a single grid interval. For $a,u\in\R$ and $R,T>0$, let $J_\Xi(a,R)=[a-R,a+R]$ if $\Xi=\R$ and $J_\Xi(a,R)=[0,2\pi]$ if $\Xi=\T$, and define
    \[
        \Omega_n(a,u;R,T):= \max_{\substack{k\in\Z: [k/n,(k+1)/n]\cap[u,u+T]\ne\varnothing}} \sup_{x\in J_\Xi(a,R)}\sup_{0\le v\le 1/n} \bigl|\Phi^n_{k/n,k/n+v}(x)-x\bigr|.
    \]
    The following lemma shows that, with high probability,  $\Omega_n$ is small.

    \begin{lemma}\label{lem:automatic-clock-oscillation}
    For every $R,T,\varepsilon>0$, there is $C_{R,T,\varepsilon}<\infty$ such that
    \[
        \sup_{a,u}\P\{\Omega_n(a,u;R,T)>\varepsilon\} \le C_{R,T,\varepsilon}(1/n+j_n^2)\rightarrow0.
    \]
    \end{lemma}

    \begin{proof}
    Choose an $\varepsilon/4$-mesh of $J_\Xi(a,R)$, including the endpoints. Equation \eqref{eq:clock-cell-fourth-moment}, Markov's inequality, and a union bound over the $T\cdot O(n)$ grid intervals give the required bound at the mesh points. Monotonicity then extends the bound to all $x\in J_\Xi(a,R)$, with an additional error no larger than the mesh spacing. Periodicity gives the circular case.
    \end{proof}
    Together with Proposition~\ref{prop:abstract-fdd}, the preceding estimate yields finite-dimensional convergence for the grid interpolations.

    \begin{corollary}\label{cor:interpolated-fdd}
    If $s_i^n\in n^{-1}\Z$ for all $i,n$, then $\bigl(I_n\Phi^n_{s_i^n,\cdot}(x_i^n)\bigr)_{i=1}^m$ converges jointly in distribution to the coalescing Brownian motions identified in Proposition~\ref{prop:abstract-fdd}.
    \end{corollary}

    \begin{proof}
    Fix a bounded time interval. On $\R$, since $K_n(0)$ is bounded, Doob's inequality gives
    \[
        \P\!\left( \max_{1\le i\le m} \sup_{s_i^n\le t\le T} \bigl|\Phi^n_{s_i^n,t}(x_i^n)-x_i^n\bigr|>R_\varepsilon \right) \le \frac{C_T}{R_\varepsilon^2}.
    \]
    Hence, for every $\varepsilon>0$, taking $R_\varepsilon=\sqrt{C_T/\varepsilon}$,  and choosing $B<\infty$ so that $\max_{1\le i\le m}|x_i^n|\le B$ for all sufficiently large $n$, we obtain
    \[
        \P\!\left( \Phi^n_{s_i^n,t}(x_i^n)\in [-B-R_\varepsilon,B+R_\varepsilon] \ \text{for all }1\le i\le m,\ s_i^n\le t\le T \right) \ge 1-\varepsilon .
    \]
    On this event, the cocycle property and the interpolation defining $I_n\Phi^n$ give
    \[
        \max_{1\le i\le m}\sup_{t\le T} d_\Xi\!\left( I_n\Phi^n_{s_i^n,t}(x_i^n), \Phi^n_{s_i^n,t}(x_i^n) \right) \le 2\Omega_n.
    \]
    By Lemma~\ref{lem:automatic-clock-oscillation}, the interpolation error vanishes in probability, and Proposition~\ref{prop:abstract-fdd} gives the joint coalescing Brownian limit on every bounded time interval. Since the path metrics compactify time, this gives the convergence.
    \end{proof}

    \subsection{Three-point survival}\label{SS.three-point survival}
    The goal of this subsection is the following three-point survival estimate under Assumption~\ref{ass:abstract-no-extra-path}(I). We prove it at the end of the subsection, after establishing one-gap and two-gap estimates. For a starting time $s$, a stopping time $\tau$, and $t>0$, define the non-coalescence event of $x_1\le x_2\le x_3$ (the three trajectories started from these points at time $s$ are still distinct at time $s+\tau+t$) by
    \[
        E_{\tau}(t) := \left\{ \#\{\Phi^n_{s,s+\tau+t}(x_1),\Phi^n_{s,s+\tau+t}(x_2),\Phi^n_{s,s+\tau+t}(x_3)\}=3 \right\}.
    \]
    Let $\ell=x_3-x_1$ be the diameter of the starting triple. In the rest of this section, the sequences $a_n,r_n,q_n,M_n$ are understood to satisfy \eqref{parameters-control}, and we write $\beta_n=a_n/r_n$.
    \begin{thm}\label{thm:abstract-three-point}
    There exists $c>0$ such that, for every $t>0$, there is $C(t)<\infty$ for which, for all sufficiently large $n$,
    \begin{equation}\label{eq:abstract-three-point-bound}
        \mathbb P\bigl(E_0(t)\bigr)\le C(t)\left\{\ell^2+\beta_n\ell+\frac{j_n}{r_n}\exp\!\left(-\frac{ct}{r_n^2}\right)\right\}.
    \end{equation}
    \end{thm}
    \begin{remark}\label{rmk:constants-chosen}
    The time-dependent constants  $C(t), C_i(t), 1\le i\le 3$ in this section can be chosen non-increasing in $t$ and locally bounded on $(0,\infty)$.
    \end{remark}
    We now prepare the proof, starting with a single gap.

    \medskip

    \paragraph{\textbf{One-gap process}}    For $d>0$, let $D^{n,d}_{t}$ be the length of the gap between the trajectories started from $x$ and $x+d$ at time $s$:
    \[
        D_t^{n,d}=\Phi^n_{s,s+t}(x+d)-\Phi^n_{s,s+t}(x).
    \]
    For test functions $f$, the generator $\mathcal L_n$ is given by
    \begin{align}\label{one-gap-generator}
        \mathcal L_nf(d) =\rho_n\int_\Xi \Big[f\bigl(d+\Delta_n(z+d)-\Delta_n(z)\bigr)-f(d)\Big] \mu_\Xi(dz),
    \end{align}
    where $\Delta_n$ is defined in (\ref{map-deviation}). On $\R$, all generator integrals below are understood as symmetric principal values. Equivalently,
    \begin{align}\label{one-gap-generator-martingale}
        f(D^{n,d}_r)-f(D^{n,d}_0) -\int_0^r \mathcal L_n f(D^{n,d}_q)\,dq
    \end{align}
    is a local martingale. For $R>0$, define the hitting times by
    \[
        \tau^n:=\inf\{t:D_t^{n,d}=0\},\qquad \tau_R^n:=\inf\{t:D_t^{n,d}=0\ \hbox{or}\ D_t^{n,d}\ge R\}.
    \]
    Note that $D^{n,d}_{t}\equiv0$ for $t\ge\tau^n$.  The next lemma converts  Assumption~\ref{ass:abstract-no-extra-path}(I) into hitting-time estimates for the one-gap process.
    \begin{lemma}\label{lem:abstract-one-gap-estimates}
    For every $t>0$, there exists $C_1(t)<\infty$ such that for all sufficiently large $n$,
    \begin{equation}\label{eq:abstract-one-gap-survival}
        \begin{aligned}
            \mathbb P(\tau^n>t)\le C_1(t)(d+j_n),\quad 0<d<1,\qquad 
            \mathbb P(D_{\tau_{r_n}^n}^{n,d}\ge r_n)\le\frac{d}{r_n},\quad0<d<r_n,
        \end{aligned}
    \end{equation}
    where $j_n, r_n$ are as in Definition~\ref{def:abstract-localized} and \eqref{parameters-control}.
    Moreover, there exist $C,c>0$ and $R_0\in(0,1]$ such that, for all sufficiently large $n$, integers $p\ge1$, $r_n\le R\le R_0$, and $0<d<R$,
    \begin{equation}\label{eq:abstract-one-gap-moment-tail}
        \mathbb E[(\tau_R^n)^p]\le p!C^pR^{2p-1}(d+j_n),\qquad \mathbb P(\tau_R^n>t)\le C\frac{d+j_n}{R}\exp\!\left(-\frac{ct}{R^2}\right)\quad(t\ge R^2).
    \end{equation}
    \end{lemma}
    \begin{proof}
    The identity $\mathcal L_n\id(d)=0$ follows from the zero-drift condition in \eqref{drift-covariance}. Fix $\mathfrak b>2$ and choose $R_0\in(0,1]$.  For $R\in[r_n,R_0]$, define the test function
    \begin{align*}
        F_{n,R}(d):=\ind_{\{d>0\}}(d+\mathfrak b j_n)(R+\mathfrak b j_n-d)=\ind_{\{d>0\}}(Rd-d^2+C_{n,R}),
    \end{align*}
    where $C_{n,R}:=\mathfrak b j_nR+\mathfrak b^2j_n^2$.  For $0<d<R$,
    \begin{equation}\label{eq:abstract-LF-from-Q-Lambda}
        \begin{aligned}
        \mathcal L_nF_{n,R}(d) \ &\overset{\mathclap{\eqref{one-gap-generator}}}{=}\ \rho_n\int_{\Xi} \bigl(F_{n,R}(d+\Delta_n(z+d)-\Delta_n(z))-F_{n,R}(d)\bigr)\,\mu_\Xi(dz) \\
        &=\  \rho_n\int_{\Xi} \bigl(-(\Delta_n(z+d)-\Delta_n(z))^2 -C_{n,R}\ind_{\{g_n(z+d)=g_n(z)\}}\bigr)\,\mu_\Xi(dz) \\
        &\overset{\mathclap{\eqref{eq:abstract-Q-covariance-identity},\,    \eqref{eq:abstract-Q-coalescence-rates}}}{=}\quad\,-Q_n(d)-C_{n,R}\Lambda_n(d).
        \end{aligned}
    \end{equation}
    For all sufficiently  large $n$, Definition~\ref{def:abstract-localized}\textup{(L2)--(L3)}, \eqref{eq:abstract-Q-covariance-identity},   and $j_n=o(a_n)$ in \eqref{parameters-control} give $Q_n(d)\ge \sigma^2/2$  uniformly in $d\in[a_n,R_0]$.  If $0<d<a_n$, then $C_{n,R}\ge j_nr_n$ and \eqref{eq:abstract-micro-coalescence} in Assumption~\ref{ass:abstract-no-extra-path}(I)  give
    \[
        Q_n(d)+C_{n,R}\Lambda_n(d)\ge Q_n(d)+j_nr_n\Lambda_n(d)\ge c_0 .
    \]
    Hence there exists $c>0$ such that, for $R\in[r_n,R_0]$ and $0<d<R$,
    \begin{equation}\label{eq:abstract-one-gap-scale}
        \mathcal L_nF_{n,R}(d)\overset{\eqref{eq:abstract-LF-from-Q-Lambda}}{=} -Q_n(d)-C_{n,R}\Lambda_n(d)\le -c.
    \end{equation}
    On $\{D_{\tau_R^n}^{n,d}\ge R\}$, the bound    $|\Delta_n(z+d)-\Delta_n(z)|\le 2j_n$ and $\mathfrak b>2$ imply
    \[
        R\le D_{\tau_R^n}^{n,d}\le R+2j_n \ \Rightarrow\ R+\mathfrak b j_n-D_{\tau_R^n}^{n,d}\ge(\mathfrak b-2)j_n>0 \ \Rightarrow\ F_{n,R}(D_{\tau_R^n}^{n,d})\ge0 \ \text{on }\{D_{\tau_R^n}^{n,d}\ge R\}.
    \]
    For $T>0$, the expected number of jumps from a positive gap to zero in $(0,T\wedge\tau_R^n]$ is $\E\bigl[\int_0^{T\wedge\tau_R^n}\Lambda_n(D_{q-}^{n,d})\,dq\bigr]$, which is at most $\P(D_{T\wedge\tau_R^n}^{n,d}=0)$. Applying \eqref{one-gap-generator-martingale} to  $f(d)=Rd-d^2+C_{n,R}$ gives
    \[
        \begin{aligned}
        0&\le \E\bigl[F_{n,R}(D^{n,d}_{T\wedge\tau_R^n})\bigr] =F_{n,R}(d)-\E\!\left[\int_0^{T\wedge\tau_R^n}Q_n(D_q^{n,d})\,dq\right] -C_{n,R}\P(D^{n,d}_{T\wedge\tau_R^n}=0) \\
        &\le F_{n,R}(d)-\E\!\left[\int_0^{T\wedge\tau_R^n} \bigl\{Q_n(D_q^{n,d})+C_{n,R}\Lambda_n(D_{q-}^{n,d})\bigr\}\,dq\right]\overset{\eqref{eq:abstract-one-gap-scale}}{\le} F_{n,R}(d)-c\E[T\wedge\tau_R^n].
        \end{aligned}
    \]

    Letting $T\to\infty$, we obtain, for some $C>0$ and all $0<d<R$,
    \begin{equation}\label{eq:abstract-one-gap-mean-exit}
        \mathbb E[\tau_R^n] \le c^{-1}F_{n,R}(d) \le C R(d+j_n).
    \end{equation}
    Since $\tau_R^n<\infty$ by \eqref{eq:abstract-one-gap-mean-exit}, optional stopping for the bounded martingale $D^{n,d}_{\cdot\wedge\tau_R^n}$ gives $\mathbb P(D_{\tau_R^n}^{n,d}\ge R)\le d/R$, proving the second bound in \eqref{eq:abstract-one-gap-survival}.
    For the first bound in \eqref{eq:abstract-one-gap-survival}, fix $R_*\in(0,R_0]$.  For $0<d<R_*$,
    \[
        \{\tau^n>t\}\subset \{D_{\tau_{R_*}^n}^{n,d}\ge R_*\}\cup\{\tau_{R_*}^n>t\}.
    \]
    The preceding estimates with $R=R_*$ give $\mathbb P(\tau^n>t)\le C_1(t)(d+j_n)$.  Enlarging $C_1(t)$ covers $R_*\le d<1$.

    Let $m_p(d)=\mathbb E[(\tau_R^n)^p]$. By \eqref{eq:abstract-one-gap-mean-exit} and $j_n=o(r_n)$, $\sup_{0<d<R}m_1(d)\le C_0R^2$ uniformly for $r_n\le R\le R_0$. The Markov property and induction give
    \[
        m_p(d)=p\,\mathbb E\left[\int_0^{\tau_R^n}m_{p-1}(D_u^{n,d})\,du\right]\le p!(C_0R^2)^{p-1}m_1(d),\qquad p\ge2.
    \]
    Together with \eqref{eq:abstract-one-gap-mean-exit}, this proves the moment bound. Summing the exponential series gives $\mathbb E[e^{c\tau_R^n/R^2}-1]\le C(d+j_n)/R$ for sufficiently small $c>0$, and Markov's inequality yields the exponential bound.
    \end{proof}
    We next turn to two adjacent gaps.

    \medskip

    \paragraph{\textbf{Two-gap process}} For starting positions $x_1\le x_2\le x_3$, define the adjacent gap lengths
    \[
        U_r:=\Phi^n_{s,s+r}(x_2)-\Phi^n_{s,s+r}(x_1),\qquad V_r:=\Phi^n_{s,s+r}(x_3)-\Phi^n_{s,s+r}(x_2),
    \]
    and the test functions
    \[
        H^{\rm lin}(u,v):=uv,\qquad H^{\rm cyl}(u,v):=uv(2\pi-u-v).
    \]
    These functions are motivated by the case of independent Brownian motions.  For three independent standard Brownian motions, the adjacent gaps before collision have generator
    \[
        \mathscr L^{\rm BM} =(\partial_{uu}+\partial_{vv}-\partial_{uv}),
    \]
    and one has
    \[
        \mathscr L^{\rm BM}H^{\rm lin}=-1,\qquad \mathscr L^{\rm BM}H^{\rm cyl}=-2\pi.
    \]
    Hence $H^{\rm lin}(U_t,V_t)+t$ and $H^{\rm cyl}(U_t,V_t)+2\pi t$ are local martingales up to collision.

    Returning to the gap process $(U_r,V_r)$, for the test functions $f$ used below, its generator $\mathscr L_n$ is given by
    \begin{equation}\label{two-gap-generator}
        \mathscr L_nf(u,v)=\rho_n\int_\Xi        \Big[ f\bigl(u+\Delta_n(z+u)-\Delta_n(z),v+\Delta_n(z+u+v)-\Delta_n(z+u)\bigr)-f(u,v) \Big]\mu_\Xi(dz),
    \end{equation}
    where $\Delta_n$ is defined in (\ref{map-deviation}). Equivalently,
    \begin{align}\label{two-gap-generator-martingale}
        f(U_r,V_r)-f(U_0,V_0)-\int_0^r \mathscr L_nf(U_q,V_q)\,dq
    \end{align}
    is a local martingale.  The next lemma shows that, as in the Brownian case, $\mathscr L_n$ applied to these functions is uniformly negative as long as all relevant gaps are at least $3a_n$.

    \begin{lemma}\label{lem:abstract-two-gap-drift}
    There exists $c_\Xi>0$  such that for all sufficiently large $n$
    \[
        \mathscr L_nH^{\rm lin}(u,v) \le -c_{\mathbb R}\ \text{on}\ u,v\ge 3a_n,\qquad \mathscr L_nH^{\rm cyl}(u,v) \le -c_{\mathbb T}\ \text{on}\ u,v,2\pi-u-v\ge 3a_n .
    \]
    \end{lemma}

    \begin{proof}
    On $\mathbb R$, the zero-drift condition in \eqref{drift-covariance} gives
    \begin{equation}\label{eq:abstract-line-barrier-generator}
        \begin{aligned}
        \mathscr L_nH^{\rm lin}(u,v) \ &\overset{\mathclap{\eqref{two-gap-generator}}}{=}\ \rho_n\int_\Xi        \Big[ H^{\rm lin}\bigl(u+\Delta_n(z+u)-\Delta_n(z),v+\Delta_n(z+u+v) -\Delta_n(z+u)\bigr) \\
        &\hspace{19.7em}-H^{\rm lin}(u,v) \Big]\mu_\Xi(dz) \\
        &=\ \rho_n\int_{\Xi} (\Delta_n(z+u)-\Delta_n(z))(\Delta_n(z+u+v)-\Delta_n(z+u)) \,\mu_\Xi(dz) \\
        &\overset{\mathclap{\eqref{drift-covariance}}}{=}\ K_n(u)+K_n(v)-K_n(u+v)-K_n(0).
        \end{aligned}
    \end{equation}
    Using $u,v\ge3a_n$,  $j_n=o(a_n)$ in (\ref{parameters-control}),  and Definition~\ref{def:abstract-localized}\textup{(L2)--(L3)}, we obtain
    \[
        \mathscr L_nH^{\rm lin}(u,v) \overset{\eqref{eq:abstract-line-barrier-generator}}{=} K_n(u)+K_n(v)-K_n(u+v)-K_n(0) \le 3\kappa_n(3a_n)-K_n(0) \le -\frac{\sigma^2}{2}.
    \]
    On $\mathbb T$, let $w=2\pi-u-v$.  Similarly, using the zero-drift condition and $K_n(u+v)=K_n(w)$, we obtain
    \begin{equation}
        \mathscr L_nH^{\rm cyl}(u,v) =(2\pi-2u)K_n(u)+(2\pi-2v)K_n(v) +(2\pi-2w)K_n(w)-2\pi K_n(0)+R_n(u,v), \label{eq:abstract-circular-barrier-generator}
    \end{equation}
    where
    \[
        \begin{aligned}
        R_n(u,v):=\rho_n\int_\Xi \bigl(\Delta_n(z+u)-\Delta_n(z)\bigr) &\bigl(\Delta_n(z+u+v)-\Delta_n(z+u)\bigr) \\
        \times& \bigl(\Delta_n(z)-\Delta_n(z+u+v)\bigr)\,\mu_\Xi(d z).
        \end{aligned}
    \]
    Since $|\Delta_n(z)-\Delta_n(z+u+v)|\le 2j_n,$ by $2|ab|\le a^2+b^2$ and $\max_{h\in\R}|K_n(h)|\le K_n(0)$,
    \begin{equation}\label{eq:abstract-circular-remainder}
        \begin{aligned}
        |R_n(u,v)| &\le 2j_n\rho_n\int_\Xi \bigl|\Delta_n(z+u)-\Delta_n(z)\bigr| \bigl|\Delta_n(z+u+v)-\Delta_n(z+u)\bigr| \,\mu_\Xi(dz) \\
        &\le j_n\rho_n\int_\Xi \Bigl[ \bigl(\Delta_n(z+u)-\Delta_n(z)\bigr)^2 +\bigl(\Delta_n(z+u+v)-\Delta_n(z+u)\bigr)^2 \Bigr]\,\mu_\Xi(dz) \\
        &\overset{\mathclap{\eqref{eq:abstract-Q-covariance-identity}}}{=}\ 2j_n\bigl(2K_n(0)-K_n(u)-K_n(v)\bigr) \le 8j_nK_n(0).
        \end{aligned}
    \end{equation}
    Together with $u,v,w\ge3a_n$,  $j_n=o(a_n)$ in (\ref{parameters-control}),  Definition~\ref{def:abstract-localized}\textup{(L2)--(L3)}, and \eqref{eq:abstract-circular-remainder}, the identity \eqref{eq:abstract-circular-barrier-generator} yields
    \[
        \mathscr L_nH^{\rm cyl}(u,v) \le -2\pi K_n(0)+6\pi\kappa_n(3a_n)+8j_nK_n(0) \le -\pi\sigma^2.
    \]
    This completes the proof.
    \end{proof}

    We now stop the three trajectories at the first time one of the gaps falls to the cutoff scale $3a_n$:
    \[
        \tau_{\rm cut}^n:=\inf\{r\ge0:m_r\le3a_n\}, \qquad m_r:=
        \begin{cases}
        \min(U_r,V_r), & \Xi=\mathbb R, \\
        \min(U_r,V_r,2\pi-U_r-V_r), & \Xi=\mathbb T.
        \end{cases}
    \]

    Let $S_r:=U_r+V_r$ be the total length of the two adjacent gaps, so that $S_0=x_3-x_1=\ell$. The next lemma bounds the cutoff time and the total gap length at that time.
    \begin{lemma}\label{lem:abstract-two-gap}
    For every $T>0$, there exists $C_2(T)<\infty$ such that for all sufficiently large $n$
    \begin{equation}\label{two-gap-eq}
        \mathbb P(\tau_{\rm cut}^n>T)\le C_2(T)\ell^2, \qquad \mathbb E[S_{\tau_{\rm cut}^n\wedge T}]\le\ell.
    \end{equation}
    \end{lemma}

    \begin{proof}
    Assume that $m_0>3a_n$ (otherwise $\tau_{\rm cut}^n=0$ and there is nothing to prove). Let
    \[
        \tau_{\rm cut}^{n,R}:=
        \begin{cases}
        \tau_{\rm cut}^n\wedge T\wedge\inf\{r: S_r\ge R\}, & \Xi=\mathbb R, \\
        \tau_{\rm cut}^n\wedge T, & \Xi=\mathbb T
        \end{cases}
        , \quad  H_\Xi(u,v):=
        \begin{cases}
        H^{\rm lin}(u,v), & \Xi=\mathbb R, \\
        H^{\rm cyl}(u,v), & \Xi=\mathbb T.
        \end{cases}
    \]
    Since $H_\Xi$ is non-negative,  Lemma~\ref{lem:abstract-two-gap-drift} gives
    \[
        \begin{aligned}
        0&\le \mathbb E\bigl[H_\Xi(U_{\tau_{\rm cut}^{n,R}},V_{\tau_{\rm cut}^{n,R}})\bigr] \\
        &\overset{\eqref{two-gap-generator-martingale}}{=} H_\Xi(U_0,V_0) +\mathbb E\Big[\int_0^{\tau_{\rm cut}^{n,R}} \mathscr L_nH_\Xi(U_q,V_q)\,dq\Big] \le H_\Xi(U_0,V_0)-c_\Xi\mathbb E[\tau_{\rm cut}^{n,R}] .
        \end{aligned}
    \]
    Letting $R\to\infty$ gives
    \begin{equation}\label{eq:abstract-two-gap-mean-cutoff}
        \mathbb E[\tau_{\rm cut}^n\wedge T] \le c_\Xi^{-1}H_\Xi(U_0,V_0)\le
        \begin{cases}
        c_\Xi^{-1}\ell^2/4, & \Xi=\mathbb R, \\
        c_\Xi^{-1}\pi\ell^2/2, & \Xi=\mathbb T.
        \end{cases}
    \end{equation}
    Therefore
    \[
        \mathbb P(\tau_{\rm cut}^n>T) \le T^{-1}\mathbb E[\tau_{\rm cut}^n\wedge T] \overset{\eqref{eq:abstract-two-gap-mean-cutoff}}{\le} C_2(T)\ell^2.
    \]
    Before relabelling, $S_t=U_t+V_t$ is a nonnegative martingale. On $\mathbb T$, cyclically relabel the points at $\tau_{\rm cut}^n$ so that a largest circular gap is the third gap. One of the two adjacent gaps is then at most $3a_n$, and their sum does not increase: if the three gaps before relabelling are $u,v,w$, the new sum is $2\pi-\max\{u,v,w\}\le u+v$. Thus optional stopping before relabelling gives
    \[
        \mathbb E[S_{\tau_{\rm cut}^n\wedge T}]\le\ell.
    \]
    This completes the proof.
    \end{proof}
    We now combine the one-gap and two-gap estimates, via a multiscale stopping-time decomposition, to prove the three-point survival estimate.

    \begin{proof}[Proof of Theorem~\ref{thm:abstract-three-point}]
    We first prove that there exist $c>0$ and $C_3(t)<\infty$ such that
    \begin{equation}\label{eq:abstract-post-cutoff-control}
        \mathbb P\bigl(E_{\tau_{\rm cut}^n}(t)\mid\mathcal F_{s+\tau_{\rm cut}^n}\bigr)\le C_3(t)\left\{\beta_n S_{\tau_{\rm cut}^n}+\frac{j_n}{r_n}\exp\!\left(-\frac{ct}{r_n^2}\right)\right\}\quad\text{on }\{\tau_{\rm cut}^n<\infty\}.
    \end{equation}
    We restart the process at time $s+\tau_{\rm cut}^n$. After the possible cyclic relabelling when $\Xi=\mathbb T$, a gap of size at most $3a_n$ is one of the two adjacent gaps.   If $S_{\tau_{\rm cut}^n}\ge r_n$, then
    \[
        \mathbb P\bigl(E_{\tau_{\rm cut}^n}(t)\mid \mathcal F_{s+\tau_{\rm cut}^n}\bigr) \overset{\eqref{eq:abstract-one-gap-survival}}{\le} C_1(t)(3a_n+j_n) \le C_3(t)\beta_n S_{\tau_{\rm cut}^n}.
    \]
    It remains to consider the case $S_{\tau_{\rm cut}^n}<r_n$, which we treat in three steps.

    \smallskip
    \noindent\textbf{Step 1: one-gap estimate up to mesoscopic scale $r_n$.} Define
    \[
        T_1:=\inf\{u\ge0:U_u=0\ \hbox{or}\ V_u=0\ \hbox{or}\ U_u\vee V_u\ge r_n\}, \qquad E_1:=\{U_{T_1}\vee V_{T_1}\ge r_n\}.
    \]
    The exponential bound in \eqref{eq:abstract-one-gap-moment-tail} and optional stopping for $U+V$ give, conditionally on $\mathcal F_{s+\tau_{\rm cut}^n}$ and for all sufficiently large $n$,
    \begin{equation}\label{eq:abstract-post-cutoff-step1}
        \begin{aligned}
        \mathbb P(T_1>t/3\mid\mathcal F_{s+\tau_{\rm cut}^n})\le C\frac{S_{\tau_{\rm cut}^n}+j_n}{r_n}\exp\!\left(-\frac{ct}{r_n^2}\right),\quad 
        \mathbb P(E_1\mid\mathcal F_{s+\tau_{\rm cut}^n})\le\frac{S_{\tau_{\rm cut}^n}}{r_n}.
        \end{aligned}
    \end{equation}

    \smallskip
    \noindent\textbf{Step 2:  two-gap estimate up to the next cutoff scale $3a_n$.} On $E_1$, by the jump bound the two adjacent gaps at time $T_1$ are bounded by $r_n+Cj_n$.  Starting from $T_1$, let $\tau'_{\rm cut}$ be the next cutoff time at scale $3a_n$.  Lemma~\ref{lem:abstract-two-gap} and the strong Markov property give
    \begin{equation}\label{eq:abstract-post-cutoff-step2}
        \mathbb P(\tau'_{\rm cut}>t/3\mid \mathcal F_{s+\tau_{\rm cut}^n+T_1}) \le 4C_2(t/3)(r_n+Cj_n)^2  \le 9C_2(t/3)r_n^2.
    \end{equation}

    \smallskip
    \noindent\textbf{Step 3: one-gap estimate from cutoff scale $3a_n$.} On $\{\tau'_{\rm cut}\le t/3\}$, one relevant gap is at most $3a_n$.  Thus
    \begin{equation}\label{eq:abstract-post-cutoff-step3}
        \mathbb P\!\left( E_{\tau_{\rm cut}^n+T_1+\tau'_{\rm cut}}       \bigl(t/3\bigr) \,\middle|\,  \mathcal F_{s+\tau_{\rm cut}^n+T_1+\tau'_{\rm cut}} \right) \overset{\eqref{eq:abstract-one-gap-survival}}{\le} C_1(t/3)(3a_n+j_n).
    \end{equation}
    Combining \eqref{eq:abstract-post-cutoff-step1}--\eqref{eq:abstract-post-cutoff-step3}, and using $r_n^2=O(a_n)$, $j_n=o(a_n)$, and $e^{-ct/r_n^2}\le C(t)r_n^2$, we obtain
    \[
        \begin{aligned}
        \mathbb P\bigl(E_{\tau_{\rm cut}^n}(t)\mid\mathcal F_{s+\tau_{\rm cut}^n}\bigr)&\le C\frac{S_{\tau_{\rm cut}^n}+j_n}{r_n}e^{-ct/r_n^2}+\frac{S_{\tau_{\rm cut}^n}}{r_n}\left\{9C_2(t/3)r_n^2+C_1(t/3)(3a_n+j_n)\right\}\\
        &\le C_3(t)\left\{\beta_n S_{\tau_{\rm cut}^n}+\frac{j_n}{r_n}e^{-ct/r_n^2}\right\},
        \end{aligned}
    \]
    which proves \eqref{eq:abstract-post-cutoff-control}. It remains to derive \eqref{eq:abstract-three-point-bound}. Since coalesced paths cannot separate again,
    \[
        E_0(t) \subseteq \{\tau_{\rm cut}^n>t/2\}\cup\bigl(\{\tau_{\rm cut}^n\le t/2\}\cap E_{\tau_{\rm cut}^n}(t/2)\bigr).
    \]
    Lemma~\ref{lem:abstract-two-gap} bounds the probability of the first event by $C_2(t/2)\ell^2$. For the second event,
    \[
        \begin{aligned}
        \mathbb P\bigl(\{\tau_{\rm cut}^n\le t/2\}\cap E_{\tau_{\rm cut}^n}(t/2)\bigr)&\overset{\eqref{eq:abstract-post-cutoff-control}} {\le} C_3(t/2)\left\{\beta_n\mathbb E\!\left[S_{\tau_{\rm cut}^n}\ind_{\{\tau_{\rm cut}^n\le t/2\}}\right]+\frac{j_n}{r_n}e^{-ct/(2r_n^2)}\right\}\\
        &\overset{\eqref{two-gap-eq}} {\le} C(t)\left\{\beta_n\ell+\frac{j_n}{r_n}e^{-ct/(2r_n^2)}\right\}.
        \end{aligned}
    \]
    This proves \eqref{eq:abstract-three-point-bound} after decreasing $c$.
    \end{proof}

    \subsection{Proof of Theorem~\ref{thm:abstract-boundary}}\label{SS.abstract-boundary-proof}
    We now deduce the counting condition (B2) from Theorem~\ref{thm:abstract-three-point}, using a dyadic cutting argument over finite ordered meshes together with  Assumption~\ref{ass:abstract-no-extra-path}(II), and thereby complete the proof of Theorem~\ref{thm:abstract-boundary}.

    We first establish compactness, measurability, and tightness of the interpolated path collections.

    \begin{lemma}\label{lem:abstract-tightness}
    For localized disturbance flows, the collections $\mathcal X_n^\Phi$ are measurable random compact sets whose laws are tight in the path-Hausdorff topology.
    \end{lemma}

    \begin{proof}
    Fix $n$. On $\R$, the martingale identity in the proof of Proposition~\ref{prop:abstract-fdd} gives, for $s\le t$ in $n^{-1}\Z$, $m\ge1$, and $e\in\{-1,1\}$,
    \[
        \P\bigl\{|\Phi^n_{s,t}(em)-em|>m/2\bigr\}  \le 4K_n(0)(t-s)m^{-2}.
    \]
    Since the set of $s,t$ is countable, the Borel--Cantelli lemma and     monotonicity imply, almost surely, $\Phi^n_{s,t}(x)\rightarrow\pm\infty$ as $x\rightarrow\pm\infty$  for every $s\le t$ in $n^{-1}\Z$.

    Fix a bounded time interval $[a,b]$ with $a,b\in n^{-1}\Z$, and let $a=t_0<\cdots<t_J=b$ be its grid times. By the cocycle property, a path started before time $a$, restricted to $[a,b]$, is determined by its value $y$ at time $a$. If $|y|\le M$, monotonicity gives
    \[
        \Phi^n_{a,t_j}(-M)\le \Phi^n_{a,t_j}(y) \le \Phi^n_{a,t_j}(M), \qquad 0\le j\le J,
    \]
    whereas $y\to\pm\infty$ implies $\Phi^n_{a,t_j}(y)\to\pm\infty$ for every $j$. On $\T$, for $y\in[0,2\pi]$,
    \[
        \Phi^n_{a,t_j}(0)\le \Phi^n_{a,t_j}(y) \le \Phi^n_{a,t_j}(0)+2\pi.
    \]
    Now consider a sequence of generating paths. Since their starting times lie in $n^{-1}\Z$, after passing to a subsequence they are either constant or tend to $\pm\infty$. For fixed starting times, the preceding bounds and linear interpolation yield uniform subsequential convergence in the compactified spatial coordinate on each bounded interval. The same argument and time compactification apply when starting times tend to $-\infty$. When they tend to $+\infty$, spatial compactification and the time weights in \eqref{eq:path-metric} give convergence. Thus the generating paths are relatively compact in the path space, and their closure $\mathcal X_n^\Phi$ is compact.

    For measurability, let
    \[
        \mathcal G_n:=\left\{t\mapsto I_n\Phi^n_{s,t}(p_\Xi(q)):s\in n^{-1}\Z,\ q\in\mathbb Q     \right\}.
    \]
    Right-continuity in the initial label, together with the linear interpolation, gives $\mathcal X_n^\Phi=\overline{\mathcal G_n}$. Since $\mathcal G_n$ is a countable collection of measurable paths and $\mathcal X_n^\Phi$ is compact, it is a measurable random compact set. Monotonicity at the grid times and the linear interpolation preserve spatial order, so $\mathcal X_n^\Phi$ is non-crossing.

    It remains to prove tightness as $n\to\infty$. Let $D\subset\Xi\times\R$ be a deterministic countable dense set. For $z=(x,s)\in D$, set $s_n=n^{-1}\lfloor ns\rfloor$ and let $\pi_n^z:=I_n\Phi^n_{s_n,\cdot}(x)\in\mathcal X_n^\Phi.$ For every finite collection $z_1,\ldots,z_k\in D$, Corollary~\ref{cor:interpolated-fdd} gives $(\pi_n^{z_i})_{i=1}^k \Rightarrow (\pi^\sigma_{z_i})_{i=1}^k,$ where the limit consists of coalescing Brownian motions with variance $\sigma^2$. Together with non-crossing, this gives tightness on $\R$ by \cite[Proposition~B.2]{fontes2004brownian}.

    The same argument gives tightness on every bounded time interval in $\T$. The time compactification in Section~\ref{S.BW} then yields tightness in the full path-Hausdorff topology.
    \end{proof}

    We next use a dyadic cutting argument, whose basic idea goes back to \cite[Lemma~9.4]{toth1998true}, to pass from a single triple to a finite ordered set of starting points of controlled cardinality.    Let $t>0$ and let $Y_n=\{y_1<\cdots<y_{m_n}\}$ be a finite ordered set of $m_n\le M_n$ labels contained in an interval or arc of length $\eps$.  For $1\le i<j<k\le m_n$, let $E_{i,j,k}(0,t)$ denote the event $E_0(t)$ for the triple started from $y_i,y_j,y_k$, and let $E(Y_n)$ be the event that the trajectories started from $Y_n$ occupy at least three distinct positions at time $s+t$, that is,
    \[
        E(Y_n):=\bigcup_{1\le i<j<k\le m_n}E_{i,j,k}(0,t).
    \]
    \begin{proposition}\label{prop:abstract-binary-reduction}
    For all sufficiently large $n$, with $c>0$ as in Theorem~\ref{thm:abstract-three-point},
    \begin{equation}\label{eq:abstract-binary-bound}
        \P(E(Y_n))\le C(t)\left\{\eps^2+\beta_n\eps\log(eM_n)+M_n\frac{j_n}{r_n}\exp\!\left(-\frac{ct}{r_n^2}\right)\right\}.
    \end{equation}
    \end{proposition}
    \begin{proof}
    We may assume $0<\eps<1$, since for larger $\eps$ the bound follows from $\P(E(Y_n))\le1$ after enlarging $C(t)$. Let $[m_n]:=\{1,\ldots,m_n\}$. For an index block $B:=[p,q]\subset [m_n]$, define
    \[
        E_B:=\{\#\{\Phi^n_{s,s+t}(y_i):i\in B\}\ge 3\},\qquad D(B):=y_q-y_p.
    \]
    Then $E(Y_n)=E_{[m_n]}$.  When $|B|\ge3$, define
    \begin{equation}\label{midpoint}
        m(B):=\max\{i\in B:\ y_i\le (y_p+y_q)/2\}, \quad B_-:=[p,m(B)],\quad B_+:=[m(B)+1,q],
    \end{equation}
    and
    \[
        C_B:= \bigcup_{\substack{i<j<k \\
        i,j,k\in\{p,m(B),m(B)+1,q\}}} E_{i,j,k}(0,t).
    \]
    We first prove
    \begin{equation}\label{event-monotonous}
        E_B\subset E_{B_-}\cup E_{B_+}\cup C_B .
    \end{equation}
    On $E_B\setminus(E_{B_-}\cup E_{B_+})$, each subblock has at most two distinct images; on $\R$, these are the images of its endpoints, by monotonicity.  On $\T$, represent the images by real numbers $z_p\le\cdots\le z_q$, with $z_q-z_p\le2\pi$. If either $z_m-z_p=2\pi$ or $z_q-z_{m+1}=2\pi$, the other subblock has a single image, so $B$ has at most two distinct images, a contradiction. Thus both differences are less than $2\pi$, and the images of each subblock are again represented by its endpoints. Consequently, at least three of $p,m,m+1,q$ have distinct images, proving \eqref{event-monotonous}.

    If $m_n<3$, then $E(Y_n)=\varnothing$. Assume that $m_n\ge3$, and build the dyadic tree of blocks with at least three elements as follows.  Start with $\mathcal T_0=\{[m_n]\}$.  Given $\mathcal T_k$, let $\mathcal T_{k+1}:=\{B_\sigma:\, B\in\mathcal T_k,\ \sigma\in\{-,+\},\ |B_\sigma|\ge 3\}$.  Let $K_*:=\inf\{k\ge0:\mathcal T_k=\varnothing\}<\infty$ and
    \[
        \mathcal T:=\bigcup_{k\ge0}\mathcal T_k=\bigcup_{0\le k<K_*}\mathcal T_k .
    \]
    Since $E_B=\varnothing$ whenever $|B|<3$, iterating \eqref{event-monotonous} gives
    \[
        E(Y_n)=E_{[m_n]}\subseteq \bigcup_{B\in\mathcal T} C_B .
    \]
    Each $C_B$ is a union of at most four triple events, each with initial diameter at most $D(B)$.  Hence Theorem~\ref{thm:abstract-three-point} gives
    \begin{align}\label{pre-abstract-binary-reduction}
        \P(E(Y_n))\le C(t)\sum_{B\in\mathcal T}\left\{D(B)^2+\beta_nD(B)+\frac{j_n}{r_n}\exp\!\left(-\frac{ct}{r_n^2}\right)\right\}.
    \end{align}
    For each $k$, the block intervals in $\mathcal T_k$ are pairwise disjoint. Since $y_{m_n}-y_1\le\varepsilon$, \eqref{midpoint} gives $D(B)\le2^{-k}\eps$ for any $B\in\mathcal T_k$. The binary tree has at most $m_n$ non-terminal blocks, so $|\mathcal T|\le m_n\le M_n$ and
    \[
        \sum_{B\in\mathcal T_k}D(B)\le\eps\min\{1,m_n2^{-k}\}.
    \]
    Consequently,
    \[
        \begin{aligned}
        \sum_{B\in\mathcal T}D(B)^2\le\sum_{k\ge0}2^{-k}\eps^2\le2\eps^2,\quad 
        \sum_{B\in\mathcal T}D(B)\le\eps\sum_{k\ge0}\min\{1,m_n2^{-k}\}\le C\eps\log(em_n).
        \end{aligned}
    \]
    Substitution into \eqref{pre-abstract-binary-reduction} proves \eqref{eq:abstract-binary-bound}.
    \end{proof}
    Recall from Subsection~\ref{SS.convergence criteria} that, for a compact path collection $K$, $\eta_K(t_0,t;A)$ is the number of distinct positions at time $t_0+t$ reached by paths in $K$ meeting $A\times\{t_0\}$. For $r>0$, define the corresponding $r$-separated count (the largest number of such positions that are pairwise at distance at least $r$) by
    \begin{align*}
        \eta_K^{(r)}(t_0,t;A):=\max\Bigl\{|F|:\;& F\subset\Xi\cap\{\pi(t_0+t):\pi\in K,\ \sigma_\pi\le t_0,\ \pi(t_0)\in A\}, \\
        &d_\Xi(x,y)\ge r\ \forall\,x,y\in F,\ x\ne y\Bigr\}.
    \end{align*}
    For a finite set $R\subset\Xi$, let  $\eta_{\Phi^n}(t_0,t;R):=|\mathsf D^{\Phi^n}(t_0,t;R)|$ be the number of distinct images of $R$ at time $t_0+t$. Finally, let
    \[
        G_n(t_0,t;A):=\left\{ d_H^\to\left( \mathsf D_{\mathrm{traj}}^{\Phi^n}(t_0,t;A), \mathsf D^{\Phi^n}(t_0,t;\mathcal R_n(t_0,A)) \right)\le q_n \right\},
    \]
    where $q_n \downarrow 0$ is as in Assumption~\ref{ass:abstract-no-extra-path}\textup{(II)}.   The next lemma uses Assumption~\ref{ass:abstract-no-extra-path}\textup{(II)} to reduce counts for the path collections $\mathcal X_n^\Phi$ to counts for the finite label sets $\mathcal R_n$, up to a vanishing error.

    Fix $r,t,\delta>0$. In the next lemma, $A$ ranges over intervals or arcs with $|A|+2\delta<1$.
    \begin{lemma}\label{lem:abstract-shadow-count-reduction}
    Uniformly in $t_0$ and $A$, as $n\to\infty$,
    \[
        \begin{aligned}
        \P\{\eta_{\mathcal X_n^\Phi}^{(r)}(t_0,t;A)\ge3\} \le \P\{\eta_{\Phi^n}(t_0,t;\mathcal R_n(t_0,A^\delta))\ge3\}+\P\bigl(G_n(t_0,t;A^\delta)^c\bigr)+o(1).
        \end{aligned}
    \]
    \end{lemma}

    \begin{proof}
    We give the proof for $\Xi=\R$; the case $\Xi=\T$ follows by the same argument.  Let $u=n^{-1}\lfloor nt_0\rfloor$ and choose $0<\eta<\min\{1,\delta,r\}/16$.  Write $A=[a,a+\eps]$ and fix $R\ge4$. By Doob's inequality, with probability at least $1-C_tR^{-2}$, the trajectories started from $(a-1,u)$ and $(a+\eps+1,u)$ remain in $[a-R,a+R]$ up to time $t_0+t+1/n$. Moreover, Lemma~\ref{lem:automatic-clock-oscillation} gives
    \[
        \P\{\Omega_n(a,u;R,t+2)\le\eta\} \ge1-C_{R,t,\eta}(n^{-1}+j_n^2)=1-o(1).
    \]
    On the intersection of these events, their interpolated values at time $t_0$ lie on the two sides of $A^\eta$. Since every generating path started no later than $t_0$ starts no later than $u$ (starting times are grid times), monotonicity shows that any such path with $\pi(t_0)\in A^\eta$ remains between the two bounding trajectories for all times after $u$. The cocycle property and $\Omega_n\le\eta$ then imply that, at both $t_0$ and $t_0+t$, its interpolated value differs from the corresponding flow value by at most $2\eta$. These inequalities extend to the closure.

    Now suppose $\eta_{\mathcal X_n^\Phi}^{(r)}(t_0,t;A)\ge3$. Since the starting times of generating paths lie in $n^{-1}\Z$, the three corresponding paths may be approximated, at both times and within $\eta$, by generating paths started no later than $t_0$. Hence their flow values at time $t_0$ lie in $A^\delta$, while their flow values at time $t_0+t$ are pairwise separated by more than $ r-6\eta>r/2$. On $G_n(t_0,t;A^\delta)$, each of these three terminal values lies within $q_n$ of some $\Phi^n_{t_0,t_0+t}(y_i)$, $y_i\in\mathcal R_n(t_0,A^\delta)$. If $2q_n<r/2$, these three images are distinct. Consequently,
    \[
        \P\!\left\{ \eta_{\mathcal X_n^\Phi}^{(r)}(t_0,t;A)\ge3,\, \eta_{\Phi^n}(t_0,t;\mathcal R_n(t_0,A^\delta))<3,\, G_n(t_0,t;A^\delta) \right\} \le C_tR^{-2}+o(1)
    \]
    uniformly in $t_0$ and $A$ on $\R$. Letting first $n\to\infty$ and then $R\to\infty$ proves the claim.
    \end{proof}
    Now we derive the $r$-separated count estimate needed for condition $(B2)$.
    \begin{proposition}\label{prop:abstract-B2}
    Let $I\subset(0,\infty)$ be a compact interval.  There exists $C_I<\infty$ such that, for every $t\in I$,  $r>0$, and  $0<\eps<1$,
    \[
        \limsup_{n\to\infty} \sup_{t_0,\ |A|=\eps} \mathbb P\{\eta_{\mathcal X_n^\Phi}^{(r)}(t_0,t;A)\ge3\} \le C_I\eps^2 .
    \]
    \end{proposition}
    \begin{proof}
    Write $I=[\hat{t},\tilde{t}]$ and fix $t\in I$, $r>0$ and $\eps$ as above. Choose $\delta>0$ with $\eps+2\delta<1$. By \eqref{eq:abstract-shadow-size} and \eqref{eq:abstract-shadow-cover},
    \[
        \P\{|\mathcal R_n(t_0,A^\delta)|>M_n\} +\P\bigl(G_n(t_0,t;A^\delta)^c\bigr)\rightarrow0
    \]
    uniformly in $t_0$ and $|A|=\eps$. Conditioning on $\mathcal F_{t_0}$, Proposition~\ref{prop:abstract-binary-reduction} applies to $\mathcal R_n(t_0,A^\delta)$, since the flow after $t_0$ is independent of $\mathcal F_{t_0}$.  On $\{|\mathcal R_n(t_0,A^\delta)|\le M_n\}$, this set has at most $M_n$ labels in an interval or arc of length at most $\eps+2\delta+2q_n$. Therefore, Lemma~\ref{lem:abstract-shadow-count-reduction} and Proposition~\ref{prop:abstract-binary-reduction} give
    \[
        \begin{aligned}
        \P\{\eta_{\mathcal X_n^\Phi}^{(r)}(t_0,t;A)\ge3\}\le o(1)+C_I\Big [(\eps+2\delta+2q_n)^2+\beta_n(\eps+2\delta+2q_n)\log(eM_n)+M_n\frac{j_n}{r_n}\text{e}^{-\tfrac{c\hat{t}}{r_n^2}} \Big ],
        \end{aligned}
    \]
    where $o(1)$ is uniform in $t_0,A$ for fixed $t,r,\delta$, and we may take $C_I=C(\hat{t})$ by Remark~\ref{rmk:constants-chosen}. Since $r_n^2=O(a_n)$, \eqref{parameters-control} gives $\beta_n\log(eM_n)\to0$, $r_n^2\log M_n\to0$, and $j_n/r_n\to0$. Consequently,
    \[
        M_n\frac{j_n}{r_n}e^{-c\hat{t}/r_n^2}=\frac{j_n}{r_n}\exp\!\left(-\frac{c\hat{t}-r_n^2\log M_n}{r_n^2}\right)\longrightarrow0.
    \]
    Taking $n\to\infty$ yields
    \[
        \limsup_{n\to\infty}\sup_{t_0,\ |A|=\eps} \P\{\eta_{\mathcal X_n^\Phi}^{(r)}(t_0,t;A)\ge3\} \le C_I(\eps+2\delta)^2.
    \]
    Letting $\delta\downarrow0$ completes the proof.
    \end{proof}
    We can now complete the proof of Theorem~\ref{thm:abstract-boundary}. The main step is to verify condition $(B2)$ for subsequential limits, using the $r$-separated count estimate.
    \begin{proof}[Proof of Theorem~\ref{thm:abstract-boundary}]
    Lemma~\ref{lem:abstract-tightness} gives tightness, and monotonicity gives non-crossing.  Let $\mathcal X$ be a subsequential path-Hausdorff limit of $\mathcal X_{n_k}^\Phi$ with $n_k\to\infty$.  The non-crossing property of $\mathcal X_{n_k}^\Phi$ ensures that $\mathcal X$ is non-crossing.

    Consider the countable collection $(\pi_n^z)_{z\in D}$ from the proof of Lemma~\ref{lem:abstract-tightness}. By Corollary~\ref{cor:interpolated-fdd} and tightness, there exists a further subsequence admitting a joint Skorokhod representation in which
    \[
        \mathcal X_{n_k}^\Phi\rightarrow\mathcal X, \qquad \pi_{n_k}^z\rightarrow\pi^z\quad(z\in D).
    \]
    The limits are coalescing Brownian motions, and since $\pi_n^z\in\mathcal X_n^\Phi$ and $(\mathcal X_n^\Phi,\pi_n^z)\to(\mathcal X,\pi^z)$, we have $\pi^z\in\mathcal X$ for every $z\in D$. Thus $\mathcal X$, viewed as a constant sequence, satisfies \textup{(I)} of Theorem~\ref{thm:BW-CBW-convergence-criterion}. On $\R$, non-crossing implies that paths meeting a bounded interval at a fixed time remain spatially bounded on every subsequent bounded time interval.

    It remains to verify condition $(B2)$ for $\mathcal X$.  Fix $t>0$,  $0<\eps<1/3$, $a\in\Xi$, and $t_0\in\mathbb R$.  For $0<\alpha<t/2$, $\gamma>0$ with $\eps+2\gamma<1$, and $r>0$, define the open event $\mathcal O_{\alpha,\gamma,r}(a,t_0,\eps)$ by
    \[
        \begin{aligned}
        \mathcal O_{\alpha,\gamma,r}(a,t_0,\eps):=\Bigl\{K:\;&\exists\,\pi_1,\pi_2,\pi_3\in K,  \sigma_{\pi_l}<t_0+\alpha,\ \pi_l(t_0+\alpha)\in A_{\eps+2\gamma}^{\circ}(a-\gamma),\ 1\le l\le3, \\
        &\pi_l(t_0+t)\in\Xi,\ 1\le l\le3,\   d_\Xi\bigl(\pi_i(t_0+t),\pi_j(t_0+t)\bigr)>2r, 1\le i<j\le3 \Bigr\} .
        \end{aligned}
    \]
    On $\{\eta_{\mathcal X}(t_0,t;A_\eps(a))\ge3\}$, choose $\pi_1,\pi_2,\pi_3\in\mathcal X$ such that
    \[
        \sigma_{\pi_i}\le t_0, \qquad \pi_i(t_0)\in A_\eps(a), \qquad \pi_i(t_0+t)\ne \pi_j(t_0+t)\quad (i\ne j).
    \]
    The terminal positions are distinct, so their pairwise distances are larger than $2r$ for all sufficiently small $r>0$.  For the fixed $\gamma>0$, continuity gives that, for all sufficiently small $\alpha>0$, the positions at time $t_0+\alpha$ lie in $A_{\eps+2\gamma}^{\circ}(a-\gamma)$.  Consequently,
    \begin{equation*}
        \ind_{\{\eta_{\mathcal X}(t_0,t;A_\eps(a))\ge3\}} \le \liminf_{\alpha\downarrow0} \liminf_{r\downarrow0} \ind_{\{\mathcal X\in\mathcal O_{\alpha,\gamma,r}(a,t_0,\eps)\}} .
    \end{equation*}
    By Fatou's lemma,
    \begin{equation}\label{eq:abstract-limit-event-cover}
        \mathbb P\{\eta_{\mathcal X}(t_0,t;A_\eps(a))\ge3\} \le \liminf_{\alpha\downarrow0} \liminf_{r\downarrow0} \mathbb P\{\mathcal X\in\mathcal O_{\alpha,\gamma,r}(a,t_0,\eps)\}.
    \end{equation}
    Since the event  $\mathcal O_{\alpha,\gamma,r}(a,t_0,\eps)$ is open,  the portmanteau  theorem and the definition of $\eta^{(r)}$ give
    \[
        \mathbb P\{\mathcal X\in\mathcal O_{\alpha,\gamma,r}(a,t_0,\eps)\} \le \limsup_{k\to\infty} \mathbb P\{\eta_{\mathcal X_{n_k}^{\Phi}}^{(r)} (t_0+\alpha,t-\alpha;A_{\eps+2\gamma}(a-\gamma))\ge3\}.
    \]
    Combining this with \eqref{eq:abstract-limit-event-cover}, we obtain
    \begin{equation}
        \mathbb P\{\eta_{\mathcal X}(t_0,t;A_\eps(a))\ge3\} \le \liminf_{\alpha\downarrow0} \liminf_{r\downarrow0} \limsup_{k\to\infty} \mathbb P\{\eta_{\mathcal X_{n_k}^{\Phi}}^{(r)} (t_0+\alpha,t-\alpha;A_{\eps+2\gamma}(a-\gamma))\ge3\}. \label{eq:abstract-portmanteau-separated-bound}
    \end{equation}
    Let $I_t=[t/2,t]$ and let $C_t:=C_{I_t}$ be the constant of  Proposition~\ref{prop:abstract-B2}. Choose $\gamma=\eps$ with $0<\eps<1/3$. Taking the supremum over $(a,t_0)\in\Xi\times\R$ and using $t-\alpha\in I_t$ and Proposition~\ref{prop:abstract-B2}, we obtain
    \[
        \sup_{(a,t_0)\in\Xi\times\R}\mathbb P\{\eta_{\mathcal X}(t_0,t;A_\eps(a))\ge3\} \le C_t(\eps+2\gamma)^2= 9C_t\eps^2 .
    \]
    Thus
    \[
        \lim_{\eps\downarrow0}\eps^{-1} \sup_{(a,t_0)\in\Xi\times\R} \mathbb P\{\eta_{\mathcal X}(t_0,t;A_\eps(a))\ge3\}=0 .
    \]
    Hence $\mathcal X$ satisfies condition $(B2)$ of Theorem~\ref{thm:BW-CBW-convergence-criterion}, and therefore has the law of $\sigma\mathcal W_\Xi$.    Since every subsequential limit has this law, $\mathcal X_n^\Phi\dto\sigma\mathcal W_\Xi$.
    \end{proof}

    \section{\texorpdfstring{Hastings--Levitov}{Hastings--Levitov}(0): Definitions and geometric preliminaries}\label{S.HL-application}
    In this section, we introduce the three HL(0) geometries and their observables, develop the estimates needed to control their boundary and geometric representations, and prove the required geometric tightness. These results are used in Section~\ref{S.HL-observables-BW} to establish the corresponding BW convergence.
    \subsection{HL\texorpdfstring{$(0)$}{(0)} geometries and observables}\label{S.HL-models}
\begin{figure}[H]
		\begin{center}

			\definecolor{myred}{RGB}{255, 0, 0}
			\definecolor{myblue}{RGB}{0, 0, 255}

			\tikzset{
				redcirc/.style={circle, draw=black, fill=myred, thin, inner sep=0pt, minimum size=3mm},
				bluecirc/.style={circle, draw=black, fill=myblue, thin, inner sep=0pt, minimum size=5mm},
			}

			\begin{tikzpicture}[scale=0.62]
				\node (v1) at (-1,-0.3) {\tiny{$-N\pi$}};
				\node (v2) at (1,-0.3) {\tiny{$N\pi$}};
				\draw [ultra thick] (-1,0) to (1,0);
				\draw[pattern=north west lines, pattern color=blue, opacity=0.2] (-1,0) rectangle (1,3);
				\node (v4) at (0,0) [redcirc] {\tiny{$0$}};
				\node (v3) at (1.5,1.5) {\tiny{$\overset{f_N}{\to}$}};
				\draw[pattern=north west lines, pattern color=blue, opacity=0.2] (1.8,0) rectangle (4.2,3);
				\filldraw[color=black!160, fill=white!5, ultra thick](3,1.5) circle (1);
				\node (v6) at (4,1.5) [redcirc] {\tiny{$1$}};
				\node (v5) at (4.6,1.5) {\tiny{$\overset{g}{\to}$}};
				\draw[pattern=north west lines, pattern color=blue, opacity=0.2] (5,0) rectangle (8,2);
				\draw [ultra thick] (5,0) to (8,0);
				\node (v7) at (6.5,0) [redcirc] {\tiny{$0$}};
				\node (v8) at (8.3,1.5) {\tiny{$\overset{\phi_\delta}{\to}$}};
				\draw[pattern=north west lines, pattern color=blue, opacity=0.2] (8.7,0) rectangle (11.7,2);
				\draw [ultra thick] (8.7,0) to (11.7,0);
				\draw [ultra thick] (10.2,0) to (10.2,0.7);
				\node (v8) at (10.2,0.7) [redcirc] {\tiny{$i\delta$}};

				\node (v8) at (10.2,-1) {\tiny{${\downarrow}{g^{-1}}$}};
				\draw[pattern=north west lines, pattern color=blue, opacity=0.2] (8.7,-5) rectangle (11.7,-2);
				\filldraw[color=black!160, fill=white!5, ultra thick](10,-3.5) circle (1);
				\draw [ultra thick] (11,-3.5) to (11.4,-3.5);
				\node (v9) at (11.4,-3.5) [redcirc] {};
				\node (v10) at (11.4,-3.1) {\tiny{$1+\eta$}};
				\node (v8) at (7.9,-3.5) {\tiny{$\overset{f_N^{-1}}{\leftarrow}$}};
				\begin{scope}[shift={(6,-5)}]
					\node (v11) at (-1,-0.3) {\tiny{$-N\pi$}};
					\node (v12) at (1,-0.3) {\tiny{$N\pi$}};
					\draw [ultra thick] (-1,0) to (1,0);
					\draw[pattern=north west lines, pattern color=blue, opacity=0.2] (-1,0) rectangle (1,3);
					\draw [ultra thick] (0,0) to (0,1.3);
					\node (v14) at (0,1.3) [redcirc] {\tiny{$i \lambda$}};
				\end{scope}

				\draw [-stealth](4,-0.5) -- (8.5,-2.5);
				\node (v15) at (6.7,-1.1) {$H^{\eta}$};

				\draw [-stealth](0,-0.5) -- (4.5,-2.5);
				\node (v15) at (2.7,-1.1) {$\slitN$};

				\draw[pattern=north west lines, pattern color=blue, opacity=0.2] (12.7,-2) rectangle (15.7,0);
				\draw [ultra thick] (12.7,-2) to (15.7,-2);
				\node (v7) at (14.2,-2) [redcirc] {\tiny{$0$}};

				\draw[pattern=north west lines, pattern color=blue, opacity=0.2] (16.7,-2) rectangle (19.7,0);
				\draw [ultra thick] (16.7,-2) to (19.7,-2);
				\draw [ultra thick] (18.2,-2) to (18.2,-1.3);
				\node (v8) at (18.2,-1.3) [redcirc] {\tiny{$i\lambda$}};

				\draw [-stealth](15.8,-0.9) -- (16.6,-0.9);
				\node (v15) at (16.2,-0.5) {$S^{\lambda}$};

			\end{tikzpicture}
		\end{center}
		\caption{Basic slit maps $H^{\eta}$, $S^{\lambda}$ and $C^{N,\lambda}$.}
		\label{basic_slit_map}
	\end{figure}
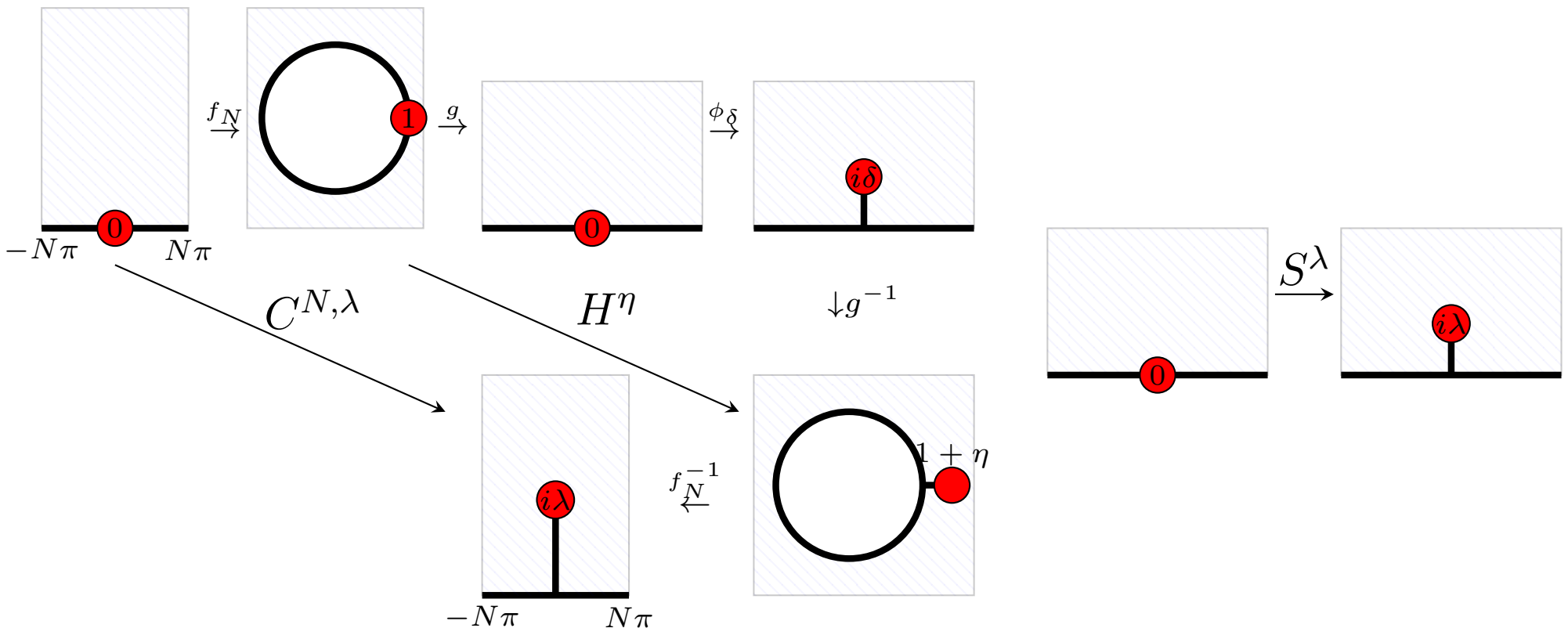

    We recall the radial, stationary, and cylindrical HL(0) processes, specify the normalizations, and define the boundary and geometric observables whose scaling limits are described by the BW or CBW.  We mainly follow the constructions and slit normalizations of \cite{hastings1998laplacian,norris2012hastings,berger2022growth,procaccia2023cylindrical, chen2025one}.  Let
    \begin{align*}
        &\mathbb T^N:=\{z\in\mathbb C:\operatorname{Im}z>0\}/(z\sim z+2\pi N), &\mathbb D_0&:=\{z\in\mathbb C:|z|>1\}, &\mathbb H&:=\{z\in\mathbb C:\operatorname{Im}z>0\}.
    \end{align*}
    We introduce the conformal maps $f_N:\mathbb  T^N\to\mathbb D_0$, $g:\mathbb D_0\to\mathbb H$, and $\phi_\delta:\mathbb H\to\mathbb H\setminus[0,\mathrm{i}\delta]$, given by
    \begin{align*}
        f_N(z)&=e^{-\mathrm i z/N}, &g(z)&=\mathrm i\frac{z-1}{z+1}, &\phi_\delta(z)&=\sqrt{z^2(1-\delta^2)-\delta^2}.
    \end{align*}
    For $\eta,\lambda>0$, the corresponding slit heights in $\mathbb H$ for the HL(0) and CHL(0) maps are given by
    \begin{equation*}
        \delta_{\mathcal H}(\eta):=\frac{\eta}{\eta+2}, \qquad \delta_N(\lambda):=1-\frac{2}{1+e^{\lambda/N}} =\tanh\!\left(\frac{\lambda}{2N}\right).
    \end{equation*}
    The three basic slit maps $H^\eta:\mathbb D_0\to\mathbb D_0\setminus[1,1+\eta]$, $S^\lambda:\mathbb H\to\mathbb H\setminus[0,\mathrm i\lambda]$, and $C^{N,\lambda}:\mathbb T^N\to\mathbb T^N\setminus[0,\mathrm i\lambda]$, illustrated in Figure~\ref{basic_slit_map}, are given by
    \[
        \resizebox{\linewidth}{!}{$\displaystyle    H^\eta(z)=g^{-1}\circ\phi_{\delta_{\mathcal H}(\eta)}\circ g(z), \quad S^\lambda(z)=\sqrt{z^2-\lambda^2}, \quad C^{N,\lambda}(z)=f_N^{-1}\circ g^{-1}   \circ\phi_{\delta_N(\lambda)}\circ g\circ f_N(z). $}
    \]
    For $\theta\in\mathbb T$, $a\in\mathbb R$, and $x\in\mathbb T_N$, the translated slit maps are defined by
    \begin{equation*}
        \begin{aligned}
        H_\theta^\eta(z)=e^{\mathrm i\theta}H^\eta(e^{-\mathrm i\theta}z),\qquad S_a^\lambda(z)=a+S^\lambda(z-a),\qquad C_x^{N,\lambda}(z)=x+C^{N,\lambda}(z-x).
        \end{aligned}
    \end{equation*}
    The three HL$(0)$ processes are constructed from these slit maps as follows.
    \begin{definition}\label{def:HL-geometries}
    \begin{enumerate}[label=\textup{(\roman*)},leftmargin=2.7em]
    \item \emph{HL$(0)$.}  Let $\mathcal P^{\mathcal H}$ be a PPP on $\mathbb R_+\times\mathbb T$ with intensity $\mathrm dt\,\mathrm dx/2\pi$, and write
    \[
        \mathcal P^{\mathcal H}\cap([0,t]\times\mathbb T) =\{(t_1,x_1),\ldots,(t_n,x_n)\}, \qquad 0<t_1<\cdots<t_n\le t.
    \]
    The HL(0) process is defined by
    \[
        \mathcal H_t^{\eta} :=H_{x_1}^{\eta}\circ\cdots\circ H_{x_n}^{\eta}, \qquad \mathcal H_0^{\eta}=\operatorname{id}.
    \]
    The logarithmic capacity increment of one particle is
    \begin{equation}\label{logarithmic capacity}
        c_{\mathcal H}(\eta):=\log (H^{\eta})'(\infty) =-\log\bigl(1-\delta_{\mathcal H}(\eta)^2\bigr).
    \end{equation}

    \item \emph{SHL$(0)$.} Let $\mathcal P^{\mathcal S}$ be a PPP on $\mathbb R_+\times\mathbb R$ with intensity $dt\,dx$, and let $A^{\mathcal S}$ be its set of Poisson points. Let $A_{t,m}:=A^{\mathcal S}\cap\{(s,x):s\le t,\ |x|\le m\}$.  The SHL(0) process is defined by
    \[
        \mathcal S_t^\lambda(z) =z+\lim_{m\to\infty}\sum_{i:(t_i,x_i)\in A_{t,m}} \left[ \mathcal S_{t_i-}^\lambda(S_{x_i}^\lambda(z)) -\mathcal S_{t_i-}^\lambda(z) \right], \qquad t\ge0,\ z\in\mathbb H,
    \]
    where the limit exists almost surely; see the remark below.

    \item \emph{CHL$(0)$.}  For each $N$, let $\mathcal P^{\mathcal C,N}$ be a PPP on $\mathbb R_+\times\mathbb T_N$ with intensity $\mathrm dt\,\mathrm dx$, and write
    \[
        \mathcal P^{\mathcal C,N}\cap([0,t]\times\mathbb T_N)=\{(t_1,x_1),\ldots,(t_n,x_n)\}, \qquad 0<t_1<\cdots<t_n\le t.
    \]
    The CHL(0) process is defined by
    \[
        \mathcal C_t^{N,\lambda} :=C_{x_1}^{N,\lambda}\circ C_{x_2}^{N,\lambda} \circ\cdots\circ C_{x_n}^{N,\lambda}, \qquad \mathcal C_0^{N,\lambda}=\operatorname{id}.
    \]
    \end{enumerate}
    \end{definition}
    \begin{remark}[{\cite[Remark~3.2] {procaccia2023cylindrical}}]\label{rmk:radial-from-cylindrical}
    For fixed $\eta>0$, set $\lambda_N^{\mathcal H}:=N\log(1+\eta/N)$.	Under the coupling of the driving PPPs via the map $(t,x)\mapsto(2\pi Nt,-x/N)$, one has
    \[
        \mathcal H_{N^3t}^{N^{-1}\eta}		=f_N\circ\mathcal C_{N^2t/(2\pi)}^{N,\lambda_N^{\mathcal H}}	\circ f_N^{-1}.
    \]
    Thus radial HL$(0)$ with particle size $N^{-1}\eta$ is a conformal image of CHL$(0)$ with particle height $\lambda_N^{\mathcal H}$.
    \end{remark}

    \begin{remark}
    The existence of SHL$(0)$ and uniqueness in distribution were proved in \cite[Theorems~3.1 and~4.3]{berger2022growth}.
    \end{remark}
    \paragraph{\textbf{Geometric coordinates.}} We use $\mathcal A\in\{\mathcal H,\mathcal S,\mathcal C\}$ to index radial HL$(0)$, SHL$(0)$, and CHL$(0)$, respectively. For the observables below, we use the following notation for the set of geometries, its subsets of models with circular boundary and with flat geometry, and their spatial domains:
    \[
        \mathfrak M:=\{\mathcal H,\mathcal S,\mathcal C\}, \qquad \mathfrak M_{\rm c}:=\{\mathcal H,\mathcal C\}, \qquad \mathfrak M_{\rm f}:=\{\mathcal S,\mathcal C\}, \qquad E^{\mathcal H}=E^{\mathcal C}:=\mathbb T, \qquad E^{\mathcal S}:=\mathbb R.
    \]
    The geometric projections are defined by
    \[
        \begin{aligned}
        \mathfrak p_N^{\mathcal H}(z)&:=\exp\!\left(N^{-1}\log|z|+\mathrm i\arg z\right),\qquad 	\mathfrak p_N^{\mathcal S}(z):=N^{-1}\rea z+\mathrm iN^{-2}\im z, \\
        \mathfrak p_N^{\mathcal C}(z)&:=(N^{-1}\rea z\ {\rm mod}\ 2\pi)+\mathrm iN^{-2}\im z.
        \end{aligned}
    \]
    The geometric spaces and their coordinate maps are
    \[
        \mathbb G_{\mathcal H}=\bar{\mathbb D}_0, \chi_{\mathcal H}(z)=(\arg z,\log|z|);\ \mathbb G_{\mathcal S}=\bar{\mathbb H}, \chi_{\mathcal S}(z)=(\rea z,\im z);\ \mathbb G_{\mathcal C}=\bar{\mathbb T}^1, \chi_{\mathcal C}(z)=(\rea z,\im z).
    \]

    \medskip
    \paragraph{\textbf{Harmonic measure and Root paths.}}

    For $0\le s\le t$, the rescaled boundary transition, which sends the harmonic measure coordinate of a boundary point at time $s$ to its coordinate at time $t$, is
    \begin{equation}\label{eq:rescaled-boundary-transition}
        U_{s,t}^{\mathcal A,N}(x)=
        \begin{cases}
        \arg[(\mathcal H_{N^3t}^{N^{-1}\eta})^{-1} (\mathcal H_{N^3s}^{N^{-1}\eta}(e^{\mathrm ix}))],&\mathcal A=\mathcal H, \\
        N^{-1}(\mathcal S_{N^2t}^{\lambda})^{-1} (\mathcal S_{N^2s}^{\lambda}(Nx)),&\mathcal A=\mathcal S, \\
        N^{-1}(\mathcal C_{N^2t}^{N,\lambda})^{-1} (\mathcal C_{N^2s}^{N,\lambda}(Nx)),&\mathcal A=\mathcal C.
        \end{cases}
    \end{equation}
    Boundary values are taken to be non-decreasing and right-continuous. At each discontinuity $x$, the completed graph contains the vertical segment $\{x\}\times[U_{s,t}^{\mathcal A,N}(x-),U_{s,t}^{\mathcal A,N}(x)]$. For a monotone boundary map $f$, define its generalized inverse by
    \begin{align}\label{generalized-inverse}
        f^\leftarrow(y):=\inf\{x:f(x)>y\},\qquad R_{s,t}^{\mathcal A,N}:=(U_{s,t}^{\mathcal A,N})^\leftarrow.
    \end{align}
    To interpolate the paths, define the grid times
    \[
        \Delta_N^{\mathcal H}=2\pi N^{-1},\qquad \Delta_N^{\mathcal S}=\Delta_N^{\mathcal C}=N^{-1},    \qquad \Gamma_N^{\mathcal A}=\Delta_N^{\mathcal A}\mathbb Z_+.
    \]
    For harmonic measure paths started at time $s\ge0$, we use the interpolation times $\{s\}\cup(\Gamma_N^{\mathcal A}\cap(s,\infty))$, and for Root paths ending at time $t\ge0$, we use $(\Gamma_N^{\mathcal A}\cap[0,t])\cup\{t\}$. For each path, let $a<b$ be two of its interpolation times with no other interpolation time between them. For $u=(1-\theta)a+\theta b$, $0\le\theta\le1$, define the harmonic measure paths and Root paths by
    \[
        \begin{aligned}
        \mathtt{HM}_{s,u}^{\mathcal A,N}(x):=(1-\theta)U_{s,a}^{\mathcal A,N}(x)+\theta U_{s,b}^{\mathcal A,N}(x),\qquad \mathtt R_{u,t}^{\mathcal A,N}(x):=(1-\theta)R_{a,t}^{\mathcal A,N}(x)+\theta R_{b,t}^{\mathcal A,N}(x).
        \end{aligned}
    \]
    In the circular cases, the interpolation is carried out for the lifts and then projected to $\mathbb T$. Thus  $\mathtt{HM}_{s,s}^{\mathcal A,N}(x)   =\mathtt R_{t,t}^{\mathcal A,N}(x)=x$.   Note that the generalized inverse is taken before time interpolation.

    For the $i$-th particle, with arrival time and position $(t_i,x_i)$ as in Definition~\ref{def:HL-geometries}, define its rescaled coordinates and the labels of the two sides of its base by
    \begin{equation}\label{eq:particle-side-labels}
        \begin{aligned}
        (\bar t_i,\bar x_i)=
        \begin{cases}
        (N^{-3}t_i,x_i),&\mathcal A=\mathcal H, \\
        (N^{-2}t_i,N^{-1}x_i),&\mathcal A\in\mathfrak M_{\rm f},
        \end{cases}
        \qquad \bar x_{i,\pm}=\bar x_i\pm
        \begin{cases}
        2\arcsin\delta_{\mathcal H}(N^{-1}\eta), &\mathcal A=\mathcal H, \\
        \lambda/N,&\mathcal A=\mathcal S, \\
        2\arcsin\delta_N(\lambda),&\mathcal A=\mathcal C.
        \end{cases}
        \end{aligned}
    \end{equation}
    The tip of the particle has label $\bar x_i$. These identities follow from the one-slit maps; for example, in SHL$(0)$,
    \[
        S_{x_i}^{\lambda}(N\bar x_{i,\pm})=x_i, \qquad S_{x_i}^{\lambda}(N\bar x_i)=x_i+\mathrm i\lambda.
    \]
    The labels $\bar x_{i,-}$ and $\bar x_{i,+}$ correspond to the two sides of the slit in boundary order.

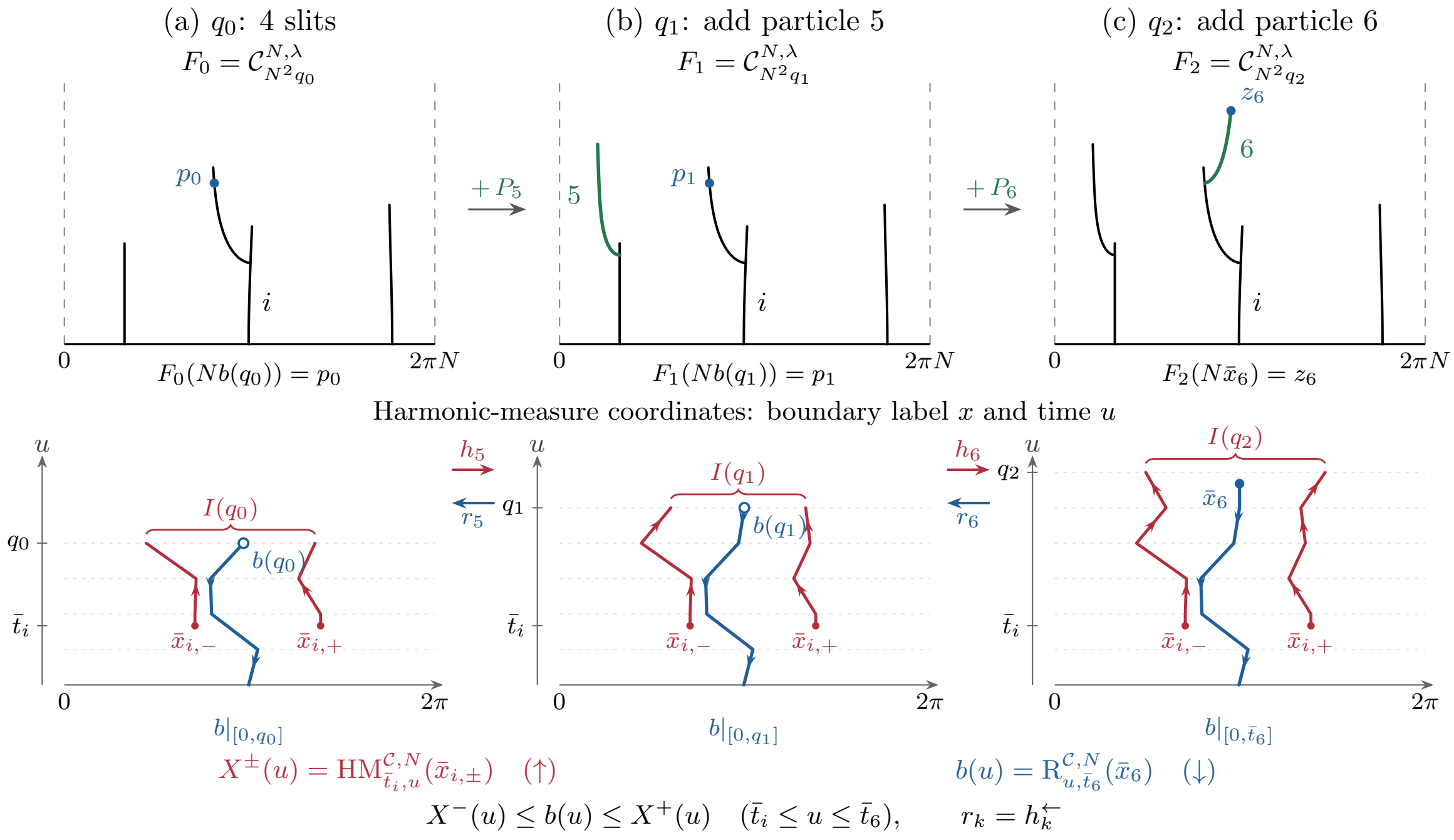
\begin{figure}[H]
        \centering
        \resizebox{\linewidth}{!}{%
%
%
\begingroup
\definecolor{CHLhm}{HTML}{B42E3A}
\definecolor{CHLroot}{HTML}{205F9A}
\definecolor{CHLnew}{HTML}{24794F}
\definecolor{CHLgray}{HTML}{7A7A7A}
\tikzset{
  chl/old/.style={draw=black,line width=.83pt,line cap=round,line join=round},
  chl/new/.style={draw=CHLnew,line width=1.15pt,line cap=round,line join=round},
  chl/hm/.style={draw=CHLhm,line width=1.0pt,line cap=round,line join=round},
  chl/root/.style={draw=CHLroot,line width=1.08pt,line cap=round,line join=round},
  chl/rail/.style={draw=black!45,line width=.45pt,dash pattern=on 3pt off 3pt},
  chl/grid/.style={draw=black!14,line width=.3pt,dash pattern=on 1pt off 2pt},
  chl/text/.style={font=\fontsize{9}{10.5}\selectfont,inner sep=1pt},
  chl/small/.style={font=\fontsize{8}{9.5}\selectfont,inner sep=1pt},
}
\begin{tikzpicture}[x=1cm,y=1cm,every node/.style={chl/text}]

\def\CHLedgeOne{%
  (1.02,0) -- (1.02,0.991094) -- (1.02,1.12)%
}
\def\CHLedgeTwo{%
  (3.125171,0) -- (3.125284,0.0502) -- (3.125626,0.100686) -- (3.126197,0.151438) -- (3.126998,0.202436) --
  (3.128033,0.253963) -- (3.1293,0.305997) -- (3.130809,0.35882) -- (3.132556,0.412406) -- (3.134539,0.466731) --
  (3.136764,0.522069) -- (3.139253,0.578993) -- (3.142019,0.637762) -- (3.145072,0.698635) -- (3.148439,0.762156) --
  (3.152179,0.829441) -- (3.156353,0.90158) -- (3.165286,1.049393) -- (3.181435,1.30829)%
}
\def\CHLedgeThree{%
  (5.559843,0) -- (5.559735,0.052063) -- (5.559408,0.104647) -- (5.558866,0.157128) -- (5.558102,0.210591) --
  (5.557112,0.264975) -- (5.555901,0.320215) -- (5.554447,0.377361) -- (5.552734,0.436874) -- (5.54958,0.533678) --
  (5.545631,0.641867) -- (5.541535,0.746407) -- (5.529999,1.032239) -- (5.526605,1.12016) -- (5.5236,1.201736) --
  (5.520504,1.290933) -- (5.517721,1.377208) -- (5.515189,1.462552) -- (5.512892,1.547495)%
}
\def\CHLedgeFour{%
  (3.156353,0.90158) -- (3.141874,0.902665) -- (3.127442,0.904211) -- (3.113079,0.906214) -- (3.10593,0.907385) --
  (3.098334,0.908757) -- (3.084181,0.911669) -- (3.069699,0.915131) -- (3.055842,0.918908) -- (3.041705,0.923238) --
  (3.027772,0.927985) -- (3.013619,0.933305) -- (2.999714,0.939032) -- (2.98607,0.945144) -- (2.972699,0.951622) --
  (2.959204,0.958662) -- (2.946014,0.966046) -- (2.932751,0.97399) -- (2.919821,0.982251) -- (2.906862,0.991062) --
  (2.894258,1.000162) -- (2.881663,1.009799) -- (2.86944,1.019694) -- (2.85726,1.03011) -- (2.845461,1.040755) --
  (2.833732,1.051902) -- (2.822096,1.06355) -- (2.810857,1.075385) -- (2.79973,1.0877) -- (2.788999,1.100176) --
  (2.778395,1.113112) -- (2.767936,1.126503) -- (2.757636,1.140342) -- (2.747508,1.154622) -- (2.737998,1.16868) --
  (2.728647,1.183157) -- (2.719467,1.198047) -- (2.710658,1.213009) -- (2.701838,1.228704) -- (2.69339,1.244456) --
  (2.685129,1.260595) -- (2.67706,1.277113) -- (2.66919,1.294005) -- (2.661522,1.311266) -- (2.65406,1.32889) --
  (2.64694,1.346533) -- (2.640017,1.364529) -- (2.633172,1.383212) -- (2.62654,1.402239) -- (2.620119,1.421606) --
  (2.61391,1.441308) -- (2.607911,1.461342) -- (2.60212,1.481705) -- (2.596534,1.502392) -- (2.591068,1.52374) --
  (2.585811,1.545407) -- (2.580759,1.567389) -- (2.575836,1.590024) -- (2.571118,1.612969) -- (2.566536,1.636561) --
  (2.562154,1.66046) -- (2.557911,1.685) -- (2.553811,1.710179) -- (2.549905,1.735661) -- (2.54614,1.761777) --
  (2.542518,1.788528) -- (2.539894,1.809014) -- (2.537347,1.829889) -- (2.534876,1.851152) -- (2.532483,1.872802) --
  (2.527927,1.917264) -- (2.523654,1.963502)%
}
\def\CHLedgeFive{%
  (1.02,0.991094) -- (1.014756,0.991162) -- (1.009078,0.99139) -- (1.003408,0.991777) -- (0.997749,0.992323) --
  (0.992538,0.992967) -- (0.98691,0.993817) -- (0.981304,0.994823) -- (0.975723,0.995984) -- (0.970171,0.997299) --
  (0.964651,0.998766) -- (0.959166,1.000383) -- (0.953719,1.002148) -- (0.948313,1.004059) -- (0.942952,1.006114) --
  (0.937638,1.008311) -- (0.932373,1.010647) -- (0.92716,1.013119) -- (0.922002,1.015725) -- (0.916899,1.018462) --
  (0.91147,1.021553) -- (0.906491,1.024553) -- (0.901574,1.027676) -- (0.89672,1.030917) -- (0.891566,1.034537) --
  (0.886849,1.038016) -- (0.881844,1.041885) -- (0.876918,1.045878) -- (0.872073,1.049991) -- (0.867309,1.054219) --
  (0.862626,1.058559) -- (0.858026,1.063008) -- (0.853509,1.067562) -- (0.849074,1.072216) -- (0.844722,1.076968) --
  (0.840453,1.081815) -- (0.835971,1.087107) -- (0.831873,1.092137) -- (0.827573,1.097618) -- (0.819525,1.108474) --
  (0.815497,1.114214) -- (0.81156,1.120035) -- (0.803956,1.131906) -- (0.796467,1.144469) -- (0.789341,1.157305) --
  (0.782354,1.170814) -- (0.775728,1.184559) -- (0.769256,1.198956) -- (0.762957,1.213998) -- (0.757013,1.229228) --
  (0.751248,1.245078) -- (0.74567,1.26154) -- (0.74029,1.278604) -- (0.735642,1.294395) -- (0.731136,1.310744) --
  (0.726778,1.327643) -- (0.722574,1.345085) -- (0.71863,1.362587) -- (0.714734,1.381084) -- (0.710999,1.400095) --
  (0.707422,1.419609) -- (0.704003,1.439615) -- (0.700665,1.460578) -- (0.697488,1.48201) -- (0.694465,1.5039) --
  (0.691474,1.527183) -- (0.688642,1.550897) -- (0.685908,1.575499) -- (0.683274,1.600974) -- (0.680826,1.626366) --
  (0.678423,1.65307) -- (0.676108,1.680597) -- (0.673842,1.709395) -- (0.671627,1.739438) -- (0.669432,1.771159) --
  (0.667261,1.804525) -- (0.665114,1.839502) -- (0.66084,1.915039) -- (0.656396,2.001416) -- (0.651632,2.101597) --
  (0.646337,2.220019)%
}
\def\CHLedgeSix{%
  (2.542518,1.788528) -- (2.552394,1.790005) -- (2.562478,1.79189) -- (2.57246,1.794137) -- (2.582328,1.796738) --
  (2.592075,1.799684) -- (2.601973,1.803068) -- (2.611448,1.806686) -- (2.621052,1.81074) -- (2.630498,1.815119) --
  (2.63978,1.81981) -- (2.649153,1.82495) -- (2.658346,1.830395) -- (2.667355,1.836135) -- (2.676424,1.842328) --
  (2.685297,1.848806) -- (2.694207,1.855743) -- (2.702912,1.862954) -- (2.711414,1.870429) -- (2.719716,1.878154) --
  (2.72803,1.886331) -- (2.736345,1.894968) -- (2.744449,1.903846) -- (2.752346,1.912956) -- (2.760232,1.922521) --
  (2.767911,1.932308) -- (2.77557,1.942553) -- (2.783025,1.953009) -- (2.790452,1.963925) -- (2.797679,1.975045) --
  (2.804872,1.986624) -- (2.811868,1.9984) -- (2.818827,2.010637) -- (2.825594,2.023063) -- (2.83232,2.035952) --
  (2.83886,2.049022) -- (2.845357,2.062557) -- (2.851674,2.076267) -- (2.857945,2.090441) -- (2.864167,2.105085) --
  (2.870213,2.119898) -- (2.876209,2.135181) -- (2.882151,2.150936) -- (2.887924,2.166854) -- (2.893643,2.183244) --
  (2.899304,2.200109) -- (2.904803,2.217129) -- (2.910244,2.234623) -- (2.915529,2.252264) -- (2.920758,2.270378) --
  (2.925927,2.288965) -- (2.931035,2.308026) -- (2.936079,2.327563) -- (2.941057,2.347576) -- (2.945888,2.367722) --
  (2.950655,2.388341) -- (2.955356,2.409431) -- (2.959991,2.430993) -- (2.964485,2.452675) -- (2.968984,2.475176) --
  (2.973345,2.497791) -- (2.977703,2.521227) -- (2.981926,2.54477) -- (2.986081,2.568775) -- (2.990165,2.593239)%
}

\begin{scope}[shift={(0.66,0)}]
\node[font=\fontsize{10}{12}\selectfont] at (2.17,8.64) {(a) $q_0$: 4 slits};
\node at (2.17,8.15) {$F_0=\mathcal C_{N^2q_0}^{N,\lambda}$};
\begin{scope}[shift={(0,4.87)},x=0.690732cm,y=1.055cm]
\draw[chl/rail] (0,0) -- (0,2.94);
\draw[chl/rail] (6.283185,0) -- (6.283185,2.94);
\draw[black,line width=.65pt] (0,0) -- (6.283185,0);
\draw[chl/old] \CHLedgeOne;
\draw[chl/old] \CHLedgeTwo;
\draw[chl/old] \CHLedgeThree;
\draw[chl/old] \CHLedgeFour;
\node[anchor=west] at (3.29471,0.471175) {$i$};
\fill[CHLroot] (2.542518,1.788528) circle[radius=1.55pt];
\node[text=CHLroot,anchor=east] at (2.412518,1.848528) {$p_0$};
\end{scope}
\node[chl/small,anchor=north] at (0,4.82) {$0$};
\node[chl/small,anchor=north] at (4.34,4.82) {$2\pi N$};
\node[chl/small] at (2.17,4.51) {$F_0(Nb(q_0))=p_0$};
\begin{scope}[shift={(0,0.88)}]
\draw[chl/grid] (0,0.415) -- (4.34,0.415);
\draw[chl/grid] (0,0.83) -- (4.34,0.83);
\draw[chl/grid] (0,1.245) -- (4.34,1.245);
\draw[chl/grid] (0,1.66) -- (4.34,1.66);
\draw[black!60,line width=.45pt,-{Stealth[length=1.5mm]}] (-.26,0) -- (-0.26,2.69);
\node[anchor=south,text=black!75] at (-0.26,2.69) {$u$};
\draw[black!60,line width=.45pt,-{Stealth[length=1.5mm]}] (0,0) -- (4.44,0);
\node[chl/small,anchor=north] at (0,-.06) {$0$};
\node[chl/small,anchor=north] at (4.34,-0.06) {$2\pi$};
\draw[black!60,line width=.4pt] (-0.32,0.69305) -- (-0.2,0.69305);
\node[chl/small,anchor=east] at (-0.35,0.69305) {$\bar t_i$};
\draw[black!60,line width=.4pt] (-0.32,1.66) -- (-0.2,1.66);
\node[chl/small,anchor=east] at (-0.35,1.66) {$q_0$};
\draw[chl/hm] (1.529509,0.69305) -- (1.529509,0.83) -- (1.543325,1.245) -- (0.958088,1.66);
\draw[chl/hm,-{Stealth[length=1.55mm,width=1.05mm]}] (1.538075,1.0873) -- (1.541667,1.1952);
\draw[chl/hm] (3.001695,0.69305) -- (3.001695,0.83) -- (2.74312,1.245) -- (2.93607,1.66);
\draw[chl/hm,-{Stealth[length=1.55mm,width=1.05mm]}] (2.841379,1.0873) -- (2.774149,1.1952);
\draw[chl/root] (2.158657,0) -- (2.265602,0.415) -- (1.724556,0.83) -- (1.711296,1.245) -- (2.098279,1.66);
\draw[chl/root,-{Stealth[length=1.55mm,width=1.05mm]}] (2.236727,0.30295) -- (2.208921,0.19505);
\draw[chl/root,-{Stealth[length=1.55mm,width=1.05mm]}] (1.712224,1.21595) -- (1.715672,1.10805);
\fill[CHLhm] (1.529509,0.69305) circle[radius=1.25pt];
\node[chl/small,text=CHLhm,anchor=north] at (1.529509,0.64305) {$\bar x_{i,-}$};
\fill[CHLhm] (3.001695,0.69305) circle[radius=1.25pt];
\node[chl/small,text=CHLhm,anchor=north] at (3.001695,0.64305) {$\bar x_{i,+}$};
\draw[CHLroot,line width=.75pt,fill=white] (2.098279,1.66) circle[radius=1.6pt];
\node[chl/small,text=CHLroot,anchor=north west] at (2.158279,1.62) {$b(q_0)$};
\draw[CHLhm,line width=.6pt,decorate,decoration={brace,amplitude=3.0pt}] (0.958088,1.765) -- (2.93607,1.765);
\node[chl/small,text=CHLhm,anchor=south] at (1.947079,1.875) {$I(q_0)$};
\end{scope}
\node[chl/small,text=CHLroot] at (2.17,0.34) {$b|_{[0,q_0]}$};
\end{scope}

\begin{scope}[shift={(6.46,0)}]
\node[font=\fontsize{10}{12}\selectfont] at (2.17,8.64) {(b) $q_1$: add particle $5$};
\node at (2.17,8.15) {$F_1=\mathcal C_{N^2q_1}^{N,\lambda}$};
\begin{scope}[shift={(0,4.87)},x=0.690732cm,y=1.055cm]
\draw[chl/rail] (0,0) -- (0,2.94);
\draw[chl/rail] (6.283185,0) -- (6.283185,2.94);
\draw[black,line width=.65pt] (0,0) -- (6.283185,0);
\draw[chl/old] \CHLedgeOne;
\draw[chl/old] \CHLedgeTwo;
\draw[chl/old] \CHLedgeThree;
\draw[chl/old] \CHLedgeFour;
\draw[chl/new] \CHLedgeFive;
\node[anchor=west] at (3.29471,0.471175) {$i$};
\fill[CHLroot] (2.542518,1.788528) circle[radius=1.55pt];
\node[text=CHLroot,anchor=east] at (2.412518,1.848528) {$p_1$};
\node[text=CHLnew,anchor=east] at (0.429039,1.656052) {$5$};
\end{scope}
\node[chl/small,anchor=north] at (0,4.82) {$0$};
\node[chl/small,anchor=north] at (4.34,4.82) {$2\pi N$};
\node[chl/small] at (2.17,4.51) {$F_1(Nb(q_1))=p_1$};
\begin{scope}[shift={(0,0.88)}]
\draw[chl/grid] (0,0.415) -- (4.34,0.415);
\draw[chl/grid] (0,0.83) -- (4.34,0.83);
\draw[chl/grid] (0,1.245) -- (4.34,1.245);
\draw[chl/grid] (0,1.66) -- (4.34,1.66);
\draw[chl/grid] (0,2.075) -- (4.34,2.075);
\draw[black!60,line width=.45pt,-{Stealth[length=1.5mm]}] (-.26,0) -- (-0.26,2.69);
\node[anchor=south,text=black!75] at (-0.26,2.69) {$u$};
\draw[black!60,line width=.45pt,-{Stealth[length=1.5mm]}] (0,0) -- (4.44,0);
\node[chl/small,anchor=north] at (0,-.06) {$0$};
\node[chl/small,anchor=north] at (4.34,-0.06) {$2\pi$};
\draw[black!60,line width=.4pt] (-0.32,0.69305) -- (-0.2,0.69305);
\node[chl/small,anchor=east] at (-0.35,0.69305) {$\bar t_i$};
\draw[black!60,line width=.4pt] (-0.32,2.075) -- (-0.2,2.075);
\node[chl/small,anchor=east] at (-0.35,2.075) {$q_1$};
\draw[chl/hm] (1.529509,0.69305) -- (1.529509,0.83) -- (1.543325,1.245) -- (0.958088,1.66) -- (1.305438,2.075);
\draw[chl/hm,-{Stealth[length=1.55mm,width=1.05mm]}] (1.538075,1.0873) -- (1.541667,1.1952);
\draw[chl/hm,-{Stealth[length=1.55mm,width=1.05mm]}] (1.103975,1.8343) -- (1.194286,1.9422);
\draw[chl/hm] (3.001695,0.69305) -- (3.001695,0.83) -- (2.74312,1.245) -- (2.93607,1.66) -- (2.883781,2.075);
\draw[chl/hm,-{Stealth[length=1.55mm,width=1.05mm]}] (2.841379,1.0873) -- (2.774149,1.1952);
\draw[chl/hm,-{Stealth[length=1.55mm,width=1.05mm]}] (2.914109,1.8343) -- (2.900514,1.9422);
\draw[chl/root] (2.158657,0) -- (2.265602,0.415) -- (1.724556,0.83) -- (1.711296,1.245) -- (2.098279,1.66) -- (2.16519,2.075);
\draw[chl/root,-{Stealth[length=1.55mm,width=1.05mm]}] (2.236727,0.30295) -- (2.208921,0.19505);
\draw[chl/root,-{Stealth[length=1.55mm,width=1.05mm]}] (1.712224,1.21595) -- (1.715672,1.10805);
\draw[chl/root,-{Stealth[length=1.55mm,width=1.05mm]}] (2.150469,1.9837) -- (2.133072,1.8758);
\fill[CHLhm] (1.529509,0.69305) circle[radius=1.25pt];
\node[chl/small,text=CHLhm,anchor=north] at (1.529509,0.64305) {$\bar x_{i,-}$};
\fill[CHLhm] (3.001695,0.69305) circle[radius=1.25pt];
\node[chl/small,text=CHLhm,anchor=north] at (3.001695,0.64305) {$\bar x_{i,+}$};
\draw[CHLroot,line width=.75pt,fill=white] (2.16519,2.075) circle[radius=1.6pt];
\node[chl/small,text=CHLroot,anchor=north west] at (2.22519,2.035) {$b(q_1)$};
\draw[CHLhm,line width=.6pt,decorate,decoration={brace,amplitude=3.0pt}] (1.305438,2.18) -- (2.883781,2.18);
\node[chl/small,text=CHLhm,anchor=south] at (2.094609,2.29) {$I(q_1)$};
\end{scope}
\node[chl/small,text=CHLroot] at (2.17,0.34) {$b|_{[0,q_1]}$};
\end{scope}

\begin{scope}[shift={(12.26,0)}]
\node[font=\fontsize{10}{12}\selectfont] at (2.17,8.64) {(c) $q_2$: add particle $6$};
\node at (2.17,8.15) {$F_2=\mathcal C_{N^2q_2}^{N,\lambda}$};
\begin{scope}[shift={(0,4.87)},x=0.690732cm,y=1.055cm]
\draw[chl/rail] (0,0) -- (0,2.94);
\draw[chl/rail] (6.283185,0) -- (6.283185,2.94);
\draw[black,line width=.65pt] (0,0) -- (6.283185,0);
\draw[chl/old] \CHLedgeOne;
\draw[chl/old] \CHLedgeTwo;
\draw[chl/old] \CHLedgeThree;
\draw[chl/old] \CHLedgeFour;
\draw[chl/old] \CHLedgeFive;
\draw[chl/new] \CHLedgeSix;
\node[anchor=west] at (3.29471,0.471175) {$i$};
\fill[CHLroot] (2.990165,2.593239) circle[radius=1.55pt];
\node[text=CHLroot,anchor=south west] at (3.090165,2.643239) {$z_6$};
\node[text=CHLnew,anchor=west] at (3.083411,2.174356) {$6$};
\end{scope}
\node[chl/small,anchor=north] at (0,4.82) {$0$};
\node[chl/small,anchor=north] at (4.34,4.82) {$2\pi N$};
\node[chl/small] at (2.17,4.51) {$F_2(N\bar x_6)=z_6$};
\begin{scope}[shift={(0,0.88)}]
\draw[chl/grid] (0,0.415) -- (4.34,0.415);
\draw[chl/grid] (0,0.83) -- (4.34,0.83);
\draw[chl/grid] (0,1.245) -- (4.34,1.245);
\draw[chl/grid] (0,1.66) -- (4.34,1.66);
\draw[chl/grid] (0,2.075) -- (4.34,2.075);
\draw[chl/grid] (0,2.49) -- (4.34,2.49);
\draw[black!60,line width=.45pt,-{Stealth[length=1.5mm]}] (-.26,0) -- (-0.26,2.69);
\node[anchor=south,text=black!75] at (-0.26,2.69) {$u$};
\draw[black!60,line width=.45pt,-{Stealth[length=1.5mm]}] (0,0) -- (4.44,0);
\node[chl/small,anchor=north] at (0,-.06) {$0$};
\node[chl/small,anchor=north] at (4.34,-0.06) {$2\pi$};
\draw[black!60,line width=.4pt] (-0.32,0.69305) -- (-0.2,0.69305);
\node[chl/small,anchor=east] at (-0.35,0.69305) {$\bar t_i$};
\draw[black!60,line width=.4pt] (-0.32,2.49) -- (-0.2,2.49);
\node[chl/small,anchor=east] at (-0.35,2.49) {$q_2$};
\draw[chl/hm] (1.529509,0.69305) -- (1.529509,0.83) -- (1.543325,1.245) -- (0.958088,1.66) -- (1.305438,2.075) -- (1.06613,2.49);
\draw[chl/hm,-{Stealth[length=1.55mm,width=1.05mm]}] (1.538075,1.0873) -- (1.541667,1.1952);
\draw[chl/hm,-{Stealth[length=1.55mm,width=1.05mm]}] (1.103975,1.8343) -- (1.194286,1.9422);
\draw[chl/hm,-{Stealth[length=1.55mm,width=1.05mm]}] (1.192963,2.27005) -- (1.130743,2.37795);
\draw[chl/hm] (3.001695,0.69305) -- (3.001695,0.83) -- (2.74312,1.245) -- (2.93607,1.66) -- (2.883781,2.075) -- (3.168951,2.49);
\draw[chl/hm,-{Stealth[length=1.55mm,width=1.05mm]}] (2.841379,1.0873) -- (2.774149,1.1952);
\draw[chl/hm,-{Stealth[length=1.55mm,width=1.05mm]}] (2.914109,1.8343) -- (2.900514,1.9422);
\draw[chl/hm,-{Stealth[length=1.55mm,width=1.05mm]}] (3.017811,2.27005) -- (3.091955,2.37795);
\draw[chl/root] (2.158657,0) -- (2.265602,0.415) -- (1.724556,0.83) -- (1.711296,1.245) -- (2.098279,1.66) -- (2.16519,2.075) -- (2.16519,2.3572);
\draw[chl/root,-{Stealth[length=1.55mm,width=1.05mm]}] (2.236727,0.30295) -- (2.208921,0.19505);
\draw[chl/root,-{Stealth[length=1.55mm,width=1.05mm]}] (1.712224,1.21595) -- (1.715672,1.10805);
\draw[chl/root,-{Stealth[length=1.55mm,width=1.05mm]}] (2.150469,1.9837) -- (2.133072,1.8758);
\fill[CHLhm] (1.529509,0.69305) circle[radius=1.25pt];
\node[chl/small,text=CHLhm,anchor=north] at (1.529509,0.64305) {$\bar x_{i,-}$};
\fill[CHLhm] (3.001695,0.69305) circle[radius=1.25pt];
\node[chl/small,text=CHLhm,anchor=north] at (3.001695,0.64305) {$\bar x_{i,+}$};
\fill[CHLroot] (2.16519,2.3572) circle[radius=1.6pt];
\node[chl/small,text=CHLroot,anchor=north east] at (2.08519,2.3122) {$\bar x_6$};
\draw[CHLhm,line width=.6pt,decorate,decoration={brace,amplitude=3.0pt}] (1.06613,2.595) -- (3.168951,2.595);
\node[chl/small,text=CHLhm,anchor=south] at (2.117541,2.705) {$I(q_2)$};
\end{scope}
\node[chl/small,text=CHLroot] at (2.17,0.34) {$b|_{[0,\bar t_6]}$};
\end{scope}
\draw[black!65,line width=.65pt,-{Stealth[length=2mm]}] (5.39,6.45) -- (6.07,6.45);
\node[chl/small,text=CHLnew,anchor=south] at (5.73,6.54) {$+\,P_5$};
\draw[CHLhm,line width=.85pt,-{Stealth[length=1.7mm]}] (5.2,3.4) -- (5.7,3.4);
\node[chl/small,text=CHLhm,anchor=south] at (5.45,3.47) {$h_5$};
\draw[CHLroot,line width=.85pt,-{Stealth[length=1.7mm]}] (5.7,3.01) -- (5.2,3.01);
\node[chl/small,text=CHLroot,anchor=north] at (5.45,2.93) {$r_5$};
\draw[black!65,line width=.65pt,-{Stealth[length=2mm]}] (11.19,6.45) -- (11.87,6.45);
\node[chl/small,text=CHLnew,anchor=south] at (11.53,6.54) {$+\,P_6$};
\draw[CHLhm,line width=.85pt,-{Stealth[length=1.7mm]}] (11,3.4) -- (11.5,3.4);
\node[chl/small,text=CHLhm,anchor=south] at (11.25,3.47) {$h_6$};
\draw[CHLroot,line width=.85pt,-{Stealth[length=1.7mm]}] (11.5,3.01) -- (11,3.01);
\node[chl/small,text=CHLroot,anchor=north] at (11.25,2.93) {$r_6$};
\node[font=\fontsize{9}{11}\selectfont] at (8.63,4.065) {Harmonic-measure coordinates: boundary label $x$ and time $u$};
\node[text=CHLhm] at (4.45,-.14) {$X^{\pm}(u)=\mathrm{HM}^{\mathcal C,N}_{\bar t_i,u}(\bar x_{i,\pm})\quad(\uparrow)$};
\node[text=CHLroot] at (12.63,-.14) {$b(u)=\mathrm{R}^{\mathcal C,N}_{u,\bar t_6}(\bar x_6)\quad(\downarrow)$};
\node at (8.63,-.65) {$X^-(u)\le b(u)\le X^+(u)\quad (\bar t_i\le u\le\bar t_6),\qquad r_k=h_k^{\leftarrow}$};
\end{tikzpicture}
\endgroup
        }
        \caption{Evolution of harmonic measure and Root paths in CHL$(0)$. }
        \label{fig:chl-hm-root-evolution}
    \end{figure}

    With particle $i$ we associate two forward harmonic measure paths, started from the two sides of its base, and one backward Root path, ending at its tip:
    \[
        \mathtt{HM}_{\bar t_i,\cdot}^{\mathcal A,N} (\bar x_{i,\varepsilon}),\ \varepsilon\in\{-,+\}, \qquad \mathtt R_{\cdot,\bar t_i}^{\mathcal A,N}(\bar x_i).
    \]
    At time $\bar t_i$ and at subsequent grid times, the length of the interval between the two side paths is the harmonic measure of particle $i$ together with all particles later attached to it, directly or indirectly; lengths are measured by Lebesgue measure on $\R$ and normalized arclength on $\T$.  Figure~\ref{fig:chl-hm-root-evolution} illustrates CHL$(0)$ slit attachments and the non-crossing duality between the two forward harmonic measure paths and a backward Root path.

\begin{figure}[H]
        \centering
        \resizebox{0.96\linewidth}{!}{%
\begingroup
\definecolor{FGfin}{HTML}{205F9A}
\definecolor{FGgap}{HTML}{24794F}
\tikzset{
  fg/tree/.style={draw=black,line width=.72pt,line cap=round,line join=round},
  fg/fin/.style={draw=FGfin,line width=1.12pt,line cap=round,line join=round},
  fg/gap/.style={draw=FGgap,line width=1.05pt,line cap=round,line join=round},
  fg/hood/.style={draw=black!15,line width=4.0pt,line cap=round,line join=round},
  fg/guide/.style={draw=black!45,line width=.4pt},
  fg/text/.style={font=\fontsize{9}{10.5}\selectfont,inner sep=1pt},
  fg/tiny/.style={font=\fontsize{8}{9.4}\selectfont,inner sep=1pt},
  fg/arrow/.style={-{Stealth[length=1.45mm,width=.95mm]}}
}
\def\FGedgeA{%
  (0.000000,0.000000) -- (0.000000,0.720000)%
}
\def\FGedgeB{%
  (0.000000,0.417600) -- (-0.016123,0.418069) -- (-0.035320,0.419854) -- (-0.054191,0.422912) -- (-0.072590,0.427148) -- (-0.090405,0.432445) -- (-0.110351,0.439794) -- (-0.126682,0.446950) -- (-0.144820,0.456132) -- (-0.161983,0.466055) -- (-0.178205,0.476582) -- (-0.195651,0.489203) -- (-0.212007,0.502326) -- (-0.227355,0.515844) -- (-0.243517,0.531425) -- (-0.258610,0.547315) -- (-0.274245,0.565257) -- (-0.288799,0.583435) -- (-0.302377,0.601804) -- (-0.315071,0.620328) -- (-0.328108,0.640854) -- (-0.340265,0.661513) -- (-0.352620,0.684182) -- (-0.364117,0.706975) -- (-0.375697,0.731792) -- (-0.386453,0.756731) -- (-0.396461,0.781785) -- (-0.406481,0.808888) -- (-0.415786,0.836116) -- (-0.425037,0.865421) -- (-0.433612,0.894862) -- (-0.442084,0.926412) -- (-0.449923,0.958111)%
}
\def\FGedgeC{%
  (0.000000,0.252000) -- (0.021778,0.253054) -- (0.043038,0.256101) -- (0.063383,0.260849) -- (0.082588,0.266929) -- (0.100580,0.273985) -- (0.122739,0.284413) -- (0.142989,0.295553) -- (0.165986,0.309994) -- (0.186816,0.324661) -- (0.209463,0.342277) -- (0.233255,0.362636) -- (0.254759,0.382678) -- (0.277067,0.405152) -- (0.299810,0.429895) -- (0.320511,0.454109) -- (0.341535,0.480462) -- (0.360817,0.506305) -- (0.380332,0.534226) -- (0.398333,0.561689) -- (0.416494,0.591205) -- (0.433325,0.620326) -- (0.450256,0.651495) -- (0.467172,0.684688) -- (0.482896,0.717554) -- (0.498582,0.752459) -- (0.513211,0.787093) -- (0.527784,0.823783) -- (0.541410,0.860248) -- (0.554962,0.898784) -- (0.568380,0.939385) -- (0.580938,0.979809) -- (0.593355,1.022316)%
}
\def\FGedgeD{%
  (-0.396758,0.782557) -- (-0.384600,0.787689) -- (-0.372975,0.793515) -- (-0.361929,0.799951) -- (-0.349478,0.808361) -- (-0.337924,0.817385) -- (-0.327257,0.826889) -- (-0.317445,0.836757) -- (-0.307018,0.848598) -- (-0.297606,0.860678) -- (-0.289115,0.872902) -- (-0.281453,0.885207) -- (-0.273604,0.899308) -- (-0.266613,0.913411) -- (-0.259649,0.929247) -- (-0.253522,0.945035) -- (-0.248130,0.960771) -- (-0.243388,0.976456) -- (-0.238793,0.993830) -- (-0.234827,1.011160) -- (-0.231108,1.030184) -- (-0.227986,1.049182) -- (-0.225182,1.069895) -- (-0.221174,1.111341) -- (-0.218786,1.156378) -- (-0.217958,1.206902) -- (-0.218715,1.263132) -- (-0.221052,1.327111)%
}
\def\FGedgeE{%
  (0.552947,0.892903) -- (0.541570,0.897292) -- (0.528489,0.903584) -- (0.516126,0.910897) -- (0.504558,0.919113) -- (0.493837,0.928110) -- (0.483989,0.937765) -- (0.475019,0.947963) -- (0.466911,0.958595) -- (0.459639,0.969566) -- (0.452161,0.982686) -- (0.445699,0.996047) -- (0.439459,1.011498) -- (0.434914,1.025103) -- (0.430664,1.040679) -- (0.427318,1.056233) -- (0.424512,1.073657) -- (0.422771,1.089054) -- (0.421556,1.106255) -- (0.421031,1.123319) -- (0.421107,1.140246) -- (0.421802,1.158898) -- (0.423051,1.177398) -- (0.424971,1.197590) -- (0.427612,1.219453) -- (0.434236,1.262775) -- (0.443559,1.312883) -- (0.455484,1.369935) -- (0.475711,1.461148)%
}
\def\FGedgeF{%
  (0.460398,1.392409) -- (0.452588,1.394414) -- (0.444915,1.396957) -- (0.437403,1.400029) -- (0.429176,1.404106) -- (0.422085,1.408263) -- (0.414383,1.413517) -- (0.407800,1.418678) -- (0.400707,1.425008) -- (0.394692,1.431075) -- (0.388260,1.438364) -- (0.382188,1.446112) -- (0.376476,1.454285) -- (0.366122,1.471781) -- (0.356678,1.491701) -- (0.348653,1.512887) -- (0.341592,1.536351) -- (0.335791,1.560845) -- (0.330674,1.588827) -- (0.326531,1.619117) -- (0.323307,1.651748) -- (0.320855,1.688187) -- (0.319232,1.727149)%
}
\def\FGStem{%
  (0.000000,0.000000) -- (0.000000,0.417600)%
}
\def\FGPartOne{%
  (0.000000,0.417600) -- (-0.014936,0.418003) -- (-0.029787,0.419202) -- (-0.044471,0.421174) -- (-0.058916,0.423881) -- (-0.073060,0.427273) -- (-0.086855,0.431295) -- (-0.100262,0.435888) -- (-0.116442,0.442342) -- (-0.128906,0.448005) -- (-0.143878,0.455623) -- (-0.158183,0.463753) -- (-0.171838,0.472314) -- (-0.184867,0.481239) -- (-0.199716,0.492345) -- (-0.213758,0.503808) -- (-0.227047,0.515561) -- (-0.239635,0.527551) -- (-0.253501,0.541786) -- (-0.264736,0.554159) -- (-0.277145,0.568761) -- (-0.288857,0.583511) -- (-0.299929,0.598385) -- (-0.310409,0.613365) -- (-0.321720,0.630596) -- (-0.332375,0.647929) -- (-0.342427,0.665352) -- (-0.351923,0.682856) -- (-0.361995,0.702635) -- (-0.371469,0.722499) -- (-0.380393,0.742444) -- (-0.396758,0.782557)%
}
\def\FGPartTwo{%
  (-0.396758,0.782557) -- (-0.384266,0.787844) -- (-0.372338,0.793861) -- (-0.361025,0.800519) -- (-0.350357,0.807725) -- (-0.340346,0.815389) -- (-0.328752,0.825484) -- (-0.318146,0.836013) -- (-0.308473,0.846855) -- (-0.299665,0.857911) -- (-0.290139,0.871354) -- (-0.281633,0.884903) -- (-0.274033,0.898494) -- (-0.267236,0.912088) -- (-0.260199,0.927920) -- (-0.254009,0.943706) -- (-0.248561,0.959439) -- (-0.243768,0.975121) -- (-0.239001,0.992987) -- (-0.234902,1.010805) -- (-0.231396,1.028588) -- (-0.228078,1.048565) -- (-0.225349,1.068529) -- (-0.223139,1.088494) -- (-0.221222,1.110696) -- (-0.218750,1.157445) -- (-0.217958,1.206705) -- (-0.218715,1.263161) -- (-0.221052,1.327111)%
}
\begin{tikzpicture}[x=1cm,y=1cm,every node/.style={fg/text}]
\path[use as bounding box] (0,-.88) rectangle (14.5,4.50);
\node[anchor=south] at (3.55,4.19)
  {(a) Finger $\mathtt{Fin}_i^{\mathcal C,N}$};
\node[anchor=south] at (10.95,4.19)
  {(b) Gap $\mathtt{Gap}_i^{\mathcal C,N}$};

\begin{scope}[shift={(3.55,0)},scale=2]
  \draw[fg/tree] (-1.04,0) -- (1.04,0);
  \foreach \ed in {\FGedgeA,\FGedgeB,\FGedgeC,\FGedgeD,\FGedgeE,\FGedgeF}
    \draw[fg/tree] \ed;
  \draw[fg/fin] \FGStem;
  \draw[fg/fin] \FGPartOne;
  \draw[fg/fin] \FGPartTwo;
  \draw[fg/fin,fg/arrow] (0.000000,0.244800) -- (0.000000,0.228343) -- (0.000000,0.211886) -- (0.000000,0.195429) -- (0.000000,0.178971) -- (0.000000,0.162514) -- (0.000000,0.146057) -- (0.000000,0.129600);
  \draw[fg/fin,fg/arrow] (-0.317753,0.624421) -- (-0.312395,0.616308) -- (-0.306881,0.608221) -- (-0.301205,0.600161) -- (-0.295358,0.592131) -- (-0.289335,0.584134) -- (-0.283126,0.576172) -- (-0.276725,0.568250) -- (-0.270121,0.560370) -- (-0.263305,0.552539) -- (-0.256269,0.544761) -- (-0.249000,0.537043);
  \draw[fg/fin,fg/arrow] (-0.218098,1.185951) -- (-0.218265,1.175816) -- (-0.218502,1.165696) -- (-0.218813,1.155589) -- (-0.219202,1.145493) -- (-0.219675,1.135410) -- (-0.220236,1.125336) -- (-0.220891,1.115271) -- (-0.221646,1.105213) -- (-0.222506,1.095162) -- (-0.223479,1.085116) -- (-0.224570,1.075073);
  \fill[FGfin] (0.000000,0.417600) circle[radius=1.12pt];
  \fill[FGfin] (-0.396758,0.782557) circle[radius=1.12pt];
  \fill[FGfin] (-0.221052,1.327111) circle[radius=1.12pt];
\end{scope}
\node[text=FGfin,anchor=west] at (3.66,.83) {$a$};
\node[text=FGfin,anchor=east] at (2.59,1.58) {$\dfrac{a+b}{2}$};
\node[text=FGfin,anchor=south east] at (3.01,2.82) {$b=\bar t_i$};
\node[text=FGfin,anchor=east] at (3.07,0.9)
  {$\mathfrak p_N^{\mathcal C}(e_1)$};
\node[text=FGfin,anchor=east] at (2.83,2.22)
  {$\mathfrak p_N^{\mathcal C}(e_2)$};
\node[fg/tiny,anchor=north] at (3.55,-.12) {initial boundary};
\node[fg/tiny] at (3.55,-.59)
  {$m=2$: each piece uses time $(b-a)/2$.};

\begin{scope}[shift={(10.95,0)},scale=2]
  \fill[black!15] (-0.305000,0.116978) -- (-0.317213,0.120684) -- (-0.320659,0.121858) -- (-0.331924,0.126119) -- (-0.335273,0.127513) --
    (-0.349201,0.133859) -- (-0.352334,0.135409) -- (-0.365106,0.142239) -- (-0.367909,0.143845) -- (-0.384942,0.154379) --
    (-0.400316,0.165065) -- (-0.414668,0.176044) -- (-0.433638,0.192054) -- (-0.448866,0.206175) -- (-0.464428,0.221817) --
    (-0.478579,0.237169) -- (-0.493024,0.254022) -- (-0.507513,0.272226) -- (-0.520699,0.290038) -- (-0.532824,0.307577) --
    (-0.545016,0.326459) -- (-0.556250,0.345098) -- (-0.566703,0.363638) -- (-0.577119,0.383409) -- (-0.586729,0.402946) --
    (-0.596309,0.423834) -- (-0.605801,0.446114) -- (-0.614549,0.468249) -- (-0.622612,0.490216) -- (-0.630557,0.513552) --
    (-0.638232,0.537930) -- (-0.645379,0.562513) -- (-0.651964,0.587026) -- (-0.658158,0.611966) -- (-0.659432,0.619047) --
    (-0.659771,0.626233) -- (-0.659168,0.633402) -- (-0.657636,0.640431) -- (-0.655198,0.647200) -- (-0.651898,0.653593) --
    (-0.647792,0.659501) -- (-0.642950,0.664822) -- (-0.637455,0.669465) -- (-0.631401,0.673352) -- (-0.624891,0.676415) --
    (-0.618037,0.678602) -- (-0.610956,0.679875) -- (-0.603770,0.680214) -- (-0.596601,0.679612) -- (-0.589572,0.678079) --
    (-0.586423,0.676945) -- (-0.584755,0.684601) -- (-0.583636,0.691068) -- (-0.581734,0.706328) -- (-0.580539,0.723307) --
    (-0.580038,0.743503) -- (-0.580322,0.765894) -- (-0.581418,0.791601) -- (-0.583338,0.820671) -- (-0.586297,0.855919) --
    (-0.586429,0.863113) -- (-0.585620,0.870261) -- (-0.583886,0.877244) -- (-0.581255,0.883939) -- (-0.577772,0.890235) --
    (-0.573498,0.896022) -- (-0.568505,0.901201) -- (-0.562878,0.905684) -- (-0.556714,0.909395) -- (-0.550119,0.912269) --
    (-0.543205,0.914258) -- (-0.536091,0.915327) -- (-0.528898,0.915459) -- (-0.521749,0.914651) -- (-0.514767,0.912916) --
    (-0.512515,0.912031) -- (-0.509014,0.921889) -- (-0.505644,0.933982) -- (-0.502586,0.948226) -- (-0.499969,0.964588) --
    (-0.497878,0.983096) -- (-0.496427,1.002704) -- (-0.495528,1.024324) -- (-0.495179,1.049659) -- (-0.494609,1.056831) --
    (-0.493109,1.063867) -- (-0.490702,1.070647) -- (-0.487432,1.077055) -- (-0.483353,1.082981) -- (-0.478535,1.088324) --
    (-0.473061,1.092993) -- (-0.467025,1.096907) -- (-0.460529,1.100000) -- (-0.453686,1.102218) -- (-0.446611,1.103524) --
    (-0.439426,1.103896) -- (-0.432254,1.103326) -- (-0.425218,1.101826) -- (-0.419470,1.099786) -- (-0.416451,1.106223) --
    (-0.412369,1.116849) -- (-0.408384,1.129576) -- (-0.404856,1.143518) -- (-0.401629,1.159492) -- (-0.398937,1.176492) --
    (-0.396543,1.196410) -- (-0.394660,1.218045) -- (-0.393569,1.225156) -- (-0.391559,1.232064) -- (-0.388665,1.238650) --
    (-0.384936,1.244803) -- (-0.380435,1.250415) -- (-0.375241,1.255393) -- (-0.369441,1.259650) -- (-0.363135,1.263113) --
    (-0.356431,1.265724) -- (-0.349444,1.267437) -- (-0.342293,1.268224) -- (-0.335100,1.268071) -- (-0.327989,1.266980) --
    (-0.321081,1.264970) -- (-0.319048,1.264076) -- (-0.316450,1.269323) -- (-0.312298,1.279594) -- (-0.308236,1.291959) --
    (-0.304643,1.305533) -- (-0.301365,1.321112) -- (-0.298644,1.337692) -- (-0.296241,1.357121) -- (-0.294371,1.378271) --
    (-0.293270,1.385380) -- (-0.291250,1.392286) -- (-0.288346,1.398868) -- (-0.284608,1.405015) -- (-0.280100,1.410621) --
    (-0.274898,1.415592) -- (-0.269093,1.419840) -- (-0.262782,1.423295) -- (-0.256074,1.425896) -- (-0.249085,1.427599) --
    (-0.241932,1.428376) -- (-0.234740,1.428212) -- (-0.227630,1.427111) -- (-0.220725,1.425092) -- (-0.218794,1.424239) --
    (-0.217053,1.427519) -- (-0.212857,1.437099) -- (-0.208676,1.448697) -- (-0.204902,1.461485) -- (-0.201371,1.476228) --
    (-0.198346,1.491992) -- (-0.195687,1.509595) -- (-0.193352,1.529856) -- (-0.192063,1.536934) -- (-0.189861,1.543783) --
    (-0.186784,1.550286) -- (-0.182885,1.556332) -- (-0.178229,1.561817) -- (-0.172898,1.566647) -- (-0.166981,1.570741) --
    (-0.160581,1.574027) -- (-0.153807,1.576449) -- (-0.146775,1.577967) -- (-0.139604,1.578554) -- (-0.132419,1.578199) --
    (-0.125341,1.576910) -- (-0.119529,1.575042) -- (-0.118289,1.577484) -- (-0.114349,1.586990) -- (-0.110471,1.598503) --
    (-0.107028,1.611179) -- (-0.105672,1.617450) -- (-0.107193,1.618622) -- (-0.113530,1.624109) -- (-0.117517,1.627927) --
    (-0.123533,1.634298) -- (-0.126579,1.637809) -- (-0.132302,1.645000) -- (-0.134669,1.648206) -- (-0.140445,1.656654) --
    (-0.142312,1.659577) -- (-0.149555,1.671904) -- (-0.156222,1.684706) -- (-0.161884,1.696747) -- (-0.168835,1.713038) --
    (-0.174942,1.728750) -- (-0.181072,1.745832) -- (-0.187261,1.764377) -- (-0.193268,1.783582) -- (-0.194961,1.790574) --
    (-0.195728,1.797727) -- (-0.195555,1.804919) -- (-0.194444,1.812027) -- (-0.192415,1.818930) -- (-0.189502,1.825508) --
    (-0.185756,1.831650) -- (-0.181240,1.837250) -- (-0.176032,1.842214) -- (-0.170220,1.846454) -- (-0.163905,1.849900) --
    (-0.157194,1.852492) -- (-0.150202,1.854186) -- (-0.143048,1.854953) -- (-0.135856,1.854780) -- (-0.128748,1.853669) --
    (-0.121846,1.851640) -- (-0.115267,1.848727) -- (-0.109125,1.844981) -- (-0.103525,1.840465) -- (-0.098562,1.835256) --
    (-0.094321,1.829445) -- (-0.090875,1.823129) -- (-0.088283,1.816418) -- (-0.082436,1.797725) -- (-0.076946,1.781303) --
    (-0.071682,1.766671) -- (-0.066660,1.753791) -- (-0.061910,1.742640) -- (-0.057197,1.732619) -- (-0.054712,1.727890) --
    (-0.051598,1.728636) -- (-0.044443,1.729386) -- (-0.037251,1.729194) -- (-0.030146,1.728066) -- (-0.023249,1.726020) --
    (-0.016677,1.723092) -- (-0.010545,1.719330) -- (-0.004955,1.714801) -- (-0.000005,1.709580) -- (0.004222,1.703758) --
    (0.007652,1.697434) -- (0.010228,1.690717) -- (0.011904,1.683721) -- (0.012654,1.676565) -- (0.012463,1.669374) --
    (0.010535,1.648482) -- (0.007724,1.626251) -- (0.003906,1.603742) -- (-0.000496,1.583738) -- (-0.004746,1.568090) --
    (-0.005701,1.564946) -- (-0.010701,1.550103) -- (-0.012015,1.546601) -- (-0.017483,1.533407) -- (-0.019256,1.529557) --
    (-0.025459,1.517346) -- (-0.027200,1.514180) -- (-0.031984,1.506246) -- (-0.037175,1.498524) -- (-0.042797,1.491049) --
    (-0.048874,1.483864) -- (-0.055426,1.477017) -- (-0.062468,1.470566) -- (-0.069823,1.464710) -- (-0.076928,1.459826) --
    (-0.080945,1.457372) -- (-0.084646,1.455302) -- (-0.089232,1.452965) -- (-0.094385,1.450656) -- (-0.098260,1.434476) --
    (-0.103634,1.416005) -- (-0.109999,1.398067) -- (-0.111359,1.394655) -- (-0.117122,1.381496) -- (-0.118924,1.377770) --
    (-0.125370,1.365627) -- (-0.127116,1.362578) -- (-0.131989,1.354777) -- (-0.135632,1.349436) -- (-0.141043,1.342185) --
    (-0.146845,1.335207) -- (-0.153056,1.328541) -- (-0.159691,1.322228) -- (-0.166761,1.316318) -- (-0.173594,1.311321) --
    (-0.176803,1.309208) -- (-0.181315,1.306457) -- (-0.185066,1.304367) -- (-0.189641,1.302047) -- (-0.193339,1.300358) --
    (-0.197291,1.281499) -- (-0.202310,1.262253) -- (-0.208366,1.243517) -- (-0.209628,1.240070) -- (-0.215248,1.226167) --
    (-0.216951,1.222374) -- (-0.223324,1.209503) -- (-0.228469,1.200377) -- (-0.233479,1.192473) -- (-0.240709,1.182415) --
    (-0.246820,1.174965) -- (-0.253403,1.167858) -- (-0.260478,1.161147) -- (-0.267871,1.155030) -- (-0.275008,1.149912) --
    (-0.279159,1.147259) -- (-0.282756,1.145147) -- (-0.287503,1.142598) -- (-0.293400,1.139832) -- (-0.297298,1.120429) --
    (-0.302121,1.101106) -- (-0.307915,1.082317) -- (-0.309060,1.079029) -- (-0.314476,1.064932) -- (-0.316023,1.061299) --
    (-0.322159,1.048219) -- (-0.327023,1.039096) -- (-0.331779,1.031141) -- (-0.338576,1.021078) -- (-0.344369,1.013527) --
    (-0.350614,1.006286) -- (-0.357331,0.999404) -- (-0.364364,0.993073) -- (-0.371160,0.987720) -- (-0.378546,0.982669) --
    (-0.383211,0.979850) -- (-0.388009,0.977254) -- (-0.388454,0.971801) -- (-0.391149,0.948470) -- (-0.394196,0.929418) --
    (-0.398444,0.909269) -- (-0.399235,0.906055) -- (-0.403563,0.890522) -- (-0.404718,0.886873) -- (-0.409899,0.872289) --
    (-0.414029,0.862308) -- (-0.419410,0.851086) -- (-0.424244,0.842291) -- (-0.429579,0.833716) -- (-0.433435,0.828136) --
    (-0.438948,0.820924) -- (-0.445424,0.813456) -- (-0.447949,0.810793) -- (-0.451932,0.806862) -- (-0.454813,0.804214) --
    (-0.459482,0.800220) -- (-0.462734,0.797644) -- (-0.467039,0.794490) -- (-0.471320,0.791646) -- (-0.470342,0.768113) --
    (-0.470046,0.741808) -- (-0.470603,0.719331) -- (-0.472109,0.697131) -- (-0.474641,0.676168) -- (-0.476920,0.662861) --
    (-0.478725,0.654314) -- (-0.481351,0.643797) -- (-0.484899,0.631999) -- (-0.488290,0.622493) -- (-0.492207,0.613036) --
    (-0.496706,0.603663) -- (-0.501842,0.594419) -- (-0.507527,0.585568) -- (-0.511058,0.580667) -- (-0.513379,0.577669) --
    (-0.517741,0.572420) -- (-0.520437,0.569407) -- (-0.524539,0.565150) -- (-0.529919,0.560151) -- (-0.526205,0.548351) --
    (-0.519099,0.527438) -- (-0.511969,0.507971) -- (-0.504283,0.488477) -- (-0.495972,0.468924) -- (-0.487649,0.450734) --
    (-0.479390,0.433903) -- (-0.470483,0.416952) -- (-0.460837,0.399843) -- (-0.451311,0.384080) -- (-0.440957,0.368085) --
    (-0.430803,0.353437) -- (-0.419765,0.338567) -- (-0.407686,0.323436) -- (-0.395828,0.309648) -- (-0.384390,0.297282) --
    (-0.371814,0.284688) -- (-0.359911,0.273697) -- (-0.346778,0.262605) -- (-0.334556,0.253255) -- (-0.323700,0.245766) --
    (-0.311807,0.238477) -- (-0.305000,0.234837) -- (-0.305000,0.388300) -- (-0.311113,0.392102) -- (-0.316964,0.396294) --
    (-0.322528,0.400861) -- (-0.327782,0.405782) -- (-0.332702,0.411035) -- (-0.337269,0.416600) -- (-0.341462,0.422451) --
    (-0.343412,0.425476) -- (-0.347011,0.431710) -- (-0.350195,0.438166) -- (-0.352950,0.444816) -- (-0.355263,0.451632) --
    (-0.357126,0.458585) -- (-0.357886,0.462103) -- (-0.359059,0.469205) -- (-0.359764,0.476369) -- (-0.359941,0.479964) --
    (-0.359941,0.487162) -- (-0.359764,0.490758) -- (-0.359059,0.497921) -- (-0.357886,0.505023) -- (-0.357126,0.508542) --
    (-0.355263,0.515495) -- (-0.352950,0.522311) -- (-0.350195,0.528961) -- (-0.348656,0.532215) -- (-0.345263,0.538563) --
    (-0.341462,0.544676) -- (-0.337269,0.550527) -- (-0.335031,0.553347) -- (-0.332702,0.556091) -- (-0.327782,0.561345) --
    (-0.322528,0.566266) -- (-0.316964,0.570832) -- (-0.314073,0.572976) -- (-0.308087,0.576976) -- (-0.305000,0.578826) --
    (-0.298652,0.582219) -- (-0.292095,0.585190) -- (-0.288748,0.586513) -- (-0.281931,0.588827) -- (-0.274978,0.590690) --
    (-0.267919,0.592094) -- (-0.264358,0.592622) -- (-0.257194,0.593328) -- (-0.253599,0.593504) -- (-0.246401,0.593504) --
    (-0.239218,0.593034) -- (-0.235642,0.592622) -- (-0.228540,0.591450) -- (-0.221530,0.589815) -- (-0.214642,0.587726) --
    (-0.207905,0.585190) -- (-0.201348,0.582219) -- (-0.195000,0.578826) -- (-0.191913,0.576976) -- (-0.188887,0.575025) --
    (-0.183036,0.570832) -- (-0.177472,0.566266) -- (-0.172218,0.561345) -- (-0.167298,0.556091) -- (-0.164969,0.553347) --
    (-0.160587,0.547636) -- (-0.158538,0.544676) -- (-0.156588,0.541651) -- (-0.152989,0.535417) -- (-0.149805,0.528961) --
    (-0.147050,0.522311) -- (-0.144737,0.515495) -- (-0.143748,0.512033) -- (-0.142114,0.505023) -- (-0.141469,0.501482) --
    (-0.140530,0.494345) -- (-0.140236,0.490758) -- (-0.140059,0.487162) -- (-0.140059,0.479964) -- (-0.140530,0.472781) --
    (-0.140941,0.469205) -- (-0.142114,0.462103) -- (-0.142874,0.458585) -- (-0.143748,0.455093) -- (-0.145838,0.448205) --
    (-0.147050,0.444816) -- (-0.149805,0.438166) -- (-0.152989,0.431710) -- (-0.156588,0.425476) -- (-0.158538,0.422451) --
    (-0.162731,0.416600) -- (-0.164969,0.413780) -- (-0.169715,0.408368) -- (-0.174805,0.403278) -- (-0.177472,0.400861) --
    (-0.180217,0.398532) -- (-0.185927,0.394150) -- (-0.191913,0.390151) -- (-0.195000,0.388300) -- (-0.195000,0.317642) --
    (-0.192587,0.318558) -- (-0.182920,0.322896) -- (-0.171506,0.328845) -- (-0.160765,0.335252) -- (-0.150632,0.342025) --
    (-0.139267,0.350468) -- (-0.128707,0.359132) -- (-0.117238,0.369462) -- (-0.106554,0.379978) -- (-0.096556,0.390638) --
    (-0.085916,0.402898) -- (-0.076010,0.415207) -- (-0.065626,0.429097) -- (-0.055922,0.443050) -- (-0.046815,0.457061) --
    (-0.037367,0.472603) -- (-0.027681,0.489674) -- (-0.018605,0.506796) -- (-0.010070,0.523971) -- (-0.001336,0.542704) --
    (0.006876,0.561464) -- (0.015231,0.581783) -- (0.023096,0.602131) -- (0.031059,0.624038) -- (0.038563,0.645976) --
    (0.046124,0.669459) -- (0.053832,0.694922) -- (0.056362,0.701657) -- (0.059750,0.708004) -- (0.063937,0.713854) --
    (0.068852,0.719108) -- (0.074410,0.723675) -- (0.080517,0.727478) -- (0.087069,0.730451) -- (0.093952,0.732543) --
    (0.101050,0.733720) -- (0.108240,0.733959) -- (0.115400,0.733258) -- (0.122407,0.731629) -- (0.129142,0.729099) --
    (0.135489,0.725711) -- (0.141339,0.721524) -- (0.146593,0.716609) -- (0.151161,0.711050) -- (0.154770,0.705254) --
    (0.159251,0.708820) -- (0.164678,0.713580) -- (0.170182,0.718855) -- (0.175733,0.724627) -- (0.181304,0.730872) --
    (0.186840,0.737530) -- (0.193055,0.745540) -- (0.198502,0.753031) -- (0.204573,0.761904) -- (0.210539,0.771168) --
    (0.216415,0.780828) -- (0.222179,0.790830) -- (0.228306,0.802051) -- (0.239967,0.825109) -- (0.251531,0.850256) --
    (0.263371,0.878457) -- (0.266584,0.884894) -- (0.270610,0.890856) -- (0.275379,0.896243) -- (0.280811,0.900960) --
    (0.286812,0.904928) -- (0.293279,0.908079) -- (0.300103,0.910359) -- (0.307165,0.911729) -- (0.314347,0.912165) --
    (0.321523,0.911660) -- (0.328572,0.910223) -- (0.335374,0.907878) -- (0.341811,0.904665) -- (0.347774,0.900639) --
    (0.353160,0.895870) -- (0.357877,0.890438) -- (0.361481,0.884988) -- (0.367393,0.889129) -- (0.372601,0.893181) --
    (0.377927,0.897684) -- (0.383345,0.902628) -- (0.388827,0.908000) -- (0.394353,0.913784) -- (0.399874,0.919932) --
    (0.405852,0.927010) -- (0.417749,0.942376) -- (0.429214,0.958843) -- (0.425812,0.964327) -- (0.423186,0.969039) --
    (0.421654,0.971997) -- (0.417944,0.979948) -- (0.414491,0.988617) -- (0.411779,0.996679) -- (0.409950,1.002992) --
    (0.409015,1.006694) -- (0.406304,1.019253) -- (0.405549,1.023587) -- (0.403664,1.037719) -- (0.403296,1.041448) --
    (0.402220,1.059526) -- (0.402151,1.078705) -- (0.402866,1.095971) -- (0.404545,1.117087) -- (0.407117,1.139553) --
    (0.410545,1.163440) -- (0.414983,1.189907) -- (0.420455,1.219116) -- (0.422242,1.226085) -- (0.424923,1.232761) --
    (0.428453,1.239030) -- (0.432771,1.244785) -- (0.437802,1.249926) -- (0.443462,1.254368) -- (0.449654,1.258032) --
    (0.456270,1.260857) -- (0.463199,1.262794) -- (0.470321,1.263810) -- (0.477515,1.263888) -- (0.484657,1.263026) --
    (0.491626,1.261239) -- (0.498302,1.258558) -- (0.504571,1.255029) -- (0.507675,1.252700) -- (0.515057,1.260245) --
    (0.523387,1.269771) -- (0.532006,1.280629) -- (0.541344,1.293465) -- (0.547306,1.302338) -- (0.543992,1.308445) --
    (0.541788,1.313113) -- (0.537959,1.322316) -- (0.535985,1.327890) -- (0.533025,1.337933) -- (0.531799,1.342944) --
    (0.529559,1.354417) -- (0.528884,1.358816) -- (0.527506,1.371084) -- (0.527216,1.374786) -- (0.526590,1.390964) --
    (0.526885,1.408128) -- (0.528023,1.425576) -- (0.530000,1.444204) -- (0.532562,1.462454) -- (0.536042,1.482999) --
    (0.540494,1.505817) -- (0.545755,1.530108) -- (0.547738,1.537024) -- (0.550606,1.543621) -- (0.554311,1.549788) --
    (0.558789,1.555419) -- (0.563964,1.560417) -- (0.569747,1.564697) -- (0.576039,1.568185) -- (0.582733,1.570822) --
    (0.589713,1.572563) -- (0.596862,1.573378) -- (0.604055,1.573253) -- (0.611170,1.572190) -- (0.618086,1.570207) --
    (0.624684,1.567339) -- (0.630356,1.563931) -- (0.634218,1.567534) -- (0.642013,1.575625) -- (0.650696,1.585523) --
    (0.659697,1.596667) -- (0.663231,1.601373) -- (0.660564,1.607236) -- (0.657401,1.615442) -- (0.655655,1.620761) --
    (0.653224,1.629682) -- (0.652147,1.634465) -- (0.650327,1.644640) -- (0.649742,1.648845) -- (0.648653,1.659728) --
    (0.648414,1.663275) -- (0.647986,1.677878) -- (0.648249,1.690763) -- (0.649479,1.709131) -- (0.651342,1.725995) --
    (0.653958,1.744014) -- (0.657179,1.762460) -- (0.660888,1.781138) -- (0.665774,1.803406) -- (0.667776,1.810316) --
    (0.670663,1.816906) -- (0.674385,1.823063) -- (0.678878,1.828681) -- (0.684067,1.833665) -- (0.689862,1.837928) --
    (0.696163,1.841399) -- (0.702864,1.844018) -- (0.709850,1.845739) -- (0.717000,1.846535) -- (0.724193,1.846390) --
    (0.731305,1.845307) -- (0.738215,1.843305) -- (0.744805,1.840419) -- (0.750962,1.836697) -- (0.756580,1.832203) --
    (0.761564,1.827014) -- (0.765828,1.821220) -- (0.769299,1.814918) -- (0.771917,1.808217) -- (0.773639,1.801232) --
    (0.774434,1.794081) -- (0.774289,1.786889) -- (0.773206,1.779776) -- (0.768671,1.759158) -- (0.765200,1.741681) --
    (0.762472,1.725977) -- (0.760358,1.711278) -- (0.758957,1.698367) -- (0.758199,1.687214) -- (0.758025,1.678682) --
    (0.761739,1.677804) -- (0.768464,1.675249) -- (0.774799,1.671837) -- (0.780633,1.667628) -- (0.785869,1.662694) --
    (0.790415,1.657118) -- (0.794195,1.650997) -- (0.797144,1.644435) -- (0.799210,1.637543) -- (0.800360,1.630442) --
    (0.800573,1.623250) -- (0.799845,1.616093) -- (0.798189,1.609092) -- (0.795634,1.602366) -- (0.792222,1.596032) --
    (0.781333,1.578644) -- (0.770204,1.561794) -- (0.758444,1.545022) -- (0.747071,1.529834) -- (0.735567,1.515535) --
    (0.723856,1.502118) -- (0.712421,1.490161) -- (0.700670,1.479093) -- (0.698032,1.476785) -- (0.688430,1.468911) --
    (0.685178,1.466441) -- (0.676126,1.460079) -- (0.672126,1.457518) -- (0.663115,1.452285) -- (0.658246,1.449772) --
    (0.649777,1.445924) -- (0.643964,1.443670) -- (0.640863,1.442666) -- (0.639068,1.429880) -- (0.637548,1.415302) --
    (0.636755,1.402612) -- (0.636601,1.391082) -- (0.636743,1.387872) -- (0.642596,1.386843) -- (0.649458,1.384682) --
    (0.655980,1.381644) -- (0.662049,1.377781) -- (0.667562,1.373159) -- (0.672424,1.367856) -- (0.676553,1.361965) --
    (0.679878,1.355585) -- (0.682341,1.348825) -- (0.683901,1.341802) -- (0.684531,1.334635) -- (0.684220,1.327448) --
    (0.682973,1.320362) -- (0.680813,1.313500) -- (0.677775,1.306978) -- (0.667037,1.287388) -- (0.662216,1.279030) --
    (0.654591,1.266068) -- (0.643616,1.248494) -- (0.630899,1.229583) -- (0.620278,1.214984) -- (0.608730,1.200357) --
    (0.596848,1.186675) -- (0.594759,1.184417) -- (0.584505,1.173936) -- (0.581869,1.171415) -- (0.571568,1.162197) --
    (0.568311,1.159500) -- (0.558512,1.152004) -- (0.554495,1.149205) -- (0.544662,1.142985) -- (0.539763,1.140226) --
    (0.530454,1.135593) -- (0.526087,1.133652) -- (0.521155,1.131712) -- (0.516889,1.130221) -- (0.516297,1.126123) --
    (0.514105,1.107206) -- (0.512726,1.090322) -- (0.512112,1.075486) -- (0.512207,1.061869) -- (0.512830,1.052202) --
    (0.519869,1.052010) -- (0.526974,1.050877) -- (0.533869,1.048826) -- (0.540439,1.045892) -- (0.546569,1.042127) --
    (0.552155,1.037593) -- (0.557101,1.032369) -- (0.561323,1.026544) -- (0.564749,1.020218) -- (0.567320,1.013498) --
    (0.568992,1.006501) -- (0.569736,0.999345) -- (0.569540,0.992154) -- (0.568407,0.985049) -- (0.566356,0.978153) --
    (0.563422,0.971584) -- (0.550762,0.947584) -- (0.545269,0.937797) -- (0.536564,0.922830) -- (0.522274,0.900091) --
    (0.505579,0.876135) -- (0.490645,0.856929) -- (0.483374,0.848319) -- (0.475629,0.839656) -- (0.467736,0.831356) --
    (0.459632,0.823375) -- (0.449745,0.814357) -- (0.442768,0.808458) -- (0.433994,0.801574) -- (0.422852,0.793670) --
    (0.415733,0.789099) -- (0.406237,0.783575) -- (0.396496,0.778576) -- (0.393626,0.777232) -- (0.386510,0.774146) --
    (0.376287,0.770332) -- (0.365840,0.767183) -- (0.362111,0.766243) -- (0.355190,0.764755) -- (0.351177,0.764046) --
    (0.344369,0.763102) -- (0.340094,0.762679) -- (0.331441,0.762244) -- (0.325266,0.750095) -- (0.317802,0.736456) --
    (0.310733,0.724220) -- (0.304151,0.713399) -- (0.296638,0.701700) -- (0.287931,0.688976) -- (0.280487,0.678784) --
    (0.273178,0.669365) -- (0.264029,0.658363) -- (0.257135,0.650635) -- (0.248686,0.641795) -- (0.238123,0.631682) --
    (0.230795,0.625255) -- (0.221320,0.617634) -- (0.211471,0.610506) -- (0.201237,0.603927) -- (0.191645,0.598494) --
    (0.188649,0.596947) -- (0.180837,0.593213) -- (0.177530,0.591764) -- (0.170555,0.588976) -- (0.166929,0.587672) --
    (0.159073,0.585154) -- (0.155133,0.584051) -- (0.148167,0.582378) -- (0.144185,0.581575) -- (0.137261,0.580444) --
    (0.132034,0.579836) -- (0.125899,0.562992) -- (0.117189,0.540491) -- (0.107893,0.517921) -- (0.098635,0.496810) --
    (0.088741,0.475630) -- (0.078922,0.455916) -- (0.068381,0.436077) -- (0.057897,0.417598) -- (0.046748,0.399212) --
    (0.035839,0.382382) -- (0.024143,0.365517) -- (0.011504,0.348557) -- (-0.000861,0.333136) -- (-0.014177,0.317736) --
    (-0.028608,0.302353) -- (-0.042717,0.288542) -- (-0.056033,0.276549) -- (-0.070519,0.264588) -- (-0.086130,0.252911) --
    (-0.103200,0.241506) -- (-0.119267,0.232027) -- (-0.133511,0.224604) -- (-0.149028,0.217537) -- (-0.152027,0.216296) --
    (-0.168341,0.210201) -- (-0.179470,0.206697) -- (-0.182790,0.205765) -- (-0.195000,0.202748) -- (-0.195000,0.055000) --
    (1.040000,0.055000) -- (1.047179,0.054529) -- (1.054235,0.053126) -- (1.061048,0.050813) -- (1.067500,0.047631) --
    (1.073482,0.043634) -- (1.078891,0.038891) -- (1.083634,0.033482) -- (1.087631,0.027500) -- (1.090813,0.021048) --
    (1.093126,0.014235) -- (1.094529,0.007179) -- (1.095000,0.000000) -- (1.094529,-0.007179) -- (1.093126,-0.014235) --
    (1.090813,-0.021048) -- (1.087631,-0.027500) -- (1.083634,-0.033482) -- (1.078891,-0.038891) -- (1.073482,-0.043634) --
    (1.067500,-0.047631) -- (1.061048,-0.050813) -- (1.054235,-0.053126) -- (1.047179,-0.054529) -- (1.040000,-0.055000) --
    (-1.040000,-0.055000) -- (-1.047179,-0.054529) -- (-1.054235,-0.053126) -- (-1.061048,-0.050813) -- (-1.067500,-0.047631) --
    (-1.073482,-0.043634) -- (-1.078891,-0.038891) -- (-1.083634,-0.033482) -- (-1.087631,-0.027500) -- (-1.090813,-0.021048) --
    (-1.093126,-0.014235) -- (-1.094529,-0.007179) -- (-1.095000,0.000000) -- (-1.094529,0.007179) -- (-1.093126,0.014235) --
    (-1.090813,0.021048) -- (-1.087631,0.027500) -- (-1.083634,0.033482) -- (-1.078891,0.038891) -- (-1.073482,0.043634) --
    (-1.067500,0.047631) -- (-1.061048,0.050813) -- (-1.054235,0.053126) -- (-1.047179,0.054529) -- (-1.040000,0.055000) --
    (-0.305000,0.055000) -- (-0.305000,0.116978) -- cycle;
  \draw[fg/tree] (-1.04,0) -- (1.04,0);
  \draw[fg/tree] (-0.250000,0.000000) -- (-0.250000,0.483563);
  \draw[fg/tree] (-0.250000,0.164412) -- (-0.263442,0.165032) -- (-0.276595,0.166834) -- (-0.289228,0.169659) -- (-0.301204,0.173303) --
    (-0.312469,0.177563) -- (-0.326397,0.183909) -- (-0.339169,0.190739) -- (-0.353709,0.199650) -- (-0.366900,0.208750) --
    (-0.381252,0.219729) -- (-0.396329,0.232464) -- (-0.409945,0.245036) -- (-0.424052,0.259164) -- (-0.436879,0.273031) --
    (-0.450041,0.288336) -- (-0.463351,0.305008) -- (-0.475497,0.321371) -- (-0.486654,0.337466) -- (-0.497944,0.354905) --
    (-0.508340,0.372108) -- (-0.518792,0.390649) -- (-0.528431,0.408993) -- (-0.537352,0.427174) -- (-0.546296,0.446719) --
    (-0.555184,0.467629) -- (-0.563383,0.488424) -- (-0.570967,0.509131) -- (-0.578481,0.531246) -- (-0.585425,0.553308) --
    (-0.592268,0.576805) -- (-0.598586,0.600283) -- (-0.604780,0.625223);
  \draw[fg/tree] (-0.250000,0.251453) -- (-0.237651,0.251856) -- (-0.222390,0.253466) -- (-0.210448,0.255585) -- (-0.195986,0.259159) --
    (-0.184857,0.262662) -- (-0.171548,0.267715) -- (-0.158930,0.273378) -- (-0.144686,0.280801) -- (-0.131376,0.288741) --
    (-0.118930,0.297060) -- (-0.105405,0.307108) -- (-0.092841,0.317416) -- (-0.079525,0.329410) -- (-0.067189,0.341552) --
    (-0.055714,0.353787) -- (-0.043709,0.367619) -- (-0.032547,0.381489) -- (-0.021011,0.396920) -- (-0.010276,0.412355) --
    (-0.000249,0.427783) -- (0.010060,0.444740) -- (0.020545,0.463219) -- (0.030327,0.481674) -- (0.039488,0.500107) --
    (0.048788,0.520053) -- (0.057508,0.539975) -- (0.066321,0.561408) -- (0.074597,0.582820) -- (0.082929,0.605742) --
    (0.090764,0.628646) -- (0.098625,0.653061) -- (0.106473,0.678988);
  \draw[fg/tree] (-0.594696,0.585611) -- (-0.589174,0.587373) -- (-0.583838,0.589601) -- (-0.578711,0.592275) -- (-0.573811,0.595372) --
    (-0.569155,0.598866) -- (-0.564146,0.603312) -- (-0.560044,0.607569) -- (-0.555682,0.612819) -- (-0.552150,0.617720) --
    (-0.548436,0.623642) -- (-0.545058,0.629874) -- (-0.542005,0.636385) -- (-0.539264,0.643150) -- (-0.536819,0.650145) --
    (-0.534656,0.657347) -- (-0.532539,0.665677) -- (-0.530734,0.674224) -- (-0.529219,0.682969) -- (-0.526973,0.700993) --
    (-0.525587,0.720694) -- (-0.525029,0.743171) -- (-0.525337,0.767414) -- (-0.526495,0.794585) -- (-0.528490,0.824785) --
    (-0.531490,0.860520);
  \draw[fg/tree] (0.093769,0.637804) -- (0.100842,0.635882) -- (0.107317,0.634744) -- (0.114778,0.634116) -- (0.121541,0.634144) --
    (0.128397,0.634725) -- (0.135321,0.635856) -- (0.142286,0.637529) -- (0.150142,0.640047) -- (0.157117,0.642835) --
    (0.164929,0.646569) -- (0.171822,0.650400) -- (0.179497,0.655258) -- (0.187071,0.660669) -- (0.194529,0.666604) --
    (0.201857,0.673031) -- (0.209046,0.679922) -- (0.216090,0.687245) -- (0.222984,0.694974) -- (0.229724,0.703081) --
    (0.237033,0.712500) -- (0.243449,0.721323) -- (0.250398,0.731480) -- (0.257162,0.741982) -- (0.263744,0.752803) --
    (0.270148,0.763917) -- (0.276994,0.776454) -- (0.289509,0.801200) -- (0.301886,0.828115) -- (0.314083,0.857166);
  \draw[fg/tree] (-0.528761,0.828317) -- (-0.523310,0.829009) -- (-0.518508,0.830073) -- (-0.513767,0.831560) -- (-0.509104,0.833463) --
    (-0.504539,0.835778) -- (-0.499540,0.838861) -- (-0.495235,0.842014) -- (-0.490565,0.846008) -- (-0.486582,0.849939) --
    (-0.482301,0.854769) -- (-0.478681,0.859406) -- (-0.474825,0.864986) -- (-0.471200,0.870929) -- (-0.467810,0.877210) --
    (-0.464653,0.883807) -- (-0.461726,0.890699) -- (-0.456545,0.905283) -- (-0.452218,0.920816) -- (-0.448506,0.938103) --
    (-0.445459,0.957154) -- (-0.443106,0.977977) -- (-0.441512,0.999531) -- (-0.440544,1.022802) -- (-0.440184,1.048901);
  \draw[fg/tree] (0.299689,0.823154) -- (0.305328,0.820936) -- (0.311211,0.819216) -- (0.317322,0.818008) -- (0.323641,0.817328) --
    (0.330147,0.817183) -- (0.336816,0.817581) -- (0.343623,0.818525) -- (0.350544,0.820014) -- (0.357554,0.822043) --
    (0.364628,0.824606) -- (0.371744,0.827692) -- (0.378880,0.831288) -- (0.386016,0.835380) -- (0.393135,0.839951) --
    (0.400220,0.844983) -- (0.407257,0.850458) -- (0.414235,0.856357) -- (0.421141,0.862661) -- (0.427968,0.869349) --
    (0.434707,0.876404) -- (0.441354,0.883806) -- (0.448625,0.892415) -- (0.462091,0.909807) -- (0.475755,0.929434) --
    (0.488890,0.950256) -- (0.502115,0.973246) -- (0.514776,0.997246);
  \draw[fg/tree] (-0.440637,1.019590) -- (-0.435795,1.019665) -- (-0.430930,1.020195) -- (-0.426062,1.021185) -- (-0.421214,1.022638) --
    (-0.416405,1.024551) -- (-0.411659,1.026921) -- (-0.406994,1.029740) -- (-0.402932,1.032614) -- (-0.398963,1.035827) --
    (-0.394623,1.039831) -- (-0.390427,1.044233) -- (-0.386388,1.049015) -- (-0.382516,1.054159) -- (-0.378816,1.059646) --
    (-0.375294,1.065458) -- (-0.371953,1.071575) -- (-0.365818,1.084656) -- (-0.360402,1.098752) -- (-0.355440,1.114600) --
    (-0.351209,1.131319) -- (-0.347487,1.149741) -- (-0.344453,1.168907) -- (-0.341828,1.190743) -- (-0.339867,1.213278);
  \draw[fg/tree] (0.498395,0.966569) -- (0.494544,0.968932) -- (0.490383,0.971958) -- (0.486956,0.974877) -- (0.483286,0.978505) --
    (0.479881,0.982429) -- (0.476739,0.986632) -- (0.473857,0.991097) -- (0.471231,0.995810) -- (0.468856,1.000756) --
    (0.466473,1.006584) -- (0.464606,1.011986) -- (0.462777,1.018299) -- (0.460066,1.030858) -- (0.458181,1.044990) --
    (0.457219,1.059906) -- (0.457105,1.076432) -- (0.457819,1.093698) -- (0.459363,1.112611) -- (0.461751,1.133222) --
    (0.464980,1.155577) -- (0.469221,1.180782) -- (0.474512,1.208970);
  \draw[fg/tree] (-0.342363,1.185715) -- (-0.337638,1.185426) -- (-0.332850,1.185600) -- (-0.328020,1.186245) -- (-0.323173,1.187369) --
    (-0.318334,1.188971) -- (-0.313529,1.191050) -- (-0.308782,1.193600) -- (-0.304631,1.196253) -- (-0.300560,1.199262) --
    (-0.296094,1.203057) -- (-0.291765,1.207270) -- (-0.287589,1.211883) -- (-0.283577,1.216874) -- (-0.279739,1.222223) --
    (-0.276083,1.227908) -- (-0.272612,1.233910) -- (-0.266239,1.246781) -- (-0.260619,1.260683) -- (-0.255480,1.276326) --
    (-0.251113,1.292825) -- (-0.247290,1.310990) -- (-0.244193,1.329860) -- (-0.241539,1.351322) -- (-0.239585,1.373426);
  \draw[fg/tree] (0.469031,1.179716) -- (0.473289,1.179128) -- (0.478123,1.178859) -- (0.482450,1.178965) -- (0.487345,1.179464) --
    (0.491711,1.180239) -- (0.496635,1.181476) -- (0.501018,1.182893) -- (0.505950,1.184833) -- (0.515259,1.189466) --
    (0.525093,1.195686) -- (0.534891,1.203183) -- (0.545192,1.212400) -- (0.555446,1.222881) -- (0.565651,1.234551) --
    (0.575803,1.247340) -- (0.586423,1.261938) -- (0.596967,1.277630) -- (0.607942,1.295203) -- (0.618807,1.313823) --
    (0.629544,1.333414);
  \draw[fg/tree] (-0.242078,1.346385) -- (-0.237577,1.346088) -- (-0.233013,1.346215) -- (-0.228404,1.346776) -- (-0.223770,1.347775) --
    (-0.219134,1.349216) -- (-0.214518,1.351099) -- (-0.209943,1.353419) -- (-0.205430,1.356170) -- (-0.201488,1.358969) --
    (-0.197146,1.362503) -- (-0.192921,1.366433) -- (-0.188825,1.370742) -- (-0.184871,1.375413) -- (-0.181069,1.380428) --
    (-0.177426,1.385768) -- (-0.173949,1.391415) -- (-0.167503,1.403559) -- (-0.161740,1.416718) -- (-0.156385,1.431572) --
    (-0.151747,1.447287) -- (-0.147593,1.464633) -- (-0.144126,1.482698) -- (-0.141159,1.502338) -- (-0.138714,1.523561);
  \draw[fg/tree] (0.615715,1.308401) -- (0.612135,1.310662) -- (0.608756,1.313204) -- (0.605581,1.316017) -- (0.602614,1.319087) --
    (0.599855,1.322404) -- (0.597305,1.325957) -- (0.594961,1.329734) -- (0.592569,1.334239) -- (0.588741,1.343442) --
    (0.585780,1.353484) -- (0.583540,1.364957) -- (0.582162,1.377225) -- (0.581585,1.390229) -- (0.581778,1.404697) --
    (0.582727,1.419872) -- (0.584466,1.436557) -- (0.587028,1.454807) -- (0.590259,1.473753) -- (0.594469,1.495246) --
    (0.599503,1.518441);
  \draw[fg/tree] (-0.141807,1.497604) -- (-0.137321,1.497195) -- (-0.132757,1.497231) -- (-0.128139,1.497721) -- (-0.123491,1.498673) --
    (-0.118837,1.500089) -- (-0.114205,1.501968) -- (-0.109620,1.504305) -- (-0.105603,1.506759) -- (-0.101659,1.509559) --
    (-0.097327,1.513109) -- (-0.093123,1.517066) -- (-0.089064,1.521410) -- (-0.085163,1.526123) -- (-0.081430,1.531182) --
    (-0.077872,1.536566) -- (-0.074495,1.542255) -- (-0.068292,1.554466) -- (-0.062823,1.567660) -- (-0.057824,1.582503) --
    (-0.053574,1.598151) -- (-0.049852,1.615366) -- (-0.046831,1.633233) -- (-0.044232,1.653536) -- (-0.042305,1.674427);
  \draw[fg/tree] (0.593883,1.492400) -- (0.597677,1.491775) -- (0.601992,1.491419) -- (0.605862,1.491405) -- (0.610247,1.491725) --
    (0.618601,1.493269) -- (0.627022,1.495996) -- (0.635492,1.499845) -- (0.644502,1.505078) -- (0.653554,1.511440) --
    (0.663156,1.519314) -- (0.672809,1.528317) -- (0.682512,1.538391) -- (0.692779,1.550093) -- (0.703093,1.562863) --
    (0.713444,1.576645) -- (0.724337,1.592146) -- (0.734719,1.607834) -- (0.745608,1.625223);
  \draw[fg/tree] (-0.044761,1.648867) -- (-0.051456,1.650148) -- (-0.058335,1.652536) -- (-0.064917,1.655911) -- (-0.071192,1.660202) --
    (-0.077529,1.665689) -- (-0.083545,1.672059) -- (-0.089268,1.679250) -- (-0.095044,1.687698) -- (-0.100867,1.697486) --
    (-0.106451,1.708113) -- (-0.112113,1.720153) -- (-0.117593,1.733019) -- (-0.123189,1.747368) -- (-0.128910,1.763272) --
    (-0.134768,1.780795) -- (-0.140775,1.800000);
  \draw[fg/tree] (0.731605,1.603036) -- (0.725713,1.607573) -- (0.723034,1.610169) -- (0.720233,1.613337) -- (0.715812,1.619581) --
    (0.711883,1.627019) -- (0.708719,1.635225) -- (0.706288,1.644146) -- (0.704469,1.654321) -- (0.703380,1.665203) --
    (0.702975,1.676756) -- (0.703238,1.689640) -- (0.704158,1.703198) -- (0.705782,1.718163) -- (0.708147,1.734602) --
    (0.711125,1.751746) -- (0.714835,1.770424) -- (0.719490,1.791591);
  \draw[fg/guide,dash pattern=on 2.5pt off 2.2pt] (-1.04,1.94) -- (1.04,1.94);
  \draw[fg/gap] (-0.250000,0.483563) -- (-0.250000,0.593563);
  \draw[fg/gap,fg/arrow] (-0.250000,0.593563) -- (-0.000496,1.583738) -- (0.003906,1.603742) -- (0.007724,1.626251) -- (0.010535,1.648482) --
    (0.012463,1.669374) -- (0.012654,1.676565) -- (0.012654,1.940000);
  \draw[fg/gap,fg/arrow] (-0.144840,1.010900) -- (-0.125287,1.088500);
  \fill[black] (-0.250000,0.483563) circle[radius=1.10pt];
  \draw[draw=FGgap,fill=white,line width=.65pt] (-0.250000,0.593563) circle[radius=.9pt];
\end{scope}
\node[fg/tiny,anchor=south east] at (13.03,3.94) {$\mathbb T\times\{L\}$};
\node[text=FGgap,anchor=west] at (11.26,3.22) {$\gamma_i^L$};
\node[text=FGgap,anchor=west] at (10.61,1.23) {$z_i$};
\node[text=FGgap,anchor=east] at (10.23,1.10) {$\ell_i$};
\draw[draw=FGgap,line width=.4pt] (10.27,1.09) -- (10.42,1.06);
\node[anchor=west] at (10.70,.35) {$(\xi_i,h_i)$};
\draw[fg/guide] (10.75,.51) -- (10.56,.64) -- (10.48,.93);
\node[fg/tiny,anchor=west] at (8.86,2.50)
  {$\mathscr T_N^{\mathcal C,\varepsilon_N}$};
\draw[fg/guide] (9.84,2.49) -- (10.19,2.58);
\node[fg/tiny,anchor=north] at (10.95,-.12) {initial boundary};
\node[fg/tiny] at (10.95,-.59)
  {$\gamma_i^{L_k}\to\gamma_i,\quad u=\theta_i(r),\quad \sigma_i=h_i/v_{\mathcal C}$.};
\end{tikzpicture}
\endgroup
        }
        \caption{Illustration of Finger and Gap constructions in rescaled CHL$(0)$.}
        \label{fig:chl-finger-gap}
    \end{figure}
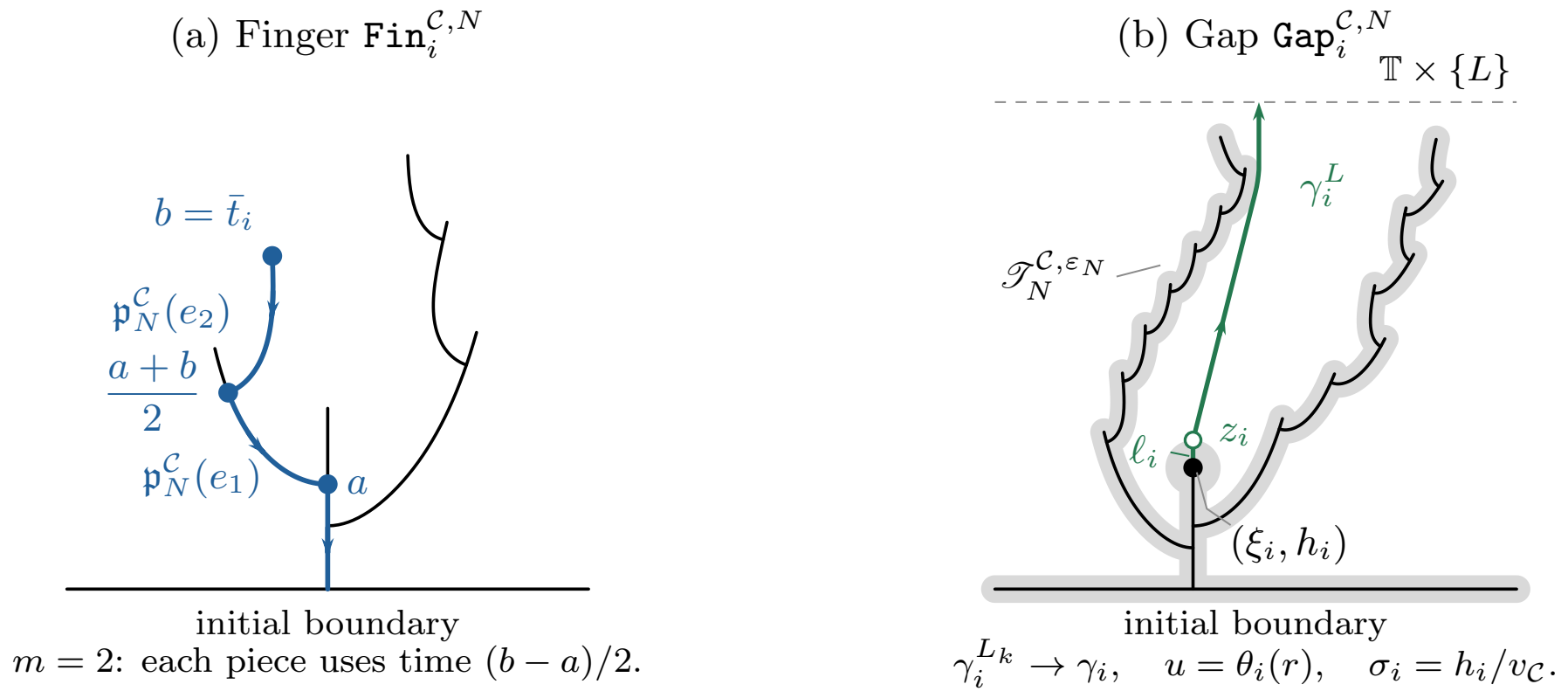

    \medskip
    \paragraph{\textbf{Fingers.}} For particle $i$, consider the unique path in the particle tree connecting its tip to the initial boundary. We regard it as a backward path, traced from the tip to the initial boundary as $u$ decreases from $\bar t_i$ to $0$. Define its rescaled parametrization $\mathtt{Fin}_i^{\mathcal A,N}:[0,\bar t_i]\to\mathbb G_{\mathcal A}$ as follows.
    \begin{itemize}
    \item     At each $u\in\Gamma_N^{\mathcal A}\cap[0,\bar t_i)$, choose the point of this path that is furthest from the initial boundary (along the path) among those already present at time $u$; for $u=\bar t_i$, choose the tip.
    \item For two of these times $a<b$ with no other such time between them, let $e_1,\ldots,e_m$ be the parts of particle edges between the points chosen at $a$ and $b$, listed in order from the initial boundary to the tip. For each particle $k$, parametrize its edge by
    \[
        \beta_k^{\mathcal A}(r):=
        \begin{cases}
        \mathcal H_{t_k-}^{N^{-1}\eta} \bigl(e^{\mathrm i x_k}(1+rN^{-1}\eta)\bigr), &\mathcal A=\mathcal H, \\
        \mathcal S_{t_k-}^{\lambda} \bigl(x_k+\mathrm i\lambda r\bigr), &\mathcal A=\mathcal S, \\
        \mathcal C_{t_k-}^{N,\lambda} \bigl(x_k+\mathrm i\lambda r\bigr), &\mathcal A=\mathcal C,
        \end{cases}
        \qquad 0\le r\le1.
    \]
    If $m\ge1$, write $e_j=\beta_{k_j}^{\mathcal A}([0,r_j]), 0\le r_j\le1,$    where $k_j$ is the index of the particle containing $e_j$. For $0\le\theta\le1,\ 1\le j\le m$, define
    \[
        \mathtt{Fin}_i^{\mathcal A,N} \left(a+\frac{j-1+\theta}{m}(b-a)\right) := \mathfrak p_N^{\mathcal A} \left( \beta_{k_j}^{\mathcal A} \bigl(\theta r_j\bigr) \right).
    \]
    If $m=0$, keep the path constant on $[a,b]$.
    \end{itemize}
    Figure~\ref{fig:chl-finger-gap}(a) illustrates the tree path and its interpolation along particle-edge pieces.   In SHL$(0)$, the number of particle-edge pieces between any two selected points is finite; see \cite[Sections~6 and 7.1]{berger2022growth}. The resulting path is continuous and contains the full rescaled edge of particle $i$.

    \medskip
    \paragraph{\textbf{Gaps.}} Gaps are paths that start at particle tips and escape to infinity in the complement of the complete particle tree. Let
    \begin{equation}\label{eq:gap-complete-tree}
        \mathscr T_N^{\mathcal A} :=(E^{\mathcal A}\times\{0\})\cup \overline{\bigcup_i \chi_{\mathcal A}\bigl( \mathtt{Fin}_i^{\mathcal A,N}([0,\bar t_i])\bigr)}, \qquad \mathscr D_N^{\mathcal A} :=(E^{\mathcal A}\times(0,\infty)) \setminus\mathscr T_N^{\mathcal A}.
    \end{equation}
    Set $\varepsilon_N=N^{-1}$. Distances and lengths are Euclidean in the $\chi_{\mathcal A}$-coordinates. By the construction in the proof of Lemma~\ref{lem:gap-exterior-geometry}, choose a closed neighbourhood $\mathscr T_N^{\mathcal A,\varepsilon_N}$ of $\mathscr T_N^{\mathcal A}$ by slightly thickening its edges and vertices, without creating new connections or enclosed regions. The neighbourhood is periodic in the compact cases, contains the initial boundary in its relative interior, and satisfies
    \[
        \mathscr T_N^{\mathcal A} \subset \operatorname{int} \mathscr T_N^{\mathcal A,\varepsilon_N} \subset \mathscr T_N^{\mathcal A,\varepsilon_N} \subset \{z:\dist(z,\mathscr T_N^{\mathcal A})\le\varepsilon_N\}.
    \]
    Moreover, every $z\in\operatorname{int}\mathscr T_N^{\mathcal A,\varepsilon_N}$ can be joined to $\mathscr T_N^{\mathcal A}$ by a curve contained in $\operatorname{int}\mathscr T_N^{\mathcal A,\varepsilon_N}$ of diameter at most $\varepsilon_N$.

    Let $(\xi_i,h_i):=\chi_{\mathcal A}\bigl(\mathtt{Fin}_i^{\mathcal A,N}(\bar t_i)\bigr)$ be the coordinates of the tip of particle $i$. Choose a segment $\ell_i$ of length at most $\varepsilon_N$ from $(\xi_i,h_i)$, stopped at its first exit from $\operatorname{int}\mathscr T_N^{\mathcal A,\varepsilon_N}$, with $\ell_i\cap\mathscr T_N^{\mathcal A}=\{(\xi_i,h_i)\}$. Let $z_i\in\partial\mathscr T_N^{\mathcal A,\varepsilon_N}$ be its other endpoint.

    For $L>h_i+\varepsilon_N$, let $\gamma_i^L$ be a shortest curve from $z_i$ to $E^{\mathcal A}\times\{L\}$ in $(E^{\mathcal A}\times[0,\infty))\setminus\operatorname{int}\mathscr T_N^{\mathcal A,\varepsilon_N}$, parametrized by arc length. By Lemma~\ref{lem:gap-exterior-geometry}, choose measurably a locally uniform subsequential limit $\gamma_i$ as $L\to\infty$ with $\sup_{r\ge0}\gamma_i^{(2)}(r)=\infty$. Let $g_i=(g_i^{(1)},g_i^{(2)})$ be the curve obtained by following $\ell_i$ and then $\gamma_i$, parametrized by arc length from $(\xi_i,h_i)$.  We now define the gap path starting at the tip of particle $i$ by reparametrizing $g_i$, so that its height at time $u$ is at most $v_{\mathcal A}u$:
    \begin{equation}\label{eq:gap-height-clock}
        \begin{split}
        \sigma_i:=h_i/v_{\mathcal A},\quad \theta_i(r):=v_{\mathcal A}^{-1}\max_{0\le a\le r}g_i^{(2)}(a)  +\varepsilon_N\frac{r}{1+r},\quad \mathtt{Gap}_i^{\mathcal A,N}(u):=\chi_{\mathcal A}^{-1}\bigl(g_i(\theta_i^{-1}(u))\bigr),
        \end{split}
    \end{equation}
    where $u\ge\sigma_i$.  Then $\theta_i$ is a continuous, strictly increasing bijection from $[0,\infty)$ onto $[\sigma_i,\infty)$. Moreover,
    \begin{equation}\label{eq:gap-strict-exterior}
        \begin{gathered}
        \chi_{\mathcal A}(\mathtt{Gap}_i^{\mathcal A,N}(\sigma_i)) =(\xi_i,h_i),\quad \chi_{\mathcal A}(\mathtt{Gap}_i^{\mathcal A,N}(u)) \in\mathscr D_N^{\mathcal A},\quad \chi_{\mathcal A}^{(2)} (\mathtt{Gap}_i^{\mathcal A,N}(u))\le v_{\mathcal A}u \quad(u>\sigma_i).
        \end{gathered}
    \end{equation}

    Figure~\ref{fig:chl-finger-gap}(b) illustrates the exterior shortest-curve construction from a particle tip.

    \medskip
    \paragraph{\textbf{Path collections and topologies.}}Equip $\mathbb  G_{\mathcal A}$ with the following metric:
    \[
        d^{\mathcal A}(z,z')= |\tanh\chi_{\mathcal A}^{(2)}(z)-\tanh \chi_{\mathcal A}^{(2)}(z')| \vee
        \begin{cases}
        d_{\mathbb T}(\chi_{\mathcal A}^{(1)}(z),\chi_{\mathcal A}^{(1)}(z')), &\mathcal A\in\mathfrak M_{\rm c}, \\
        |\tanh\rea z-\tanh\rea z'|,&\mathcal A=\mathcal S,
        \end{cases}
    \]
    where $\tanh$ induces a compactification of the unbounded coordinates. For a backward path  $\gamma:[0,t_\gamma]\to\mathbb G_{\mathcal A}$ and $T>0$,  restrict $\gamma$ to $[0,t_\gamma\wedge T]$ and extend it to $[0,T]$ by its terminal value. For a boundary path $\gamma:[0,t_\gamma]\to E^{\mathcal A}$, we first embed it in $\mathbb G_{\mathcal A}$ via $u\mapsto\mathfrak e_{\mathcal A}^{v_{\mathcal A}}(u,\gamma(u))$, and then proceed in the same way. For $\mathcal A\in\mathfrak M_{\rm c}$, define the distance between paths up to time $T$ by
    \[
        d_T^{\mathcal A}(\gamma,\gamma')= |\tanh(t_\gamma\wedge T)-\tanh(t_{\gamma'}\wedge T)| \vee\sup_{0\le u\le T}d^{\mathcal A}(\gamma(u),\gamma'(u)).
    \]
    For $\mathcal A=$ $\mathcal S$, we also localize in space: for each $m\ge1$, let
    \[
        \zeta_m(x):=1\wedge (m+1-|x|)_+, \quad    \kappa_m(x):=\int_0^x\zeta_m(q)\,dq,\quad \chi_{\mathcal S,m}(x+\mathrm ih):=(\kappa_m(x),\zeta_m(x)\tanh h),
    \]
    and define the distance between paths up to time $T$ by
    \[
        d_T^{\mathcal S,m}(\gamma,\gamma')= |\tanh(t_\gamma\wedge T)-\tanh(t_{\gamma'}\wedge T)| \vee\sup_{0\le u\le T} \|\chi_{\mathcal S,m}(\gamma(u))-\chi_{\mathcal S,m}(\gamma'(u))\|_\infty.
    \]
    For compact path collections, let $d_{H,T}^{\mathcal A}$ and $d_{H,T}^{\mathcal S,m}$ denote the Hausdorff distances induced by $d_T^{\mathcal A}$ and $d_T^{\mathcal S,m}$, respectively. Define the path metric
    \begin{equation*}
        D_{\mathcal A}(\gamma,\gamma') =
        \begin{cases}
        \displaystyle \sum_{k\ge1}2^{-k} \bigl(1\wedge d_k^{\mathcal A}(\gamma,\gamma')\bigr), & \mathcal A\in\mathfrak M_{\rm c}, \\[2mm]
        \displaystyle \sum_{k,m\ge1}2^{-(k+m)} \bigl(1\wedge d_k^{\mathcal S,m}(\gamma,\gamma')\bigr), & \mathcal A=\mathcal S.
        \end{cases}
    \end{equation*}
    The corresponding Hausdorff metric on compact path collections is denoted by $D_{\mathcal A,H}$.

    For forward paths $\gamma:[\sigma_\gamma,\infty)\to\mathbb G_{\mathcal A}$, use the same construction with $\gamma(u\vee\sigma_\gamma)$ and $\tanh(\sigma_\gamma\wedge T)$ in place of the extended backward path $\gamma(u)$ and $\tanh(t_\gamma\wedge T)$, respectively. Denote the resulting path and Hausdorff metrics by $D_{\mathcal A}^\uparrow$ and $D_{\mathcal A,H}^\uparrow$. Closures below are taken in the metric completions of the respective path spaces.

    For $\ell\ge m+1$, the map $(x,y)\mapsto(\kappa_m(x),\zeta_m(x)y)$, which expresses the coordinates $\chi_{\mathcal S,m}$ in terms of $\chi_{\mathcal S,\ell}$, is $2$-Lipschitz on the compactified coordinate range. Since restriction in time does not increase distance, for every $h\in\mathbb N$ and compact path collections $K,L$,
    \begin{equation}\label{eq:geometric-path-metric-comparison}
        D_{\mathcal A,H}(K,L)\le 2^{1-h}+2
        \begin{cases}
        d_{H,h}^{\mathcal A}(K,L), &\mathcal A\in\mathfrak M_{\rm c}, \\
        d_{H,h}^{\mathcal S,h+1}(K,L), &\mathcal A=\mathcal S.
        \end{cases}
    \end{equation}
    The same estimate holds for forward paths.

    We now define the full harmonic measure, Root, Finger, and Gap path collections in their respective path spaces.

    \begin{definition}\label{def:process-generated-observables}
    For each $\mathcal A\in\mathfrak M$, define the full path collections by
    \begin{equation*}
        \begin{aligned}
        &\mathbf{HM}^{\mathcal A,N} :=\overline{\bigl\{ \mathtt{HM}_{\bar t_i,\cdot}^{\mathcal A,N} (\bar x_{i,\varepsilon}): i\in\mathbb N^*,\ \varepsilon\in\{-,+\} \bigr\}},\quad &\mathbf R^{\mathcal A,N} :=\overline{\bigl\{ \mathtt R_{\cdot,\bar t_i}^{\mathcal A,N} (\bar x_i):i\in\mathbb N^* \bigr\}}, \\
        &\mathbf{Fin}^{\mathcal A,N} :=\overline{\bigl\{ \mathtt{Fin}_i^{\mathcal A,N}:i\in\mathbb N^* \bigr\}},\quad &\mathbf{Gap}^{\mathcal A,N} :=\overline{\bigl\{ \mathtt{Gap}_i^{\mathcal A,N}:i\in\mathbb N^* \bigr\}}.
        \end{aligned}
    \end{equation*}
    \end{definition}

    \medskip
    \paragraph{\textbf{Backbones.}}

    For $\mathcal A\in\{\mathcal H,\mathcal C\}$, define the backbone $\mathtt{B}^{\mathcal A,N,\downarrow}:[0,\infty)\to\mathbb G_{\mathcal A}$ by
    \[
        \mathtt{B}^{\mathcal A,N,\downarrow}|_{[0,T]} :=\lim_{i\to\infty}\mathtt{Fin}_{i}^{\mathcal A,N}|_{[0,T]}, \qquad T<\infty.
    \]
    The restrictions on the right are eventually identical, as shown in the proof of Proposition~\ref{prop:process-root-finger-saturation}. In words, up to any fixed time $T$, the fingers of all sufficiently late particles coincide, and the backbone is this common path. Let $U$ be uniform on $\mathbb T$, independent of the driving PPP.  Among particles with $0<\bar t_i\le\Delta_N^{\mathcal A}$, choose $i_N$ with $\bar x_{i_N}$ nearest to $U$ in $d_{\mathbb T}$, taking the smaller index in a tie; if there is no such particle,  set $i_N=1$. Define the selected dual backbone by
    \[
        \mathtt{B}^{\mathcal A,N,\uparrow}:=\mathtt{Gap}_{i_N}^{\mathcal A,N}.
    \]

    \medskip
    \paragraph{\textbf{One-arm domination time.}}

    For $\mathcal A\in \mathfrak M_{\rm c}$, using the continuous boundary extensions of the finite-time aggregate maps, define the arms at time $t$ as the connected components of
    \[
        \mathcal C_t^{N,\lambda}(\mathbb T_N)\setminus\mathbb T_N \quad\text{for CHL}(0),\qquad \mathcal H_t^{N^{-1}\eta}(\partial\mathbb D_0) \setminus\partial\mathbb D_0 \quad\text{for HL}(0).
    \]
    A particle is attached to the component containing its edge with the attachment point removed. An arm is infinite if it receives infinitely many   particles. The almost sure existence and uniqueness  of the infinite arm follow from \cite{rohde2005some,norris2012hastings,chen2025one}.
    \begin{definition}\label{def:one-arm-domination-time}
    Let $T_\infty$ denote the unique infinite arm. Let $t_j^{\mathcal A}$, $\mathcal A\in\{\mathcal H,\mathcal C\}$, be the particle arrival times of Definition~\ref{def:HL-geometries}, which occur at rates $1$ and $2\pi N$, respectively. Define the one-arm domination time by
    \[
        \upsilon_{N,\lambda}^{\mathcal C} := \sup\left\{ t_j^{\mathcal C}: \text{the particle arriving at }t_j^{\mathcal C} \text{ is not attached to }T_\infty \right\},
    \]
    \[
        \upsilon_{N,\eta}^{\mathcal H} := \sup\left\{ t_j^{\mathcal H}: \text{the particle arriving at }t_j^{\mathcal H} \text{ is not attached to }T_\infty \right\},
    \]
    with $\sup\varnothing=0$. Their rescaled versions are
    \[
        \tau_{\mathcal C}^{N} := N^{-2}\upsilon_{N,\lambda}^{\mathcal C}, \qquad \tau_{\mathcal H}^{N} := N^{-3}\upsilon_{N,\eta}^{\mathcal H}.
    \]
    \end{definition}

    \subsection{Preliminaries}\label{S.preliminaries}
    In this subsection, we relate the boundary transitions of Subsection~\ref{S.HL-models} to one-slit maps and time truncation along the corresponding tree paths. For $\mathcal A\in\mathfrak M_{\rm f}$, the growth rate $v_{\mathcal A}$ and variance $\sigma_{\mathcal A}^2$ are defined in \eqref{growth-rates-variance}; since they depend only on $\lambda$, we also write $\sigma_\lambda^2:=\sigma_{\mathcal S}^2=\sigma_{\mathcal C}^2=4\lambda^3/3$. For these models, write
    \[
        (A_\xi,\mathcal A_t)=
        \begin{cases}
        (S_\xi^\lambda,\mathcal S_t^\lambda), &\mathcal A=\mathcal S,\quad \xi\in\mathbb R, \\
        (C_\xi^{N,\lambda},\chl_t^{N,\lambda}), &\mathcal A=\mathcal C,\quad \xi\in\mathbb T_N.
        \end{cases}
    \]

    For $a\in\mathbb R$, define the boundary maps of a single slit and their generalized inverses (see \eqref{generalized-inverse}) by
    \begin{equation}\label{eq:one-slit-boundary-maps}
        \begin{aligned}
        r_x^{\mathcal S}(a)&:=\rea S_x^\lambda(a), & h_x^{\mathcal S}&:=(r_x^{\mathcal S})^\leftarrow, &&x\in\mathbb R, \\
        r_x^{\mathcal C,N}(a)&:=\rea C_x^{N,\lambda}(a), & h_x^{\mathcal C,N}&:=(r_x^{\mathcal C,N})^\leftarrow, &&x\in\mathbb T_N.
        \end{aligned}
    \end{equation}
    The maps $h_x$ are the real boundary values of the inverse slit maps, with the right value taken at each jump.  For example,
    \[
        r_x^{\mathcal S}(a)=x\ \Longleftrightarrow\ |a-x|\le\lambda, \qquad h_x^{\mathcal S}(x-)=x-\lambda,\quad h_x^{\mathcal S}(x)=x+\lambda.
    \]
    Thus $r_x^{\mathcal S}$ collapses the interval $[x-\lambda,x+\lambda]$ to the attachment coordinate $x$, while $S_x^\lambda$ maps only its endpoints to the base point $x$ of the slit. The following lemma gives the moment estimates and uniform bounds for the one-slit boundary maps.
    \begin{lemma}\label{lem:one-slit-moments}
    For every $a\in\mathbb R$,
    \begin{align}
        \int_{\mathbb R}(r_x^{\mathcal S}(a)-a)\,dx =\int_{\mathbb R}(h_x^{\mathcal S}(a)-a)\,dx &=\int_{\mathbb T_N}(r_x^{\mathcal C,N}(a)-a)\,dx =\int_{\mathbb T_N}(h_x^{\mathcal C,N}(a)-a)\,dx=0,\notag \\
        \int_{\mathbb R}(r_x^{\mathcal S}(a)-a)^2\,dx &=\int_{\mathbb R}(h_x^{\mathcal S}(a)-a)^2\,dx =\sigma_\lambda^2,\notag \\
        \int_{\mathbb T_N}(r_x^{\mathcal C,N}(a)-a)^2\,dx &=\int_{\mathbb T_N}(h_x^{\mathcal C,N}(a)-a)^2\,dx =\sigma_\lambda^2+O_\lambda(N^{-1}). \label{eq:CHL-moments}
    \end{align}
    For all sufficiently large $N$, there exists a constant $c_\lambda>0$ such that
    \begin{equation}\label{eq:bounded-one-slit-jumps}
        \sup_{a,x}|r_x^{\mathcal C,N}(a)-a| +\sup_{a,x}|h_x^{\mathcal C,N}(a)-a|\le c_\lambda.
    \end{equation}
    For every compact $\Lambda\Subset(0,\infty)$, $q\ge2$, $h>0$, and $N\ge1$, there exist constants $C_{q,\Lambda}, C_{q,h,\Lambda}>0$ such that
    \begin{align}
        \sup_{\lambda\in\Lambda}\sup_{w\in\overline{\mathbb H}} \int|A_\xi(w)-w|^q\,d\xi &\le C_{q,\Lambda},\label{eq:complex-slit-global-moment} \\
        \sup_{\lambda\in\Lambda}\sup_{\im w\ge hN} \left\{N\int|A_\xi(w)-w|^2\,d\xi+ N^{q-1}\int|A_\xi(w)-w|^q\,d\xi\right\} &\le C_{q,h,\Lambda}.\label{eq:complex-slit-high-moment}
    \end{align}
    For $w\in\overline{\mathbb H}$,
    \[
        \int_{\mathbb R}(S_\xi^\lambda(w)-w)\,d\xi =\mathrm iv_{\mathcal S},\qquad \int_{\mathbb T_N}(C_\xi^{N,\lambda}(w)-w)\,d\xi =\mathrm i\{v_{\mathcal C}+O_\Lambda(N^{-2})\},\qquad \im A_\xi(w)\ge\im w.
    \]
    \end{lemma}
    The proof is given in Appendix~\ref{appendix:one-slit-moments}.

    We now express the maps $U^{\mathcal A,N},R^{\mathcal A,N}$ in terms of the Poisson points $(t_i,x_i)$ of Definition~\ref{def:HL-geometries}. In CHL(0), let $s<\bar t_{i_1}<\cdots<\bar t_{i_k}\le t$ be the rescaled arrival times in $(s,t]$. Then
    \begin{equation}\label{eq:CHL-root-flow-def}
        \begin{aligned}
        U_{s,t}^{\mathcal C,N}(x) =N^{-1}(h_{x_{i_k}}^{\mathcal C,N}\circ\cdots\circ h_{x_{i_1}}^{\mathcal C,N})(Nx),\quad   R_{s,t}^{\mathcal C,N}(x) =N^{-1}(r_{x_{i_1}}^{\mathcal C,N}\circ\cdots\circ r_{x_{i_k}}^{\mathcal C,N})(Nx).
        \end{aligned}
    \end{equation}
    Thus $U$ composes the inverse slit maps forward in time, while $R$ composes the slit boundary maps backward in time. In SHL(0), for $M>0$, define
    \begin{equation}\label{eq:SHL-root-cutoff}
        R_{s,t}^{\mathcal S,N;M}(x):=N^{-1}(r_{x_{i_1}}^{\mathcal S}\circ\cdots\circ r_{x_{i_k}}^{\mathcal S})(Nx),
    \end{equation}
    where $\{i_1,\ldots,i_k\}=\{i:s<\bar t_i\le t,\ |x_i|\le M\}, \bar t_{i_1}<\cdots<\bar t_{i_k}.$    Here $M$ restricts the spatial coordinates $x_i$ to $[-M,M]$, and the limit $M\to\infty$ is taken with $N$ fixed.

    For SHL$(0)$, the map in \eqref{eq:rescaled-boundary-transition} is defined  for $x\in N^{-1}(\mathcal S_{N^2s}^\lambda)^{-1}(\mathcal S_{N^2t}^\lambda(\mathbb H))$, with the right-continuous boundary values taken in     Subsection~\ref{S.HL-models}; see \cite[Theorem~4.3 and Proposition~4.5]{berger2022growth}.  The next proposition shows that these boundary values equal $(R_{s,t}^{\mathcal S,N})^\leftarrow$, including when $s$ is an arrival time.
    \begin{proposition}\label{prop:SHL-boundary-cocycle}
    For each fixed $N\in\mathbb N$, the map $ (\omega,s,t,a)\longmapsto R_{s,t}^{\mathcal S,N}(a)(\omega), 0\le s\le t,\ a\in\mathbb R $, is jointly measurable. Moreover, almost surely, for every $0<T,L<\infty$,
    \begin{equation}\label{eq:SHL-root-cutoff-limit}
        \lim_{M\to\infty}\sup_{0\le s\le t\le T}\sup_{|a|\le L}   \bigl|R_{s,t}^{\mathcal S,N;M}(a)-R_{s,t}^{\mathcal S,N}(a)\bigr|=0.
    \end{equation}
    In addition, $R_{s,t}^{\mathcal S,N}$ is non-decreasing and
    \begin{equation}\label{eq:SHL-root-uniform-continuity}
        \lim_{\delta\downarrow0}\sup_{0\le s\le t\le T}    \sup_{\substack{a,b\in[-L,L],\ |a-b|\le\delta}}    \bigl|R_{s,t}^{\mathcal S,N}(a)-R_{s,t}^{\mathcal S,N}(b)\bigr|=0,
    \end{equation}
    while
    \begin{equation}\label{eq:SHL-root-uniform-endpoints}
        \lim_{a\to+\infty}\inf_{0\le s\le t\le T} R_{s,t}^{\mathcal S,N}(a)=+\infty, \qquad \lim_{a\to-\infty}\sup_{0\le s\le t\le T} R_{s,t}^{\mathcal S,N}(a)=-\infty.
    \end{equation}
    If $s$ is a nonnegative rational number or an arrival time $s=\bar t_i$, then for $t\ge s$, $U_{s,t}^{\mathcal S,N}=(R_{s,t}^{\mathcal S,N})^\leftarrow$. Consequently, for $a,y\in\mathbb R$,
    \begin{equation}\label{eq:SHL-grid-inverse}
        R_{s,t}^{\mathcal S,N}(a)=y \quad\Longleftrightarrow\quad U_{s,t}^{\mathcal S,N}(y-)\le a\le U_{s,t}^{\mathcal S,N}(y).
    \end{equation}
    For $s\le t\le u$,
    \begin{equation}\label{eq:SHL-grid-cocycle}
        R_{s,u}^{\mathcal S,N}=R_{s,t}^{\mathcal S,N}\circ R_{t,u}^{\mathcal S,N}, \qquad U_{s,u}^{\mathcal S,N}=U_{t,u}^{\mathcal S,N}\circ U_{s,t}^{\mathcal S,N},
    \end{equation}
    where the second identity holds when $s,t$ are nonnegative rational times or arrival times.  For fixed $0\le s\le t$, $(R_{s,t}^{\mathcal S,N},U_{s,t}^{\mathcal S,N})$ is measurable with respect to $\sigma\!\left( \mathcal P^{\mathcal S}\cap((N^2s,N^2t]\times\mathbb R) \right).$
    \end{proposition}
    We prove this proposition in Appendix~\ref{appendix:SHL-boundary-cocycle}.  For each fixed $c\in\mathbb R$, almost surely, a shift $x_i\mapsto x_i+Nc$ of all locations induces
    \[
        R_{s,t}^{\mathcal S,N}(a)\mapsto R_{s,t}^{\mathcal S,N}(a-c)+c, \qquad U_{s,t}^{\mathcal S,N}(a)\mapsto U_{s,t}^{\mathcal S,N}(a-c)+c .
    \]
    By \eqref{eq:CHL-root-flow-def}, the same inverse and composition relations hold for CHL(0). These statements concern the maps themselves, before the time interpolation defining $\mathtt{HM}$ and $\mathtt R$.

    We now turn to the geometry of particle-tree paths. For $\mathcal A\in\mathfrak M_{\rm f}$, let $e_i^{\mathcal A}:=\beta_i^{\mathcal A}([0,1])$ be the closed edge of particle $i$. Let $E_0^{\mathcal S}=\mathbb R$ and  $E_0^{\mathcal C}=\mathbb T_N$ denote the unrescaled initial boundaries, and define the particle tree at macroscopic time $u$ by
    \[
        \mathscr T_u^{\mathcal A,N}:=E_0^{\mathcal A}\cup\bigcup_{i:\,t_i\le N^2u}e_i^{\mathcal A}.
    \]
    To trace a point back along its tree path, we define the projection $P_{s,t}^{\mathcal A,N}$ as follows. For $0\le s\le t$ and $z\in\mathscr T_t^{\mathcal A,N}$, let  $\gamma_z:[0,1]\to\mathscr T_t^{\mathcal A,N}$ be a simple path, oriented so that $\gamma_z(0)\in E_0^{\mathcal A}$ and $\gamma_z(1)=z$. Set
    \[
        \theta_s(z):=\sup\{q\in[0,1]:      \gamma_z([0,q])\subset\mathscr T_s^{\mathcal A,N}\}, \qquad P_{s,t}^{\mathcal A,N}(z):=\gamma_z(\theta_s(z)).
    \]
    Since $\mathscr T_s^{\mathcal A,N}\subset\mathscr T_t^{\mathcal A,N}$ and its intersection with a path is an initial segment, $P_{s,t}^{\mathcal A,N}(z)$ is independent of the choice and parametrization of $\gamma_z$. In words, $P_{s,t}^{\mathcal A,N}(z)$ is the last point on the tree path from the initial boundary to $z$ that was already present at time $s$. The following lemma shows that, in boundary coordinates, this projection is given by the Root map $R_{s,t}^{\mathcal A,N}$. This relation is used below to compare Fingers with Root paths.

    \begin{lemma}\label{lem:root-finger-identity}
    For every $x\in\mathbb R$ and $0\le s<t$,
    \[
        \begin{aligned}
        P_{s,t}^{\mathcal S,N} \bigl(\mathcal S_{N^2t}^\lambda(Nx)\bigr) =\mathcal S_{N^2s}^\lambda \bigl(NR_{s,t}^{\mathcal S,N}(x)\bigr),\qquad P_{s,t}^{\mathcal C,N} \bigl(\mathcal C_{N^2t}^{N,\lambda}(Nx)\bigr) =\mathcal C_{N^2s}^{N,\lambda} \bigl(NR_{s,t}^{\mathcal C,N}(x)\bigr).
        \end{aligned}
    \]
    \end{lemma}
    \begin{proof}
    For SHL(0), the first identity follows from \cite[Lemma~7.5]{berger2022growth} and the scaling relation $S_x^\lambda(z)=\lambda S_{x/\lambda}^1(z/\lambda)$. The CHL(0) case follows in the same way.
    \end{proof}
    \subsection{Geometric tightness}\label{SS.geometric-tightness}
    In this subsection, we establish the geometric estimates needed to compare Fingers with Root paths and Gaps with harmonic measure paths. For $\mathcal A\in\mathfrak M$, let $I_{k,N}^{\mathcal A}=[k\Delta_N^{\mathcal A},(k+1)\Delta_N^{\mathcal A}]$. Write $\mathcal P_j^{\mathcal A,N}$ for the entire particle attached at time $t_j$ (as a subset of the unrescaled geometric domain), and write $j\prec i$ if particle $j$ lies on the tree path from the initial boundary to particle $i$. The values $R_{u,\bar t_i}^{\mathcal A,N}(\bar x_i)$ refer to the maps before interpolation in $u$, and $R_{u\pm,\bar t_i}^{\mathcal A,N}(\bar x_i)$ to their one-sided limits. Let $\widetilde{\mathcal A}_s$ be obtained by reversing the time order of the slit maps in $\mathcal A_s$, and write its martingale decomposition as
    \begin{equation}\label{eq:centered-reverse-composition}
        \widetilde{\mathcal A}_s(w)=w+\mathrm{i}v_Ns+M_s^{\mathcal A}(w),\qquad v_N=
        \begin{cases}
        v_{\mathcal A}, & \mathcal A=\mathcal S, \\[1mm]
        -2\pi N^2\log\bigl(1-\delta_N(\lambda)^2\bigr), & \mathcal A=\mathcal C.
        \end{cases}
    \end{equation}

    The following lemma controls the deviations of the rescaled aggregate maps from linear growth at rate $v_{\mathcal A}$.
    \begin{lemma}\label{lem:future-label-flatness}
    Fix $T,R_0>0$, $\Lambda\Subset(0,\infty)$, $p\ge8$, and $\mathcal A\in\mathfrak M_{\rm f}$. There exist $N_0,C<\infty$, depending only on $p,T,R_0,\Lambda$, such that, for $N\ge N_0$, $\lambda\in\Lambda$, $0\le u\le T$, $y\in E^{\mathcal A}$, $0<\varepsilon\le1$, and $w\in\overline{\mathbb H}$ satisfying
    \[
        |w-Ny|\le R_0,\quad \mathcal A=\mathcal S, \qquad \inf_{k\in2\pi N\mathbb Z}|w-Ny-k|\le R_0,\quad \mathcal A=\mathcal C,
    \]
    we have
    \begin{equation}\label{eq:future-label-flatness}
        \mathbb P\!\left\{   d^{\mathcal A}\!\left(\mathfrak p_N^{\mathcal A}    \bigl(\mathcal A_{N^2u}(w)\bigr), \mathfrak e_{\mathcal A}^{v_{\mathcal A}}(u,y) \right)>\varepsilon \right\} \le C\varepsilon^{-p}N^{-p/2}.
    \end{equation}

    For $u\le t\le T$, \eqref{eq:future-label-flatness} also holds conditionally on the Poisson points in $(N^2u,N^2t]$ for any    $(Y,W)$ measurable with respect to them and satisfying the same conditions as $(y,w)$. Moreover, for SHL(0),
    \begin{equation}\label{eq:SHL-whole-slit-moment}
        \sup_{\lambda\in\Lambda}\sup_{x\in\mathbb R} \mathbb E\!\left[ \sup_{0\le s\le N^2T} \sup_{0<\vartheta\le R_0} \bigl|M_s^{\mathcal S}(x+\mathrm i\vartheta)\bigr|^p \right] \le CN^{p/2}.
    \end{equation}
    \end{lemma}
    \begin{proof}
    For CHL(0), let $\Delta_\xi=A_\xi-\mathrm{id}$. By \cite[Lemmas~5.5 and~6.1]{procaccia2023cylindrical},
    \[
        M_s^{\mathcal C}(w)=\int_{(0,s]}\!\int_{\mathbb T_N} \Delta_\xi(\widetilde{\mathcal A}_{r-}(w))\, \widetilde{\mathcal P}(dr,d\xi),\qquad \widetilde{\mathcal P}(dr,d\xi) =\mathcal P^{\mathcal C,N}(dr,d\xi)-dr\,d\xi,
    \]
    and $v_N=v_{\mathcal A}+O_\Lambda(N^{-2})$. Continuity extends this identity to $w\in\overline{\mathbb H}$. For $q\ge2$, \cite[Theorem~3.2, $\mathsf{BJ}_{2,q}$]{marinelli2014maximal} and Lemma~\ref{lem:one-slit-moments} give
    \begin{equation}\label{eq:SHL-BDG-Rosenthal}
        \begin{aligned}
        \mathbb E\left[\sup_{r\le s}|M_r^{\mathcal C}(w)|^q\right] &\le C_q\mathbb E\!\left[ \left(\int_0^s\!\int_{\mathbb T_N} |\Delta_\xi(\widetilde{\mathcal A}_{r-}(w))|^2\,d\xi\,dr\right)^{q/2} +\int_0^s\!\int_{\mathbb T_N} |\Delta_\xi(\widetilde{\mathcal A}_{r-}(w))|^q\,d\xi\,dr\right] \\
        &\le C_{q,\Lambda}(s^{q/2}+s).
        \end{aligned}
    \end{equation}
    As in \cite[proof of Theorem~7.1]{berger2022growth}, let $h=\tfrac14\inf_{\lambda\in\Lambda}v_{\mathcal A}>0$, $L=N\vee N^2T$, and $E_N(w)=\{\im\widetilde{\mathcal A}_N(w)\ge hN\}$. On $E_N(w)$, the image of $w$ has reached height at least $hN$ by time $N$; since heights do not decrease, the later slit maps then act on points far from the boundary, where their deviations are small by \eqref{eq:complex-slit-high-moment}. For all sufficiently large $N$, $v_N\ge2h$, so
    \begin{equation*}
        E_N(w)^c\subset\{|M_N^{\mathcal C}(w)|\ge hN\},\qquad \mathbb P(E_N(w)^c) \le(hN)^{-4q}\mathbb E\left[\sup_{r\le N}|M_r^{\mathcal C}(w)|^{4q}\right] \le C_{q,\Lambda}N^{-2q}.
    \end{equation*}
    Define
    \[
        J_s(w)=\int_{(0,s]}\!\int_{\mathbb T_N} \mathbf1_{\{r>N\}}\mathbf1_{E_N(w)} \Delta_\xi(\widetilde{\mathcal A}_{r-}(w))\, \widetilde{\mathcal P}(dr,d\xi).
    \]
    The event $E_N(w)$ is measurable at time $N$, so the integrand is predictable. Height monotonicity and \eqref{eq:complex-slit-high-moment} imply, for $a\in\{2,q\}$,
    \[
        \int_0^L\!\int_{\mathbb T_N} \mathbf1_{\{r>N\}}\mathbf1_{E_N(w)} |\Delta_\xi(\widetilde{\mathcal A}_{r-}(w))|^a\,d\xi\,dr \le C_{a,h,\Lambda}LN^{1-a}.
    \]
    Applying the same maximal inequality to $J$ gives
    \begin{equation*}
        \mathbb E\left[\sup_{s\le L}|J_s(w)|^q\right] \le C_{q,h,\Lambda}\bigl((L/N)^{q/2}+LN^{1-q}\bigr) \le C_{q,T,\Lambda}N^{q/2}.
    \end{equation*}
    On $E_N(w)$, $M_s^{\mathcal C}(w)=M_{s\wedge N}^{\mathcal C}(w)+J_s(w)$. Hence
    \[
        \begin{aligned}
        \mathbb E\!\left[\mathbf1_{E_N(w)}\sup_{s\le L}|M_s^{\mathcal C}(w)|^q\right] &\le2^{q-1}\!\left( \mathbb E\left[\sup_{s\le N}|M_s^{\mathcal C}(w)|^q\right]+ \mathbb E\left[\sup_{s\le L}|J_s(w)|^q\right]\right) \le C_{q,T,\Lambda}N^{q/2}, \\
        \mathbb E\!\left[\mathbf1_{E_N(w)^c}\sup_{s\le L}|M_s^{\mathcal C}(w)|^q\right] &\le\mathbb P(E_N(w)^c)^{1/2} \left(\mathbb E\left[\sup_{s\le L}|M_s^{\mathcal C}(w)|^{2q}\right]\right)^{1/2} \le C N^{-q}(L^q+L)^{1/2}\le C.
        \end{aligned}
    \]
    Hence
    \begin{equation}\label{eq:SHL-fixed-input-moment}
        \sup_{\lambda\in\Lambda}\sup_{w\in\overline{\mathbb H}} \mathbb E\left[\sup_{s\le N^2T}|M_s^{\mathcal C}(w)|^q\right] \le C_{q,T,\Lambda}N^{q/2},\qquad q\ge2.
    \end{equation}
    By the definitions of $\mathfrak p_N^{\mathcal C}$ and  $\mathfrak e_{\mathcal C}^{v_{\mathcal C}}$,
    \begin{equation}\label{eq:future-label-chart-bound}
        d^{\mathcal C}\!\left( \mathfrak p_N^{\mathcal C} \bigl(\widetilde{\mathcal C}_{N^2u}(w)\bigr), \mathfrak e_{\mathcal C}^{v_{\mathcal C}}(u,y) \right) \le \frac{C_{T,R_0,\Lambda}}{N} \left(1+\bigl|M_{N^2u}^{\mathcal C}(w)\bigr|\right).
    \end{equation}
    Since $\widetilde{\mathcal C}_s\stackrel{\rm d}{=}     \mathcal C_s$ at each fixed time $s$,   \eqref{eq:SHL-fixed-input-moment} and   \eqref{eq:future-label-chart-bound} yield
    \begin{equation*}
        \mathbb E\!\left[ d^{\mathcal C}\!\left( \mathfrak p_N^{\mathcal C}(\mathcal C_{N^2u}(w)), \mathfrak e_{\mathcal C}^{v_{\mathcal C}}(u,y) \right)^p\right]\le CN^{-p/2}.
    \end{equation*}
    Markov's inequality proves  \eqref{eq:future-label-flatness} for CHL(0).

    For SHL(0), \cite[proof of Lemma~5.1]{berger2022growth} gives the  corresponding integral identity on $\mathbb R$ with $\mathrm i v_{\mathcal A}$.    Lemma~\ref{lem:one-slit-moments} and the preceding argument give  \eqref{eq:SHL-fixed-input-moment} with $\mathcal C$ replaced by $\mathcal S$, for every $w\in\mathbb H$.   By \cite[Theorem~4.3 and Proposition~4.7]{berger2022growth}, $s\mapsto M_s^{\mathcal S}$ is c\`adl\`ag in  $\operatorname{Hol}(\mathbb H)$, and its law is translation invariant. It remains to prove \eqref{eq:SHL-whole-slit-moment}, which requires a bound that is uniform down to the boundary; for this we use the Poisson representation of the holomorphic function $M_s^{\mathcal S}$ and a maximal inequality. Let $Q\subset[0,N^2T]$ be countable and dense with $N^2T\in Q$.  For $\delta>0$, let $ B_\delta(x):=\sup_{s\in Q}    |M_s^{\mathcal S}(x+\mathrm i\delta)|$.  Then
    \begin{equation}\label{eq:SHL-strip-bound}
        \sup_{\delta>0}\sup_{x\in\mathbb R} \mathbb E [B_\delta(x)^p] \le CN^{p/2}.
    \end{equation}
    Tonelli's theorem and monotonicity imply that, for each $s$,
    \[
        \zeta\longmapsto M_s^{\mathcal S}\!\left( \mathrm i\delta+\mathrm i\frac{1-\zeta}{1+\zeta} \right)
    \]
    belongs almost surely to $H^2(\mathbb D)$. Define the  Hardy--Littlewood maximal operator and its truncated version by
    \[
        \mathcal Mf(x):=\sup_{r>0}\frac1{2r}\int_{x-r}^{x+r}|f(a)|\,da,    \qquad    \mathcal M_\rho f(x):=\sup_{0<r\le\rho}\frac1{2r}\int_{x-r}^{x+r}|f(a)|\,da.
    \]
    Then    \cite[Theorems~I.4.2 and~II.3.1]{garnett2007bounded}, first for $s\in Q$ and then for all $s\le N^2T$ by right continuity, give
    \begin{equation*}
        \sup_{s\le N^2T} |M_s^{\mathcal S}(x+\mathrm i(\delta+\vartheta))| \le \frac1\pi\int_{\mathbb R} \frac{\vartheta B_\delta(a)} {(x-a)^2+\vartheta^2}\,da \le C\mathcal M B_\delta(x), \qquad \vartheta>0.
    \end{equation*}
    Moreover, \cite[Theorem~I.4.3(b)]{garnett2007bounded}, applied to $B_\delta\mathbf 1_{[-K-\rho,K+\rho]}$, together with translation invariance, yields
    \[
        2K\,\mathbb E\bigl[(\mathcal M_\rho B_\delta(0))^p\bigr] \le C_p(2K+2\rho)\, \mathbb E [B_\delta(0)^p].
    \]
    Letting $K\to\infty$, then $\rho\to\infty$, and finally  $\delta\downarrow0$, and using \eqref{eq:SHL-strip-bound},  we obtain \eqref{eq:SHL-whole-slit-moment}.

    By \cite[Lemma~5.4]{berger2022growth}, for fixed $s$ and $x\in\mathbb R$, the limit as $\vartheta\downarrow0$ equals the boundary value. Since $\mathcal S_s$ and  $\widetilde{\mathcal S}_s$ have the same law at each fixed time  \cite[Section~2.1]{berger2022growth}, Fatou's lemma and the preceding  bounds give
    \[
        \sup_{w\in\overline{\mathbb H}} \mathbb E\left[\left| \mathcal S_s(w)-w-\mathrm i v_Ns \right|^p\right] \le CN^{p/2},\quad  \sup_{x\in\mathbb R} \mathbb E\left[\sup_{0\le\vartheta\le R_0} \left| \mathcal S_s(x+\mathrm i\vartheta) -x-\mathrm i\vartheta-\mathrm i v_{\mathcal A} s \right|^p\right] \le CN^{p/2}.
    \]
    The SHL(0) analogue of \eqref{eq:future-label-chart-bound} then proves \eqref{eq:future-label-flatness}.

    Since the Poisson points before $N^2u$ are independent of those in $(N^2u,N^2t]$, the preceding bounds hold conditionally for $(Y,W)$, and also with one additional point at $r\in[N^2u,N^2t]$, using $\mathcal A_{r-}$ when $r=N^2u$.  By \cite[Lemma~5.4]{berger2022growth} and the Mecke identity  \cite[Theorem~4.1]{last2017lectures}, almost surely every Poisson point $(r,x)$ satisfies $ \mathcal P_{(r,x)}^{\mathcal S,N}=\mathcal S_{r-} \bigl(x+\mathrm i[0,\lambda]\bigr). $
    \end{proof}

    All particles attached to particle $j$ (directly or indirectly) remain between its two side paths. For $j\prec i$, Proposition~\ref{prop:SHL-boundary-cocycle} gives $R_{\bar t_j,\bar t_i}^{\mathcal A,N}(\bar x_i) \in[\bar x_{j,-},\bar x_{j,+}]$ and  $R_{\bar t_j-,\bar t_i}^{\mathcal A,N}(\bar x_i)=\bar x_j$. By \eqref{eq:particle-side-labels},
    \begin{equation}\label{eq:ancestral-centre-root}
        d_{E^{\mathcal A}}\!\left( R_{\bar t_j\pm,\bar t_i}^{\mathcal A,N}(\bar x_i), \bar x_j\right)\le C_{\eta,\lambda}/N, \qquad j\prec i.
    \end{equation}
    We next control the deviations of entire particles from their growth-rate embeddings, which will allow us to compare Fingers with embedded Root paths. For $\mathcal A\in\mathfrak M_{\rm c}$, define the maximal particle deviation by
    \[
        U_N^{\mathcal A}(T):=\max_{j:\bar t_j\le T} \sup_{z\in\mathcal P_j^{\mathcal A,N}} d^{\mathcal A}\!\left(\mathfrak p_N^{\mathcal A}(z), \mathfrak e_{\mathcal A}^{v_{\mathcal A}}(\bar t_j,\bar x_j)\right).
    \]
    \begin{lemma}\label{lem:completed-cell-oscillation}
    Fix $T<\infty$ and $\mathcal A\in\mathfrak M_{\rm c}$. Then
    \begin{equation}\label{eq:uniform-complete-particle-localization}
        U_N^{\mathcal A}(T)\pto0,\qquad \sup_{i:\,\bar t_i\le T}\max_k \operatorname{diam}_{d^{\mathcal A}} \mathtt{Fin}_i^{\mathcal A,N} \bigl(I_{k,N}^{\mathcal A}\cap[0,\bar t_i]\bigr)\pto0,
    \end{equation}
    and
    \begin{equation}\label{eq:finite-volume-completed-finger-transfer}
        \sup_{i:\,\bar t_i\le T} d_T^{\mathcal A}\!\left( \mathtt{Fin}_i^{\mathcal A,N}, \mathfrak e_{\mathcal A}^{v_{\mathcal A}}\bigl(\mathtt R_{\cdot,\bar t_i}^{\mathcal A,N}(\bar x_i)\bigr)\right)\pto0.
    \end{equation}
    \end{lemma}
    \begin{proof}
    First take $\mathcal A=\mathcal H$. Set $a=\eta/N$ and $c=c_{\mathcal H}(a)$ in \eqref{logarithmic capacity}.  By \cite[Theorem~3.1]{norris2012hastings}, with probability tending to one,
    \[
        \sup_{z\in\mathcal P_j^{\mathcal H,N}}|ze^{-cj-\mathrm ix_j}-1|\le a^{2/3}\{\log(1/a)\}^{8},\quad j\le\lfloor a^{-6}\rfloor.
    \]
    Let $K_N(t)$ be the number of arrival times up to time $N^3t$. Then
    \[
        \mathbb E\left[\sup_{t\le T}|K_N(t)-N^3t|^2\right]\le4N^3T, \qquad N^2c=v_{\mathcal A}+O_\eta(N^{-1}).
    \]
    Thus $K_N(T)\le\lfloor a^{-6}\rfloor$ with probability tending to one. For $|z-1|\le \tfrac12$, $|\log z|\le C|z-1|$. Hence, on these events,
    \[
        U_N^{\mathcal H}(T)\le C a^{2/3}\{\log(1/a)\}^{8}+\frac cN\sup_{t\le T}|K_N(t)-N^3t|+T|N^2c-v_{\mathcal A}|\pto0.
    \]
    The estimates are uniform for $\eta$ in compact subsets of $(0,\infty)$. For $\mathcal C$, take $\eta_N=N(e^{\lambda/N}-1)$ and apply Remark~\ref{rmk:radial-from-cylindrical}. This proves the first conclusion of \eqref{eq:uniform-complete-particle-localization}.

    For $\mathcal A=\mathcal C$, Lemma~\ref{lem:root-finger-identity} identifies, for $u\in\Gamma_N^{\mathcal C}\cap[0,\bar t_i)$, the point of the tree path of particle $i$ selected at time $u$ as $\mathcal C_{N^2u}\!(N R_{u,\bar t_i}^{\mathcal C,N}(\bar x_i))$. Hence, for any collection $\mathcal I$ of such pairs $(i,u)$, with $\bar t_i\le T$ and selected from the Poisson points after $N^2u$, Lemma~\ref{lem:future-label-flatness} and the Mecke identity \cite[Theorem~4.1]{last2017lectures} give
    \begin{equation}\label{eq:hl-geometry-birth-knot-count}
        \begin{aligned}
        &\mathbb P\Bigg\{\sup_{(i,u)\in\mathcal I}\Bigg\|    \chi_{\mathcal C}\!\left(\mathfrak p_N^{\mathcal C} \!\left( \mathcal C_{N^2u} \!\left( N R_{u,\bar t_i}^{\mathcal C,N}(\bar x_i) \right) \right) \right) -\left( R_{u,\bar t_i}^{\mathcal C,N}(\bar x_i),v_{\mathcal A} u \right) \Bigg\|_\infty >\vartheta_N \Bigg\} \le C\vartheta_N^{-p}N^{-p/2}\mathbb E[\#\mathcal I],
        \end{aligned}
    \end{equation}
    where $p\ge8$.  For the collection of all such pairs,
    \[
        \mathbb E[\#\mathcal I] \le 2\pi N^3 \sum_{u\in\Gamma_N^{\mathcal C}\cap[0,T]}(T-u) =O_T(N^4).
    \]
    Thus, taking $p=24$ and $\vartheta_N=N^{-1/4}$, \eqref{eq:hl-geometry-birth-knot-count} shows that the corresponding uniform error exceeds $\vartheta_N$ with probability $O_T(N^{-2})$.  Remark~\ref{rmk:radial-from-cylindrical} gives the corresponding estimate for $\mathcal H$, with an additional error  $O_{T,\eta}(N^{-1})$. The value at $u=\bar t_i$, including the tip, is  controlled by the first conclusion of \eqref{eq:uniform-complete-particle-localization}.

    For $\mathcal A\in\mathfrak M_{\rm c}$, let
    \[
        M_N^{\mathcal A}(T) := \max_{\substack{k: I_{k,N}^{\mathcal A}\cap[0,T]\neq\varnothing}} \; \sup_{\substack{s,t\in I_{k,N}^{\mathcal A}\cap[0,T], s\le t}} \; \sup_{x\in E^{\mathcal A}} d_{E^{\mathcal A}} \!\left( R_{s,t}^{\mathcal A,N}(x),x \right).
    \]
    For $\mathcal A=\mathcal C$, since $R_{t,b}^{\mathcal A,N}$ is onto, the cocycle property bounds the oscillation on each grid interval $[a,b]=I_{k,N}^{\mathcal A}$ by twice the maximal deviation of $R_{s,b}^{\mathcal A,N}$ from the identity over $a\le s\le b$. By \eqref{eq:clock-cell-fourth-moment}, Lemma~\ref{lem:one-slit-moments}, and the mesh argument of Lemma~\ref{lem:automatic-clock-oscillation}, applied after time reversal, for every $\varepsilon>0$,
    \[
        \mathbb P\{M_N^{\mathcal A}(T)>\varepsilon\} \le C_{\varepsilon,T}\bigl(1+T/\Delta_N^{\mathcal A}\bigr) \bigl((\Delta_N^{\mathcal A})^2+N^{-2}\Delta_N^{\mathcal A}\bigr) =O_{\varepsilon,T}(N^{-1}).
    \]
    The same bound holds for $\mathcal A=\mathcal H$ by Remark~\ref{rmk:radial-from-cylindrical} and $\lambda_N^{\mathcal H}\to\eta$. Thus $M_N^{\mathcal A}(T)\pto0$. Fix $[u,v]= I_{k,N}^{\mathcal A}\cap[0,\bar t_i]\neq\varnothing$. If a particle $j\preceq i$ contributes to  ${\mathtt{Fin}}_i^{\mathcal A,N}([u,v])$, then  $u<\bar t_j\le v$,  except  $j=i$ when $v=\bar t_i$. By the cocycle property,
    \[
        \sup_{s\in[u,v]} d_{E^{\mathcal A}} \!\left( R_{s,\bar t_i}^{\mathcal A,N}(\bar x_i), R_{u,\bar t_i}^{\mathcal A,N}(\bar x_i) \right) \le M_N^{\mathcal A}(T),
    \]
    and \eqref{eq:ancestral-centre-root} gives
    \[
        d_{E^{\mathcal A}} \!\left( \bar x_j, R_{u,\bar t_i}^{\mathcal A,N}(\bar x_i) \right) \le C_{\eta,\lambda}N^{-1}+M_N^{\mathcal A}(T), \qquad 0\le\bar t_j-u\le\Delta_N^{\mathcal A}.
    \]
    For $j=i$, the same bounds follow from  $R_{\bar t_i,\bar t_i}^{\mathcal A,N}(\bar x_i)=\bar x_i$.  Combining these estimates with the first conclusion of  \eqref{eq:uniform-complete-particle-localization} and  \eqref{eq:hl-geometry-birth-knot-count}, we obtain, uniformly in  $\bar t_i\le T$ and $k$,
    \[
        \begin{aligned}
        &\sup_{s\in[u,v]} d^{\mathcal A} \!\left( {\mathtt{Fin}}_i^{\mathcal A,N}(s), \mathfrak e_{\mathcal A}^{v_{\mathcal A}} \!\left( s,R_{s,\bar t_i}^{\mathcal A,N}(\bar x_i) \right) \right) \vee\; \operatorname{diam}_{d^{\mathcal A}} {\mathtt{Fin}}_i^{\mathcal A,N}([u,v]) \\
        &\hspace{20mm} \le C\left( \vartheta_N +U_N^{\mathcal A}(T) +M_N^{\mathcal A}(T) +\Delta_N^{\mathcal A} +N^{-1} \right) \pto0 .
        \end{aligned}
    \]
    The linear interpolation defining $\mathtt R$ differs from the Root maps by at most $2M_N^{\mathcal A}(T)$ on these intervals. Taking the supremum over $i$ and $k$ completes \eqref{eq:uniform-complete-particle-localization} and proves \eqref{eq:finite-volume-completed-finger-transfer}.
    \end{proof}

    We next extend Lemma~\ref{lem:completed-cell-oscillation}     to SHL(0) and then to the full path collections. For paths $\gamma, \gamma'$ paired with their particle arrival times $r, r'$, use the metrics
    \[
        d_{\rm path}(\gamma,\gamma')\vee |\tanh r-\tanh r'|, \qquad D_{\mathcal A}(\gamma,\gamma')\vee |\tanh r-\tanh r'|
    \]
    for boundary and geometric paths, where $d_{\rm path}$ is defined in \eqref{path-distance}. We also add an isolated point $\dagger$, which is at distance $1$ from every such pair; this ensures that the collections below are non-empty.  For $G>0$ and $K\in\{\mathbf{HM}^{\mathcal A,N},  \mathbf R^{\mathcal A,N},\mathbf{Fin}^{\mathcal A,N}\}$, set
    \[
        K\langle G\rangle :=\overline{\{(\gamma_i,\bar t_i): \bar t_i\le G\}}\cup\{\dagger\},
    \]
    where $\gamma_i$ ranges over the paths generated  by particle $i$, including both side paths for $\mathbf{HM}^{\mathcal A,N}$. We include only particles arriving by time $G$, without restricting their paths; their arrival times remain unchanged under later time restrictions. Extend the embedding by $\mathfrak e_{\mathcal A}^{v}(\gamma,r) :=(\mathfrak e_{\mathcal A}^{v}(\gamma),r)$ and $\mathfrak e_{\mathcal A}^{v}(\dagger):=\dagger$. When comparing paths alone, we ignore the arrival times and add $\dagger$ to both collections, using $1\wedge d$ between paths and distance $1$ to $\dagger$.

    \begin{proposition}\label{prop:process-tip-geometric-transfer}
    For every $T<\infty$ and $\mathcal A\in\mathfrak M$,
    \begin{equation}\label{eq:process-tip-geometric-transfer}
        d_H\!\left(\mathbf{Fin}^{\mathcal A,N}\!\langle T\rangle, \mathfrak e_{\mathcal A}^{v_{\mathcal A}}(\mathbf R^{\mathcal A,N}\!\langle T\rangle)\right) \pto0.
    \end{equation}
    For $\mathcal S$, the convergence is understood in each local spatial coordinate: for every $m\ge1$,
    \begin{equation}\label{eq:SHL-completed-cell-transfer}
        \sup_{i:\bar t_i\le T} d_T^{\mathcal S,m}\!\left( \mathtt{Fin}_i^{\mathcal S,N}, \mathfrak e_{\mathcal S}^{v_{\mathcal S}}\bigl( \mathtt R_{\cdot,\bar t_i}^{\mathcal S,N}(\bar x_i)\bigr)\right) \pto0.
    \end{equation}
    Moreover, for $p\ge8$, $R>m+1$, and all sufficiently large $N$,
    \begin{equation}\label{eq:SHL-remote-particle-summation}
        \mathbb P\!\left\{\exists j:\bar t_j\le T,\ |x_j|>NR, \inf_{z\in\mathcal P_j^{\mathcal S,N}} \frac{\operatorname{sgn}(x_j)\operatorname{Re}z}{N}\le m+1\right\} \le C_{p,T,\lambda}N^{3-p/2}(R-m-1)^{1-p}.
    \end{equation}
    \end{proposition}
    \begin{proof}
    For $\mathcal A\in\mathfrak M_{\rm c}$, \eqref{eq:process-tip-geometric-transfer} follows immediately from Lemma~\ref{lem:completed-cell-oscillation}.

    For $\mathcal A=\mathcal S$, fix $m\ge1$ and set $L=m+4$. Lemma~\ref{lem:root-finger-identity} gives, for $u\in\Gamma_N^{\mathcal S}\cap[0,\bar t_i)$,
    \[
        \mathtt{Fin}_i^{\mathcal S,N}(u) =\mathfrak p_N^{\mathcal S}\!\left( \mathcal S_{N^2u}\!\left( N R_{u,\bar t_i}^{\mathcal S,N}(\bar x_i)\right)\right).
    \]
    Proposition~\ref{prop:SHL-boundary-cocycle} and the Mecke identity \cite[Theorem~4.1]{last2017lectures} give
    \[
        \mathbb E\left[\#\!\left\{(i,u): u\in\Gamma_N^{\mathcal S}\cap[0,\bar t_i),\ \bar t_i\le T,\ \left|R_{u,\bar t_i}^{\mathcal S,N}(\bar x_i)\right|\le L\right\}\right] =2LN^3\!\sum_{u\in\Gamma_N^{\mathcal S}\cap[0,T]}(T-u) =O_{L,T}(N^4).
    \]
    The equality follows from $r_x^{\mathcal S}(x)=x$ and translation invariance of the Root flow, which gives $\int_{\mathbb R}\mathbb P\{|R_{u,t}^{\mathcal S,N}(y)|\le L\}\,dy=2L$.  Applying \eqref{eq:future-label-flatness} conditionally on the Poisson points in $(N^2u,N^2\bar t_i]$, with $p=24$ and $\vartheta_N=N^{-1/4}$, yields
    \[
        \mathbb P\!\left\{ \sup_{\substack{i:\bar t_i\le T,\ u\in\Gamma_N^{\mathcal S}\cap[0,\bar t_i),\ |R_{u,\bar t_i}^{\mathcal S,N}(\bar x_i)|\le L}} \left\|\chi_{\mathcal S}\bigl(\mathtt{Fin}_i^{\mathcal S,N}(u)\bigr) -\bigl(R_{u,\bar t_i}^{\mathcal S,N}(\bar x_i),v_{\mathcal A} u\bigr) \right\|_\infty>\vartheta_N\right\} =O_{L,T,\lambda}(N^{-2}).
    \]
    The same conditional Markov inequality, integrated over $|R_{u,\bar t_i}^{\mathcal S,N}(\bar x_i)|>L$, gives
    \[
        \begin{aligned}
        \mathbb P\!\Bigl\{\exists(i,u):\;& \bar t_i\le T,\quad u\in\Gamma_N^{\mathcal S}\cap[0,\bar t_i),\quad |R_{u,\bar t_i}^{\mathcal S,N}(\bar x_i)|>L, \\
        &\operatorname{sgn}\!\bigl( R_{u,\bar t_i}^{\mathcal S,N}(\bar x_i)\bigr) \operatorname{Re}\mathtt{Fin}_i^{\mathcal S,N}(u) \le m+1\Bigr\} \le C_{T,\lambda}N^{-8}.
        \end{aligned}
    \]
    It remains to control complete particles. For a Poisson point $(r,x)$, almost surely $    \mathcal P_{(r,x)}^{\mathcal S,N}=\mathcal S_{r-}(x+\mathrm i[0,\lambda])$. By \eqref{eq:SHL-whole-slit-moment}, Markov's inequality, and Mecke's identity, for $R>m+1$,
    \[
        \begin{aligned}
        &\mathbb P\!\left\{\exists j:\bar t_j\le T,\ |x_j|>NR, \inf_{z\in\mathcal P_j^{\mathcal S,N}} \frac{\operatorname{sgn}(x_j)\operatorname{Re}z}{N}\le m+1\right\} \\
        &\qquad\le \int_0^{N^2T}\!\int_{|x|>NR} C N^{p/2}\bigl(|x|-N(m+1)\bigr)^{-p}\,dx\,dr \le C_{p,T,\lambda}N^{3-p/2}(R-m-1)^{1-p},
        \end{aligned}
    \]
    which proves \eqref{eq:SHL-remote-particle-summation}. The same bound on $|x_j|\le N(L+2)$, with $p=24$, gives
    \[
        \mathbb P\!\left\{ \sup_{\substack{j:\bar t_j\le T,\ |x_j|\le N(L+2)}} \sup_{z\in\mathcal P_j^{\mathcal S,N}} \left\|\chi_{\mathcal S}\bigl(\mathfrak p_N^{\mathcal S}(z)\bigr) -(\bar x_j,v_{\mathcal A}\bar t_j)\right\|_\infty >C\vartheta_N\right\} =O_{L,T,\lambda}(N^{-3}).
    \]
    The fourth-moment bound for the Root flow, together with monotonicity, gives
    \[
        \max_{k:I_{k,N}^{\mathcal S}\cap[0,T]\neq\varnothing} \sup_{\substack{s,u\in I_{k,N}^{\mathcal S}\cap[0,T],\ s\le u,\ x\in\mathbb R}} \left| \tanh R_{s,u}^{\mathcal S,N}(x)-\tanh x \right| \xrightarrow{\mathbb P}0.
    \]
    Hence, uniformly for $s,u\in I_{k,N}^{\mathcal S}\cap[0,\bar t_i]$, $ \tanh R_{s,\bar t_i}^{\mathcal S,N}(\bar x_i) -\tanh R_{u,\bar t_i}^{\mathcal S,N}(\bar x_i) =o_{\mathbb P}(1). $ If $ |R_{u,\bar t_i}^{\mathcal S,N}(\bar x_i)|\le L, $ then with probability tending to one, $ R_{s,\bar t_i}^{\mathcal S,N}(\bar x_i)\in[-L-1,L+1]. $ Since $\tanh^{-1}$ is uniformly continuous on $[\tanh(-L-1),\tanh(L+1)]$, \eqref{eq:ancestral-centre-root} gives, uniformly for the relevant $s,u,j$,
    \[
        \left|R_{s,\bar t_i}^{\mathcal S,N}(\bar x_i) -R_{u,\bar t_i}^{\mathcal S,N}(\bar x_i)\right| + \left|\bar x_j-R_{u,\bar t_i}^{\mathcal S,N}(\bar x_i)\right| =o_{\mathbb P}(1).
    \]
    Together with the preceding particle estimates, this gives the required bound after applying $\chi_{\mathcal S,m}$. If no particle edge is traversed on a grid interval, the Finger is constant there, and the same bound follows from its grid value and the Root oscillation estimate. If instead $ R_{u,\bar t_i}^{\mathcal S,N}(\bar x_i)>L, $ then on the event
    \[
        \left|\tanh R_{s,\bar t_i}^{\mathcal S,N}(\bar x_i) -\tanh R_{u,\bar t_i}^{\mathcal S,N}(\bar x_i)\right| <\tanh L-\tanh(L-1),
    \]
    we have $ R_{s,\bar t_i}^{\mathcal S,N}(\bar x_i)>L-1. $ Hence \eqref{eq:ancestral-centre-root} implies $ \bar x_j>m+2 $ for all sufficiently large $N$; the case $R_{u,\bar t_i}^{\mathcal S,N}(\bar x_i)<-L$ is identical. For $|\bar x_j|\le L+2$, the preceding particle estimate keeps the     whole particle on the corresponding side of $[-m-1,m+1]$; for   $|\bar x_j|>L+2$, apply \eqref{eq:SHL-remote-particle-summation} with $R=L+2$. The same side bound holds for the selected grid points. Since $\chi_{\mathcal S,m}$ is constant on each side of $[-m-1,m+1]$, their images under $\chi_{\mathcal S,m}$ coincide there. The terminal edge satisfies the same bounds, proving \eqref{eq:SHL-completed-cell-transfer}.

    Finally, in \eqref{eq:process-tip-geometric-transfer}, each path in the two collections is indexed by the same particle and carries the same arrival time. Hence \eqref{eq:SHL-completed-cell-transfer} extends to their closures and gives \eqref{eq:process-tip-geometric-transfer} for $\mathcal S$.
    \end{proof}
    The next lemma justifies the exterior construction of the Gaps and controls  their parametrization.
    \begin{lemma}\label{lem:gap-exterior-geometry}
    Almost surely, the tree $\mathscr T_N^{\mathcal A}$ in \eqref{eq:gap-complete-tree} is locally finite, the neighbourhoods and curves of Subsection~\ref{S.HL-models} can be chosen  measurably with $\sup_{r\ge0}g_i^{(2)}(r)=\infty$, and    $\mathbf{Gap}^{\mathcal A,N}$ is a measurable compact set. Every path in $\mathbf{Gap}^{\mathcal A,N}$ with finite starting time and position satisfies \eqref{eq:gap-strict-exterior}.  For $T<\infty$,
    \begin{equation}\label{eq:gap-height-approximation}
        \sup_{i:\sigma_i\le T}\sup_{\sigma_i\le u\le T} \bigl\{v_{\mathcal A}u-\chi_{\mathcal A}^{(2)} (\mathtt{Gap}_i^{\mathcal A,N}(u))\bigr\}\pto0.
    \end{equation}
    For $\mathcal A=\mathcal S$ and each fixed $R>0$, the supremum over $i$ in \eqref{eq:gap-height-approximation} is further restricted to those satisfying $\sigma_i\le T$ and    $\inf_{\sigma_i\le u\le T}|\rea\mathtt{Gap}_i^{\mathcal S,N}(u)|\le R$. Moreover,
    \begin{equation}\label{eq:gap-spatial-localization}
        \lim_{L\to\infty}\limsup_{N\to\infty} \P\left\{
        \begin{array}
        {l@{}l} \exists i:\sigma_i\le T,\quad &\displaystyle\inf_{\sigma_i\le u\le T} |\rea\mathtt{Gap}_i^{\mathcal S,N}(u)|\le R, \\
        &\displaystyle\sup_{\sigma_i\le u\le T} |\rea\mathtt{Gap}_i^{\mathcal S,N}(u)|>L
        \end{array}
        \right\}=0.
    \end{equation}
    \end{lemma}
    The proof is in Appendix~\ref{appendix:gap-exterior-geometry}. Both estimates hold for every measurable choice of regular neighbourhoods and shortest curves satisfying the conditions of Subsection~\ref{S.HL-models}.

    \section{\texorpdfstring{BW convergence of HL$(0)$ observables}{BW Convergence of HL(0) observables}}\label{S.HL-observables-BW}
    In this section, we prove BW convergence of the path collections in the three HL(0) geometries and identify the limit of the one-arm domination time in HL(0) and CHL(0), that is, we prove Theorems~\ref{thm:harmonic-measure-webs}, \ref{thm:joint-geometric-webs}, and \ref{thm:one-arm-domination}. We first establish the limits of $\mathbf R^{\mathcal A,N}$ and $\mathbf{HM}^{\mathcal A,N}$, and then use the geometric estimates of Section~\ref{S.HL-application} to prove the convergence of $\mathbf{Fin}^{\mathcal A,N}$ and $\mathbf{Gap}^{\mathcal A,N}$.

    Since the HL$(0)$ observables are defined only for nonnegative times, we first record a version of the convergence criterion on a bounded time window. Fix $\mathcal A\in\mathfrak M$ and $J=[s_0,s_1]$, with $0\le s_0<s_1<\infty$. Forward paths are continuous on $[\sigma_\pi,s_1]$, $\sigma_\pi\in J$; backward paths are treated after time reversal. We use the path metric and spatial compactification of Subsection~\ref{SS.convergence criteria}, restricted to $J$, and the count $\eta_K$ of Theorem~\ref{thm:BW-CBW-convergence-criterion} for $s_0\le t_0<t_0+h\le s_1$. Let $\mathcal W_{\mathcal A,J}^{\uparrow}:=\overline{\{\pi|_{[\sigma_\pi,s_1]}: \pi\in\mathcal W_{\mathcal A,\ge0}^{\uparrow},\ \sigma_\pi\in J\}}$. Let $D\subset E^{\mathcal A}\times J$ be a fixed countable set such that
    \[
        D\cap\bigl(E^{\mathcal A}\times(s_0,s_1)\bigr) \quad\text{and}\quad D\cap\bigl(E^{\mathcal A}\times\{s_j\}\bigr),\ j=0,1,
    \]
    are dense in the corresponding sets.

    Let $K_N$ be random nonempty compact non-crossing path collections forming a tight sequence. Suppose \textup{(I)} of Theorem~\ref{thm:BW-CBW-convergence-criterion} holds with selections indexed by $D$ and variance $\sigma_{\mathcal A}^2$, and \textup{(B2)} holds for $t_0,t_0+h\in J$. Then $ K_N\dto \sigma_{\mathcal A}\mathcal W_{\mathcal A,J}^{\uparrow}. $ Indeed, for any subsequential limit $K$, the closure of the selected limiting paths is a BW $W$, with $W\stackrel{d}{=}\sigma_{\mathcal A}\mathcal W_{\mathcal A,J}^{\uparrow}$ and $W\subseteq K$.  The argument of \eqref{eq:abstract-portmanteau-separated-bound}, after shifting both counting times forward by the same amount and slightly enlarging the initial interval or arc, yields \textup{(B2)} for $K$ on every $[c,d]\subset(s_0,s_1)$. For $s_0<c<d<s_1$, let $K|_{[c,d]}$ be the closure of $\{\pi|_{[\sigma_\pi\vee c,d]}:\pi\in K,\ \sigma_\pi\le d\}$, and define $W|_{[c,d]}$ similarly. The proof of Theorem~\ref{thm:BW-CBW-convergence-criterion}, applied on $[c,d]$, using the local mesh argument in \cite[Lemma~5.1 and proof of Theorem~4.2]{fontes2004brownian} and its cylindrical version in \cite[Theorem~2.5]{coupier2019directed}, shows that non-crossing and (B2) imply $K|_{[c,d]}\subseteq W|_{[c,d]}$ almost surely. Since $W\subseteq K$, we obtain $K|_{[c,d]}=W|_{[c,d]}$ almost surely. Letting $c\downarrow s_0$ and $d\uparrow s_1$ through countable sequences, compactness and continuity recover all paths with starting time $<s_1$, while the selected limits on $E^{\mathcal A}\times\{s_1\}$ and their closure recover all paths with starting time $s_1$. Thus $K=W$.

    \subsection{Finite-dimensional convergence}\label{S.finite-dimensional}
    We first prove finite-dimensional convergence for the harmonic measure and Root paths. For the one-slit boundary maps in \eqref{eq:one-slit-boundary-maps}, define their deviations from the identity by
    \[
        \begin{aligned}
        \Delta_h^{\mathcal S}:=h_0^{\mathcal S}-\id,\quad \Delta_r^{\mathcal S}:=r_0^{\mathcal S}-\id,\quad \Delta_h^{\mathcal C,N}:=h_0^{\mathcal C,N}-\id,\quad \Delta_r^{\mathcal C,N}:=r_0^{\mathcal C,N}-\id.
        \end{aligned}
    \]
    For $\star\in\{h,r\}$, define the covariance rates by
    \begin{equation}\label{eq:scaled-covariance-rates-section7}
        \begin{aligned}
        \mathcal K_\star^{\mathcal S,N}(d) :=\int_{\mathbb R}\Delta_\star^{\mathcal S}(u) \Delta_\star^{\mathcal S}(u+Nd)\,du,\qquad  \mathcal K_\star^{\mathcal C,N}(d) :=\int_{\mathbb T_N}\Delta_\star^{\mathcal C,N}(u) \Delta_\star^{\mathcal C,N}(u+Nd)\,du.
        \end{aligned}
    \end{equation}

    \begin{lemma}\label{lem:covariance-decay-section7}
    Fix $\Lambda\Subset(0,\infty)$. If $c_N\downarrow0$ and $Nc_N\to\infty$, then, uniformly for $\lambda\in\Lambda$, $\mathcal A\in\mathfrak M_{\rm f}$, and $\star\in\{h,r\}$,
    \begin{equation*}
        \sup_{d_{E^{\mathcal A}}(d,0)\ge c_N} |\mathcal K_\star^{\mathcal A,N}(d)|\rightarrow0.
    \end{equation*}
    \end{lemma}
    \begin{proof}
    The one-slit formulas and \eqref{eq:Rdelta-scaling-section3}, \eqref{eq:cylinder-deviation-tail} give, for all sufficiently large $N$ and $R\ge0$,
    \begin{equation}\label{eq:C-tail-section7}
        \int_{|u|>R}|\Delta_\star^{\mathcal S}(u)|^2\,du +\int_{d_{\mathbb T_N}(u,0)>R} |\Delta_\star^{\mathcal C,N}(u)|^2\,du \le\frac{C_\Lambda}{1+R}.
    \end{equation}

    Set $R=Nc_N/3$. If $d_{E^{\mathcal A}}(d,0)\ge c_N$, the balls of radius $R$ centred at $0$ and $-Nd$ are disjoint. Thus the Cauchy--Schwarz inequality and \eqref{eq:C-tail-section7} give
    \[
        |\mathcal K_\star^{\mathcal A,N}(d)| \le C_\Lambda(1+R)^{-1/2}\longrightarrow0,
    \]
    uniformly in the stated parameters.
    \end{proof}

    For fixed space--time points, let $\pi_{x,s}^{\mathcal A,\sigma}$ and $\widehat\pi_{x,t}^{\mathcal A,\sigma}$ be the almost surely unique paths of $\sigma\mathcal W_{\mathcal A,\ge0}^{\uparrow}$ and  $\sigma\mathcal W_{\mathcal A,\ge0}^{\downarrow}$ starting at $(x,s)$  and ending at $(x,t)$, respectively. The following proposition establishes finite-dimensional convergence of harmonic measure and Root paths, allowing random spatial labels.

    \begin{proposition}\label{prop:fdd-selected-paths-section7}
    Fix $\mathcal A\in\mathfrak M$ and $m\in\mathbb N$, with $\eta>0$ if $\mathcal A=\mathcal H$, and $\lambda>0$ if $\mathcal A\in\mathfrak M_{\mathrm f}$. For $1\le i\le m$, let $x_i\in E^{\mathcal A}$ and $s_i,t_i\in[0,\infty)$ be fixed. Let $X_{N,i}\in E^{\mathcal A}$ be random labels, which may depend on all the driving Poisson points, and let $s_{N,i},t_{N,i}\in\Gamma_N^{\mathcal A}$ be deterministic, with
    \[
        s_{N,i}\to s_i,\qquad t_{N,i}\to t_i,\qquad \max_{1\le i\le m} d_{E^{\mathcal A}}(X_{N,i},x_i) \xrightarrow{\mathbb P}0.
    \]
    Then, on compact time intervals,
    \[
        \left( \mathtt{HM}^{\mathcal A,N}_{s_{N,i},\cdot}(X_{N,i}) \right)_{i=1}^m \Rightarrow \left( \pi^{\mathcal A,\sigma_{\mathcal A}}_{x_i,s_i} \right)_{i=1}^m, \qquad \left( \mathtt R^{\mathcal A,N}_{\cdot,t_{N,i}}(X_{N,i}) \right)_{i=1}^m \Rightarrow \left( \widehat\pi^{\mathcal A,\sigma_{\mathcal A}}_{x_i,t_i} \right)_{i=1}^m
    \]
    in the product path topologies of Subsection~\ref{SS.convergence criteria}.
    \end{proposition}
    \begin{proof}
    We first take $X_{N,i}=x_i$ and $\mathcal A\in\mathfrak M_{\mathrm f}$. Let the particle height lie in a  fixed $\Lambda\Subset(0,\infty)$, and choose an integer $T>\max_i(s_i\vee t_i)$. Then $s_{N,i},t_{N,i}\le T$ for all sufficiently large $N$. Let $X_i^N(u):=U_{s_{N,i},u}^{\mathcal A,N}(x_i),  u\ge s_{N,i},$ with $X_i^N(u)=x_i$ for $u<s_{N,i}$. By     Lemma~\ref{lem:one-slit-moments}, \eqref{eq:CHL-root-flow-def}, and Appendix~\ref{appendix:SHL-boundary-cocycle}, these are square-integrable martingales with
    \begin{equation}\label{eq:HL-boundary-brackets}
        \begin{aligned}
        |\Delta X_i^N(u)|&\le C_\Lambda N^{-1},\quad  \langle X_i^N,X_j^N\rangle_t=\int_{s_{N,i}\vee s_{N,j}}^t  \mathcal K_h^{\mathcal A,N} \bigl(X_j^N(u-)-X_i^N(u-)\bigr)\,du ,
        \end{aligned}
    \end{equation}
    for $t\ge s_{N,i}\vee s_{N,j}$. The c\`adl\`ag version of the reversed process $u\mapsto R_{T-u,t_{N,i}}^{\mathcal A,N}(x_i)$, $T-t_{N,i}\le u\le T$, coincides with the original process at all grid times almost surely and satisfies \eqref{eq:HL-boundary-brackets} with $h$ replaced by $r$ and $s_{N,i}$ by $T-t_{N,i}$. For $\star\in\{h,r\}$, Lemma~\ref{lem:one-slit-moments} gives
    \[
        \mathcal K_\star^{\mathcal A,N}(0) =\sigma_\lambda^2+O_\Lambda(N^{-1}).
    \]
    Together with Lemma~\ref{lem:covariance-decay-section7}, this allows the martingale argument of Proposition~\ref{prop:abstract-fdd} to be applied to both flows. Arguing as in Corollary~\ref{cor:interpolated-fdd} then gives the limits of the interpolated paths. The estimates are uniform for $\lambda\in\Lambda$, so the same conclusion holds when $\lambda_N\to\lambda>0$. For $\mathcal A=\mathcal H$, apply Remark~\ref{rmk:radial-from-cylindrical} with $\lambda_N^{\mathcal H}\to\eta$ and the time change $t\mapsto t/(2\pi)$.

    It remains to replace $x_i$ by $X_{N,i}$. On $\max_{1\le i\le m}d_{E^{\mathcal A}}(X_{N,i},x_i)\le\varepsilon$, monotonicity gives
    \[
        \mathtt{HM}^{\mathcal A,N}_{s_{N,i},\cdot}(x_i-\varepsilon) \le \mathtt{HM}^{\mathcal A,N}_{s_{N,i},\cdot}(X_{N,i}) \le \mathtt{HM}^{\mathcal A,N}_{s_{N,i},\cdot}(x_i+\varepsilon).
    \]
    Before interpolation, the gap between the two bounding paths is a nonnegative martingale with initial value $2\varepsilon$. Linear interpolation preserves order and bounds each interpolated gap by the larger endpoint gap. Doob's inequality and a union bound yield, for $a>0$ (and $a<\pi$ in the circular case),
    \begin{equation*}
        \begin{aligned}
        \mathbb P\!\left\{ \max_{1\le i\le m} d_{\rm path}\!\left( \mathtt{HM}_{s_{N,i},\cdot}^{\mathcal A,N}(X_{N,i}), \mathtt{HM}_{s_{N,i},\cdot}^{\mathcal A,N}(x_i) \right)>a \right\}\le \mathbb P\!\left\{ \max_{1\le i\le m} d_{E^{\mathcal A}}(X_{N,i},x_i)>\varepsilon \right\} +\frac{2m\varepsilon}{a}.
        \end{aligned}
    \end{equation*}
    The same estimate holds with $\mathtt{HM}$ replaced by $\mathtt R$. Letting $N\to\infty$ and then $\varepsilon\downarrow0$, we obtain the result.
    \end{proof}

    \subsection{Small-gap estimates and verification of \texorpdfstring{\textup{(B2)}}{(B2)}}\label{S.three-point-estimate-I}  For the Root paths,  we follow the proof of Theorem~\ref{thm:abstract-boundary} in     Section~\ref{S.abstract-flow}: Assumption~\ref{ass:abstract-no-extra-path}\textup{(I)}  gives the three-point survival estimate, and reducing to a finite set of labels allows the  dyadic decomposition. For paths indexed by the arriving particles, instead of Assumption~\ref{ass:abstract-no-extra-path}\textup{(II)}, we use labels chosen at most one grid interval after the given time.

    Fix $T\in\mathbb N$, $\Lambda\Subset(0,\infty)$, and $\mathcal A\in\mathfrak M_{\mathrm f}$. Let $\Phi^N$ be the flow obtained from the time-reversed Poisson points using the Root maps, with points in $(s,t]$ included in $\Phi^N_{s,t}$. Thus, at fixed endpoints,
    \[
        \Phi^N_{s,t}(x)=R^{\mathcal A,N}_{T-t,T-s}(x), \qquad 0\le s\le t\le T.
    \]
    Let $(\mathcal F_t)$ be the right-continuous filtration generated by the time-reversed points. The rates in \eqref{eq:abstract-Q-covariance-identity}--\eqref{eq:abstract-Q-coalescence-rates} become
    \begin{equation*}
        \begin{aligned}
        \mathcal Q_N^{\mathcal A}(d) =2\{\mathcal K_r^{\mathcal A,N}(0)-\mathcal K_r^{\mathcal A,N}(d)\},\quad \Lambda_N^{\mathcal S}(d)=N^2|\{u\in\mathbb R: r_0^{\mathcal S}(u+Nd)=r_0^{\mathcal S}(u)\}|.
        \end{aligned}
    \end{equation*}
    For $\Lambda_N^{\mathcal C}$, replace $\mathbb R,r_0^{\mathcal S}$ by $\mathbb T_N,r_0^{\mathcal C,N}$.  Thus $\Lambda_N^{\mathcal A}(d)$ is the rate at which a single disturbance makes the ordered gap zero. Set
    \[
        j_N=C_\Lambda/N,\qquad a_N=N^{-1}\log N,\qquad r_N=N^{-1/2},\qquad \beta_N=a_N/r_N,
    \]
    where $C_\Lambda$ bounds the increments induced by a slit map, see \eqref{eq:bounded-one-slit-jumps} in  Lemma~\ref{lem:one-slit-moments}. By Proposition~\ref{prop:SHL-boundary-cocycle} and the Poisson construction in Appendix~\ref{appendix:SHL-boundary-cocycle}, the martingale identities of Section~\ref{S.abstract-flow} hold for $\Phi^N$, also conditionally after stopping times $\tau$, for finitely many $\mathcal F_\tau$-measurable labels. We first establish the uniform small-gap bound needed to verify Assumption~\ref{ass:abstract-no-extra-path}\textup{(I)}.

    \begin{proposition}\label{lem:one-gap-section6}
    There are $a_0\in(0,1/2)$ and $c>0$, depending only on $\Lambda$, such that, for all sufficiently large $N$, uniformly for $\lambda\in\Lambda$ and $\mathcal A\in\mathfrak M_{\mathrm f}$,
    \begin{equation}\label{eq:HL-gap-lower-bounds}
        \mathcal Q_N^{\mathcal A}(d)\ge c\quad(a_0/N\le d\le1), \qquad \Lambda_N^{\mathcal A}(d)\ge cN^2\quad(0<d\le2a_0/N).
    \end{equation}
    In particular, $\Phi^N$ satisfies Assumption~\ref{ass:abstract-no-extra-path}\textup{(I)} with the above scales.
    \end{proposition}
    \begin{proof}
    The one-slit maps $r_0^{\mathcal S}$ and $r_0^{\mathcal C,N}$ are constant on intervals of lengths $2\lambda$ and $4N\arcsin\delta_N(\lambda)\to2\lambda$, respectively; see \eqref{eq:one-slit-boundary-maps} and \eqref{eq:Rdelta-formula-section3}. Choose $a_0=\min\{(\inf\Lambda)/4,1/4\}$. For $d\le2a_0/N$, the set of $u$ for which both labels $u$ and $u+Nd$ lie in the same such interval has length bounded below uniformly, proving the second bound in \eqref{eq:HL-gap-lower-bounds}. For $d>0$, $\mathcal Q_N^{\mathcal S}(d)=\|\Delta_r^{\mathcal S}(\cdot+Nd)-\Delta_r^{\mathcal S}\|_2^2>0$, since a nonzero $L^2(\mathbb R)$ function cannot be periodic. Continuity and compactness then give, for every fixed $L>a_0$,
    \[
        \inf_{\substack{\lambda\in\Lambda,\ a_0\le Nd\le L}} \mathcal Q_N^{\mathcal S}(d)>0.
    \]
    Local uniform convergence $\Delta_r^{\mathcal C,N}\to\Delta_r^{\mathcal S}$, together with \eqref{eq:C-tail-section7}, gives the same lower bound for   $\mathcal C$ when $N$ is large. For $L\le Nd\le N$, the proof of Lemma~\ref{lem:covariance-decay-section7} gives  $|\mathcal K_r^{\mathcal A,N}(d)|\le C_\Lambda L^{-1/2}$, while $\inf_{\lambda\in\Lambda}\mathcal K_r^{\mathcal A,N}(0)>c$ for large $N$. Choosing $L$ large proves the first bound in \eqref{eq:HL-gap-lower-bounds}. For $0<d<a_0/N$, $j_Nr_N\Lambda_N^{\mathcal A}(d)\ge cN^{1/2}$; together with the first bound in \eqref{eq:HL-gap-lower-bounds}, this proves \eqref{eq:abstract-micro-coalescence}. Finally, $j_N=o(a_N)$ and $r_N=O(\sqrt{a_N})$.
    \end{proof}
    The following two-point estimate will be used to prove convergence of meeting times in Subsection~\ref{SS.binary-decomposition-section6}.
    \begin{lemma}\label{lem:HL-shrinking-cell-survival}
    There are $C,c>0$, depending only on $T,\Lambda$, such that, for all sufficiently large $N$, every stopping time $\tau$, $\mathcal F_\tau$-measurable $x,y\in E^{\mathcal A}$, and $h>0$ with $\tau+h\le T$,
    \begin{equation*}
        \mathbb P\!\left\{ \Phi^N_{\tau,\tau+h}(x)\ne\Phi^N_{\tau,\tau+h}(y) \,\middle|\,\mathcal F_\tau\right\} \le C\left(1\wedge\frac{d_{E^{\mathcal A}}(x,y)}{\sqrt h}\right) +Ce^{-cN^2h}.
    \end{equation*}
    \end{lemma}
    \begin{proof}
    By the strong Markov property and translation invariance, it suffices to consider deterministic labels at time $0$. Let $D_t=\Phi^N_{0,t}(d)-\Phi^N_{0,t}(0)$. Non-coalescence on the circle implies $D_h>0$. Use $\tau_R,F_{N,R}$ from the proof of Lemma~\ref{lem:abstract-one-gap-estimates}. For $C_\Lambda/N\le R\le R_0<1$, Proposition~\ref{lem:one-gap-section6} gives, with the same fixed $b>2$,
    \[
        \mathcal Q_N^{\mathcal A}(d)+(bj_NR+b^2j_N^2)\Lambda_N^{\mathcal A}(d)\ge c, \quad 0<d<R.
    \]
    The calculation \eqref{eq:abstract-LF-from-Q-Lambda}--\eqref{eq:abstract-one-gap-mean-exit} applies for every $C_\Lambda/N\le R\le R_0$ and yields
    \[
        \mathbb E_d[\tau_R\wedge h]\le CR(d+N^{-1}),\qquad \mathbb P_d\{\tau_R\le h,\ D_{\tau_R}\ge R\}\le d/R.
    \]
    For $0<d<R$, these bounds give
    \[
        \mathbb P_d\{D_h>0\} \le \mathbb P_d\{\tau_R>h\} +\mathbb P_d\{\tau_R\le h,\ D_{\tau_R}\ge R\} \le \frac{CR(d+N^{-1})}{h}+\frac dR.
    \]
    Taking $R$ of order $\sqrt h\wedge R_0$ gives
    \[
        \mathbb P_d\{D_h>0\}\le C(1\wedge d/\sqrt h), \qquad d\ge2a_0/N;
    \]
    if $\sqrt h<C_\Lambda/N$ or $d\ge R$, enlarge $C$.

    For $0<d<2a_0/N$, let $\sigma$ be the first exit time from $(0,2a_0/N)$. Before $\sigma$, \eqref{eq:HL-gap-lower-bounds} gives a coalescence rate at least $cN^2$ for $D$. Together with the optional stopping theorem for the nonnegative martingale $D$, this yields
    \[
        \mathbb P_d\{\sigma>h/2\}\le e^{-cN^2h},\qquad \mathbb P_d\{\sigma\le h/2,\ D_\sigma\ge2a_0/N\} \le\frac{Nd}{2a_0}.
    \]
    On $\{D_\sigma\ge2a_0/N\}$, we have $D_\sigma\le2a_0/N+2j_N$. Applying the preceding bound over the remaining time interval gives
    \[
        \mathbb P_d\{D_h>0\} \le e^{-cN^2h}+\frac{Nd}{2a_0} C\left(1\wedge\frac{1}{N\sqrt h}\right) \le e^{-cN^2h}+C\frac d{\sqrt h}.
    \]
    \end{proof}

    For each time-reversed point $(b_i,\bar x_i)$ with $0\le\bar t_i\le T$, set $b_i=T-\bar t_i\in[0,T]$, set $\sigma_i^N=N^{-1}\lceil Nb_i\rceil$, and let $\gamma_i^N$ be the linear interpolation of $\gamma_i^N(u)=\Phi^N_{b_i,u}(\bar x_i)$ for  $u\in\Gamma_N^{\mathcal A}\cap[\sigma_i^N,T]$. Let $K_N$ be the  closure of an arbitrary measurable subcollection of $\{\gamma_i^N\}_i$.  We use the count $\eta_K$ from \eqref{eq:eta-BW}--\eqref{eq:eta-cyl} to verify condition \textup{(B2)}.
    \begin{proposition}\label{prop:HL-particle-B2}
    For $0<h\le T$ and $0<\varepsilon\le1/2$,
    \begin{equation}\label{eq:HL-particle-B2}
        \limsup_{N\to\infty} \sup_{\substack{0\le t_0\le T-h \\
        |A|=\varepsilon}} \mathbb P\{\eta_{K_N}(t_0,h;A)\ge3\} \le C_h\varepsilon^2.
    \end{equation}
    Here $A$ ranges over intervals or arcs in $E^{\mathcal A}$, and the bound is uniform for $\lambda\in\Lambda$, $\mathcal A\in\mathfrak M_{\mathrm f}$, and $h$ in compact subsets     of $(0,T]$. Thus $K_N$ satisfies \textup{(B2)} on $[0,T]$.
    \end{proposition}
    \begin{proof}
    For all sufficiently large $N$, set
    \[
        r=N^{-1}\lfloor Nt_0\rfloor,\qquad s=r+N^{-1},\qquad u=N^{-1}\lfloor N(t_0+h)\rfloor, \qquad q_N=N^{-1/8},\quad M_N=\lceil N^4\rceil,
    \]
    and  $Y_N=\{\Phi^N_{b_i,s}(\bar x_i):b_i\le r\}\cap A^{2q_N}$. Then $Y_N$ is $\mathcal F_s$-measurable. It is the set of positions at time $s$ of the paths started by time $r$ that lie within $2q_N$ of $A$, and it plays the role of the label set $\mathcal R_n$ in Assumption~\ref{ass:abstract-no-extra-path}\textup{(II)}. We first show that $Y_N$ is not too large.

    For both geometries, the rescaled Poisson points have intensity $N^3\,dv\,dy$ on $[0,T]\times E^{\mathcal A}$, where $dy$ denotes length measure. Translation invariance gives $\int_{E^{\mathcal A}}\mathbb P\{\Phi^N_{v,s}(y)\in B\}\,dy=|B|$ for every bounded interval or arc $B$. Since $\Phi^N_{v,s}$ uses only points in $(v,s]$, the Mecke identity \cite[Theorem~4.1]{last2017lectures} gives
    \[
        \begin{aligned}
        \mathbb E[\#Y_N]\le\mathbb E\!\left[\sum_{i:b_i\le r}\ind_{\{\Phi^N_{b_i,s}(\bar x_i)\in A^{2q_N}\}}\right]=N^3\int_0^r\int_{E^{\mathcal A}}\mathbb P\{\Phi^N_{v,s}(y)\in A^{2q_N}\}\,dy\,dv=N^3r|A^{2q_N}|.
        \end{aligned}
    \]
    Thus, since $r\le T$ and $|A^{2q_N}|=\varepsilon+4q_N=O(1)$,
    \begin{equation}\label{eq:HL-particle-label-size}
        \mathbb E[\#Y_N]\le C_TN^3,\qquad\mathbb P\{\#Y_N>M_N\}\le\frac{C_TN^3}{M_N}\le\frac{C_T}{N}\longrightarrow0.
    \end{equation}
    By \eqref{eq:clock-cell-fourth-moment} and the mesh argument of  Lemma~\ref{lem:automatic-clock-oscillation}, the probability that, on some grid interval in $[0,T]$, a path starting from a label in $A^1$ moves by more than $q_N/4$ is at most
    \[
        CNq_N^{-1}q_N^{-4}(N^{-2}+N^{-3}) \le CN^{-3/8}.
    \]
    The estimate is uniform in the location of $A$. If this event does not occur, every path $\gamma_i^N$ with $\sigma_i^N\le t_0$ and $\gamma_i^N(t_0)\in A^{q_N/2}$ satisfies $b_i\le r$ and $\gamma_i^N(s)\in Y_N$. The same holds for every path in $K_N$. Indeed, the starting times belong to $N^{-1}\mathbb Z$, and $Y_N$ is finite by \eqref{eq:HL-particle-label-size}. Hence, from any convergent sequence of generating paths, there exists a subsequence with the same starting time and the same value in $Y_N$ at time $s$. Finally, $\gamma_1(u)=\gamma_2(u)$ implies $ \gamma_1(t_0+h)=\gamma_2(t_0+h)$,  since the two paths remain equal at the next grid time and are linearly interpolated. Hence
    \begin{equation}\label{eq:HL-particle-count-reduction}
        \mathbb P\{\eta_{K_N}(t_0,h;A)\ge3\} \le \mathbb P\{\#\Phi^N_{s,u}(Y_N)\ge3\}+CN^{-3/8}.
    \end{equation}

    Conditionally on $\mathcal F_s$, apply Proposition~\ref{prop:abstract-binary-reduction} to $Y_N$, using Proposition~\ref{lem:one-gap-section6} and the independence of the Poisson points after $s$.  On $\{\#Y_N\le M_N\}$,
    \[
        \operatorname{diam}(Y_N)\le\varepsilon+4q_N, \qquad h-2N^{-1}\le u-s\le h.
    \]
    Therefore, for $h$ in a fixed compact subset of $(0,T]$,     \eqref{eq:HL-particle-label-size}--\eqref{eq:HL-particle-count-reduction}  and Proposition~\ref{prop:abstract-binary-reduction} give
    \[
        \begin{aligned}
        \mathbb P\{\eta_{K_N}(t_0,h;A)\ge3\}\le o(1)+C_h\left\{(\varepsilon+4q_N)^2+\beta_N(\varepsilon+4q_N)\log(eM_N)+M_N\frac{j_N}{r_N}\text{e}^{-\tfrac{ch}{r_N^2}} \right\}.
        \end{aligned}
    \]
    The error is uniform for $t_0$, $A$, $\lambda\in\Lambda$, and $h$ in compact subsets of $(0,T]$. Since $q_N\to0$,
    \[
        \beta_N\log(eM_N)=O\!\left(N^{-1/2}(\log N)^2\right)\longrightarrow0,\qquad M_N\frac{j_N}{r_N}e^{-ch/r_N^2}=O\!\left(N^{7/2}e^{-chN}\right)\longrightarrow0,
    \]
    \eqref{eq:HL-particle-B2} follows.
    \end{proof}

    \subsection{\texorpdfstring{Proofs of Theorems~\ref{thm:harmonic-measure-webs} and \ref{thm:joint-geometric-webs}}{Proofs of Theorems}}\label{SS.binary-decomposition-section6}

    We now use the preceding estimates to prove Theorems~\ref{thm:harmonic-measure-webs} and \ref{thm:joint-geometric-webs}. Corollary~\ref{cor:backbone-dual-backbone} follows at the end of this subsection. Fix $\mathcal A\in\mathfrak M$ and use the paths and grid times introduced in Subsection~\ref{S.HL-models}. Set
    \[
        p_i=\min\bigl(\Gamma_N^{\mathcal A}\cap[\bar t_i,\infty)\bigr), \qquad q_i=\max\bigl(\Gamma_N^{\mathcal A}\cap[0,\bar t_i]\bigr),
    \]
    that is, the first grid time after and the last grid time before the arrival of particle $i$, and define the corresponding restricted paths by
    \[
        \widehat{\mathtt{HM}}_{i,\pm}^{\mathcal A,N} :=\mathtt{HM}_{\bar t_i,\cdot}^{\mathcal A,N}(\bar x_{i,\pm}) \big|_{[p_i,\infty)}, \qquad \widehat{\mathtt R}_i^{\mathcal A,N} :=\mathtt R_{\cdot,\bar t_i}^{\mathcal A,N}(\bar x_i) \big|_{[0,q_i]}.
    \]
    In both cases the arrival time is $\bar t_i$, and $0\le p_i-\bar t_i<\Delta_N^{\mathcal A}$, $0\le\bar t_i-q_i<\Delta_N^{\mathcal A}$. Set
    \[
        \mathcal O_N[a,b]:=\sup_{a\le s\le t\le b}
        \begin{cases}
        \displaystyle\sup_{x\in\mathbb R} |\tanh R^{\mathcal A,N}_{s,t}(x)-\tanh x|,&\mathcal A=\mathcal S, \\[3pt]
        \displaystyle\sup_{x\in[0,2\pi]} |R^{\mathcal A,N}_{s,t}(x)-x|,&\mathcal A\in\mathfrak M_{\rm c}.
        \end{cases}
    \]

    The following lemma controls $\mathcal O_N$ and the error from the above path restrictions. Its conclusions hold uniformly in $\lambda\in\Lambda$ for $\mathcal A\in\mathfrak M_{\rm f}$ and in $\eta\in\Lambda$ for $\mathcal A=\mathcal H$, for every $\Lambda\Subset(0,\infty)$.
    \begin{lemma}\label{lem:particle-preparation}
    Fix $\Lambda\Subset(0,\infty)$. For every $G<\infty$,
    \begin{equation}\label{eq:actual-partial-oscillation}
        \max_{k:\,I_{k,N}^{\mathcal A}\cap[0,G]\ne\varnothing} \mathcal O_N[I_{k,N}^{\mathcal A}] \pto0.
    \end{equation}
    Moreover,
    \begin{equation}\label{eq:particle-path-restriction}
        \sup_{i:\,\bar t_i\le G} \left\{ \max_{\pm}d_{\rm path}\!\left( \mathtt{HM}_{\bar t_i,\cdot}^{\mathcal A,N}(\bar x_{i,\pm}), \widehat{\mathtt{HM}}_{i,\pm}^{\mathcal A,N} \right) \vee d_{\rm path}\!\left( \mathtt R_{\cdot,\bar t_i}^{\mathcal A,N}(\bar x_i), \widehat{\mathtt R}_i^{\mathcal A,N} \right) \right\} \pto0.
    \end{equation}
    For every $G,T<\infty$, the four path collections with $\bar t_i\le G$, restricted to $[0,T]$, are measurable random compact sets and are jointly tight.  The restricted harmonic measure and Root paths are non-crossing on the intersection of their time domains.
    \end{lemma}

    \begin{proof}
    First take $\mathcal A\in\mathfrak M_{\rm f}$. Equation~\eqref{eq:clock-cell-fourth-moment} and Lemma~\ref{lem:one-slit-moments}, applied after time reversal, give
    \[
        \sup_z\E\left[\sup_{a\le s\le b}|R^{\mathcal A,N}_{s,b}(z)-z|^4\right] \le C\bigl((b-a)^2+N^{-2}(b-a)\bigr).
    \]
    For $a\le s\le t\le b$, the map $R^{\mathcal A,N}_{t,b}$ is onto, so for every $y$ we may choose $z$ with $R^{\mathcal A,N}_{t,b}(z)=y$. Proposition~\ref{prop:SHL-boundary-cocycle} then gives $R^{\mathcal A,N}_{s,t}(y)=R^{\mathcal A,N}_{s,b}(z)$. Thus,
    \[
        \mathcal O_N[a,b] \le 2\sup_{a\le s\le b}
        \begin{cases}
        \displaystyle \sup_{x\in\mathbb R} |\tanh R^{\mathcal A,N}_{s,b}(x)-\tanh x|, & \mathcal A=\mathcal S, \\[3pt]
        \displaystyle \sup_{x\in[0,2\pi]} |R^{\mathcal A,N}_{s,b}(x)-x|, & \mathcal A\in\mathfrak M_{\rm c}.
        \end{cases}
    \]
    The mesh argument of Lemma~\ref{lem:automatic-clock-oscillation} gives, for each grid interval $I_{k,N}^{\mathcal A}$, $\P\{\mathcal O_N[I_{k,N}^{\mathcal A}]>\varepsilon\} \le C_{\varepsilon,\Lambda}N^{-2}$. A union bound over $O_G(N)$ intervals proves \eqref{eq:actual-partial-oscillation}. At grid and arrival times, \eqref{eq:SHL-grid-inverse} gives     $R^{\mathcal A,N}_{s,t}(U^{\mathcal A,N}_{s,t}(x-))=R^{\mathcal A,N}_{s,t}(U^{\mathcal A,N}_{s,t}(x))=x$, so the same bound controls both harmonic measure paths. The original and restricted paths are the same except on $[\bar t_i,p_i]$ and $[q_i,\bar t_i]$. Thus the supremum in \eqref{eq:particle-path-restriction} is bounded by
    \[
        C\left(\Delta_N^{\mathcal A} +\max_{k:\,I_{k,N}^{\mathcal A}\cap[0,G]\ne\varnothing} \mathcal O_N[I_{k,N}^{\mathcal A}]\right)\pto0.
    \]
    The same pairing gives convergence after closure. For fixed $N$, Proposition~\ref{prop:SHL-boundary-cocycle} gives compactness: on every bounded time interval and for every $R<\infty$, only finitely many particles satisfy $|\bar x_i|\le R$, while paths with $\bar x_i\to\pm\infty$ converge uniformly to $\pm\infty$ in the compactified coordinate. On the circle, the number of particles is finite. Since the collections are closures of countable measurable paths, they are measurable.

    The cocycle gives, on the respective domains,
    \[
        \begin{aligned}
        \widehat{\mathtt{HM}}_{i,\pm}^{\mathcal A,N}(u)     =\mathtt{HM}_{p_i,u}^{\mathcal A,N}\bigl(U^{\mathcal A,N}_{\bar t_i,p_i}(\bar x_{i,\pm})\bigr),\qquad   \widehat{\mathtt R}_i^{\mathcal A,N}(u) =\mathtt R_{u,q_i}^{\mathcal A,N}\bigl(R^{\mathcal A,N}_{q_i,\bar t_i}(\bar x_i)\bigr).
        \end{aligned}
    \]
    These are subcollections of the full grid-time trajectories of $U^{\mathcal A,N}$ and the time-reversed $R^{\mathcal A,N}$. By Lemmas~\ref{lem:one-slit-moments} and~\ref{lem:covariance-decay-section7}, these flows satisfy Definition~\ref{def:abstract-localized}. Therefore, Lemma~\ref{lem:abstract-tightness}  gives tightness on bounded windows; bounded arrival times and \eqref{eq:particle-path-restriction} give joint tightness.  Monotonicity preserves spatial order within each collection. At grid times, \eqref{eq:SHL-grid-inverse} also gives
    \[
        x<R^{\mathcal A,N}_{s,t}(y)\ \Longrightarrow\ U^{\mathcal A,N}_{s,t}(x)<y, \qquad x>R^{\mathcal A,N}_{s,t}(y)\ \Longrightarrow\ U^{\mathcal A,N}_{s,t}(x-)>y.
    \]
    Since the restricted paths use the same grid, linear interpolation preserves their spatial order. The case $\mathcal A=\mathcal H$ follows from Remark~\ref{rmk:radial-from-cylindrical} and $\lambda_N^{\mathcal H}\to\eta$.
    \end{proof}

    For the following proofs, let $J=[s_0,s_1]$ and let $D$ be the dense set fixed at the beginning of Section~\ref{S.HL-observables-BW}. For each $(x,t)\in D$, choose a grid interval $I_{k,N}^{\mathcal A}\subset J$ whose endpoints converge to $t$, taking the interval inside $J$ at $s_0,s_1$, and choose the particle with $\bar x_i$ nearest to $x$ among those arriving in the middle half of $I_{k,N}^{\mathcal A}$. Let $E_N$ be the event that every such middle half contains a particle; then $\P(E_N^c)\to0$.

    At the two grid endpoints, use $ U^{\mathcal A,N}_{\bar t_i,p_i}(\bar x_{i,\pm}), R^{\mathcal A,N}_{q_i,\bar t_i}(\bar x_i). $ By \eqref{eq:particle-side-labels}, Lemma~\ref{lem:particle-preparation}, and $\bar x_i\to x$, these labels converge in probability to $x$. Hence Proposition~\ref{prop:fdd-selected-paths-section7} gives the coalescing Brownian limits for the corresponding harmonic measure and Root paths. On $E_N$, these are the restricted particle paths defined above, and \eqref{eq:particle-path-restriction} gives the same limits for the original paths. The choice inside $J$ also covers $t=s_0,s_1$.

    For the Root collections, let $d_H^{\mathcal A}$ denote the Hausdorff metric on boundary paths with their arrival times, as defined in Subsection~\ref{SS.geometric-tightness}.  For $\mathcal A\in\mathfrak M_{\rm c}$, define
    \[
        d_{\mathrm{loc}}^{\mathcal A}(K,L) := \sum_{m\ge1}2^{-m} \bigl(1\wedge d_H^{\mathcal A}(K\langle m\rangle,L\langle m\rangle)\bigr).
    \]
    For $\mathcal A=\mathcal S$, let $d_{H,k}^{\mathcal S}$ be the Hausdorff distance induced by the following pseudometric on boundary paths with their arrival times  $\bigl(1\wedge d_k^{\mathcal S,k}(\mathfrak e_{\mathcal S}^{v_{\mathcal S}}(\gamma), \mathfrak e_{\mathcal S}^{v_{\mathcal S}}(\gamma'))\bigr) \vee|\tanh r-\tanh r'|$,  with $\dagger$ at distance $1$, and set
    \[
        d_{\mathrm{loc}}^{\mathcal S}(K,L) := \sum_{m,k\ge1}2^{-(m+k)} \bigl(1\wedge d_{H,k}^{\mathcal S}(K\langle m\rangle,L\langle m\rangle)\bigr).
    \]

    For the limiting backward BW, each path is paired with its terminal time, with $\dagger$ included in each restriction.

    \begin{proposition}\label{prop:process-root-webs}
    For every $\mathcal A\in\mathfrak M$,  $\mathbf R^{\mathcal A,N}\dto\sigma_{\mathcal A}\mathcal W^\downarrow_{\mathcal A,\ge0}$ in the topology $d_{\mathrm{loc}}^{\mathcal A}$ defined above.
    \end{proposition}

    \begin{proof}
    First take $\mathcal A\in\mathfrak M_{\rm f}$. Fix $G>0$ and an integer $T>G$, and use the reversed flow $\Phi^N$ of     Subsection~\ref{S.three-point-estimate-I}. For $b_i=T-\bar t_i$, at grid  times,
    \[
        T-q_i=\frac{\lceil Nb_i\rceil}{N},\qquad  \widehat{\mathtt R}_i^{\mathcal A,N}(T-u)=R^{\mathcal A,N}_{T-u,\bar t_i}(\bar x_i)  =\Phi^N_{b_i,u}(\bar x_i).
    \]
    Since the slit map at $\bar t_i$ fixes $\bar x_i$, after time reversal the restricted Root paths with $\bar t_i\le G$ are precisely the paths in  Proposition~\ref{prop:HL-particle-B2} with $b_i\ge T-G$, restricted to $[T-G,T]$. Proposition~\ref{prop:HL-particle-B2} gives (B2), Lemma~\ref{lem:particle-preparation} gives tightness and non-crossing, and the paths selected above give (I). Since $\P(E_N^c)\to0$, we may define the collection arbitrarily on $E_N^c$; the criterion at the beginning of Section~\ref{S.HL-observables-BW} then gives the backward BW on $[0,G]$.  Since  $ 0\le \bar t_i-q_i<N^{-1} $,  \eqref{eq:particle-path-restriction} allows us to replace the restricted Root paths by the original ones. Thus $\mathbf R^{\mathcal A,N}\langle G\rangle$ converges to the backward BW restricted to terminal times at most $G$.

    Fix $M\in\mathbb N$. By Lemma~\ref{lem:particle-preparation}, $\{\mathbf R^{\mathcal A,N}\langle m\rangle:1\le m\le M\}$ is jointly tight. Along any convergent subsequence, denote the limit by $(K^1,\ldots,K^M)$. For $m\le M$,  $ K^m\subset\{(\pi,t)\in K^M:t\le m\}\cup\{\dagger\}.$ For each fixed $m$, both sets have the law of the backward BW restricted to terminal times at most $m$, and hence are equal almost surely.  Thus $K^1,\ldots,K^M$ are restrictions of the same backward BW. The tails of the defining sums for $d_{\mathrm{loc}}^{\mathcal A}$ tend to zero uniformly, yielding the required convergence. The case $\mathcal A=\mathcal H$  follows from Remark~\ref{rmk:radial-from-cylindrical} and  $\lambda_N^{\mathcal H}\to\eta$.
    \end{proof}

    Now we use the convergence of $\mathbf R^{\mathcal A,N}$ and non-crossing duality to prove BW convergence of $\mathbf{HM}^{\mathcal A,N}$.

    \begin{proof}[Proof of Theorem~\ref{thm:harmonic-measure-webs}]
    First take $\mathcal A\in\mathfrak M_{\rm f}$. Let    $K_N=\mathbf{HM}^{\mathcal A,N}$, and let $K_N^G$ be the closure of its generating paths with $\bar t_i\le G$. Let $\mathcal I$ be as in Subsection~\ref{SS.convergence criteria}. For every $G>0$, the side labels of particles with $\bar t_i>G$ are dense in every spatial open set, so $\mathcal I\subset K_N$. Hence, for $G>0$,
    \[
        d_H(K_N,K_N^G\cup\mathcal I)\le \frac{2}{1+G}.
    \]
    For $T>G$, if each path in $K_N^G$ is kept fixed after time $T$, then its distance from the original path is at most $2/(1+T)$, and the same bound holds for the Hausdorff distance. Therefore, Lemma~\ref{lem:particle-preparation} gives measurable compact approximations, tight for fixed $G,T$. Letting first $T\to\infty$ and then $G\to\infty$ yields compactness,     measurability and tightness of $K_N$.

    Use the harmonic measure paths selected above, with spatially dense starting labels at times tending to infinity. Fix a countable dense set $D_0\subset E^{\mathcal A}\times[0,\infty)$. For $w=(x,t)\in D_0$, set
    \[
        t_N=\max(\Gamma_N^{\mathcal A}\cap[0,t]),  \qquad \pi_N^w=\mathtt R^{\mathcal A,N}_{\cdot,t_N}(x).
    \]
    Proposition~\ref{prop:fdd-selected-paths-section7} gives the required  finite-dimensional coalescing Brownian limits for both selected  collections.

    Fix an integer $T$ larger than the terminal times of a finite backward collection. For any two selected backward paths, in reversed time $v=T-u$, let $Y_N(v)$ be their distance in $E^{\mathcal A}$ and set $\tau_N:=\inf\{v:Y_N(v)=0\}\wedge T$. Let $\tau$ be the corresponding limiting meeting time. At a grid time $v$, Lemma~\ref{lem:HL-shrinking-cell-survival} gives, for  $v+h+\Delta_N^{\mathcal A}\le T$,
    \[
        \P\{Y_N(v)\le\gamma,\ \tau_N>v+h+\Delta_N^{\mathcal A}\} \le C(1\wedge\gamma/\sqrt h)+Ce^{-cN^2h}.
    \]
    On a coupling with uniform path convergence, positivity of the limiting distance before $\tau$ excludes $\tau_N<\tau-\varepsilon$. If $\tau<T$, choose a deterministic mesh in $\Gamma_N^{\mathcal A}\cap[0,T]$ with spacing less than $\varepsilon/4$ and cardinality bounded in $N$. On $\{\tau\le T-\varepsilon\}$, take a mesh point $v\in[\tau,\tau+\varepsilon/4]$. Uniform path convergence gives $Y_N(v)\le\gamma$ with probability tending to one as $N\to\infty$ and then $\gamma\downarrow0$. Applying the above estimate with $h=\varepsilon/2$ at the finite mesh points gives $\P(\tau_N>\tau+\varepsilon)\rightarrow0$.  If $\tau>T-\varepsilon$, the same conclusion follows from $\tau_N\le T$. Hence $\tau_N\pto\tau$, jointly for finite collections. For wedges on $\T$, use the arc whose gap vanishes at the meeting, as in Subsection~\ref{SS.convergence criteria}.

    The proof of Lemma~\ref{lem:particle-preparation} gives non-crossing between the restricted harmonic measure paths and the paths $\pi_N^w$. If two  backward paths meet at a grid time $s$ with value $x$, and have values $y_1\le y_2$ at a later grid time $u$, then the cocycle property and \eqref{eq:SHL-grid-inverse} give
    \[
        R^{\mathcal A,N}_{s,u}(y_1)=R^{\mathcal A,N}_{s,u}(y_2)=x, \qquad U^{\mathcal A,N}_{s,u}(x-)\le y_1\le y_2\le U^{\mathcal A,N}_{s,u}(x).
    \]
    For a restricted generating forward path $\pi$ active at $s$, non-crossing excludes entry when $\pi(s)\ne x$, while, if $\pi(s)=x$,
    \[
        \pi(u)\in \{U^{\mathcal A,N}_{s,u}(x-),U^{\mathcal A,N}_{s,u}(x)\} \subset\R\setminus(y_1,y_2).
    \]
    Together with the linear interpolation, this excludes entry into the wedge, also for paths starting after $s$. The same holds after closure.

    Take a joint subsequential limit $K$ of $K_N$ and the selected paths,  together with their meeting times. By \eqref{eq:particle-path-restriction}, every path of $K$ with finite starting time is also a limit of restricted harmonic measure paths. Their non-crossing and wedge exclusion pass to the limit by convergence of the backward paths and their meeting times. Therefore, the identification argument of Lemma~\ref{lem:half-axis-dual-wedge-criterion} applies. Let $\mathcal W_0$ be the closure of the selected limiting forward paths. By Theorem~\ref{thm:BW-CBW}, the paths in $\mathcal W_0$ with finite starting times form the corresponding forward BW, and $K\setminus\mathcal I\subset\mathcal W_0\subset K$. On $\R$, apply the comparison on bounded spatial intervals as in the proof of Lemma~\ref{lem:abstract-tightness}. Since the selected starting labels are spatially dense at times tending to infinity, $\mathcal I\subset\mathcal W_0$. Therefore $K=\mathcal W_0\stackrel{d}{=}\sigma_{\mathcal A}\mathcal W^\uparrow_{\mathcal A,\ge0}$. Thus every subsequential limit has the required law.

    For $\mathcal A=\mathcal H$, we use Remark~\ref{rmk:radial-from-cylindrical} and $\lambda_N^{\mathcal H}\to\eta$.
    \end{proof}

    The preceding boundary limits also hold jointly.   Lemma~\ref{lem:particle-preparation} gives joint tightness on every bounded arrival-time interval. Its non-crossing property, together with the vanishing restriction errors in \eqref{eq:particle-path-restriction}, identifies each subsequential limit by    Theorem~\ref{thm:double BWs}. Proposition~\ref{prop:process-root-webs} and the tail estimate in the proof of    Theorem~\ref{thm:harmonic-measure-webs} then yield
    \begin{equation}\label{eq:joint-root-HM-webs}
        (\mathbf R^{\mathcal A,N},\mathbf{HM}^{\mathcal A,N}) \dto (\sigma_{\mathcal A}\mathcal W_{\mathcal A,\ge0}^{\downarrow}, \sigma_{\mathcal A}\mathcal W_{\mathcal A,\ge0}^{\uparrow}),
    \end{equation}
    where convergence in the first coordinate is with respect to $d_{\mathrm{loc}}^{\mathcal A}$, and in the second with respect to the Hausdorff metric induced by \eqref{path-distance}, with the standard non-crossing coupling. We next show that $\mathbf{Gap}^{\mathcal A,N}$ is approximated by the embedded $\mathbf{HM}^{\mathcal A,N}$.

    \begin{proposition}\label{prop:gap-HM-comparison}
    For every $\mathcal A\in\mathfrak M$,
    \begin{equation}\label{eq:gap-HM-comparison}
        D_{\mathcal A,H}^{\uparrow} \bigl(\mathbf{Gap}^{\mathcal A,N}, \mathfrak e_{\mathcal A}^{v_{\mathcal A}} (\mathbf{HM}^{\mathcal A,N})\bigr)\pto0.
    \end{equation}
    \end{proposition}
    \begin{proof}
    Fix $T<\infty$ and a bounded spatial window. For $\mathcal A=\mathcal S$, \eqref{eq:gap-spatial-localization} allows us, after choosing $L<\infty$ large, to restrict to gaps contained in $[-L,L]$ up to time $T$, up to an error vanishing as $L\to\infty$; for the compact cases take $L=2\pi$.  After this localization, \eqref{eq:gap-height-approximation} holds uniformly over all relevant gaps up to time $T$.

    Let $0<\eta<1$ and set $\delta=\eta^3$ (in this proof, $\eta$ is an auxiliary small parameter, not the particle size of HL$(0)$). For each interval $[k\delta,(k+1)\delta]$, choose particles in the middle half of a grid interval whose endpoints converge to $(k+2)\delta$ near an $\eta$-mesh of $[-L,L]$, as before Proposition~\ref{prop:process-root-webs}. By Propositions~\ref{prop:fdd-selected-paths-section7} and \ref{prop:process-tip-geometric-transfer}, outside an event of probability at most
    \[
        C(1+T/\eta^3)(1+L/\eta)e^{-c/\eta}+o_N(1),
    \]
    the selected fingers remain within spatial distance $\eta/4$ of the corresponding mesh points on  $ [(k-1)\delta,(k+2)\delta]\cap[0,\infty)$,  and cross $v_{\mathcal A}(k\delta-\delta/2)_+\le h\le v_{\mathcal A}(k\delta+3\delta/2)$. By \eqref{eq:gap-height-approximation}, every relevant gap with $\sigma_i\le(k+1)\delta$ remains between these two heights on $[k\delta,(k+1)\delta]\cap[\sigma_i,\infty)$, outside an additional event of probability $o_N(1)$. Since \eqref{eq:gap-strict-exterior} prevents it from crossing the selected fingers,
    \[
        \sup_{ u,v\in[k\delta,(k+1)\delta]\cap[\sigma_i,\infty) } d_{E^{\mathcal A}} \!\left( \chi_{\mathcal A}^{(1)} \bigl(\mathtt{Gap}_i^{\mathcal A,N}(u)\bigr), \chi_{\mathcal A}^{(1)} \bigl(\mathtt{Gap}_i^{\mathcal A,N}(v)\bigr) \right) \le 4\eta .
    \]

    Letting $N\to\infty$, then $\eta\downarrow0$, and, for     $\mathcal A=\mathcal S$, finally $L\to\infty$, gives tightness on every bounded window, hence tightness in $D_{\mathcal A,H}^{\uparrow}$ by the forward version of \eqref{eq:geometric-path-metric-comparison}.

    We use the same selection of particles for a countable dense set of terminal points; this gives selected fingers. The meeting-time argument in the proof of    Theorem~\ref{thm:harmonic-measure-webs} applies jointly to every finite subcollection. Suppose two selected fingers separate at a point $b_N$ above the initial boundary, with $b_N$ in the interior of the earlier particle $j$ where the separating branch is attached. If the corresponding particles are $i_1,i_2$, then Proposition~\ref{prop:SHL-boundary-cocycle} and Lemma~\ref{lem:root-finger-identity} give
    \[
        R_{\bar t_j,\bar t_{i_1}}^{\mathcal A,N}(\bar x_{i_1}) \neq R_{\bar t_j,\bar t_{i_2}}^{\mathcal A,N}(\bar x_{i_2}), \qquad R_{\bar t_j-,\bar t_{i_1}}^{\mathcal A,N}(\bar x_{i_1}) = R_{\bar t_j-,\bar t_{i_2}}^{\mathcal A,N}(\bar x_{i_2}) = \bar x_j .
    \]
    Hence, if $\tau_N$ is the last grid time at which the two Root paths have the same value, $ |\tau_N-\bar t_j|\le \Delta_N^{\mathcal A}$. The whole-particle estimate in the proof of Proposition~\ref{prop:process-tip-geometric-transfer} gives
    \[
        \left| v_{\mathcal A}^{-1} \chi_{\mathcal A}^{(2)}(b_N)-\tau_N \right| \le \left| v_{\mathcal A}^{-1} \chi_{\mathcal A}^{(2)}(b_N)-\bar t_j \right| +\Delta_N^{\mathcal A} =o_{\mathbb P}(1).
    \]
    By tightness and \eqref{eq:joint-root-HM-webs}, along a further subsequence the gaps, the selected Root paths, and their meeting times may be coupled to converge jointly. Let $K$ denote the limiting gap collection. Equations~\eqref{eq:gap-height-approximation} and     \eqref{eq:gap-spatial-localization} imply that every path in $K$ with finite starting time and position is of the form $u\mapsto \mathfrak e_{\mathcal A}^{v_{\mathcal A}}(u,\pi(u))$.

    We claim that such a path cannot enter a wedge determined by two selected limiting  backward paths. Suppose, for some $\sigma_\pi\le s<t<r$, $(\pi(s),s)\notin\overline V, (\pi(t),t)\in V$,  where $r$ is the upper time of the wedge $V$. Choose $r_0\in(t,r)$.  Starting from their separation point, or from the initial boundary if the fingers are disjoint, follow the two finger paths towards their respective tips and stop each at its first hit of height $v_{\mathcal A}r_0$. Join the endpoints by the horizontal segment at this height, including the initial-boundary segment when needed. The resulting curve is simple.

    By convergence of the fingers and the preceding estimate for their  separation point, the approximating gap is outside this curve at time $s$ and inside it at time $t$. This is impossible: by \eqref{eq:gap-strict-exterior} the gap crosses neither finger nor the initial boundary, while for $u\le t$,
    \[
        \chi_{\mathcal A}^{(2)} \bigl(\mathtt{Gap}_i^{\mathcal A,N}(u)\bigr) \le v_{\mathcal A}u \le v_{\mathcal A}t < v_{\mathcal A}r_0,
    \]
    so it cannot cross the upper line segment before time $t$.  The wedge characterization \cite[Remark~6.3.10]{schertzer2017brownian} used in the proof of Lemma~\ref{lem:half-axis-dual-wedge-criterion}  gives $\pi\in\sigma_{\mathcal A}\mathcal W_{\mathcal A,\ge0}^{\uparrow}$,  with the coupling determined by \eqref{eq:joint-root-HM-webs}.  Conversely, let $(x,s)$ range over a countable dense set and let $i_N$ be the selected particle. Then  $(\bar x_{i_N},\bar t_{i_N})\pto(x,s)$, and Proposition~\ref{prop:process-tip-geometric-transfer}, applied to the tip, gives
    \[
        d_{E^{\mathcal A}}(\xi_{i_N},\bar x_{i_N}) +|\sigma_{i_N}-\bar t_{i_N}|\pto0.
    \]
    Hence every subsequential limit contains the embedded forward BW path  from $(x,s)$. Taking the closure over the dense set yields  $ \mathfrak E_{v_{\mathcal A}}^{\mathcal A,\uparrow}(      \sigma_{\mathcal A}\mathcal W_{\mathcal A,\ge0}^{\uparrow})\subset K$.  The remaining compactification limits agree on both sides: if $\sigma_i\to\infty$, both converge to the corresponding limits at infinite starting time; for $\mathcal A=\mathcal S$, bounded $\sigma_i$ and $|\xi_i|\to\infty$ are handled by    \eqref{eq:gap-spatial-localization}. Therefore
    \[
        K= \mathfrak E_{v_{\mathcal A}}^{\mathcal A,\uparrow} \bigl( \sigma_{\mathcal A}\mathcal W_{\mathcal A,\ge0}^{\uparrow} \bigr) \qquad\text{a.s.}
    \]

    By \eqref{eq:joint-root-HM-webs} and continuity of the embedding, every subsequential limit of $(      \mathbf{Gap}^{\mathcal A,N}, \mathfrak e_{\mathcal A}^{v_{\mathcal A}}(\mathbf{HM}^{\mathcal A,N}))$ is $(K,K)$. This proves \eqref{eq:gap-HM-comparison}.
    \end{proof}

    We next approximate the full $\mathbf{Fin}^{\mathcal A,N}$  by the subcollection of particles arriving before a finite time.
    \begin{proposition}\label{prop:process-root-finger-saturation}
    The collections $\mathbf{Fin}^{\mathcal A,N}$ and $\mathbf{Fin}^{\mathcal A,N}\langle G\rangle$, $G<\infty$, are measurable random compact sets in their respective path spaces. For every $\varepsilon>0$,
    \begin{equation}\label{eq:physical-cutoff-saturation}
        \lim_{G\to\infty}\limsup_{N\to\infty} \P\!\left\{ D_{\mathcal A,H}\bigl(\mathbf{Fin}^{\mathcal A,N}, \mathbf{Fin}^{\mathcal A,N}\langle G\rangle\bigr)> \varepsilon\right\}=0.
    \end{equation}
    For $\mathcal A=\mathcal C,\mathcal H$, the convergence is uniform in $\lambda\in\Lambda,\eta\in\Lambda$, respectively, for every $\Lambda\Subset(0,\infty)$.
    \end{proposition}

    \begin{proof}
    Set $ U_{N,\lambda}:=\inf\{t_j^{\mathcal C}:t_j^{\mathcal C}>\upsilon_{N,\lambda}^{\mathcal C}\}$, the first arrival time after the one-arm domination time. The upper bound in \cite[Theorem~1.4]{chen2025one}, together with the scaling identity $\lambda U_{N,\lambda}\stackrel{\rm d}=U_{N/\lambda,1}$, yields, for every $\Lambda\Subset(0,\infty)$ and all sufficiently large $N$,
    \begin{equation}\label{eq:hl-observables-cpz-tail}
        \sup_{\lambda\in\Lambda} \P\{N^{-2}U_{N,\lambda}>u\} \le C_\Lambda e^{-c_\Lambda u},\qquad u>0.
    \end{equation}

    Fix $T<\infty$ and $M\in\Gamma_N^{\mathcal C}$ with $M>T$.  Restart at time $N^2M$ via $(\mathcal C_{N^2M}^{N,\lambda})^{-1}$, and let $i_*$ be the index in the original process of the particle arriving immediately after the last attachment outside the infinite arm of the restarted process. By independence of the Poisson points after $N^2M$,  $\bar t_{i_*}-M\stackrel{\rm d}=N^{-2}U_{N,\lambda}$.   All particles with $i\ge i_*$ belong to the same arm of the restarted process, which meets the aggregate at time $M$ at a single point. Their paths therefore have the same selected points and the same portions of particles on every grid interval in $[0,M]$. The Finger parametrization gives
    \[
        \mathtt{Fin}_i^{\mathcal C,N}(u) =\mathtt{Fin}_{i_*}^{\mathcal C,N}(u), \qquad 0\le u\le M,\quad i\ge i_*.
    \]
    Consequently, uniformly for $\lambda\in\Lambda$,
    \[
        \begin{aligned}
        \P\!\left\{d_{H,T}^{\mathcal C} (\mathbf{Fin}^{\mathcal C,N}, \mathbf{Fin}^{\mathcal C,N}\langle M+L\rangle)>0\right\}\le\P\{\bar t_{i_*}>M+L\} =\P\{N^{-2}U_{N,\lambda}>L\} \le C_\Lambda e^{-c_\Lambda L}.
        \end{aligned}
    \]
    In particular, each such restricted collection is finite.   For $\mathcal H$, use Remark~\ref{rmk:radial-from-cylindrical}.

    Now take $\mathcal A=\mathcal S$. For $K>0$, set
    \[
        H_{N,K}:=\sup\{\bar t_i:|R^{\mathcal A,N}_{0,\bar t_i}(\bar x_i)|\le K\}, \qquad D_N(u):=U^{\mathcal A,N}_{0,u}(K)-U^{\mathcal A,N}_{0,u}(-K),
    \]
    with $\sup\varnothing=0$. For each fixed $N,K$,    \cite[Corollary~6.3 and Lemma~8.1]{berger2022growth} and the slit scaling in Lemma~\ref{lem:root-finger-identity} imply that only finitely many particles have roots in $[-K,K]$.   We need an estimate uniform in $N$. By  \eqref{eq:HL-boundary-brackets}, $D_N$ is a nonnegative martingale with $D_N(0)=2K$ and $|\Delta D_N|\le2\lambda/N$. The  proof of Proposition~\ref{lem:one-gap-section6} applied to $\Delta_h^{\mathcal S}$ gives, for $0<\vartheta<1$,
    \[
        \frac{d\langle D_N\rangle_u}{du} =\|\Delta_h^{\mathcal S}(\cdot+ND_N(u-)) -\Delta_h^{\mathcal S}\|_2^2 \ge c_{\lambda,\vartheta}>0 \quad\text{if }D_N(u-)\ge\frac{\vartheta\lambda}{N}.
    \]
    For either fixed label $x=\pm K$, the maximal inequality for the nonnegative martingale $U^{\mathcal A,N}_{0,u}(x)-U^{\mathcal A,N}_{0,u}(x-\delta)$ gives
    \[
        \P\!\left\{\sup_{u\le T} [U^{\mathcal A,N}_{0,u}(x)-U^{\mathcal A,N}_{0,u}(x-\delta)]>\varepsilon\right\} \le\frac{\delta}{\varepsilon}.
    \]
    Letting $\delta\downarrow0$ shows that, almost surely, $U^{\mathcal A,N}_{0,u}$ has no jump at the fixed labels $\pm K$, simultaneously for all $u\le T$.  The spatial coordinate of each new Poisson point has a continuous distribution and therefore almost surely avoids these endpoints. Thus Proposition~\ref{prop:SHL-boundary-cocycle} and  \eqref{eq:one-slit-boundary-maps} imply
    \[
        \begin{aligned}
        |R^{\mathcal A,N}_{0,\bar t_i}(\bar x_i)|\le K &\ \Longrightarrow\ U^{\mathcal A,N}_{0,\bar t_i-}(-K)<\bar x_i<U^{\mathcal A,N}_{0,\bar t_i-}(K) \Longrightarrow\ D_N(\bar t_i)\ge2\lambda/N.
        \end{aligned}
    \]
    For $B>2K$, let $\tau=\inf\{u\ge0:D_N(u)\notin(\vartheta\lambda/N,B)\}$. The proof of  Lemma~\ref{lem:abstract-one-gap-estimates} gives
    \[
        \E[\tau]\le\frac{2K(B+2\lambda/N-2K)}{c_{\lambda,\vartheta}}, \qquad\P\{D_N(\tau)\ge B\}\le\frac{2K}{B}.
    \]
    If $D_N$ exits through the lower boundary at time $\tau$, the conditional maximal inequality gives
    \[
        \mathbb P\!\left\{ \sup_{v\ge\tau} D_N(v)\ge \frac{2\lambda}{N} \,\middle|\,\mathcal F_\tau \right\} \le \frac{N D_N(\tau)}{2\lambda} \le \frac{\vartheta}{2}.
    \]
    Therefore
    \[
        \P\{H_{N,K}>G\} \le\frac{2K}{B} +\frac{2K(B+2\lambda/N-2K)}{c_{\lambda,\vartheta}G} +\frac{\vartheta}{2}, \qquad \lim_{G\to\infty}\sup_N\P\{H_{N,K}>G\}=0,
    \]
    where the limit follows by taking $B=\sqrt G$, then $G\to\infty$  and $\vartheta\downarrow0$.

    Fix $m\ge1$ and take $K>4(m+1)$. Let $E_{N,K}$ be the event that
    \[
        U^{\mathcal A,N}_{0,u}(K)\ge K/2,\quad U^{\mathcal A,N}_{0,u}(-K)\le-K/2 \quad(0\le u\le T+1),
    \]
    and, for every $j$ with $\bar t_j\le T+1$ and $|\bar x_j|>K/2$,
    \[
        \inf_{z\in\mathcal P_j^{\mathcal S,N}} \operatorname{sgn}(\bar x_j)N^{-1}\rea  z>m+1.
    \]
    The maximal inequality and \eqref{eq:SHL-remote-particle-summation}, with $p=8$, give    $\P(E_{N,K}^c)\le C_{T,\lambda}     \left(K^{-2}+N^{-1}K^{-7}\right)$  for all sufficiently large $N$. For $j\preceq i$ with $\bar t_j\le T+1$, \eqref{eq:ancestral-centre-root} and the cocycle property give, on $E_{N,K}$,
    \[
        \begin{aligned}
        R^{\mathcal A,N}_{0,\bar t_i}(\bar x_i) =R^{\mathcal A,N}_{0,\bar t_j-}(\bar x_j)>K\ \Longrightarrow\ \bar x_j>U^{\mathcal A,N}_{0,\bar t_j-}(K)\ge K/2\ \Longrightarrow\ \inf_{z\in\mathcal P_j^{\mathcal S,N}}N^{-1}\rea  z>m+1.
        \end{aligned}
    \]
    For $R^{\mathcal A,N}_{0,\bar t_i}(\bar x_i)<-K$, the inequalities are reversed. Every edge traversed on $[0,T]$ is present by the next point of $\Gamma_N^{\mathcal S}$, hence by time $T+N^{-1}\le T+1$, independently of $\bar t_i$. Thus, on $E_{N,K}$, if $|R^{\mathcal A,N}_{0,\bar t_i}(\bar x_i)|>K$, then
    \[
        \chi_{\mathcal S,m} \bigl(\mathtt{Fin}_i^{\mathcal S,N}(u\wedge\bar t_i)\bigr) =\bigl(\pm(m+1/2),0\bigr), \qquad 0\le u\le T.
    \]
    For fixed $N$, the proof of \eqref{eq:SHL-remote-particle-summation}, using the same moment estimate with an $N$-dependent constant, together with the maximal inequality for the two boundary paths, gives
    \[
        \mathbb P(E_{N,K}^c)\le C_{N,T,m,\lambda}K^{-2}
    \]
    for all sufficiently large $K$. By the Borel--Cantelli lemma, for every integer $T,m\ge1$, the preceding identity holds whenever the absolute value of the root exceeds an almost surely finite $K_0=K_0(N,T,m)$. Only finitely many generating paths have roots in $[-K_0,K_0]$, and the restricted terminal times lie in $[0,T]$. It follows that the collection is totally bounded in $D_{\mathcal A}$, and hence its closure is compact. The same argument for  the collections with $\bar t_i\le G$, together with their countably many  measurable generating paths, gives measurability. Both constant paths $(\pm(m+1/2),0)$ above are limits of paths with $\bar t_i\le T+2$. Indeed, choose Poisson points with
    \[
        \bar t_{i_\ell}\to T+1,\qquad \bar x_{i_\ell}\to+\infty\ \text{or }-\infty, \qquad\bar t_{i_\ell}<T+2.
    \]
    By Proposition~\ref{prop:SHL-boundary-cocycle}, their roots tend to $\pm\infty$, respectively. Applying the preceding fixed-$N$ estimate with time bound larger than $T+2$, the corresponding paths converge in every local spatial coordinate, with terminal time $T+1$. Hence the two constant limits belong to $\mathbf{Fin}^{\mathcal S,N}\langle T+2\rangle$, and their restrictions to $[0,T]$ have terminal time $T$. On $E_{N,K}\cap\{H_{N,K}\le G\}$, every path with $\bar t_i>G$ has the same image under $\chi_{\mathcal S,m}$ on $[0,T]$ as one of these two limits, while every path with $\bar t_i\le G$ belongs to both collections. Thus, for $G>T+2$,
    \[
        \begin{aligned}
        \P\!\left\{d_{H,T}^{\mathcal S,m} (\mathbf{Fin}^{\mathcal S,N}, \mathbf{Fin}^{\mathcal S,N}\langle G\rangle)>0\right\}\le\P\{H_{N,K}>G\} +C_{T,\lambda}\left(K^{-2}+N^{-1}K^{-7}\right).
        \end{aligned}
    \]
    Letting $N\to\infty$, then $G\to\infty$, and then $K\to\infty$, and using \eqref{eq:geometric-path-metric-comparison}, first at fixed $h$ and then as $h\to\infty$, we obtain \eqref{eq:physical-cutoff-saturation}.
    \end{proof}
    We now combine the preceding convergence results and geometric approximations to establish the joint convergence of $(\mathbf{Fin}^{\mathcal A,N},\mathbf{Gap}^{\mathcal A,N})$.
    \begin{proof}[Proof of Theorem~\ref{thm:joint-geometric-webs}]
    Let   $\mathcal B=\sigma_{\mathcal A}\mathcal W^\downarrow_{\mathcal A,\ge0}$  and let  $\mathcal B^G$ be the closure of its paths with terminal time at most $G$, restricted to nonnegative times. By the construction in Subsection~\ref{SS.convergence criteria},
    \[
        \mathcal B^G\subseteq\mathcal B^{G+1}, \qquad \overline{\bigcup_{G\in\N}\mathcal B^G}=\mathcal B .
    \]
    The collection $\mathcal B$ is compact. By the definition of the geometric metrics in Subsection~\ref{S.HL-models}, the embedding $\mathfrak e_{\mathcal A}^{v_{\mathcal A}}$ is continuous; for $\mathcal A=\mathcal S$, use the continuous extensions of $\kappa_m\circ\operatorname{arctanh}$ and $\zeta_m\circ\operatorname{arctanh}$ to $[-1,1]$. Hence $D_{\mathcal A,H}(\mathfrak e_{\mathcal A}^{v_{\mathcal A}}(\mathcal B^G),\mathfrak E_{v_{\mathcal A}}^{\mathcal A,\downarrow}(\mathcal B))\rightarrow0\ \text{a.s.} $ For every fixed integer $G\ge1$, Proposition~\ref{prop:process-root-webs}, Proposition~\ref{prop:process-tip-geometric-transfer}, and continuity of the embedding give $\mathbf{Fin}^{\mathcal A,N}\langle G\rangle \dto \mathfrak e_{\mathcal A}^{v_{\mathcal A}}(\mathcal B^G)$ in $D_{\mathcal A,H}$. By Theorem~\ref{thm:harmonic-measure-webs} and \eqref{eq:joint-root-HM-webs}, this convergence is joint with $\mathfrak e_{\mathcal A}^{v_{\mathcal A}}(\mathbf{HM}^{\mathcal A,N})$ in $D_{\mathcal A,H}^{\uparrow}$, with the same double BW coupling. Proposition~\ref{prop:process-root-finger-saturation} then allows us to let first $N\to\infty$ and then $G\to\infty$, while Proposition~\ref{prop:gap-HM-comparison} replaces the embedded harmonic measure collection by $\mathbf{Gap}^{\mathcal A,N}$. This proves the joint convergence.
    \end{proof}
    We finally identify the scaling limits of the two backbone paths.
    \begin{proof}[Proof of Corollary~\ref{cor:backbone-dual-backbone}]
    Let $U$ and $i_N$ be as in the definition of the selected dual backbone in Subsection~\ref{S.HL-models}.
    Proposition~\ref{prop:process-root-finger-saturation} gives $\mathtt{B}^{\mathcal A,N,\downarrow}\in\mathbf{Fin}^{\mathcal A,N}$  with terminal time $+\infty$. As in the proof of    Proposition~\ref{prop:gap-HM-comparison},     Lemma~\ref{lem:completed-cell-oscillation} gives
    \[
        d_{\mathbb T}(\xi_{i_N},U)+\sigma_{i_N}\pto0.
    \]
    Hence, by Theorem~\ref{thm:joint-geometric-webs}, every subsequential limit of the two paths belongs to the corresponding embedded double BW.  The first path has terminal time $+\infty$, and is $\mathfrak e_{\mathcal A}^{v_{\mathcal A}} (C_{\mathcal A}^{\downarrow}|_{[0,\infty)})$    by \cite[Theorem~3.1]{coupier2019directed}.  The second starts from $\mathfrak e_{\mathcal A}^{v_{\mathcal A}}(0,U)$ and is almost surely unique by Theorem~\ref{thm:BW-CBW}, since $U$ is independent of the limiting double BW.  Finally, \cite[proof of Theorem~3.6, Step~4]{coupier2019directed} identifies $C_{\mathcal A}^{\uparrow}(0)$ as uniform and independent of the double BW on nonnegative times, which gives the required joint convergence.
    \end{proof}

    \subsection{\texorpdfstring{One-arm domination time}{One-arm domination time}}\label{SS.one-arm-domination-section6}
    Write $\tau_{\mathcal C}^N(\lambda)$ for the rescaled one-arm domination time $\tau_{\mathcal C}^N$ of Definition~\ref{def:one-arm-domination-time}, to indicate its dependence on the particle height $\lambda$, and set
    \begin{equation*}
        D_N^\lambda(t):=2\pi-\sup_{y\in\mathbb T} \bigl|\{x\in\mathbb T:R_{0,t}^{\mathcal C,N}(x)=y\}\bigr|, \qquad t\ge0.
    \end{equation*}
    By Lemma~\ref{lem:root-finger-identity} with $s=0$, the sets $(R_{0,t}^{\mathcal C,N})^{-1}\{y\}$ of positive length are, up to endpoints, exactly the boundary-label intervals of the arms. They evolve by inverse images under the rescaled one-slit maps, as in \cite[Definition~3.3 and Proposition~3.4]{chen2025one}, and attachment locations almost surely avoid their endpoints. Thus $1-D_N^\lambda(t)/(2\pi)$ is the largest conditional probability, over the existing arms, that the next particle attaches to that arm; it equals $0$ before the first attachment.

    Recall that $\sigma_\lambda^2=4\lambda^3/3$ and that $\upsilon_\infty$ is the global coalescence time introduced before Theorem~\ref{thm:one-arm-domination}. The limiting analogue of $D_N^\lambda$ is defined as follows. For $t>0$, let $x\mapsto\widehat\pi_{x,t}^{\mathcal C,\sigma_\lambda}(0)$ be the map sending $x$ to the position at time $0$ of the backward path from $(x,t)$, made right-continuous by taking right limits along a fixed countable dense set of labels, and set
    \[
        D^\lambda(t):=2\pi-\sup_{y\in\mathbb T} \bigl|\{x\in\mathbb T: \widehat\pi_{x,t}^{\mathcal C,\sigma_\lambda}(0)=y\}\bigr|.
    \]
    We first establish the Brownian limit of the  $D_N^{\lambda_N}(t)$.
    \begin{proposition}\label{prop:maximal-fibre-marginal}
    If $\lambda_N\to\lambda>0$ and $t>0$, then
    \begin{equation*}
        D_N^{\lambda_N}(t)\dto D^\lambda(t),\qquad \mathbb P\{D^\lambda(t)=0\} =\mathbb P\{\upsilon_\infty\le\sigma_\lambda^2t\}.
    \end{equation*}
    \end{proposition}
    \begin{proof}
    Fix an integer $T>t$ and use the reversed Root flow $\Phi^N$ of Subsection~\ref{S.three-point-estimate-I}. Stationarity gives $R_{0,t}^{\mathcal C,N}\stackrel{\rm d}=\Phi^N_{0,t}$.  Proposition~\ref{prop:fdd-selected-paths-section7} and the meeting-time argument in the proof of Theorem~\ref{thm:harmonic-measure-webs}, applied on $[0,T]$, give joint convergence of paths and pairwise meeting times. For fixed labels $x_1,\ldots,x_m\in\mathbb T$, let $\tau_{ij}^N$ be the meeting times of the c\`adl\`ag trajectories $u\mapsto\Phi^N_{0,u}(x_i)$ and $u\mapsto\Phi^N_{0,u}(x_j)$, and let $\tau_{ij}$ be their limiting meeting times, all truncated at $T$. Since $\P\{\tau_{ij}=t\}=0$,
    \[
        \bigl(\ind_{\{\tau_{ij}^N\le t\}}\bigr)_{i<j} \dto \bigl(\ind_{\{\tau_{ij}\le t\}}\bigr)_{i<j}
    \]
    for every finite set of labels.

    For a non-decreasing circular map $f$, each set $f^{-1}\{y\}$ is a circular interval. Hence, using $m$ equally spaced labels,
    \begin{equation*}
        \left| \frac{2\pi}{m}\max_{y\in\mathbb T} \#\{0\le j<m:f(2\pi j/m)=y\} -\sup_{y\in\mathbb T}|f^{-1}\{y\}| \right| \le \frac{4\pi}{m}.
    \end{equation*}
    Hence the convergence for fixed $m$, followed by $m\to\infty$, yields $ D_N^{\lambda_N}(t)\dto D^\lambda(t). $    Since the sets $f^{-1}\{y\}$ are disjoint, at most one can have length greater than $\pi$. Thus
    \[
        D^\lambda(t)=0 \quad\Longleftrightarrow\quad \bigl|\{x\in\mathbb T: \widehat\pi_{x,t}^{\mathcal C,\sigma_\lambda}(0)=y\}\bigr|=2\pi \quad\text{for some }y\in\mathbb T.
    \]
    By right continuity, this is equivalent to $x\mapsto\widehat\pi_{x,t}^{\mathcal C,\sigma_\lambda}(0)$ being constant, that is, all trajectories have coalesced by time $t$. This completes the proof.
    \end{proof}
    We next compare the one-arm domination time with the vanishing of $D_N^\lambda$.
    \begin{lemma}\label{lem:one-arm-two-sided}
    For every compact $\Lambda\Subset(0,\infty)$, $q\ge0$, and $s>0$,
    \begin{align}
        \lim_{\alpha\downarrow0}\limsup_{N\to\infty}\sup_{\lambda\in\Lambda} \mathbb P\{\tau_{\mathcal C}^N(\lambda)>q+s, D_N^\lambda(q)\le\alpha\}&=0, \label{eq:one-arm-lower-comparison} \\
        \lim_{N\to\infty}\sup_{\lambda\in\Lambda} \mathbb P\{\tau_{\mathcal C}^N(\lambda)\le q, D_N^\lambda(q+s)>\alpha\}&=0, \qquad 0<\alpha<\pi. \label{eq:one-arm-upper-comparison}
    \end{align}
    \end{lemma}
    \begin{proof}
    In this proof, $(\mathcal F_u)$ denotes the usual filtration generated by the Poisson points up to macroscopic time $u$. We use $j_N=C_\Lambda/N$ from Subsection~\ref{S.three-point-estimate-I} and the covariance rates in \eqref{eq:scaled-covariance-rates-section7}. Let $Y_u=U_{0,u}^{\mathcal C,N}(a+d)-U_{0,u}^{\mathcal C,N}(a)$ be the length of the evolving arc started from $[a,a+d]$, with initial length $0\le d\le1/2$. By \eqref{eq:HL-boundary-brackets}, $Y$ is a martingale with $0\le Y\le2\pi$, $Y_0=d$, and $|\Delta Y|\le2j_N$. Fix $0<\varepsilon<1$, let $a_0$ be as in Proposition~\ref{lem:one-gap-section6}, and set $\tau:=\inf\{u\ge0:Y_u\notin(\varepsilon a_0/N,1)\}$.  The proof of Proposition~\ref{lem:one-gap-section6}, with $r$ replaced by $h$, gives
    \[
        \frac{d\langle Y\rangle_u}{du} =2\{\mathcal K_h^{\mathcal C,N}(0) -\mathcal K_h^{\mathcal C,N}(Y_{u-})\} \ge c_{\Lambda,\varepsilon}>0,\qquad u<\tau.
    \]
    The optional stopping theorem for $Y$ and $Y^2-\langle Y\rangle$ gives
    \[
        \mathbb E[\tau]\le \frac{d(1+2j_N-d)}{c_{\Lambda,\varepsilon}},\qquad \mathbb P\{Y_\tau\ge1\}\le d.
    \]
    If the attachment location at time $u$ lies in the arc immediately before the attachment, then $Y_u\ge 4\arcsin\delta_N(\lambda)\ge 2a_0/N$. Hence, on $\{Y_\tau\le \varepsilon a_0/N\}$, the conditional maximal inequality gives
    \[
        \P\!\left( \text{the arc receives a particle after }\tau \,\middle|\,\mathcal F_\tau \right)\le \frac{NY_\tau}{2a_0} \le \frac{\varepsilon}{2}.
    \]
    Consequently,
    \[
        \mathbb P\{\text{the arc receives a particle after }s\} \le\mathbb P\{\tau>s\}+\mathbb P\{Y_\tau\ge1\} +\frac{\varepsilon}{2} \le C_{\Lambda,\varepsilon,s}d+\frac{\varepsilon}{2}.
    \]
    On the event $D_N^\lambda(q)\le\alpha$, where $\alpha<1/2$, the arm with label interval of length at least $2\pi-\alpha$ is unique and determined by the Poisson points up to $q$. Conditionally on these points, its complement evolves as above, with $d\le\alpha$. If the complement receives no  particle after $q+s$, this arm is the infinite arm and $\tau_{\mathcal C}^N(\lambda)\le q+s$. Taking $N\to\infty$,  then $\alpha\downarrow0$, and then $\varepsilon\downarrow0$    proves \eqref{eq:one-arm-lower-comparison}.

    Let $K=\#\{j:q<\bar t_j\le q+s\}$, so $K\sim\operatorname{Poisson}(2\pi N^3s)$. On $\{\tau_{\mathcal C}^N(\lambda)\le q\}$, every subsequent    point goes to the infinite arm, with at most one particle to arrive needed to create it. By \cite[Lemma~2.2(1)]{chen2025one},
    \[
        D_N^\lambda(q+s) \le2\pi\bigl(1-\delta_N(\lambda)^2\bigr)^{(K-1)_+/2}.
    \]
    Since $\delta_N(\lambda)\asymp_\Lambda N^{-1}$, it follows that
    \[
        \{\tau_{\mathcal C}^N(\lambda)\le q, D_N^\lambda(q+s)>\alpha\} \subseteq\{K\le C_{\Lambda,\alpha}N^2\}.
    \]
    The lower-tail bound for the Poisson distribution proves \eqref{eq:one-arm-upper-comparison}.
    \end{proof}

    \begin{proof}[Proof of Theorem~\ref{thm:one-arm-domination}]
    Fix $\lambda_N\to\lambda>0$ in a compact $\Lambda\Subset(0,\infty)$, and let $t>s>0$. Apply Lemma~\ref{lem:one-arm-two-sided} at $q=t-s$ and $q=t$, together with Proposition~\ref{prop:maximal-fibre-marginal}, and let $\alpha\downarrow0$ through continuity points of the laws of $D^\lambda(t-s)$ and $D^\lambda(t+s)$. Then
    \[
        \begin{aligned}
        \mathbb P\{\upsilon_\infty\le\sigma_\lambda^2(t-s)\}\le\liminf_{N\to\infty} \mathbb P\{\tau_{\mathcal C}^N(\lambda_N)\le t\}\le\limsup_{N\to\infty} \mathbb P\{\tau_{\mathcal C}^N(\lambda_N)\le t\} \le\mathbb P\{\upsilon_\infty\le\sigma_\lambda^2(t+s)\}.
        \end{aligned}
    \]
    The law of $\upsilon_\infty$ is continuous \cite[Section~4, eq.~(17)]{norris2015weak}; letting $s\downarrow0$  gives $\tau_{\mathcal C}^N(\lambda_N)\dto\sigma_\lambda^{-2}\upsilon_\infty$.     Moreover, $0\le\tau_{\mathcal C}^N(\lambda_N)\le N^{-2}U_{N,\lambda_N}$, so \eqref{eq:hl-observables-cpz-tail} gives uniformly bounded second  moments and hence   $\mathbb E[\tau_{\mathcal C}^N(\lambda_N)]     \rightarrow\sigma_\lambda^{-2}\mathbb E[\upsilon_\infty]$.    For $\mathcal A=\mathcal H$, Remark~\ref{rmk:radial-from-cylindrical} gives $\tau_{\mathcal H}^N=2\pi\,\tau_{\mathcal C}^N(\lambda_N^{\mathcal H}), \lambda_N^{\mathcal H}\to\eta, $  and the conclusion follows.
    \end{proof}
    
    \appendix
    \section{Complementary proofs}\label{appendix:complementary-proofs}
    \subsection{Proof of Lemma~\ref{lem:one-slit-moments}}\label{appendix:one-slit-moments}
    \begin{proof}
    Fix $\Lambda\Subset(0,\infty)$; all bounds below are uniform for $\lambda\in\Lambda$. The SHL(0) identities in \eqref{eq:CHL-moments} follow from \cite[proofs of Theorem~7.2 and Lemma~7.6]{berger2022growth} and $S_\xi^\lambda(z)=\lambda S_{\xi/\lambda}^1(z/\lambda)$.

    For CHL(0), let $\delta_N=\tanh(\lambda/(2N))$. Then
    \begin{equation*}
        N\delta_N=\frac\lambda2+O_\Lambda(N^{-2}),\qquad 2N\arcsin\delta_N\le\lambda,
    \end{equation*}
    since $(\arcsin\tanh)'=\operatorname{sech}\le1$. For $0<\delta<1$, define $R_\delta:\mathbb R\to\mathbb R$ on $[0,\pi]$ by
    \begin{equation}\label{eq:Rdelta-formula-section3}
        R_\delta(u)=
        \begin{cases}
        0,&0\le u\le 2\arcsin\delta, \\[1mm]
        2\arctan\sqrt{(1-\delta^2)\tan^2(u/2)-\delta^2}, &2\arcsin\delta\le u<\pi
        \end{cases}.
    \end{equation}
    Set $R_\delta(\pi)=\pi$ by continuity and extend it by $R_\delta(2\pi-u)=2\pi-R_\delta(u)$ and $R_\delta(u+2\pi)=R_\delta(u)+2\pi$. With \eqref{eq:one-slit-boundary-maps},
    \begin{equation}\label{eq:Rdelta-scaling-section3}
        r_0^{\mathcal C,N}(Nu)=N R_{\delta_N}(u),\qquad h_0^{\mathcal C,N}(Nu)=N R_{\delta_N}^{\leftarrow}(u).
    \end{equation}
    By symmetry,
    \[
        \int_0^{2\pi}\bigl(R_\delta(u)-u\bigr)\,du = \int_0^{2\pi}\bigl(R_\delta^{\leftarrow}(u)-u\bigr)\,du =0.
    \]
    Since $R_\delta$ is continuous and non-decreasing, with
    \[
        \int_0^{2\pi} \bigl(R_\delta(u)-u\bigr)^2 d\bigl(R_\delta(u)-u\bigr) = \frac13 \Bigl[ \bigl(R_\delta(u)-u\bigr)^3 \Bigr]_{0}^{2\pi} =0,
    \]
    we have
    \begin{equation}\label{eq:Rdelta-second-moment-identity}
        \begin{aligned}
        \int_0^{2\pi} \bigl(R_\delta^{\leftarrow}(y)-y\bigr)^2\,dy= \int_0^{2\pi} \bigl(u-R_\delta(u)\bigr)^2\,dR_\delta(u)= \int_0^{2\pi} \bigl(R_\delta(u)-u\bigr)^2\,du ,
        \end{aligned}
    \end{equation}
    where $R_\delta(0)=0$ and $R_\delta(2\pi)=2\pi$. The estimate \cite[(2.18)]{chen2025one} and the same reflection imply
    \begin{equation}\label{eq:cylinder-deviation-tail}
        |R_\delta(u)-u|+|R_\delta^{\leftarrow}(u)-u| \le\frac{C\delta^2}{\delta+d_{\mathbb T}(u,0)}, \qquad 0<\delta\le\frac12.
    \end{equation}
    The expansion \cite[Lemma~2.3]{chen2025one} gives
    \begin{equation*}
        \int_0^{2\pi}\bigl(R_\delta^{\leftarrow}(u)-u\bigr)^2\,du =\frac{32}{3}\delta^3-2\pi\delta^4+O(\delta^5).
    \end{equation*}
    By \eqref{eq:Rdelta-scaling-section3} and \eqref{eq:Rdelta-second-moment-identity},
    \[
        N^3\int_0^{2\pi}\bigl(R_{\delta_N}(u)-u\bigr)^2\,du =N^3\int_0^{2\pi}\bigl(R_{\delta_N}^{\leftarrow}(u)-u\bigr)^2\,du =\sigma_\lambda^2+O_\Lambda(N^{-1}).
    \]
    Also, the two deviations from the identity are each bounded by $2N\arcsin\delta_N\le\lambda$. This proves \eqref{eq:CHL-moments} and \eqref{eq:bounded-one-slit-jumps}.

    For the SHL(0) estimates, on $\mathbb R$, $|S^\lambda(z)-z|\le\lambda$ and $\im (S^\lambda(z)-z)\ge0$, while $S^\lambda(z)-z=-\lambda^2/(2z)+O_\Lambda(|z|^{-3})$ at infinity. The maximum-modulus and minimum principles extend both inequalities to $\overline{\mathbb H}$. Since $\rea(S^\lambda(z)\overline z)\ge0$,
    \[
        S^\lambda(z)-z=-\frac{\lambda^2}{S^\lambda(z)+z},\quad |S^\lambda(z)+z|^2\ge|z|^2, \quad |S^\lambda(x+\mathrm iy)-(x+\mathrm iy)|\le\frac{C_\Lambda}{1+|x|+y} \quad(y\ge0).
    \]
    Hence, for $x\in\mathbb R$, $y\ge0$, and $q\ge2$,
    \[
        \int_{\mathbb R}|S_\xi^\lambda(x+\mathrm iy)-(x+\mathrm iy)|^q\,d\xi \le C_{q,\Lambda}\int_{\mathbb R} (1+|x-\xi|+y)^{-q}\,d\xi \le C_{q,\Lambda}(1+y)^{1-q}.
    \]
    This proves \eqref{eq:complex-slit-global-moment} and \eqref{eq:complex-slit-high-moment} for SHL(0).

    For CHL(0), the boundary values satisfy $|C^{N,\lambda}(z)-z|\le2\lambda$ and $\im(C^{N,\lambda}(z)-z)\ge0$. After the change of variables $\zeta=e^{\mathrm iz/N}$, the periodic deviation $C^{N,\lambda}(z)-z$ is analytic in $\{|\zeta|<1\}$ (including at $\zeta=0$ by \cite[Lemma~5.4]{procaccia2023cylindrical}) and continuous on its closure. Hence its modulus is bounded by its boundary supremum and its imaginary part by its boundary infimum, so the boundary bounds extend to $\overline{\mathbb H}$. Thus \cite[Lemma~5.2]{procaccia2023cylindrical} and Fatou's lemma give
    \[
        \int_{\mathbb T_N}|C_\xi^{N,\lambda}(z)-z|^q\,d\xi \le(2\lambda)^{q-2} \int_{\mathbb T_N}|C_\xi^{N,\lambda}(z)-z|^2\,d\xi \le C_{q,\Lambda}.
    \]
    Note that for $\operatorname{Im}z\ge hN$,
    \[
        \bigl|C^{N,\lambda}(z)-z\bigr|\le-N\log(1-\delta_N^2)\coth(h/2) \le\frac{C_{h,\Lambda}}{N} \implies   \int_{\mathbb T_N} \bigl|C_\xi^{N,\lambda}(z)-z\bigr|^q\,d\xi \le C_{q,h,\Lambda}N^{1-q}.
    \]
    This proves \eqref{eq:complex-slit-global-moment} and \eqref{eq:complex-slit-high-moment} for CHL(0).

    Finally, \cite[Lemma~3.3(ii)]{berger2022growth} and \cite[Lemma~5.5]{procaccia2023cylindrical} give, for $\im z>0$,
    \[
        \int_{\mathbb R}(S_\xi^\lambda(z)-z)\,d\xi =\mathrm iv_{\mathcal S},\qquad \int_{\mathbb T_N}(C_\xi^{N,\lambda}(z)-z)\,d\xi =-2\pi\mathrm iN^2\log(1-\delta_N^2) =\mathrm i\{v_{\mathcal C}+O_\Lambda(N^{-2})\}.
    \]
    Both identities extend to the boundary; for the principal value, the symmetric tails cancel since $S_\xi^\lambda(z)-z=\lambda^2/(2\xi)+O_{\Lambda,K}(|\xi|^{-2})$, uniformly for $z$ in a compact set $K\subset\overline{\mathbb H}$.
    \end{proof}
    \subsection{Proof of Proposition~\ref{prop:SHL-boundary-cocycle}}\label{appendix:SHL-boundary-cocycle}
    \begin{proof}
    We follow the proof of  \cite[Theorem~4.3 and Section~7.1]{berger2022growth}. Since the truncated maps in \eqref{eq:SHL-root-cutoff} satisfy
    \[
        R_{s,t}^{\mathcal S,N;M}(a) =N^{-1}R_{N^2s,N^2t}^{\mathcal S,1;M}(Na),
    \]
    it suffices to take $N=1$. Let $R^M=R^{\mathcal S,1;M}$,  $R=R^{\mathcal S,1}$, $r_x=r_x^{\mathcal S}$, $h_x=h_x^{\mathcal S}$, and  $v_x=r_x-\id$, $g_x=h_x-\id$. Let $\Pi$ be the counting measure of the Poisson points and set $\widetilde\Pi=\Pi-du\,dx$; all sums below are over $(u,x)\in\Pi$. Fix integers $T,K\ge1$ and take $B=2K+2\lambda+1$.  By the explicit one-slit formulas \cite[Section~7.1]{berger2022growth}, for $|x|>B$, define
    \[
        w_x(z):=v_x(z)-\frac{\lambda^2}{2x},\qquad  \widetilde w_x(z):=g_x(z)+\frac{\lambda^2}{2x}.
    \]
    Then
    \begin{equation}\label{eq:SHL-tail-splitting}
        \sup_{|z|\le K}\bigl(|w_x(z)|+|\widetilde w_x(z)| +|w'_x(z)|+|\widetilde w'_x(z)|\bigr) \le C_{K,\lambda}x^{-2}.
    \end{equation}
    For $M>0$, let $A_B^M(t)=\sum_{u\le t,\,B<|x|\le M}\lambda^2/(2x)$.  Applying Doob's inequality first in time and then as $M$ varies yields
    \[
        \E\left[\sup_{\mu\le M\le\nu}\sup_{t\le T} |A_B^M(t)-A_B^\mu(t)|^2\right] \le16T\int_{\mu<|x|\le\nu}\frac{\lambda^4}{4x^2}\,dx, \qquad B\le\mu<\nu.
    \]
    For the events
    \[
        E_j:= \left\{ \sup_{2^j\le M\le 2^{j+1}} \sup_{t\le T} |A_B^M(t)-A_B^{2^j}(t)|>2^{-j/4} \right\},
    \]
    the preceding estimate and Markov's inequality give
    \[
        \mathbb P(E_j)\le C_T2^{-j/2}, \qquad \sum_{j\ge1}\mathbb P(E_j)<\infty.
    \]
    Hence, by the Borel--Cantelli lemma and
    \[
        \E\left[\sum_{\substack{u\le T, |x|>B}}x^{-2}\right] =T\int_{|x|>B}x^{-2}\,dx<\infty,
    \]
    almost surely, for all integers $T,K$,
    \[
        A_B^M\rightarrow A_B \quad\text{uniformly on }[0,T], \qquad \sum_{u\le T,  |x|>B}x^{-2}<\infty .
    \]
    Since $\Pi([0,T]\times\{x:M_1<|x|\le M_2\})<\infty$ for $B\le M_1<M_2<\infty$, the convergence holds for every $M\to\infty$.

    Define $\chi_K(z):=(-K)\vee(z\wedge K), r_x^K(z):=z+v_x(\chi_K(z)).$  Then $r_x^K$ is continuous and non-decreasing, and for $|x|\le B$,
    \[
        |r_x^K(a)-r_x^K(b)| \le \omega_K(|a-b|), \qquad \omega_K(d):=d+C_{K,\lambda}\sqrt d .
    \]
    Let $\widehat R^{M,K}_{s,t}$ be the composition of the maps $r_x^K$ over $(u,x)\in\Pi\cap((s,t]\times[-M,M])$ in increasing order of $u$. The time order is reversed from that of $R^{M}_{s,t}$. Define $\widehat R^M_{s,t}$ similarly with $r_x$.

    On any time interval containing no point with $|x|\le B$, suppose the two trajectories have values $y$ and $y'$ at the beginning of the interval, and let $d:=|y-y'|$. Then, for $M'\ge M\ge B$,
    \begin{equation}\label{eq:SHL-cutoff-error}
        \left|\widehat R_{s,t}^{M',K}(y')-\widehat R_{s,t}^{M,K}(y)\right|\le(d+\varepsilon_M)\exp\!\left(C_{K,\lambda}\sum_{u\le T,\,|x|>B}x^{-2}\right),
    \end{equation}
    where
    \[
        \varepsilon_M:= 2\sup_{M'\ge M}\|A_B^{M'}-A_B^M\|_\infty + C_{K,\lambda} \sum_{u\le T, |x|>M}x^{-2} \rightarrow0 .
    \]
    There are only finitely many points with $|x|\le B$ before time $T$; at each such point the difference is replaced by at most $\omega_K$ of its previous value. Hence
    \[
        \sup_{0\le s\le t\le T}\sup_{a\in\mathbb R} |\widehat R^{M',K}_{s,t}(a)-\widehat R^{M,K}_{s,t}(a)| \rightarrow0 \qquad (M,M'\to\infty).
    \]
    With $\varepsilon_M=0$, the same argument shows that $\widehat R^{M,K}_{s,t}(a)$ is uniformly continuous in $a$, uniformly in $M,s,t$. The same bounds hold for $[s,t)$ with the endpoint convention reversed. Lemma~\ref{lem:one-slit-moments} and its proof give
    \[
        \int_{\R}v_x(z)^2\,dx=\frac43\lambda^3, \qquad \sup_{|z|\le K}\left|\int_{-M}^M v_x(z)\,dx\right| \le\frac{2K\lambda^2}{M-K},\quad M>K+\lambda.
    \]
    Since $\chi_K(\R)\subset[-K,K]$, the integral is square-integrable. For the limit $\widehat R^K$, we have
    \[
        \widehat R^K_{0,t}(a)-a = \int_{(0,t]\times\R} v_x\!\left(\chi_K(\widehat R^K_{0,u-}(a))\right)\,\widetilde\Pi(du\,dx),
    \]
    and Doob's inequality and the preceding second-moment identity give
    \[
        \mathbb E\!\left[\sup_{0\le t\le T} \bigl|\widehat R^K_{0,t}(a)-a\bigr|^2\right] \le 4\mathbb E\!\left[\int_0^T\!\int_{\R} v_x\!\left(\chi_K(\widehat R^K_{0,u-}(a))\right)^2\,dx\,du\right] =\frac{16}{3}\lambda^3T.
    \]
    Thus, for $K>1+|a|$,
    \[
        \P\!\left\{ \sup_{t\le T}|\widehat R^K_{0,t}(a)|\ge K-1 \right\} \le \frac{16\lambda^3T} {3(K-1-|a|)^2}.
    \]
    For every $a\in\mathbb Q$, the Borel--Cantelli lemma gives, almost surely, $\widehat R^{M,K}_{0,\cdot}(a)=\widehat R^M_{0,\cdot}(a)$ for all sufficiently large $K$ and  $M$. By monotonicity, the same  holds uniformly for $a$ in every compact interval bounded by rational endpoints. Consequently, $\widehat R^M_{0,t}(a)\rightarrow \widehat R_{0,t}(a)$  locally and uniformly in $a$ and $t\le T$. Moreover, $\widehat R_{0,t}$ is uniformly continuous in $a$, uniformly for $t\le T$.

    By Fatou's lemma, the preceding second moment bound also holds for $\widehat R$. Hence the Borel--Cantelli lemma, applied to $a=\pm m$, $m\in\N$, and monotonicity give
    \[
        \lim_{a\to+\infty}\inf_{t\le T}\widehat R_{0,t}(a)=+\infty, \qquad \lim_{a\to-\infty}\sup_{t\le T}\widehat R_{0,t}(a)=-\infty.
    \]
    Fix $L<\infty$. Choose $a_0\in\N$ such that
    \[
        \widehat R_{0,s}(-a_0)<-L-1<L+1<\widehat R_{0,s}(a_0), \qquad 0\le s\le T.
    \]
    Then, for all sufficiently large $M$,
    \[
        \widehat R^M_{0,t}(-a_0) \le \widehat R^M_{s,t}(a) \le \widehat R^M_{0,t}(a_0), \qquad |a|\le L,\quad 0\le s\le t\le T.
    \]
    Since $\widehat R^M_{0,t}(\pm a_0)$ are uniformly bounded for $t\le T$ and all sufficiently large $M$, we can fix $K$ large enough so that $\widehat R^{M,K}_{s,t}(a)=\widehat R^M_{s,t}(a)$ for all sufficiently large $M$, whenever $0\le s\le t\le T$ and $|a|\le L$. Thus  $\widehat R^M_{s,t}(a)\rightarrow \widehat R_{s,t}(a) $  uniformly on $\{(s,t,a):0\le s\le t\le T,\ |a|\le L\}$,  and $\widehat R_{s,t}$ is uniformly continuous in $a$ there. By time reflection, this proves \eqref{eq:SHL-root-cutoff-limit} and \eqref{eq:SHL-root-uniform-continuity}. Since $\widehat R^M_{s,u}=\widehat R^M_{t,u}\circ\widehat R^M_{s,t}$, letting $M\to\infty$ gives $\widehat R_{s,u}=    \widehat R_{t,u}\circ\widehat R_{s,t},\ 0\le s\le t\le u\le T$. By time reflection, this proves the first identity in \eqref{eq:SHL-grid-cocycle}. Finally, for $b_0>0$, choose $a_1,a_2$ so that
    \[
        \inf_{t\le T}\widehat R_{0,t}(a_1)>b_0, \qquad a_2>\sup_{s\le T}\widehat R_{0,s}(a_1).
    \]
    Then
    \[
        \widehat R_{s,t}(a_2) \ge \widehat R_{s,t}\!\left(\widehat R_{0,s}(a_1)\right) = \widehat R_{0,t}(a_1) >b_0,
    \]
    and the negative bound follows similarly. By time reflection, this proves \eqref{eq:SHL-root-uniform-endpoints}.

    For continuous non-decreasing functions $f,g:\R\to\R$ with $f(\pm\infty)=g(\pm\infty)=\pm\infty$, \cite[Section~7]{norris2015weak} gives
    \[
        f^{-1}(\{y\}) = [f^\leftarrow(y-),f^\leftarrow(y)], \qquad (f^\leftarrow)^\leftarrow=f, \qquad (g\circ f)^\leftarrow = f^\leftarrow\circ g^\leftarrow.
    \]
    Hence, with  $Q_{s,t}:=R_{s,t}^\leftarrow,$    the map $Q_{s,t}$ is strictly increasing and right-continuous, and $   Q_{s,u}=Q_{t,u}\circ Q_{s,t}, 0\le s\le t\le u$. Moreover, letting $\Gamma(f):=\overline{\{(f(a),a):a\in\R\}}$ and
    \[
        \Gamma_{\mathrm{comp}}(f^\leftarrow) :=\{(x,y):f^\leftarrow(y-)\le x\le f^\leftarrow(y)\},
    \]
    one has $\Gamma_{\mathrm{comp}}(f^\leftarrow) =\{(y,x):(x,y)\in\Gamma(f)\}$. The local uniform convergence of $R^M_{s,t}$, together with
    \[
        \lim_{a\to+\infty}\inf_{0\le s\le t\le T}R_{s,t}(a)=+\infty, \qquad \lim_{a\to-\infty}\sup_{0\le s\le t\le T}R_{s,t}(a)=-\infty,
    \]
    gives, uniformly for $0\le s\le t\le T$,
    \begin{equation}\label{eq:SHL-inverse-graph-limit}
        d_H\!\left( \overline{(\tanh\times\tanh) \{(R^M_{s,t}(a),a):a\in\R\}}, \overline{(\tanh\times\tanh) \{(R_{s,t}(a),a):a\in\R\}} \right) \rightarrow 0 .
    \end{equation}
    The same convergence holds for the compositions over $[s,t)$, with the endpoint convention reversed.

    To define the maps simultaneously for all $s\le t$ on a single measurable event, set
    \[
        R_{s,t} :=
        \begin{cases}
        \displaystyle \lim_{m\to\infty}R^m_{s,t} & \text{in }C_{\mathrm{loc}}(\R), \quad\text{if the limit exists and } \displaystyle\lim_{a\to\pm\infty}R_{s,t}(a)=\pm\infty, \\[1ex]
        \mathrm{id} & \text{otherwise},
        \end{cases}
        \quad Q_{s,t}:=R_{s,t}^\leftarrow .
    \]
    The first case is determined by the countable conditions
    \[
        \lim_{m,n\to\infty}\sup_{|a|\le k}|R^m_{s,t}(a)-R^n_{s,t}(a)|=0, \qquad k\in\N,
    \]
    together with the corresponding endpoint conditions at rational labels. Thus  $(\omega,s,t,a)\mapsto\bigl(R_{s,t}(a)(\omega),Q_{s,t}(a)(\omega)\bigr)$  is jointly measurable, and, for each deterministic $s\le t$,
    \[
        R_{s,t},\,Q_{s,t} \quad\text{are }\; \sigma\!\left(\Pi|_{(s,t]\times\R}\right)\text{-measurable}.
    \]

    Fix $c\in\R$. For $j\le M\le j+1$, shifting the symmetric window by $c$ changes it only on
    \[
        S_j(c) := [j-|c|-1,j+|c|+1] \cup [-j-|c|-1,-j+|c|+1].
    \]
    Let $Z_j:=\Pi\bigl([0,T]\times S_j(c)\bigr). $ Since $ \sup_{j\ge1}\E[Z_j^2]<\infty, $ Markov's inequality and the Borel--Cantelli lemma give $ Z_j/j\rightarrow 0$ almost surely. The comparison in \cite[proof of Proposition~4.7]{berger2022growth} and \eqref{eq:SHL-cutoff-error} makes the error tend uniformly to zero on every interval containing no point with $|x|\le B$. There are only finitely many such points before time $T$.  At each such point, an error $d$ is bounded by $\omega_K(d)$. Since $\omega_K(0)=0$ and $\omega_K$ is continuous at $0$, after finitely many such points the resulting error still tends to $0$ as $d\to0$. Hence
    \[
        \sup_{\substack{j\le M\le j+1,\ 0\le s\le t\le T,\ |a|\le K}} \bigl|R^{M,c}_{s,t}(a+c)-c-R^M_{s,t}(a)\bigr|\rightarrow0 \qquad\text{a.s.}
    \]
    where $R^{M,c}_{s,t}$ is obtained after shifting each point $(u,x)$ to $(u,x+c)$. Letting $M\to\infty$, the corresponding limit map is $    a\mapsto R_{s,t}(a-c)+c.$ Taking generalized inverses gives  $ a\mapsto Q_{s,t}(a-c)+c.$

    Let  $\mathcal F_s:=\sigma\!\left(\Pi|_{[0,s]\times\R}\right), \Pi_s:=\{(u-s,x):(u,x)\in\Pi,\ u>s\}$, and $    G_s(u):=(\mathcal S_s^\lambda)^{-1}\circ\mathcal S_{s+u}^\lambda, u\ge0$.  For fixed $s$, the identity for $|x|\le M$ in \cite[proof of Proposition~4.5]{berger2022growth}, together with \cite[Theorem~4.3]{berger2022growth} and the time reversal in \cite[Section~2.1]{berger2022growth}, gives
    \[
        (\Pi_s,G_s)\stackrel{d}{=}(\Pi,\mathcal S^\lambda), \qquad (\Pi_s,G_s)\ \text{is independent of }\mathcal F_s.
    \]

    Let $Y_{s,t}$ denote the right-continuous boundary values in  \eqref{eq:rescaled-boundary-transition} for $N=1$.  For fixed $s,a$, \cite[proof of Theorem~7.2]{berger2022growth}    gives
    \[
        Y_{s,t}(a) = a+\int_{(s,t]\times\R} g_x\!\left(Y_{s,u-}(a)\right)\,\widetilde\Pi(du\,dx).
    \]
    For every $B<\infty$,
    \[
        \E\left[\sum_{\substack{(u,x)\in\Pi,\ s<u\le T,\ |x|\le B}} \ind_{\{Y_{s,u-}(a)=x\}}\right]=0.
    \]
    Hence, almost surely, $Y_{s,u-}(a)\neq x$, for every $(u,x)\in\Pi,\ s<u\le T,\ |x|\le B$.

    Choose $K>1+\sup_{s\le t\le T}|Y_{s,t}(a)|, B:=2K+2\lambda+1$, and set
    \[
        \delta:= 1\wedge \inf_{\substack{(u,x)\in\Pi \\
        s<u\le T,\ |x|\le B}} |Y_{s,u-}(a)-x|>0, \qquad \inf\varnothing:=\infty.
    \]
    By \eqref{eq:SHL-tail-splitting},
    \[
        \begin{aligned}
        Y_{s,t}(a)-a=\sum_{\substack{s<u\le t, \ |x|\le B}} g_x\!\left(Y_{s,u-}(a)\right) -A_B(t)+A_B(s)+\sum_{\substack{s<u\le t, \ |x|>B}} \widetilde w_x\!\left(Y_{s,u-}(a)\right).
        \end{aligned}
    \]
    Compare this with $(R^M_{s,t})^\leftarrow(a)$. On each interval containing no Poisson point with $|x|\le B$, the same argument as in \eqref{eq:SHL-cutoff-error}, using \eqref{eq:SHL-tail-splitting} for $\widetilde w_x$, controls the error up to the first exit from $[-K,K]$. At each such point, an error smaller than $\delta/2$ keeps the two trajectories on the same side of $x$, where $h_x$ is 1-Lipschitz. Since there are only finitely many such points, the error before the first exit is bounded by $C(\Pi,K,T)\varepsilon_M$, which is smaller than $1\wedge\delta/2$ for all sufficiently large $M$. Since $\sup_{s\le t\le T}|Y_{s,t}(a)|<K-1$, no exit occurs, and the comparison holds throughout $[s,T]$. Hence
    \[
        \sup_{s\le t\le T} \left|(R^M_{s,t})^\leftarrow(a)-Y_{s,t}(a)\right| \rightarrow0.
    \]

    Intersecting over rational $s,a$ and integer $T$, and using \eqref{eq:SHL-inverse-graph-limit}, gives $Q_{s,t}(a-)\le Y_{s,t}(a)\le Q_{s,t}(a)$. For rational $a_n\downarrow y$, $a_n>y$,
    \[
        Q_{s,t}(y) \le Q_{s,t}(a_n-) \le Y_{s,t}(a_n) \le Q_{s,t}(a_n) \rightarrow Q_{s,t}(y).
    \]
    By right continuity, $Y_{s,t}(y)=Q_{s,t}(y)$  for every $y$ and $t\ge s$; the left limits give the two values at each jump. Thus, for every fixed $s$, $Q_{s,\cdot}$ is a version of $U^{\mathcal S,1}_{s,\cdot}$.

    Now let $\tau$ be a finite stopping time for $(\mathcal F_s)$ and set $\tau_m:=2^{-m}\lceil2^m\tau\rceil\downarrow\tau$. Since $\tau_m$ takes countable values and $\mathcal F_\tau\subseteq\mathcal F_{\tau_m}$, the preceding increment identity applies at $\tau_m$. For every $u\in\mathbb Q_{>0}$,
    \[
        G_{\tau_m}(u) = (\mathcal S_{\tau_m}^\lambda)^{-1} \circ\mathcal S_{\tau_m+u}^\lambda \rightarrow (\mathcal S_\tau^\lambda)^{-1} \circ\mathcal S_{\tau+u}^\lambda = G_\tau(u)
    \]
    locally uniformly, while $\Pi_{\tau_m}\rightarrow\Pi_\tau$  against compactly supported continuous functions at positive times.  Therefore, for bounded $\mathcal F_\tau$-measurable $Z$ and bounded  continuous function $\Psi$ depending on finitely many such map values and  integrals,
    \begin{equation*}
        \E\!\left[Z\Psi(\Pi_\tau,G_\tau)\right] = \lim_{m\to\infty} \E\!\left[Z\Psi(\Pi_{\tau_m},G_{\tau_m})\right] = \E[Z]\,\E[\Psi(\Pi,\mathcal S^\lambda)].
    \end{equation*}
    Hence  $(\Pi_\tau,G_\tau)\stackrel{d}{=}(\Pi,\mathcal S^\lambda)\text{, and } (\Pi_\tau,G_\tau)$ is independent of $\mathcal F_\tau$. Rational $u>0$ and right continuity determine all $u\ge0$.  Therefore, the identity proved above for fixed $s$ applies after time $\tau$, giving $Y_{\tau,\tau+u}=Q_{\tau,\tau+u}, u\ge0$.

    Every Poisson point is the $n$-th arrival point in $\mathbb R_+\times[j,j+1)$, for some $j\in\Z$ and $n\ge1$, and its arrival time is a finite stopping time. Taking the countable intersection over  $j\in\Z$, $n\ge1$, and rational start times proves the    identity at every required arrival time. The interval $(s,t]$ excludes the point at the starting time. Rescaling gives the result for each fixed $N$, including the shift $Nc$ and the arrival times $\bar t_i$. Thus, \eqref{eq:SHL-grid-inverse} and the second identity in \eqref{eq:SHL-grid-cocycle} follow.
    \end{proof}


    \subsection{Proof of Lemma~\ref{lem:gap-exterior-geometry}}
    \label{appendix:gap-exterior-geometry}
    \begin{proof}
    For $\mathcal A=\mathcal S$, the slit moment estimate in the proof of Lemma~\ref{lem:future-label-flatness}, applied at scale $\max\{N_0,\lceil\sqrt r\rceil\}$ with $p=8$, gives
    \[
        \sup_{x\in\mathbb R} \mathbb E\!\left[ \sup_{0\le \vartheta\le\lambda} \left| \mathcal S^\lambda_{r-}(x+\mathrm i\vartheta) -x-\mathrm i\vartheta-\mathrm i v_{\mathcal S}r \right|^8 \right] \le C(1+r^2).
    \]
    Fix $R,H>0$. If $r\ge N^2G$, $G>2H/v_{\mathcal S}$, and $\mathcal S^\lambda_{r-}\bigl(x+\mathrm i[0,\lambda]\bigr) \cap\bigl([-NR,NR]+\mathrm i[0,N^2H]\bigr)\ne\varnothing$, then, for some $\vartheta\in[0,\lambda]$,
    \[
        \left| \mathcal S^\lambda_{r-}(x+\mathrm i\vartheta) -x-\mathrm i\vartheta-\mathrm i v_{\mathcal S}r \right| \ge c_\lambda\bigl\{r+(|x|-NR)_+\bigr\}.
    \]
    Hence Mecke's identity gives, for all sufficiently large $N$,
    \[
        \E\left[\#\{i:\bar t_i\ge G,\ e_i^{\mathcal S} \cap([-NR,NR]+\mathrm i[0,N^2H])\ne\varnothing\}\right] \le C\int_{N^2G}^{\infty}(NRr^{-6}+r^{-5})\,dr=o_N(1).
    \]
    For fixed $N$, the same integral is finite after increasing its lower endpoint. On bounded arrival-time intervals, the spatial tail is integrable, as in \eqref{eq:SHL-remote-particle-summation}. Thus only finitely many particles meet each compact window almost surely.  For $\mathcal C$ and fixed $N$, \eqref{eq:SHL-BDG-Rosenthal} and the proof of Lemma~\ref{lem:future-label-flatness} give
    \[
        \sup_{x\in\mathbb T_N}\E\!\left[ \sup_{0\le\vartheta\le\lambda} |\mathcal C^{N,\lambda}_{r-}(x+\mathrm i\vartheta) -x-\mathrm i\vartheta-\mathrm i v_Nr|^8 \right]\le C_{N,\lambda}(1+r^4),
    \]
    where $v_N>0$ is defined in \eqref{eq:centered-reverse-composition}. The probability that a slit arriving at time $r$ meets a fixed bounded-height window is therefore $O_{N,\lambda,H}(r^{-4})$. Mecke's identity proves local finiteness, and Remark~\ref{rmk:radial-from-cylindrical} gives the radial case. By \cite[Theorem~3.1]{norris2012hastings} and Remark~\ref{rmk:radial-from-cylindrical}, for each fixed bounded-height window and sufficiently large fixed $G$, with probability tending to one as $N\to\infty$, no rescaled particle with $\bar t_i\ge G$ meets the window.

    Recall that $\varepsilon_N=N^{-1}$. By local finiteness, divide each edge into arcs of diameter less than $\varepsilon_N/10$. Choose compact sets $K_1\subset K_2\subset\cdots$ covering the domain. On each $K_m$, surround these arcs, the vertices, and the initial boundary by thin regions bounded by finite straight line segments, keeping the choices on $K_{m-1}$ unchanged. Choose the regions small enough that those around disjoint parts of the tree remain disjoint, and every interior point can be joined to the tree by a curve in the interior of diameter at most $\varepsilon_N$. Their union is a closed neighbourhood with locally piecewise-linear boundary and introduces no new connections or holes. Near each tip, choose the region to allow the segment $\ell_i$ specified in Subsection~\ref{S.HL-models}. Make the choices periodic in the compact cases. Measurability follows by selecting the first suitable choice from a fixed countable collection of regions with rational vertices.

    For a continuous curve $\eta:[a,b]\to E^{\mathcal A}\times[0,\infty)$, define its length by
    \[
        \operatorname{len}(\eta):=\sup_{a=t_0<\cdots<t_m=b}\sum_{j=1}^m |\eta(t_j)-\eta(t_{j-1})|.
    \]
    Write
    \[
        F_N=(E^{\mathcal A}\times[0,\infty))          \setminus\operatorname{int}\mathscr T_N^{\mathcal A,\varepsilon_N}, \qquad d_{\mathrm{ext}}(z,w) =\inf_{\eta\subset F_N:\,z\to w}\operatorname{len}(\eta)
    \]
    for $z,w$ in the same component of $F_N$. Its boundary is locally formed by finite straight line segments, so sufficiently close points in a component can be joined inside  $F_N$ by short curves. Thus $d_{\mathrm{ext}}$ gives the usual topology there. Since each component is closed and $d_{\mathrm{ext}}$ is at least the usual distance, its closed    $d_{\mathrm{ext}}$-balls are compact. The Arzel\`a--Ascoli theorem then gives shortest curves between points and from a point to any nonempty closed set in the same component.

    Fix $h<\infty$. In the stationary case, the moment bound above and the bounded spatial range of each tree path give positive   probability of a finger reaching above $h$ within a fixed bounded interval. Horizontal ergodicity \cite[Proposition~4.7]{berger2022growth} gives such fingers     arbitrarily far to the left and right. In the compact cases, local finiteness and infinitely many arrival points give tips of     unbounded height; repeating their finite tree paths with the spatial period on $\mathbb R\times[0,\infty)$ gives the same     conclusion. Hence a component of $F_N$ with bounded height would also be bounded horizontally, leaving a bounded hole in the neighbourhood. This contradicts its construction,  so every component reaches arbitrarily large heights.

    For each integer $L>h_i+\varepsilon_N$, choose a shortest curve $\gamma_i^L$ from $z_i$ to height $L$, parametrized by arc length.  Its length is at least $L-z_i^{(2)}$ and tends to infinity.  Since the curves have unit speed and closed $d_{\mathrm{ext}}$-balls are compact, a subsequence converges locally uniformly to $\gamma_i$ satisfying
    \[
        d_{\mathrm{ext}}\bigl(\gamma_i(a),\gamma_i(b)\bigr)=b-a,  \qquad 0\le a\le b.
    \]
    If $\gamma_i$ had bounded height, the same fingers would confine it to a compact set, contradicting $d_{\mathrm{ext}}(z_i,\gamma_i(r))=r$. Hence    $\sup_{r\ge0}\gamma_i^{(2)}(r)=\infty$.

    The distances $d_{\mathrm{ext}}$ are measurable by approximation with curves having finite straight segments and rational  intermediate vertices. If $\eta_n\to\eta$ uniformly, then $\operatorname{len}(\eta)\le\liminf_n\operatorname{len}(\eta_n)$. Since the condition that a curve stays in $F_N$ is closed and measurable, the shortest curves to each height $L$ form a     measurable nonempty compact set, so $\gamma_i^L$ can be chosen measurably. After extending these curves constantly beyond their endpoints, their locally uniform subsequential limits also form a measurable nonempty compact set; choose $\gamma_i$ from this set measurably. Let $g_i$ follow $\ell_i$ and then $\gamma_i$.     Equation~\eqref{eq:gap-height-clock} gives \eqref{eq:gap-strict-exterior}.

    Fix $N,T$ and a bounded spatial interval. Fingers on either side reaching above $v_{\mathcal A}T$ keep any gap entering the interval before $T$ between them: by \eqref{eq:gap-strict-exterior}, the gap     cannot meet the tree or rise above these fingers. In the compact cases, use periodic copies on $\mathbb R\times[0,\infty)$.     Thus its starting tip lies in a compact set, which contains only finitely many particle tips by local finiteness. Every convergent    sequence with finite limiting starting time and position has a subsequence with fixed particle index; its limit satisfies     \eqref{eq:gap-strict-exterior}. If starting times tend to infinity, the extended paths are constant on bounded time intervals. If    starting positions tend to spatial infinity, the same fingers prevent entry into bounded spatial intervals before $T$. The compactification in $D_{\mathcal A}^{\uparrow}$ then gives subsequential convergence. Hence $\mathbf{Gap}^{\mathcal A,N}$, the closure of a countable collection of measurable paths, is measurable and compact.

    To prove \eqref{eq:gap-spatial-localization}, fix $L>4R$ and choose particles near $(\pm L/2,T+1)$ as in Subsection~\ref{SS.binary-decomposition-section6}.  By Propositions~\ref{prop:fdd-selected-paths-section7} and \ref{prop:process-tip-geometric-transfer}, with probability at least  $1-C_TL^{-2}-o_N(1)$ their fingers stay respectively in $[-3L/4,-L/4]$ and $[L/4,3L/4]$ and reach above $v_{\mathcal S}T$. On this event, any gap entering $[-R,R]$ by time $T$ stays in $[-L,L]$ up to time $T$, by \eqref{eq:gap-strict-exterior}. Thus the probability in \eqref{eq:gap-spatial-localization} has limit superior at most     $C_TL^{-2}$. Letting $L\to\infty$ proves the claim.

    For \eqref{eq:gap-height-approximation}, use the same event in the  stationary case and set $W=[-L,L]\times[0,v_{\mathcal S}T+1]$.    Any segment of a relevant exterior curve that stays at height at most $v_{\mathcal S}T$ and contains a point traced by the gap by time $T$ lies in $W$, since it cannot meet either selected finger. In the compact cases set $W=E^{\mathcal A}\times[0,v_{\mathcal A}T+1]$ and use periodic copies in $\mathbb R\times[0,\infty)$. The estimates above and     \eqref{eq:SHL-remote-particle-summation} show that, for suitable fixed $G,B<\infty$, with probability $1-o_N(1)$ every particle     meeting $W$ has $\bar t_j\le G$ and, if $\mathcal A=\mathcal S$, $|\bar x_j|\le B$. By Lemma~\ref{lem:completed-cell-oscillation} and the proof of Proposition~\ref{prop:process-tip-geometric-transfer}, there is $\alpha_N\pto0$ such that, for all these particles,
    \[
        \left|\bigl(\chi_{\mathcal A}\mathfrak p_N^{\mathcal A}z\bigr)^{(2)}               -v_{\mathcal A}\bar t_j\right|\le\alpha_N, \quad z\in\mathcal P_j^{\mathcal A,N}.
    \]
    In the compact cases, $\bar t_j\le G$ allows us to pass from the $\tanh$ estimate to this height bound. Thus, if $z,w\in W$ lie on edges $j\preceq k$, with $z$ on the tree path from $w$ to the initial boundary,
    \[
        z^{(2)}-w^{(2)} \le v_{\mathcal A}(\bar t_j-\bar t_k)+2\alpha_N\le2\alpha_N.
    \]
    Following \cite[proof of Theorem~8.2]{norris2012hastings}, suppose a part of $\gamma_i$ containing a point traced by the gap by time $T$ starts and ends at height $H\le v_{\mathcal A}T$, stays below $H$ between its endpoints, and reaches below $H-2\alpha_N-2\varepsilon_N$.  A shortest curve has no self-intersections, so this part and the horizontal segment joining its endpoints enclose a region $Q$. Choose $h<H-2\alpha_N-2\varepsilon_N$ and a horizontal interval     $I$ at height $h$ with endpoints on $\gamma_i$ and interior in $Q$. Since $\gamma_i$ is shortest,  $I$ must meet $\operatorname{int}\mathscr T_N^{\mathcal A,\varepsilon_N}$; otherwise $I$ would shorten the curve. By the construction of the neighbourhood, a point $w$ of this intersection can be joined inside it to $w'\in\mathscr T_N^{\mathcal A}$ by a curve of diameter at most $\varepsilon_N$. This curve cannot meet $\gamma_i$ or reach height $H$, so $w'\in Q$ and $(w')^{(2)}\le h+\varepsilon_N$. The tree path from $w'$ to the initial boundary must leave $Q$ through the horizontal segment at height $H$, at some point $z$. Since $z,w'\in W$ (or its periodic copy), the preceding bound gives
    \[
        2\alpha_N\ge z^{(2)}-(w')^{(2)} \ge H-h-\varepsilon_N>2\alpha_N,
    \]
    a contradiction.

    At any point of $\gamma_i$ traced by the gap by time $T$, its maximum height up to that point is at most $v_{\mathcal A}T$ by \eqref{eq:gap-height-clock}. If the point lies below this maximum, the segment from the last preceding time at that height to the first subsequent return is one of the segments considered above; the return exists since $\sup_{r\ge0}\gamma_i^{(2)}(r)=\infty$. Thus the height difference is at most $2\alpha_N+2\varepsilon_N$. The initial segment $\ell_i$ has length at most $\varepsilon_N$. Hence, for $r=\theta_i^{-1}(u)$ and $\sigma_i\le u\le T$,
    \[
        \max_{0\le a\le r}g_i^{(2)}(a)-g_i^{(2)}(r) \le2\alpha_N+3\varepsilon_N.
    \]
    Equation~\eqref{eq:gap-height-clock} now yields
    \[
        0\le v_{\mathcal A}u-g_i^{(2)}(r)       =\max_{0\le a\le r}g_i^{(2)}(a)-g_i^{(2)}(r)+v_{\mathcal A}\varepsilon_N\frac{r}{1+r}\le2\alpha_N+(3+v_{\mathcal A})\varepsilon_N.
    \]
    This bound holds for all the paths in     \eqref{eq:gap-height-approximation}. Let $N\to\infty$ and, in the stationary case, then $L\to\infty$ to remove the error     $C_TL^{-2}$. This proves \eqref{eq:gap-height-approximation}.
    \end{proof}

    \bibliographystyle{alpha}
    \bibliography{bibliography}

@phdthesis{arratia1979coalescing,
  author = {Arratia, Richard A.}, title = {Coalescing {B}rownian {M}otions on the {L}ine},
  school = {University of Wisconsin--Madison}, year = {1979}
}

@article{berestycki2015coalescing,
  author = {Berestycki, Nathana{\"e}l and Garban, Christophe and Sen, Arnab},
  title = {Coalescing {B}rownian flows: a new approach}, journal = {Annals of Probability},
  volume = {43}, number = {6}, pages = {3177--3215}, year = {2015},
  doi = {10.1214/14-AOP957}, eprint = {1307.4313v3}, archivePrefix = {arXiv}
}

@article{coletti2009scaling,
  author = {Coletti, C. F. and Fontes, L. R. G. and Dias, E. S.},
  title = {Scaling limit for a drainage network model}, journal = {Journal of Applied Probability},
  volume = {46}, number = {4}, pages = {1184--1197}, year = {2009},
  doi = {10.1239/jap/1261670696}, eprint = {0809.3454}, archivePrefix = {arXiv}
}

@article{coletti2014convergence,
  author = {Coletti, Cristian and Valle, Glauco},
  title = {Convergence to the {B}rownian web for a generalization of the drainage network model},
  journal = {Annales de l'Institut Henri Poincar\'e, Probabilit\'es et Statistiques},
  volume = {50}, number = {3}, pages = {899--919}, year = {2014}
}

@article{coupier2021directedspanningforest,
  author = {Coupier, David and Saha, Kumarjit and Sarkar, Anish and Tran, Viet Chi},
  title = {The 2{D}-directed spanning forest converges to the {B}rownian web},
  journal = {Annals of Probability}, volume = {49}, number = {1}, pages = {435--484}, year = {2021}
}

@article{ferrari2005two,
  author = {Ferrari, Pablo A. and Fontes, Luiz R. G. and Wu, Xian-Yuan},
  title = {Two-dimensional {P}oisson trees converge to the {B}rownian web},
  journal = {Annales de l'Institut Henri Poincar\'e, Probabilit\'es et Statistiques},
  volume = {41}, number = {5}, pages = {851--858}, year = {2005}
}

@article{newman2005convergence,
  author = {Newman, Charles M. and Ravishankar, Krishnamurthi and Sun, Rongfeng},
  title = {Convergence of coalescing nonsimple random walks to the {B}rownian web},
  journal = {Electronic Journal of Probability}, volume = {10}, pages = {21--60}, year = {2005},
  doi = {10.1214/EJP.v10-235}
}

@article{berger2022growth,
  author = {Berger, Noam and Procaccia, Eviatar B. and Turner, Amanda},
  title = {Growth of {S}tationary {H}astings--{L}evitov}, journal = {Annals of Applied Probability},
  volume = {32}, number = {5}, pages = {3331--3360}, year = {2022},
  doi = {10.1214/21-AAP1761}, eprint = {2008.05792v1}, archivePrefix = {arXiv}
}

@book{last2017lectures,
  author = {Last, G{\"u}nter and Penrose, Mathew},
  title = {Lectures on the {P}oisson {P}rocess},
  series = {Institute of Mathematical Statistics Textbooks}, volume = {7},
  publisher = {Cambridge University Press}, address = {Cambridge}, year = {2017},
  doi = {10.1017/9781316104477}, isbn = {9781107088016}
}

@misc{chen2025one,
  author = {Chen, Guanyi and Procaccia, Eviatar B. and Zong, Yuxuan},
  title = {One-arm domination time in {C}ylindrical {H}astings--{L}evitov(0)},
  year = {2025}, eprint = {2507.11028v2}, archivePrefix = {arXiv}, note = {{arXiv}:2507.11028v2. Version 2, revised 21 June 2026}
}

@article{coupier2019directed,
  author = {Coupier, David and Marckert, Jean-Fran\c{c}ois and Tran, Viet Chi},
  title = {Directed, cylindric and radial {B}rownian webs}, journal = {Electronic Journal of Probability},
  volume = {24}, number = {20}, pages = {1--48}, year = {2019},
  doi = {10.1214/18-EJP255}, eprint = {1707.08336}, archivePrefix = {arXiv}
}

@incollection{coupier2020collision,
  author = {Coupier, David and Saha, Kumarjit and Sarkar, Anish and Tran, Viet Chi},
  title = {Collision {T}imes of {R}andom {W}alks and {A}pplications to the {B}rownian {W}eb},
  booktitle = {Genealogies of {I}nteracting {P}article {S}ystems},
  series = {Lecture Notes Series, Institute for Mathematical Sciences, National University of Singapore},
  volume = {38}, publisher = {World Scientific}, pages = {267--293}, year = {2020},
  doi = {10.1142/9789811206092_0007}, eprint = {1508.06919v4}, archivePrefix = {arXiv}
}

@article{fontes2004brownian,
  author = {Fontes, Luiz Renato G. and Isopi, Marco and Newman, Charles M. and Ravishankar, Krishnamurthi},
  title = {The {B}rownian web: characterization and convergence}, journal = {Annals of Probability},
  volume = {32}, number = {4}, pages = {2857--2883}, year = {2004}, doi = {10.1214/009117904000000568},
  eprint = {math/0311254v2}, archivePrefix = {arXiv}
}

@incollection{schertzer2017brownian,
  author = {Schertzer, Emmanuel and Sun, Rongfeng and Swart, Jan M.},
  title = {The {B}rownian web, the {B}rownian net, and their universality},
  booktitle = {Advances in {D}isordered {S}ystems, {R}andom {P}rocesses and {S}ome {A}pplications},
  publisher = {Cambridge University Press}, pages = {270--368}, year = {2017},
  doi = {10.1017/9781316403877.007},
  eprint = {1506.00724v3}, archivePrefix = {arXiv}
}

@article{sun2008brownian,
  author = {Sun, Rongfeng and Swart, Jan M.},
  title = {The {B}rownian net}, journal = {Annals of Probability},
  volume = {36}, number = {3}, pages = {1153--1208}, year = {2008},
  doi = {10.1214/07-AOP357}, eprint = {math/0610625v5}, archivePrefix = {arXiv}
}

@article{hastings1998laplacian,
  author = {Hastings, Matthew B. and Levitov, Leonid S.}, title = {Laplacian growth as one-dimensional turbulence},
  journal = {Physica D: Nonlinear Phenomena}, volume = {116}, number = {1--2}, pages = {244--252}, year = {1998},
  doi = {10.1016/S0167-2789(97)00244-3}, eprint = {cond-mat/9607021v1}, archivePrefix = {arXiv}
}

@article{norris2012hastings,
  author = {Norris, James and Turner, Amanda}, title = {Hastings--{L}evitov aggregation in the small-particle limit},
  journal = {Communications in Mathematical Physics}, volume = {316}, number = {3}, pages = {809--841}, year = {2012},
  doi = {10.1007/s00220-012-1552-6}, eprint = {1106.3546v2}, archivePrefix = {arXiv}
}

@article{norris2015weak,
  author = {Norris, James and Turner, Amanda}, title = {Weak convergence of the localized disturbance flow to the coalescing {B}rownian flow},
  journal = {Annals of Probability}, volume = {43}, number = {3}, pages = {935--970}, year = {2015},
  doi = {10.1214/13-AOP845}, eprint = {1106.3252v4}, archivePrefix = {arXiv}
}

@misc{berger2025logarithmic,
  author = {Berger, Noam and Procaccia, Eviatar B.}, title = {Logarithmic fluctuations of {S}tationary {H}astings--{L}evitov},
  year = {2025}, eprint = {2502.03554}, archivePrefix = {arXiv},
  note = {{arXiv}:2502.03554}
}

@misc{procaccia2023cylindrical,
  author = {Procaccia, Eviatar B. and Zhuchenko, Anna}, title = {Cylindrical {H}astings {L}evitov},
  year = {2023}, eprint = {2301.12737v1}, archivePrefix = {arXiv},
  note = {{arXiv}:2301.12737v1. Version 1, 30 January 2023}
}

@article{kesten1987hitting,
  author = {Kesten, Harry}, title = {Hitting probabilities of random walks on $\mathbb{Z}^d$},
  journal = {Stochastic Processes and their Applications}, volume = {25}, pages = {165--184}, year = {1987},
  doi = {10.1016/0304-4149(87)90196-7}
}

@article{kesten1987long,
  author = {Kesten, Harry}, title = {How long are the arms in {DLA}?}, journal = {Journal of Physics A: Mathematical and General},
  volume = {20}, number = {1}, pages = {L29--L33}, year = {1987}, doi = {10.1088/0305-4470/20/1/007}
}

@article{kol1998solution,
  author = {Kol, Boaz and Aharony, Amnon}, title = {Solution of diffusion limited aggregation in a narrow cylindrical geometry},
  journal = {Physical Review E}, volume = {58}, number = {4}, pages = {4716--4729}, year = {1998},
  doi = {10.1103/PhysRevE.58.4716}, eprint = {cond-mat/9802048}, archivePrefix = {arXiv}
}

@article{benjamini2008diffusion,
  author = {Benjamini, Itai and Yadin, Ariel}, title = {Diffusion limited aggregation on a cylinder},
  journal = {Communications in Mathematical Physics}, volume = {279}, number = {1}, pages = {187--223}, year = {2008},
  doi = {10.1007/s00220-008-0424-6}, eprint = {math/0701201v2}, archivePrefix = {arXiv}
}

@article{marchetti2012stationary,
  author = {Marchetti, Riccardo and Taloni, Alessandro and Caglioti, Emanuele and Loreto, Vittorio and Pietronero, Luciano},
  title = {Stationary growth and unique invariant harmonic measure of cylindrical diffusion limited aggregation},
  journal = {Physical Review Letters}, volume = {109}, number = {6}, pages = {065501}, year = {2012},
  doi = {10.1103/PhysRevLett.109.065501}
}

@article{procaccia2019stationary,
  author = {Procaccia, Eviatar B. and Zhang, Yuan}, title = {Stationary harmonic measure and {DLA} in the upper half plane},
  journal = {Journal of Statistical Physics}, volume = {176}, number = {4}, pages = {946--980}, year = {2019},
  doi = {10.1007/s10955-019-02327-y}, eprint = {1711.01011v2}, archivePrefix = {arXiv}
}

@article{procaccia2020stationary,
  author = {Procaccia, Eviatar B. and Ye, Jiayan and Zhang, Yuan}, title = {Stationary {DLA} is well defined},
  journal = {Journal of Statistical Physics}, volume = {181}, number = {4}, pages = {1089--1111}, year = {2020},
  doi = {10.1007/s10955-020-02619-8}, eprint = {1907.00381v1}, archivePrefix = {arXiv}
}

@article{procaccia2021dimension,
  author = {Procaccia, Eviatar B. and Procaccia, Itamar}, title = {Dimension of diffusion-limited aggregates grown on a line},
  journal = {Physical Review E}, volume = {103}, number = {2}, pages = {L020101}, year = {2021},
  doi = {10.1103/PhysRevE.103.L020101}, eprint = {2010.05709v1}, archivePrefix = {arXiv}
}

@article{mu2022scaling,
  author = {Mu, Yingxin and Procaccia, Eviatar B. and Zhang, Yuan}, title = {Scaling limit of {DLA} on a long line segment},
  journal = {Transactions of the American Mathematical Society}, volume = {375}, number = {12}, pages = {8769--8806}, year = {2022},
  doi = {10.1090/tran/8771}, eprint = {1912.02370v1}, archivePrefix = {arXiv}
}

@article{toth1998true,
  author = {T\'oth, B\'alint and Werner, Wendelin}, title = {The true self-repelling motion},
  journal = {Probability Theory and Related Fields}, volume = {111}, number = {3}, pages = {375--452}, year = {1998},
  doi = {10.1007/s004400050172}
}

@article{witten1981diffusion,
  author = {Witten, Jr., Thomas A. and Sander, Leonard M.}, title = {Diffusion-limited aggregation, a kinetic critical phenomenon},
  journal = {Physical Review Letters}, volume = {47}, number = {19}, pages = {1400--1403}, year = {1981},
  doi = {10.1103/PhysRevLett.47.1400}
}

@book{garnett2007bounded,
  author = {Garnett, John B.},
  title = {Bounded {A}nalytic {F}unctions},
  edition = {Revised first},
  series = {Graduate Texts in Mathematics}, volume = {236},
  publisher = {Springer}, address = {New York}, year = {2007},
  doi = {10.1007/0-387-49763-3}, isbn = {9780387336213}
}

@incollection{marinelli2014maximal,
  author = {Marinelli, Carlo and R{\"o}ckner, Michael},
  title = {On {M}aximal {I}nequalities for {P}urely {D}iscontinuous {M}artingales in {I}nfinite {D}imensions},
  booktitle = {S{\'e}minaire de {P}robabilit{\'e}s {XLVI}},
  series = {Lecture Notes in Mathematics}, volume = {2123},
  publisher = {Springer International Publishing}, address = {Cham},
  pages = {293--315}, year = {2014},
  doi = {10.1007/978-3-319-11970-0_10},
  eprint = {1308.2648}, archivePrefix = {arXiv}
}

@article{rohde2005some,
  title={Some remarks on {L}aplacian growth},
  author={Rohde, Steffen and Zinsmeister, Michel},
  journal={Topology and its Applications},
  volume={152},
  number={1-2},
  pages={26--43},
  year={2005},
  publisher={Elsevier}
}

@article{zhou2008stepping,
  title={Stepping-{S}tone {M}odel with {C}ircular {B}rownian {M}igration},
  author={Zhou, Xiaowen},
  journal={Canadian Mathematical Bulletin},
  volume={51},
  number={1},
  pages={146--160},
  year={2008},
  publisher={Cambridge University Press},
  doi={10.4153/CMB-2008-017-2}
}

@article{carleson2001aggregation,
  author  = {Carleson, Lennart and Makarov, Nikolai},
  title   = {Aggregation in the {P}lane and {Loewner}'s {E}quation},
  journal = {Communications in Mathematical Physics},
  volume  = {216},
  pages   = {583--607},
  year    = {2001}
}

@article{norris2024stability,
  author  = {Norris, James and Silvestri, Vittoria and Turner, Amanda},
  title   = {Stability of {R}egularized {Hastings--Levitov} {A}ggregation in the {S}ubcritical {R}egime},
  journal = {Communications in Mathematical Physics},
  volume  = {405},
  number  = {3},
  pages   = {74},
  year    = {2024}
}

@article{berger2014stretched,
  author = {Berger, Noam and Kagan, Jacob J. and Procaccia, Eviatar B.},
  title = {Stretched {IDLA}},
  journal = {ALEA, Latin American Journal of Probability and Mathematical Statistics},
  volume = {11}, number = {2}, pages = {471--481}, year = {2014},
  eprint = {1209.3112}, archivePrefix = {arXiv}
}

@article{Johansson2012scaling,
  author  = {Johansson Viklund, Fredrik and Sola, Alan and Turner, Amanda},
  title   = {Scaling {L}imits of {A}nisotropic {Hastings--Levitov} {C}lusters},
  journal = {Annales de l'Institut Henri Poincar{\'e}, Probabilit{\'e}s et Statistiques},
  volume  = {48},
  number  = {1},
  pages   = {235--257},
  year    = {2012},
  doi     = {10.1214/10-AIHP395}
}

    \end{document}